\documentclass[a4paper,10pt]{article}
\usepackage[utf8]{inputenc}
\usepackage{amsthm}
\usepackage{amsfonts}
\usepackage{amssymb}	
\usepackage{amsmath}
\allowdisplaybreaks
\usepackage{mathtools}
\usepackage[english]{babel}
\usepackage{color}
\usepackage{cite}
\usepackage[numbers,sort]{natbib}
\usepackage{slashed}
\usepackage{enumerate}
\usepackage{graphicx}
\usepackage{systeme}
\usepackage{dsfont}
\usepackage{relsize}
\usepackage[margin=0.90in]{geometry}
\usepackage{float}
\usepackage{subcaption}
\usepackage{tikz}
\usepackage[labelformat=simple]{subcaption}

\usepackage{graphicx}
\usepackage{subcaption}
\usepackage{bmpsize}
\usepackage{epstopdf}
\usepackage{bbm}
\usepackage{xcolor,etoolbox}
\usepackage{titling}
\usepackage{authblk}
\newcommand*\samethanks[1][\value{footnote}]{\footnotemark[#1]}
\usepackage{blindtext,graphicx}
\usepackage[absolute]{textpos}
\usepackage{physics}
\usepackage[hang,flushmargin]{footmisc}

\usepackage{xfrac}
\usepackage{nicefrac}
\usepackage{esvect}
\usepackage{cases}
\usepackage{empheq}
\usepackage{stmaryrd}
\usepackage{cancel}
\usepackage[linguistics]{forest}
\usepackage{url}
\theoremstyle{definition}

\makeatletter
\newcommand{\pushright}[1]{\ifmeasuring@#1\else\omit\hfill$\displaystyle#1$\fi\ignorespaces}
\newcommand{\pushleft}[1]{\ifmeasuring@#1\else\omit$\displaystyle#1$\hfill\fi\ignorespaces}
\makeatother
\title{\vspace{-2cm}Analytical solution of the cylindrical torsion problem for the relaxed micromorphic continuum and other generalized continua (including full derivations)}
\author{
	Gianluca Rizzi\thanks{corresponding author, GEOMAS, INSA-Lyon, Universit\'e de Lyon, 20 avenue Albert Einstein,	69621, \\ \hspace*{0.55cm} Villeurbanne cedex, France, gianluca.rizzi@insa-lyon.fr},
	\quad\quad
	Geralf Hütter\thanks{Institute of Mechanics and Fluid Dynamics, Technische Universität Bergakademie Freiberg,\\ \hspace*{0.55cm} Lampadiusstr. 4, 09596 Freiberg, Germany},
	\quad\quad
	Hassam Khan\thanks{Fakultät für Mathematik, Universität Duisburg-Essen, Thea-Leymann-Straße 9, 45127 Essen, Germany},
	\quad\quad
	Ionel Dumitrel Ghiba\thanks{Department of Mathematics, Alexandru Ioan Cuza University of Ia\c si, Blvd. Carol I, no. 11, 700506 Ia\c si, Romania; and Octav Mayer Institute of Mathematics of the Romanian Academy, Ia\c si Branch, 700505 Ia\c si},
	\\
	Angela Madeo\samethanks[1]
	\quad and \quad
	Patrizio Neff\thanks{Head of Chair for Nonlinear Analysis and Modelling, Fakultät für Mathematik, Universität Duisburg-Essen, \\ \hspace*{0.55cm} Thea-Leymann-Straße 9, 45127 Essen, Germany}}

\thanksmarkseries{arabic}
\date{\today}
\begin{document}
\maketitle
\begin{abstract}
	We solve the St.Venant torsion problem for an infinite cylindrical rod whose behaviour is described by a family of isotropic generalized continua, including the relaxed micromorphic and classical micromorphic model.
	The results can be used to determine the  material parameters of these models.
	Special attention is given to the possible nonphysical stiffness singularity for a vanishing rod diameter, since slender specimens are in general described as stiffer. 
\end{abstract}

\section{Introduction }
\label{sec:intro}
In this paper we continue our endeavour to find analytical solutions to simple boundary value problems for families of generalized continua \cite{rizzi2021shear,rizzi2021bending,altenbach2019higher}.
The focus is on non-homogeneous solutions that on one side activate the additional deformations modes offered by generalized continua (curvature terms) and which may be used, on the other side, in calibrating the additional (many) material parameters.
The renewed interest in models of generalized continua comes in part from the fact that for small specimens  one may observe size-effects, not accounted for by linear Cauchy elasticity.
On the other hand, the description of man-made architecture materials/meta-materials need generalized continua to capture frequency band-gaps in the dynamic range, a prominent example being given by the relaxed micromorphic model \cite{rizzi2020towards,rizzi2020exploring,aivaliotis2021frequency,madeo2018relaxed,aivaliotis2019microstructure}.

Here, we consider the static St.~Venant torsion problem.
Since we aim at identifying material parameters, let us first review what can be said for isotropic linear elasticity.
\subsection{Material parameters in linear elasticity vs. generalized continua}
The determination of the two constitutive material parameters in isotropic linear elasticity can be achieved in several different ways.
For example the Young's modulus and Poisson's ratio
\begin{align}
E_{\tiny \mbox{macro}}
&=
\frac{9 \kappa_{\tiny \mbox{macro}}  \, \mu_{\tiny \mbox{macro}} }{3 \kappa_{\tiny \mbox{macro}} + \mu_{\tiny \mbox{macro}}} \, ,
\qquad
&\nu_{\tiny \mbox{macro}}
&=
\frac{3 \kappa_{\tiny \mbox{macro}} -2 \mu_{\tiny \mbox{macro}} }{2 (3 \kappa_{\tiny \mbox{macro}} +\mu_{\tiny \mbox{macro}} )} \, ,
\\*
\lambda_{\tiny \mbox{macro}}
&=
\frac{3\kappa_{\tiny \mbox{macro}} - 2\mu_{\tiny \mbox{macro}}}{3} \, ,
\qquad
&\kappa_{\tiny \mbox{macro}}
&=
\frac{2\mu_{\tiny \mbox{macro}} + 3\lambda_{\tiny \mbox{macro}}}{3} \, ,
\label{eq:real_parameter_macro}
\end{align}
can be uniquely determined from a \textit{homogeneous} macroscopic tension-compression test.
Moreover, the shear modulus $\mu_{\tiny \mbox{macro}}$ and the Young's modulus $E_{\tiny \mbox{macro}}$ can also be identified from the \textit{inhomogeneous} torsion and bending test, respectively.
Indeed, the classical torsional stiffness (per unit length) of a circular rod is given by
\begin{equation}
T_{\mbox{\tiny macro}} =
\mu_{\mbox{\tiny macro}} \, I_{p} =
\mu_{\mbox{\tiny macro}} \, \frac{\pi R^4}{2}
\, ,
\label{eq:torsional_stiffness_intro}
\end{equation}
and the bending stiffness (per unit length and per unit thickness) in cylindrical bending \cite{rizzi2021bending} is equivalent to 
\begin{equation}
D_{\mbox{\tiny macro}} =
\frac{h^3}{12} \frac{E_{\tiny \mbox{macro}}}{\left(1-\nu_{\tiny \mbox{macro}}^2\right)} =
\frac{h^3}{12} \frac{4 \mu _{\mbox{\tiny macro}} \left(3 \kappa _{\mbox{\tiny macro}}+\mu _{\mbox{\tiny macro}}\right)}{3 \kappa _{\mbox{\tiny macro}}+4 \mu _{\mbox{\tiny macro}}} \, .
\label{eq:bending_stiffness_intro}
\end{equation}
\begin{figure}[H]
	\begin{subfigure}{0.6\textwidth}
		\centering
		\includegraphics[height=5cm]{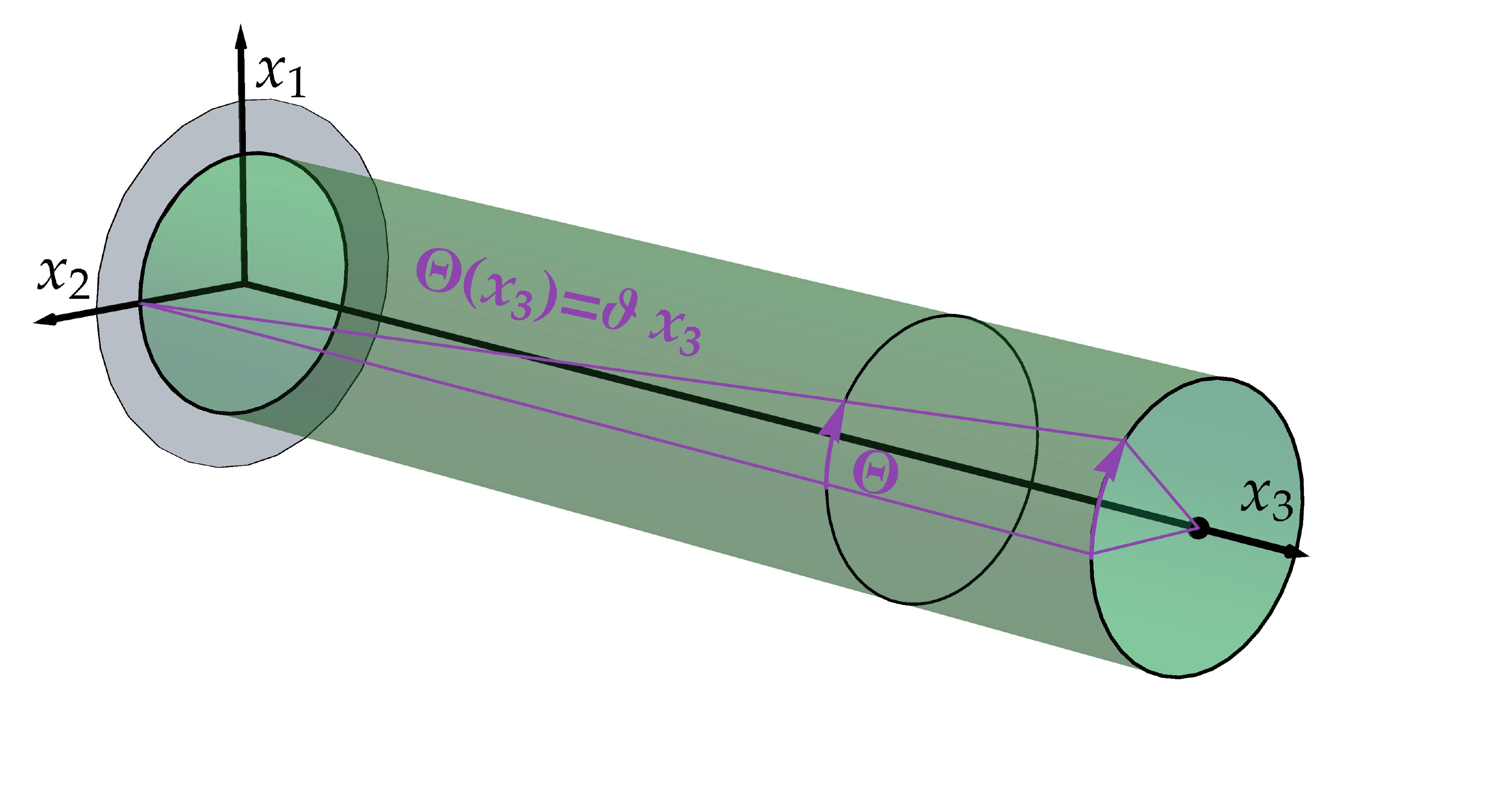}
		\caption{}
	\end{subfigure}
	\hfill
	\begin{subfigure}{0.4\textwidth}
		\centering
		\includegraphics[height=5cm]{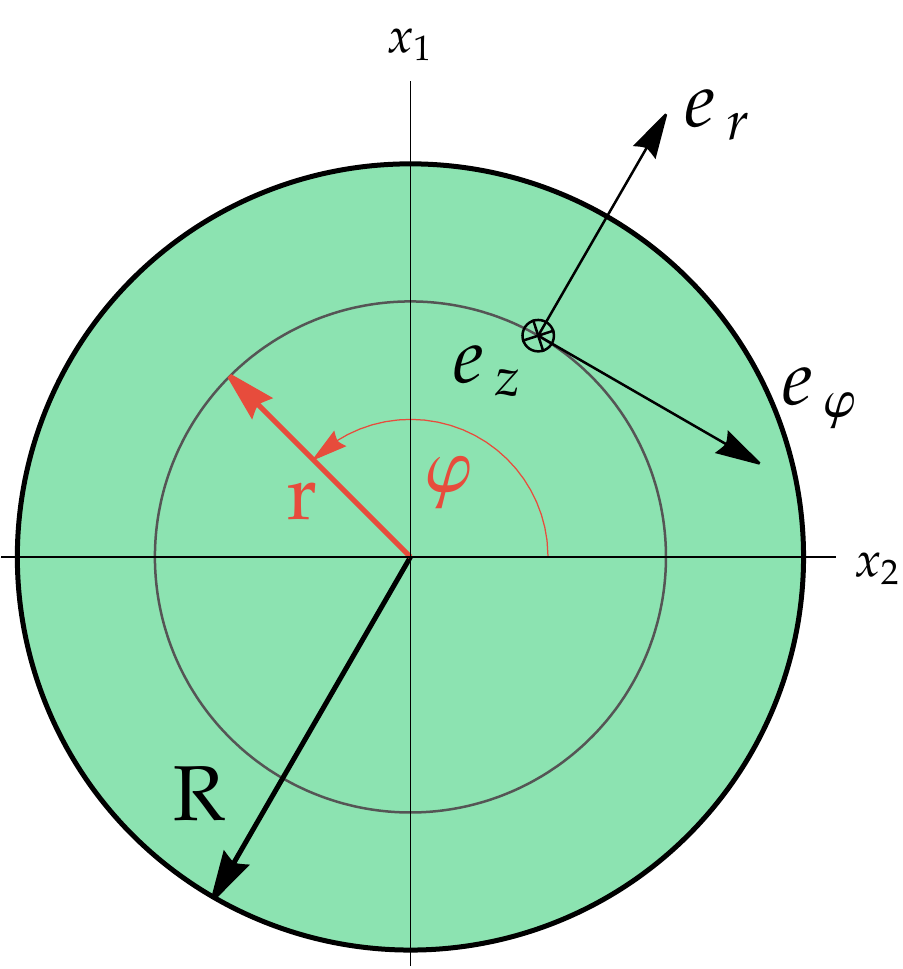}
		\caption{}
	\end{subfigure}
	\caption{
	Geometry of the torsion problem: according to the St.Venant principle, we do not consider how the resultant end torque is applied. Furthermore, we assume that each cross-section (orthogonal to $x_3$) rotates as a rigid body with a constant rate of twist $\frac{\partial \boldsymbol{\Theta}}{\partial x_3}=\boldsymbol{\vartheta}$.
	Since there is no warping, every cross-section remains in the same plane before and after the deformation.
	Note that the final solution for linear elasticity satisfies this ansatz only to within first order in the rate of  twist $\boldsymbol{\vartheta}$, see Fig. \ref{fig:Cauchy_deform}.
	}
	\label{fig:intro_1}
\end{figure}

A third independent identification can be achieved with dynamic measurements, determining the shear wave speed ($c_{s}$) and the pressure wave speed ($c_{p}$) 
\begin{equation}
    c_{s} = \sqrt{\frac{\mu_{\tiny \mbox{macro}}}{\rho}} \, ,
    \qquad\qquad
    c_{p} = \sqrt{\frac{2 \mu_{\tiny \mbox{macro}} + \lambda_{\tiny \mbox{macro}}}{\rho}} \, .
    \label{eq:wave_speed_intro}
\end{equation}

In reality, all these three methods may lead to slightly different values when used to fit real experiments due to the experimental   set up.
Nevertheless, they all are useful in complementing the identification procedure.
We note that all mentioned tests convey a precise physical meaning to the appearing material parameters and this greatly helps in the mechanical application of linear elasticity to real world structures.

The situation is much more involved when trying to determine material parameters for generalized continua.
Even when restricting the attention to linear and isotropic response, the number of additional parameters increase significantly and it is also not clear a priori what the physical meaning of the additional parameters really is.
Lakes \cite{lakes1995experimental,lakesWebsite} has prominently investigated the fitting procedure for isotropic Cosserat solids.
In the linear isotropic Cosserat model (Section  \ref{sec:Cos}) there appear already six independent parameters and a series of experiments with differently sized materials allows to determine the Cosserat constants.
A decisive tool for that purpose is the analytical solution for torsion and bending, which is already available in the literature \cite{lakes2015bending,lurie2018bending,rizzi2021bending}.
The Cosserat model allows to describe size-effects in the sense that more slender specimen have an increased apparent stiffness' in bending and torsion.
However, it is observed that the Cosserat model does have an unphysical stiffness singularity in bending \cite{rizzi2021bending} for a zero slenderness limit, the same appears in general in torsion (Section  \ref{sec:Cos}) but can be avoided upon setting to zero some curvature parameters (Sec \ref{sec:Cos_confor}).
The mentioned stiffness singularity is not only an academic issue, but it concerns  the stable identification of the material parameters \cite{neff2010linear}.
Yet, in the Cosserat theory, the Young's modulus $E_{\tiny \mbox{macro}}$ and the Poisson's ratio $\nu_{\tiny \mbox{macro}}$ can still be determined in a size-independent manner with a homogeneous tension-compression test.
In question are the so-called Cosserat couple modulus $\mu_c \geq 0$ and the three curvature parameters.

A first extension of the Cosserat model is the so-called micro-stretch model, which allows for infinitesimal rotation and volume stretch as independent kinematic fields.
For the micro-stretch model we show that the additional kinematic degree of freedom is not activated in the torsion problem.

Another extension of the Cosserat model is given by the recently introduced relaxed micromorphic model \cite{neff2015relaxed,neff2006existence,sky2020hybrid} (Section  \ref{sec:Relaxed_micro}).
In its static isotropic version it features only 8 independent material parameters comparing favourably to the large number of constitutive parameters in the classical micromorphic model.
While the kinematics of the relaxed micromorphic model coincides with the classical micromorphic model (9 additional degrees of freedom: stretch, shrink, shear, rotations) the curvature term is a direct extension of the Cosserat curvature written in terms of ${\rm Curl} \, \boldsymbol{P}$. An important advantage of the model compared to the Cosserat model is that there is no stiffness singularity in whatsoever situation and four of the eight constants ($\mu_{\tiny \mbox{macro}}$, $\lambda_{\tiny \mbox{macro}}$, $\mu_e$, and $\lambda_e$) can be determined ab initio from size-independent homogeneous tests \cite{neff2019identification}.
There remains to fit three curvature parameters and the Cosserat couple modulus $\mu_c \geq 0$ (which in some situations may be chosen to be zero since the model remains well-posed) \cite{lewintan2020korn,lewintan2019ne,neff2015poincare,bauer2016dev,lewintan2020p}. 

Another advantage of the isotropic relaxed micromorphic model is given by the fact that it can replace the isotropic Cosserat model in a straightforward manner without additional costs. Indeed, the Cosserat curvature parameters can be taken as such as well as the Cosserat couple modulus $\mu_{c}$. The only new parameter set is $\mu_{\mbox{\tiny micro}}, \lambda_{\mbox{\tiny micro}}$, an estimate of which can be inferred from the small-scale response \cite{neff2019identification}. 
Regularity and continuous dependence results for the relaxed micromorphic model have been obtained in \cite{neff2015relaxed,ghiba2014relaxed,Sebastian1,Sebastian2} and first FEM-implementations in $H(\mbox{Curl})$-space are presented in \cite{neff2019identification}.

Next, the micro-strain model \cite{forest2006nonlinear} is in a sense complementary to the Cosserat model:  it assumes an additional strain like symmetric field $S$ as extra degree of freedom.
Here, we recover a simplified micro-strain model without mixed terms and a choice for the curvature parameters, see also \cite{shaat2018reduced} who considers a degenerate micro-strain model in disguise. We recover the  analytical solution given by 
Hütter \cite{hutter2016application}   for the micro-strain torsion problem.
It turns out that for bending \cite{rizzi2021bending}, simple shear \cite{rizzi2021shear} and torsion (Section  \ref{sec:Micro_strain}) the micro-strain solution does not show a stiffness singularity either.
 However, this is not a general feature of the micro-strain model, but only related to the restricted kinematic possibilities: bending and torsion activate prominently rotations, but these are ``filtered out'' in the micro-strain model.
Therefore, bounded stiffness in bending and torsion should come as no surprise.
Next, we combine the Cosserat and the micro-strain ansatz in a novel ad-hoc model whose response is nevertheless governed by the Cosserat kinematics.

Lastly, we have the full micromorphic model \cite{mindlin1964micro,Eringen99}.
The kinematics is augmented with a non-symmetric micro-distortion tensor $\boldsymbol{P}$ (as for the relaxed micromorphic model, 
too) but the curvature energy depends on the full gradient $\mbox{D}\boldsymbol{P}$ of the micro-distortion.
For simplicity and for comparison, we consider a subclass without mixed terms and simplified curvature expression.
In general, the bending and torsion responses show a stiffness singularity, which can be avoided in torsion by a very special ansatz for the curvature energy.
However, nonphysical stiffness singularities cannot, in general, be avoided. Our investigation is complemented by considering the strain-gradient models and its couple-stress subclass. The reason for the singular stiffening behaviour in the other generalized continuum models (except the relaxed micromorphic one) can be connected to their \textit{non-redundant} formulation of the curvature measure \cite{romano2016micromorphic}.

An alternative method to study the deformation of (finite) elastic cylinders is the semi-inverse method introduced by Ie\c san in
	\cite{iecsan1986saint,iesan2006saint},  see also \cite{ghiba2014saint}. This method was also successfully used to study the deformation of elastic cylinders with microstructures, see \cite{ghiba2008semi,bulgariu2013thermal,ghiba2009deformation} and the book \cite{iesan2008classical},
	in which many of Ie\c san's results were unified.
	Regarding the semi-inverse method, all the results obtained in the classical micromorphic theory and all its subclasses (Cosserat, micro-stretch, micro-voids) are obtained by assuming that the internal energy is positive definite in terms of $\mbox{D}\boldsymbol{P}$.
	To the contrary, in the framework of the relaxed micromorphic model, the present results are valid also for internal energies which are not positive definite in terms of $\mbox{D}\boldsymbol{P}$, but rather in terms of $\mbox{Curl} \, \boldsymbol{P}$. We recall that an internal energy which is positive definite in terms of $\mbox{Curl} \, \boldsymbol{P}$ is only semi-positive definite in terms of $\mbox{D}\boldsymbol{P}$.

The paper is now structured as follows.
After fixing our notation in  Section  \ref{sec:notation} we shortly dwell on the formulation of the torsion problem in adapted variables, making it clear that we \textit{do not} revert to express stresses and moments in cylindrical coordinate but we always use a Cartesian expression written in suitable variables.
To set the stage we recall the linear isotropic torsion problem, which will then be suitably generalized.

\allowdisplaybreaks
\subsection{Notation}
\label{sec:notation}
For vectors $a,b\in\mathbb{R}^n$, we define the scalar product $\langle \boldsymbol{a},\boldsymbol{b} \rangle \coloneqq \sum_{i=1}^n a_i\,b_i \in \mathbb{R}$, the euclidean norm  $\norm{\boldsymbol{a}}^2\coloneqq\langle \boldsymbol{a},\boldsymbol{a} \rangle$ and  the dyadic product  $\boldsymbol{a}\otimes \boldsymbol{b} \coloneqq \left(a_i\,b_j\right)_{i,j=1,\ldots,n}\in \mathbb{R}^{n\times n}$. In the same way, for tensors  $\boldsymbol{P},\boldsymbol{Q}\in\mathbb{R}^{n\times n}$, we define the scalar product $\langle \boldsymbol{P},\boldsymbol{Q} \rangle \coloneqq\sum_{i,j=1}^n P_{ij}\,Q_{ij} \in \mathbb{R}$ and the Frobenius-norm $\norm{\boldsymbol{P}}^2\coloneqq\langle \boldsymbol{P},\boldsymbol{P} \rangle$.
Moreover, $\boldsymbol{P}^T\coloneqq (P_{ji})_{i,j=1,\ldots,n}$ denotes the transposition of the matrix $\boldsymbol{P}=(P_{ij})_{i,j=1,\ldots,n}$, which decomposes orthogonally into the skew-symmetric part $\mbox{skew} \, \boldsymbol{P} \coloneqq \frac{1}{2} (\boldsymbol{P}-\boldsymbol{P}^T )$ and the symmetric part $\mbox{sym} \, \boldsymbol{P} \coloneqq \frac{1}{2} (\boldsymbol{P}+\boldsymbol{P}^T)$.
The identity matrix is denoted by $\boldsymbol{\mathbbm{1}}$, so that the trace of a matrix $\boldsymbol{P}$ is given by \ $\tr \boldsymbol{P} \coloneqq \langle \boldsymbol{P},\boldsymbol{\mathbbm{1}} \rangle$, while the deviatoric component of a matrix is given by $\mbox{dev} \, \boldsymbol{P} \coloneqq \boldsymbol{P} - \frac{\tr \left( \boldsymbol{P}\right)}{3} \, \boldsymbol{\mathbbm{1}}$.
Given this, the  orthogonal decomposition possible for a matrix is $\boldsymbol{P} = \mbox{dev} \,\mbox{sym} \, \boldsymbol{P} + \mbox{skew} \, \boldsymbol{P} + \frac{\tr \left( \boldsymbol{P}\right)}{3} \, \boldsymbol{\mathbbm{1}}$.
The Lie-Algebra of skew-symmetric matrices is denoted by $\mathfrak{so}(3)\coloneqq \{\boldsymbol{A}\in\mathbb{R}^{3\times 3}\mid \boldsymbol{A}^T = -\boldsymbol{A}\}$,
while the vector space of symmetric matrices  $\mbox{Sym}(3)\coloneqq \{\boldsymbol{S}\in\mathbb{R}^{3\times 3}\mid \boldsymbol{S}^T = \boldsymbol{S}\}$.
Using the one-to-one map $\mbox{axl}:\mathfrak{so}(3)\to\mathbb{R}^3$ we have
\begin{align}
\boldsymbol{A} \, \boldsymbol{b} =\mbox{axl}(\boldsymbol{A})\times \boldsymbol{b} \quad \forall\, \boldsymbol{A}\in\mathfrak{so}(3) \, ,
\quad
\boldsymbol{b}\in\mathbb{R}^3.
\label{eq:Aanti_axl}
\end{align}
where $\times$ denotes the cross product in $\mathbb{R}^3$.
The inverse of  axl is denoted by Anti: $\mathbb{R}^3\to \mathfrak{so}(3)$.
The Jacobian matrix D$\boldsymbol{u}$ and the curl for a vector field $\boldsymbol{u}$ are defined as
\begin{equation}
\boldsymbol{\mbox{D}u}
=\!
\left(
\begin{array}{ccc}
u_{1,1} & u_{1,2} & u_{1,3} \\
u_{2,1} & u_{2,2} & u_{2,3} \\
u_{3,1} & u_{3,2} & u_{3,3}
\end{array}
\right)\, ,
\qquad
\mbox{curl} \, \boldsymbol{u} = \boldsymbol{\nabla} \times \boldsymbol{u}
=
\left(
\begin{array}{ccc}
u_{3,2} - u_{2,3}  \\
u_{1,3} - u_{3,1}  \\
u_{2,1} - u_{1,2}
\end{array}
\right)
\, .
\end{equation}
We also introduce the $\mbox{Curl}$ and the $\mbox{Div}$ operators of the $3\times 3$ matrix field $\boldsymbol{P}$ as
\begin{equation}
\mbox{Curl} \, \boldsymbol{P}
=\!
\left(
\begin{array}{c}
(\mbox{curl}\left( P_{11} , \right. P_{12} , \left. P_{13} \right)^{{T}})^T \\
(\mbox{curl}\left( P_{21} , \right. P_{22} , \left. P_{23} \right)^{T})^T \\
(\mbox{curl}\left( P_{31} , \right. P_{32} , \left. P_{33} \right)^{T})^T
\end{array}
\right) \!,
\qquad
\mbox{Div}  \, \boldsymbol{P}
=\!
\left(
\begin{array}{c}
\mbox{div}\left( P_{11} , \right. P_{12} , \left. P_{13} \right)^{T} \\
\mbox{div}\left( P_{21} , \right. P_{22} , \left. P_{23} \right)^{T} \\
\mbox{div}\left( P_{31} , \right. P_{32} , \left. P_{33} \right)^{T}
\end{array}
\right) \, .
\end{equation}
The cross product between a second order tensor and a vector is also needed and is defined row-wise as follow
\begin{equation}
\boldsymbol{m} \times \boldsymbol{b} =
\left(
\begin{array}{ccc}
(b \times (m_{11},m_{12},m_{13})^{T})^T \\
(b \times (m_{21},m_{22},m_{23})^{T})^T \\
(b \times (m_{31},m_{32},m_{33})^{T})^T \\
\end{array}
\right) =
\boldsymbol{m} \cdot \boldsymbol{\epsilon} \cdot \boldsymbol{b} =
m_{ik} \, \epsilon_{kjh} \, b_{h} 
\, ,
\end{equation}
where $\boldsymbol{m} \in \mathbb{R}^{3\times 3}$, $\boldsymbol{b} \in \mathbb{R}^{3}$, and $\boldsymbol{\epsilon}$ is the Levi-Civita tensor.
The two indices contraction $\boldsymbol{:}$ is intended as
\begin{equation}
    \boldsymbol{B} \boldsymbol{:} \boldsymbol{\nabla} \, \boldsymbol{\mathfrak{m}}
    =
    B_{ip} \, \mathfrak{m}_{ijk,p}
    =
    N_{jk}
    \,
    ,
    \qquad\qquad\qquad
    \boldsymbol{B} \boldsymbol{:} \boldsymbol{\mathfrak{m}}
    =
    B_{ij} \, \mathfrak{m}_{ijk}
    =
    b_k
    \, ,
\end{equation}
where $\boldsymbol{B}$ and $\boldsymbol{N}$ are second order tensors, $\boldsymbol{\mathfrak{m}}$ is a third order tensor, and $\boldsymbol{b}$ is a vector.

\subsection{Cartesian variables expressed through cylindrical variables}

To address the torsional problem in its natural environment but with the comfort of the classical Cartesian coordinate system, we introduce the cylindrical set of coordinates which allows us to express the classic Cartesian orthogonal set of coordinates $\boldsymbol{x} = \{x_1,x_2,x_3\}$ through a more suitable set of variables $\boldsymbol{r} = \{r,\varphi,z\}$, without switching completely to a cylindrical coordinate system, i.e., without expressing all the quantities (strains, stresses etc.)  in the basis corresponding to the cylindrical coordinates 
	\begin{equation}
	\boldsymbol{e}_{r} = 
	\left(
	\begin{array}{c}
	\cos \varphi \\
	\sin \varphi \\
	0 
	\end{array}
	\right) 
	\, ,
	\qquad
	\boldsymbol{e}_{\varphi} = 
	\left(
	\begin{array}{c}
	-\sin \varphi \\
	\cos \varphi \\
	0 
	\end{array}
	\right) 
	\, ,
	\qquad
	\boldsymbol{e}_{z} = 
	\left(
	\begin{array}{c}
	0 \\
	0 \\
	1 
	\end{array}
	\right) 
	\, .
	\label{eq:polar_unit_vector}
	\end{equation}
The relations for the coordinates are
\begin{equation}
	x_1 = r \, \cos \varphi \, ,
	\qquad
	x_2 = r \, \sin \varphi \, ,
	\qquad
	x_3 = z \, ,
	\qquad
	\label{eq:cyli_coordinate}
\end{equation}
while the relations between the first and the second derivative of a generic vector field $\boldsymbol{f}$ are
\begin{equation}
	\frac{\partial \, f_i(r,\varphi,z)}{\partial \, r_j}
	=
	\frac{\partial \, f_i(r,\varphi,z)}{\partial \, x_k} \, \frac{\partial \, x_k}{\partial \, r_j}
	\, ,
	\qquad
	\frac{\partial^2 \, f_i(r,\varphi,z)}{\partial \, r_j \partial \, r_k}
	=
	\frac{\partial^2 \, f_i(r,\varphi,z)}{\partial \, x_m \partial \, x_n}
	\frac{\partial \, x_m}{\partial \, r_j}
	\frac{\partial \, x_n}{\partial \, r_k}
	+
	\frac{\partial \, f_i(r,\varphi,z)}{\partial \, x_m}
	\frac{\partial^2 \, x_m}{\partial \, r_j \partial \, r_k}
	\, .
	\label{eq:first_second_derivative_cyli}
\end{equation}
The quantities we want to obtain are $\frac{\partial \, f_i(r,\varphi,z)}{\partial \, x_k}$ and $\frac{\partial^2 \, f_i(r,\varphi,z)}{\partial \, x_m \partial \, x_n}$, which are obtainable thanks to  (\ref{eq:first_second_derivative_cyli}) (see Appendix \ref{app:cyli_coordi} for full calculations).

It is emphasized again that \textbf{we will not represent the torsional problem in cylindrical coordinates} (namely all the differential operators, the equilibrium equation, and the kinematic fields), but we will use the classical Cartesian coordinates $\{x_1,x_2,x_3\}$ parameterized in cylindrical variables $\{r,\varphi,z\}$.

\subsection{Structure of the higher-order ansatz}
The ansatz for the displacement field for the cylindrical torsion problem, regardless of the treated model, is always given by
\begin{equation}
	\boldsymbol{u}(r,\varphi,z) =
	\boldsymbol{\vartheta}
	\left(
	\begin{array}{c}
		-x_2(r,\varphi) \, x_3(z) \\
		x_1(r,\varphi) \, x_3(z) \\
		0 
	\end{array}
	\right)
	\, .
	\label{eq:ansatz_disp}
\end{equation}
Here, $\boldsymbol{\vartheta}$ is the rate of twist per unit length.
It is highlighted that the displacement field has the third component equal to zero since we are studying a cylindrical rod, whose cross-section is not subjected to warping.
The most general ansatz for the micro-distortion tensor, which will be used for the full and the relaxed micromorphic model, is 
\begin{equation}
	\boldsymbol{P}(r,\varphi,z) =
	\boldsymbol{\vartheta}
	\left(
	\begin{array}{ccc}
		0   & -x_3(z) & - g_{2}(r) \, x_2(r,\varphi) \\
		x_3(z) &    0 &   g_{2}(r) \, x_1(r,\varphi) \\
		g_{1}(r) \, x_2(r,\varphi) & - g_{1}(r) \, x_1(r,\varphi) & 0 \\
	\end{array}
	\right)
	\, ,
	\label{eq:ansatz_P}
\end{equation}
where $g_1, g_2:[0,\infty)\to \mathbb{R}$.
Starting from  the form (\ref{eq:ansatz_P}) of the ansatz for $P$, it is possible to obtain the ansatz for the micro-stretch model \, (${\boldsymbol{A} = \mbox{skew} \, P}$ and ${\omega \boldsymbol{\mathbbm{1}}= \mbox{tr}(P) \boldsymbol{\mathbbm{1}}}$), the Cosserat model ($\boldsymbol{A} = \mbox{skew} \, P$), the micro-void model ($\omega \boldsymbol{\mathbbm{1}}= \mbox{tr}(P)\boldsymbol{\mathbbm{1}}$), and the micro-strain model ($\boldsymbol{S} = \mbox{sym} \, P$), by taking the skew-symmetric part, the trace of $\boldsymbol{P}$, or the symmetric part depending on what is needed.
Here are reported the symmetric part, the skew-symmetric part and the trace of $\boldsymbol{P}$
\begin{align}
	\boldsymbol{S}(r,\varphi,z) &=
	\mbox{sym} \, \boldsymbol{P}(r,\varphi,z) =
	\frac{\boldsymbol{\vartheta}}{2}
	\left(
	\begin{array}{ccc}
		0   & 0 & g_{m}(r) \, x_2(r,\varphi) \\
		0 &    0 &  - g_{m}(r) \, x_1(r,\varphi) \\
		g_{m}(r) \, x_2(r,\varphi) & - g_{m}(r) \, x_1(r,\varphi) & 0 \\
	\end{array}
	\right)
	\, ,
	\\
	\boldsymbol{A}(r,\varphi,z) &=
	\mbox{skew} \, \boldsymbol{P}(r,\varphi,z) =
	\frac{\boldsymbol{\vartheta}}{2}
	\left(
	\begin{array}{ccc}
		0   & -2x_3(z) & - g_{p}(r) \, x_2(r,\varphi) \\
		2x_3(z) &    0 &   g_{p}(r) \, x_1(r,\varphi) \\
		g_{p}(r) \, x_2(r,\varphi) & - g_{p}(r) \, x_1(r,\varphi) & 0 \\
	\end{array}
	\right)
	\, ,
	\label{eq:ansatz_S_A_omega}
\end{align}
were $g_p(r)=g_{1}(r) + g_{2}(r)$, $g_m(r)=g_{1}(r) - g_{2}(r)$, and $\omega$ is not reported since the ansatz  (\ref{eq:ansatz_P}) has a zero trace.

It is highlighted that each section remains ``rigid'' is not really correct, since the deformation of a cylinder section due to the displacement field  (\ref{eq:ansatz_disp}) (which is a linear approximation of a rigid rotation) looks like

\begin{equation}
    \boldsymbol{x} + \boldsymbol{u} (\boldsymbol{x})
    =
    \left(
    \begin{array}{ccc}
        x_1\\
        x_2\\
        0
    \end{array}
    \right)
    +
    \left(
    \begin{array}{ccc}
        0 & -\boldsymbol{\vartheta} \, x_3 & 0\\
         \boldsymbol{\vartheta} \, x_3 & 0 & 0\\
        0 & 0 & 0
    \end{array}
    \right)
    \left(
    \begin{array}{ccc}
        x_1\\
        x_2\\
        0
    \end{array}
    \right)
    =
    \left(
    \begin{array}{ccc}
        x_1 - x_2 \, x_3 \, \boldsymbol{\vartheta}\\
        x_2 + x_1 \, x_3 \, \boldsymbol{\vartheta}\\
        0
    \end{array}
    \right) \, .
\end{equation}

\begin{figure}[H]
	\begin{subfigure}{0.6\textwidth}
		\centering
		\includegraphics[height=5cm]{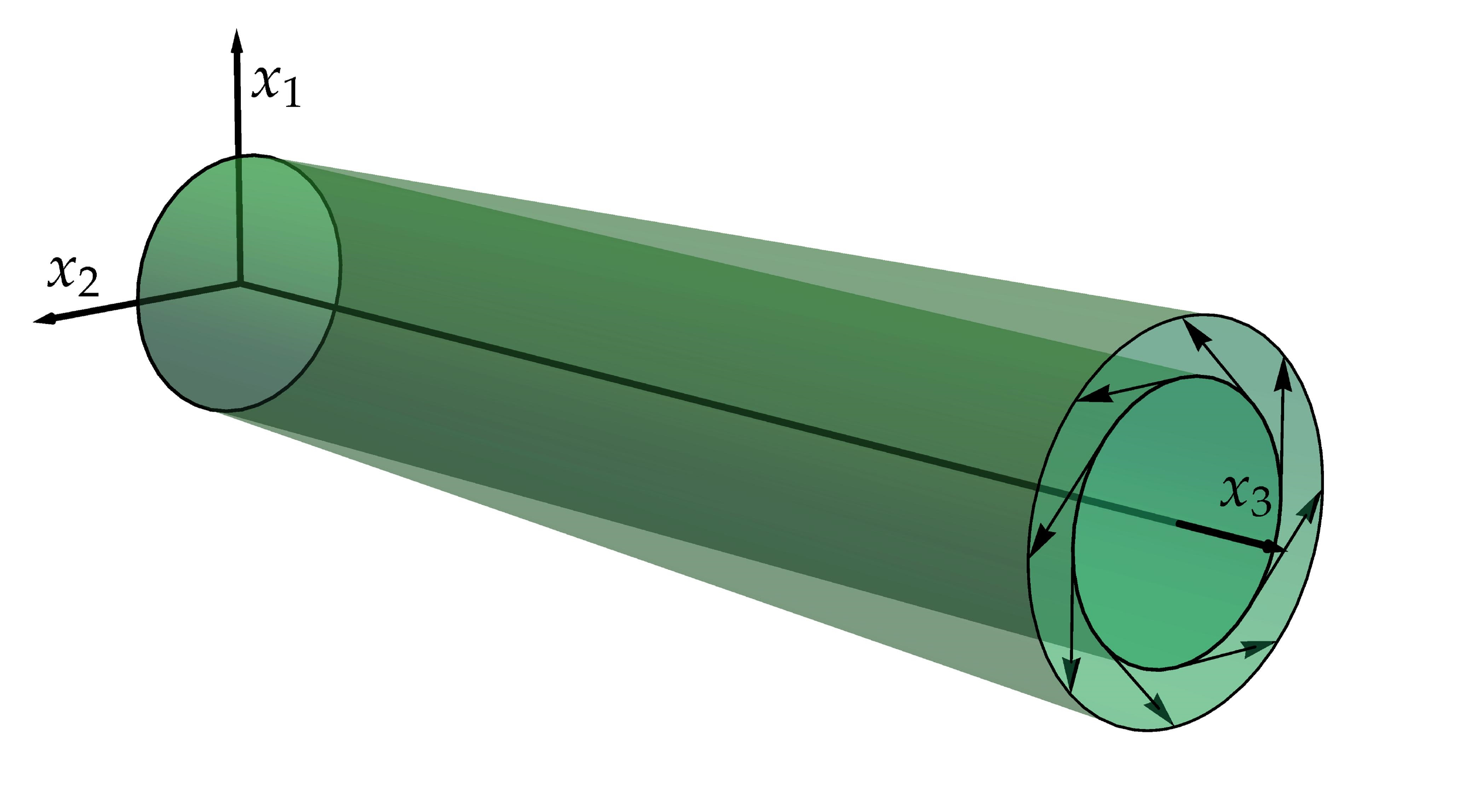}
		\caption{}
	\end{subfigure}
	\hfill
	\begin{subfigure}{0.4\textwidth}
		\centering
		\includegraphics[height=5cm]{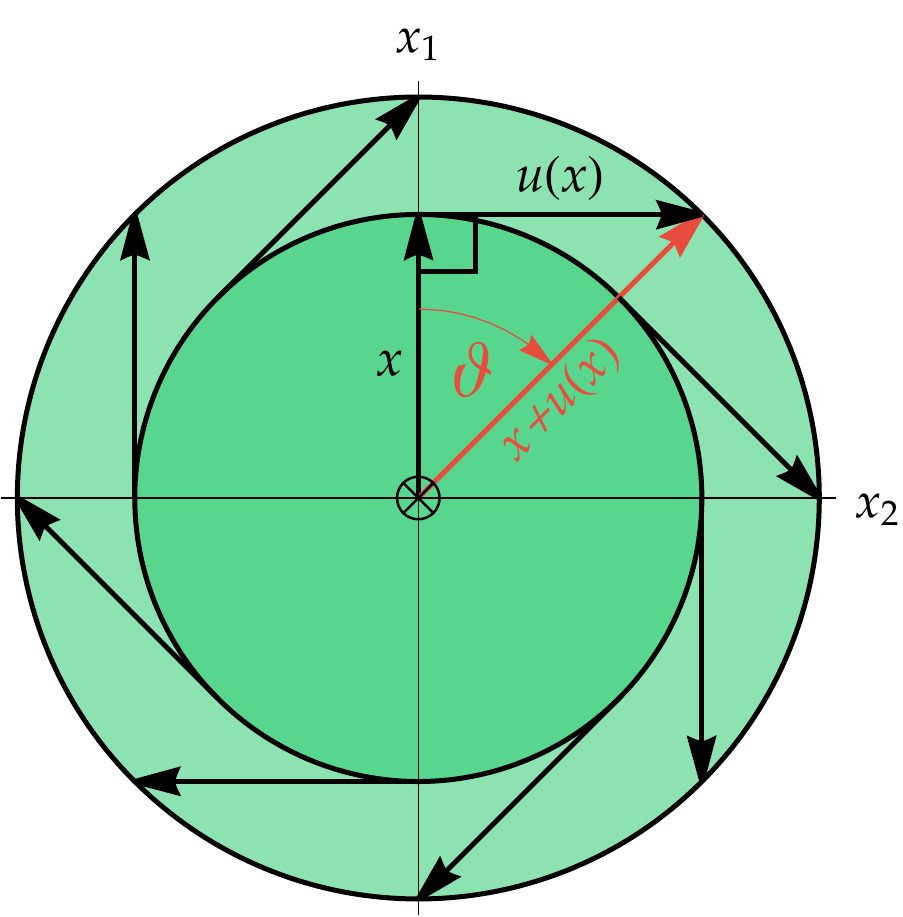}
		\caption{}
	\end{subfigure}
	\caption{
	In the linear approximation, sections of the cylindrical rod are not only rotated, but also expanded radially for non zero rate of twist $\boldsymbol{\vartheta}$. With (b) we see that the change of radius
	$\delta R =
	\frac{R}{\cos{\boldsymbol{\vartheta}}}-R= \frac{1-\cos{\boldsymbol{\vartheta}}}{\cos{\boldsymbol{\vartheta}}}R=
	\left(\frac{\boldsymbol{\vartheta}^2}{2} + \frac{5}{24}\boldsymbol{\vartheta}^4 + \mbox{h.o.t.}\right)R$, being of second order in $\boldsymbol{\vartheta}$. Thus, the linear kinematics is correct to within first order in $\boldsymbol{\vartheta}$, as is well-known.
	}
	\label{fig:Cauchy_deform}
\end{figure}

Of course this radial expansion does not contribute energetically under the small displacement hypothesis, and this can be seen from  (\ref{eq:ansatz_grad_Cau}), in which it is clear that the symmetric strain tensor $\varepsilon={\rm sym} \,\mbox{D}\boldsymbol{u}$ does not depend on $x_3 \equiv z$.

\section{Overview of some generalized continuum models and their interconnections}
\label{sec:overview}

Considering the following notations for the involved quantities:
\begin{align}
    \boldsymbol{u}&: \Omega \subset \mathbb{R}^{3} \to \mathbb{R}^{3} \, , && \mbox{displacement} \, ,
    \notag
    \\*
    \boldsymbol{P}&: \Omega \subset \mathbb{R}^{3} \to \mathbb{R}^{3\times 3} \, , && \mbox{micro-distortion} \, ,
    \notag
    \\*
    \boldsymbol{A}&: \Omega \subset \mathbb{R}^{3} \to \mathfrak{so}(3) \, , && \mbox{micro-rotation} \, ,
    \notag
    \\*
    \boldsymbol{S}&: \Omega \subset \mathbb{R}^{3} \to \mbox{Sym}(3) \, , && \mbox{micro-strain} \, ,
    \notag
    \\*
    \omega&: \Omega \subset \mathbb{R}^{3} \to \mathbb{R} \, , && \mbox{micro-dilatation} \, 
    \notag
\end{align}
and the orthogonal decomposition
\begin{align}
\boldsymbol{P}=\mbox{dev}\, \mbox{sym} \boldsymbol{P} + \mbox{skew} \, \boldsymbol{P} + \frac{1}{3}\, \tr(\boldsymbol{P})\, \boldsymbol{\mathbbm{1}}=\mbox{dev}\, \boldsymbol{S} + \boldsymbol{A} + \omega\,  \boldsymbol{\mathbbm{1}}
\end{align}
we give the following genealogy tree of the generalized continuum models:
\begin{center}
\begin{forest}
    [
    {
    \textit{classical micromorphic}
    \\
    $
    \displaystyle \min_{\boldsymbol{u},\boldsymbol{P}}
    \Big[
    W \left(\boldsymbol{\mbox{D}u},\boldsymbol{P}, \boldsymbol{\mbox{D}P} \right)
    \Big]
    $
    }
        [
        {
        \textit{micro-strain} ($\boldsymbol{P} = \boldsymbol{S}$)
        \\
        $
        \displaystyle \min_{\boldsymbol{u},\boldsymbol{S}}
        \Big[
        W \left(\boldsymbol{\mbox{D}u},\boldsymbol{S}, \boldsymbol{\mbox{D}S} \right)
        \Big]
        $
        }
            [
            {
            \textit{strain gradient} ($\boldsymbol{S} = \mbox{sym}~\boldsymbol{\mbox{D}u}$)
            \\
            $
            \displaystyle \min_{\boldsymbol{u}}
            \Big[
            W \left(\boldsymbol{\mbox{D}u}, \boldsymbol{\mbox{D} \left(\mbox{sym}~\boldsymbol{\mbox{D}u}\right)} \right)
            \Big]
            $
            }
            ]
        ]
        [
        {
        \qquad
        \textit{relaxed micromorphic}
        \qquad
        \\
        \qquad
        $
        \displaystyle \min_{\boldsymbol{u},\boldsymbol{P}}
        \Big[
        W \left(\boldsymbol{\mbox{D}u},\boldsymbol{P}, \boldsymbol{\mbox{Curl}P} \right)
        \Big]
        $
        \qquad
        }
            [
            {
            \qquad\qquad
            \textit{micro-stretch} ($\boldsymbol{P} = \boldsymbol{A} + \omega \boldsymbol{\mathbbm{1}}$)
            \qquad\qquad
            \\
            $
            \displaystyle \min_{\boldsymbol{u},\boldsymbol{A},\omega}
            \Big[
            W \left(\boldsymbol{\mbox{D}u},\boldsymbol{A} , \omega \boldsymbol{\mathbbm{1}}, \boldsymbol{\mbox{Curl} \left( \boldsymbol{A} + \omega \boldsymbol{\mathbbm{1}} \right)} \right)
            \Big]
            $
            }
                [
                {
                \qquad\qquad\qquad
                \textit{Cosserat} ($\boldsymbol{P} = \boldsymbol{A}$)
                \qquad\qquad\qquad
                \\
                $
                \displaystyle \min_{\boldsymbol{u},\boldsymbol{A}}
                \Big[
                W \left(\boldsymbol{\mbox{D}u},\boldsymbol{A}, \boldsymbol{\mbox{Curl} \, \boldsymbol{A} } \right)
                \Big]
                $
                }
                    [
                    {
                    \textit{couple stress} ($\boldsymbol{A} = \mbox{skew}~\boldsymbol{\mbox{D}u}$)
                    \\
                    $
                    \displaystyle \min_{\boldsymbol{u}}
                    \Big[
                    W \left(\boldsymbol{\mbox{D}u}, \boldsymbol{\mbox{Curl} \left(\mbox{skew}~\boldsymbol{\mbox{D}u}\right)} \right)
                    \Big]
                    $
                    }
                        [
                        {
                        \textit{skew symmetric couple stress}
                        \\
                        $
                        \displaystyle \min_{\boldsymbol{u}}
                        \Big[
                        W \left(\boldsymbol{\mbox{D}u}, \mbox{skew}~\boldsymbol{\mbox{Curl} \left(\mbox{skew}~\boldsymbol{\mbox{D}u}\right)} \right)
                        \Big]
                        $
                        }
                        ]
                        [
                        {
                        \textit{modified couple stress}
                        \\
                        $
                        \displaystyle \min_{\boldsymbol{u}}
                        \Big[
                        W \left(\boldsymbol{\mbox{D}u}, \mbox{sym}~\boldsymbol{\mbox{Curl} \left(\mbox{skew}~\boldsymbol{\mbox{D}u}\right)} \right)
                        \Big]
                        $
                        }
                        ]
                    ]
                ]
                [
                {
                \qquad\qquad\qquad
                \textit{micro-void} ($\boldsymbol{P} = \omega \boldsymbol{\mathbbm{1}}$)
                \qquad\qquad\qquad
                \\
                $
                \displaystyle \min_{\boldsymbol{u},\omega}
                \Big[
                W \left(\boldsymbol{\mbox{D}u},\omega, \boldsymbol{\mbox{Curl} \left( \omega \boldsymbol{\mathbbm{1}} \right)} \right)
                \Big]
                $
                }
                [
                {
                \qquad\qquad\qquad
                \textit{ad-hoc model}
                \qquad\qquad\qquad
                \\
                $
                \displaystyle \min_{\boldsymbol{u},\boldsymbol{A},\boldsymbol{S}}
                \Big[
                 W \left(\boldsymbol{\mbox{D}u}, \boldsymbol{A}, \boldsymbol{S}, \boldsymbol{\mbox{Curl} \, \boldsymbol{A} }, \boldsymbol{\mbox{D}S} \right)
                \Big]
                $
                }
                ]
                ]
            ]
        ]
        [
        {
        \textit{second gradient} ($\boldsymbol{P} = \boldsymbol{\mbox{D}u}$)
        \\
        $
        \displaystyle \min_{\boldsymbol{u}}
        \Big[
        W \left(\boldsymbol{\mbox{D}u}, \boldsymbol{\mbox{D}^2 u} \right)
        \Big]
        $
        }
        ]
    ]
\end{forest}
\end{center}

The strain gradient theory and second gradient theory are equivalent \cite{mindlin1964micro,altenbach2019higher}, and contain additionally the couple stress theory as a special case.
Using the ${\rm Curl}$ as primary differential operator for the curvature terms allows a neat unification of concepts.
\section{Torsional problem for the isotropic Cauchy continuum}
\label{sec:Cau}

In order to set up a comparison with the models we will present in the next sections, we start by presenting the solution of the classical cylindrical torsional problem.
The strain energy for a linear elastic isotropic Cauchy continuum is
\begin{equation}
W \left(\boldsymbol{\mbox{D}u}\right) = 
\mu_{\mbox{\tiny macro}} \left\lVert \mbox{sym} \, \boldsymbol{\mbox{D}u} \right\rVert^{2} + 
\frac{\lambda_{\mbox{\tiny macro}}}{2} \mbox{tr}^2\left(\boldsymbol{\mbox{D}u}\right)
\, ,
\label{eq:energy_Cau}
\end{equation}
where $\lambda_{\mbox{\tiny macro}}$ and $\mu_{\mbox{\tiny macro}}$ are the macroscopic Lamé constants.

In terms of the symmetric Cauchy stress tensor 
${\boldsymbol{\sigma}= 
2\,\mu_{\mbox{\tiny macro}}\,\mbox{sym} \, \boldsymbol{\mbox{D}u}
+ \lambda_{\mbox{\tiny macro}}\,\mbox{tr}\left(\boldsymbol{\mbox{D}u}\right)  \boldsymbol{\mathbbm{1}}}$,
where $\boldsymbol{\varepsilon} = \mbox{sym} \, \boldsymbol{\mbox{D}u}$ denotes the classical symmetric strain tensor, the equilibrium equation (in the absence of body forces) and the Neumann lateral boundary conditions (at the free surface) are
\begin{equation}
\mbox{Div}\, 
\boldsymbol{\sigma} 
= \boldsymbol{0}
\, ,
\qquad\qquad\qquad
\boldsymbol{t}(r = R) = \boldsymbol{\sigma}(r = R) \cdot \boldsymbol{e}_{r} = \boldsymbol{0}
\, .
\label{eq:equi_Cau}
\end{equation}
Our aim is to study a \textbf{state of uniform torsion $\boldsymbol{\vartheta}$} for an infinitely extended cylindrical rod. According to the cylindrical reference  system shown in Fig.~\ref{fig:intro_1}, the ansatz for the displacement is 
\begin{equation}
	\boldsymbol{u}(x_1,x_2,x_3) =
	\boldsymbol{u}(r,\varphi,z) =
	\boldsymbol{\vartheta}
	\left(
	\begin{array}{c}
		-x_2(r,\varphi) \, x_3(z) \\
		x_1(r,\varphi) \, x_3(z) \\
		0 
	\end{array}
	\right)
	=
	\boldsymbol{\vartheta} 
	\left(
	\begin{array}{c}
		-z \, r \, \sin \varphi \\
		 z \, r \, \cos \varphi \\
		0 
	\end{array}
	\right) 
	\, ,
	\label{eq:ansatz_disp_Cau}
\end{equation}
where $\boldsymbol{\vartheta}$ is the angle of twist per unit length.
It is underlined that the third component of the displacement is chosen equal to zero since the cross-section is circular and therefore no warping is expected.
The gradient of the displacement and its symmetric part are (the gradient is always taken with respect to the Cartesian coordinate system and then rewritten in the variables $\{r,\varphi,z\}$)
\begin{equation}
	\boldsymbol{\mbox{D}u}=
	\boldsymbol{\vartheta}
	\left(
	\begin{array}{ccc}
		0 &-z &- r \, \sin \varphi \\
		z & 0 &  r \, \cos \varphi \\
		0 & 0 & 0 
	\end{array}
	\right) \, ,
	\qquad
	\boldsymbol{\varepsilon} =
	\mbox{sym} \, \boldsymbol{\mbox{D}u}=
	\frac{\boldsymbol{\vartheta}}{2}
	\left(
	\begin{array}{ccc}
		0 & 0 & -r \, \sin \varphi \\
		0 & 0 & r \, \cos \varphi \\
		-r \, \sin \varphi & r \, \cos \varphi & 0 \\
	\end{array}
	\right) \, .
	\label{eq:ansatz_grad_Cau}
\end{equation}
Substituting the ansatz  (\ref{eq:ansatz_grad_Cau}) in the equilibrium equation  (\ref{eq:equi_Cau}), it is easy to verify that they are identically satisfied.

\begin{figure}[H]
	\centering
	\includegraphics[height=5cm]{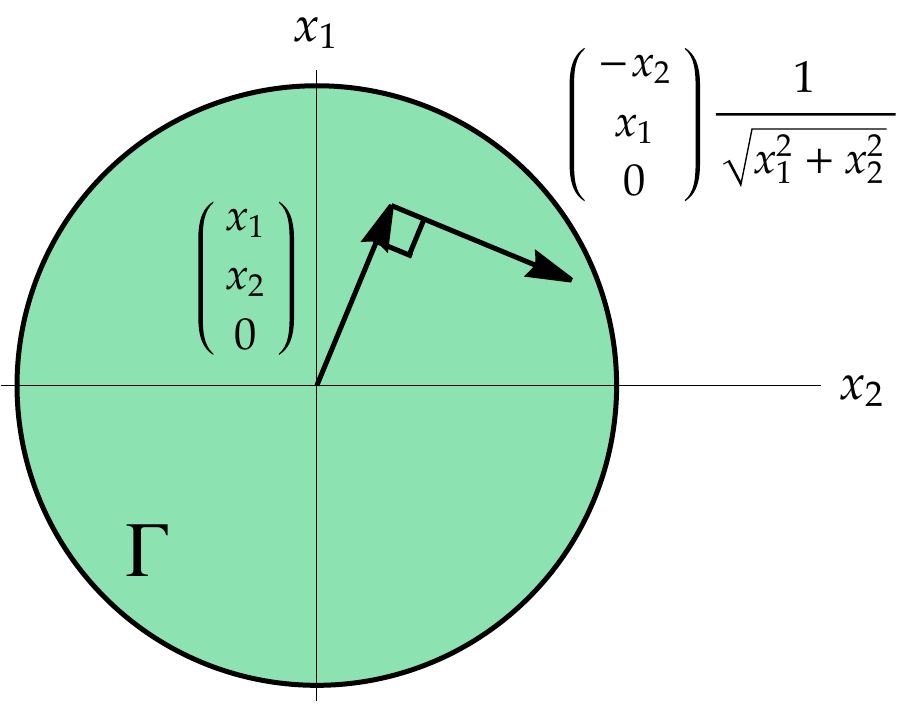}
	\caption{Calculation of the classical torque.}
	\label{fig:Cauchy_traction}
\end{figure}

In order to help the geometric interpretation of the torque, see Fig.~\ref{fig:Cauchy_traction}, we present its expression in Cartesian coordinates along with its representation in the cylindrical variables 
	\begin{align}
	M_{\mbox{c}} (\boldsymbol{\vartheta})
	\coloneqq
	&
	\iint_{\Gamma}
	\Big[
	\overbrace{
		\langle
		\overbrace{\boldsymbol{\sigma} \, \boldsymbol{e}_{3}}^{\mbox{traction}} , \left(
		\begin{array}{c} -x_2 \\ x_1 \\ 0 \\\end{array}\right)
		\frac{1}{\sqrt{x_1^2 + x_2^2}}
		\rangle
	}^{
		\mbox{twisting force per unit area}}
	\overbrace{\sqrt{x_1^2 + x_2^2}}^{\begin{array}{c} \mbox{length of} \\ \mbox{lever arm} \end{array}}
	\Big]
	\, \mbox{d}x_1 \, \mbox{d}x_2
	\label{eq:Cauchy_traction}
	\\*
	=
	&
	\iint_{\Gamma}
	\Big[
	x_1 \, \sigma_{23} - x_2 \, \sigma_{13}
	\Big]
	\, \mbox{d}x_1 \, \mbox{d}x_2
	=
	\int_{0}^{2\pi}
	\int_{0}^{R}
	\Big[
	\langle
	\boldsymbol{\sigma} \, \boldsymbol{e}_{z} , \boldsymbol{e}_{\varphi}
	\rangle
	r
	\Big] r
	\, \mbox{d}r \, \mbox{d}\varphi
	\, ,
	\end{align}
where $e_3 = e_{\widehat{z}} = \left(0,0,1\right)$ is the unit vector aligned with the mid-axis of the cylindrical rod.

The torque (or moment of torsion) about the $x_3$-axis and energy (per unit length d$x_3$) expressions are
\begin{align}
	M_{\mbox{c}} (\boldsymbol{\vartheta}) \coloneqq&
	\int_{0}^{2\pi}
	\int_{0}^{R}
	\Big[
	\langle
	\boldsymbol{\sigma} \, \boldsymbol{e}_{z} , \boldsymbol{e}_{\varphi}
	\rangle
	r
	\Big] r
	\, \mbox{d}r \, \mbox{d}\varphi
	=
	\mu_{\mbox{\tiny macro}} \, \frac{\pi R^4}{2} \, \boldsymbol{\vartheta}
	=
	\mu_{\mbox{\tiny macro}} \, I_{p} \, \boldsymbol{\vartheta}
	=
	T_{\mbox{\tiny macro}} \, \boldsymbol{\vartheta}
	\, ,
	\label{eq:torque_ene_dimensionless_Cau}
	\\
	W_{\mbox{tot}} (\boldsymbol{\vartheta}) \coloneqq&
	\int_{0}^{2\pi}
	\int_{0}^{R}
	W(\boldsymbol{\boldsymbol{\mbox{D}u}}) \, r
	\, \mbox{d}r \, \mbox{d}\varphi
	=
	\frac{1}{2} \, \mu_{\mbox{\tiny macro}} \, \frac{\pi R^4}{2} \, \boldsymbol{\vartheta}^2
	=
	\frac{1}{2} \, \mu_{\mbox{\tiny macro}} \, I_{p} \, \boldsymbol{\vartheta}^2
	=
	\frac{1}{2} \, T_{\mbox{\tiny macro}} \, \boldsymbol{\vartheta}^2
	\, ,
	\notag
\end{align}
where $\mu_{\mbox{\tiny macro}}$ is the macroscopic shear modulus, $I_{p}=\frac{\pi R^4}{2}$ is the polar moment of inertia, and $\mbox{T}_{\mbox{\tiny macro}}=\mu_{\mbox{\tiny macro}} \, I_{p}$  is the torsional stiffness.
It is also highlighted that
\begin{equation}
\frac{\mbox{d}}{\mbox{d}\boldsymbol{\vartheta}}W_{\mbox{tot}}(\boldsymbol{\vartheta}) = M_{\mbox{c}} (\boldsymbol{\vartheta})
= T_{\mbox{\tiny macro}} \, \boldsymbol{\vartheta}
\, ,
\qquad\qquad\qquad
\frac{\mbox{d}^2}{\mbox{d}\boldsymbol{\vartheta}^2}W_{\mbox{tot}}(\boldsymbol{\vartheta})
= T_{\mbox{\tiny macro}}
\, .
\end{equation}

Here and in the remainder of this work, the elastic coefficients $\mu_i,\lambda_i,\kappa_i$ are expressed in [MPa], the coefficients $a_i$ are dimensionless, the characteristic lengths $L_c$ and the radius $R$ in meter [m], the rate of twist $\boldsymbol{\vartheta}$ in [1/m].
\section{Torsional problem for the isotropic relaxed micromorphic model}
\label{sec:Relaxed_micro}
The relaxed micromorphic model, contrary to all the other proposals for generalized continua in the literature, lives on two well-defined and separated scales, each describing linear elastic response:
the classical \textit{macroscopic response} (characteristic length $L_c \to 0$, available for experiments with large specimens) is described as usual by
\begin{align}
    E_{\tiny \mbox{macro}}
    &=
    \frac{9 \kappa_{\tiny \mbox{macro}}  \, \mu_{\tiny \mbox{macro}} }{3 \kappa_{\tiny \mbox{macro}} + \mu_{\tiny \mbox{macro}}} \, ,
    \qquad
    &\nu_{\tiny \mbox{macro}}
    &=
    \frac{3 \kappa_{\tiny \mbox{macro}} -2 \mu_{\tiny \mbox{macro}} }{2 (3 \kappa_{\tiny \mbox{macro}} +\mu_{\tiny \mbox{macro}} )} \, ,
    \\*
    \lambda_{\tiny \mbox{macro}}
    &=
    \frac{3\kappa_{\tiny \mbox{macro}} - 2\mu_{\tiny \mbox{macro}}}{3} \, ,
    \qquad
    &\kappa_{\tiny \mbox{macro}}
    &=
    \frac{2\mu_{\tiny \mbox{macro}} + 3\lambda_{\tiny \mbox{macro}}}{3} \, .
    \label{eq:real_parameter_macro_2}
\end{align}
The macroscopic parameters can be uniquely determined from a \textit{homogeneous} macroscopic tension-compression test.
However, the shear modulus $\mu_{\tiny \mbox{macro}}$ and the Young's modulus $E_{\tiny \mbox{macro}}$ can also be identified from the \textit{inhomogeneous} torsion and bending test, respectively.
Indeed, the classical torsional stiffness of a circular rod is given by
\begin{equation}
    T_{\mbox{\tiny macro}} =
    \mu_{\mbox{\tiny macro}} \, I_{p} =
    \mu_{\mbox{\tiny macro}} \, \frac{\pi R^4}{2} \, .
    \label{eq:torsional_stiffness_macro}
\end{equation}

The \textit{microscopic scale} (appearing for $L_c \to \infty$), representing a \textit{surrogate stiffness} connected to the smallest meaningful scale of the material is described by the parameters
\begin{align}
    E_{\tiny \mbox{micro}}
    &=
    \frac{9 \kappa_{\tiny \mbox{micro}}  \, \mu_{\tiny \mbox{micro}} }{3 \kappa_{\tiny \mbox{micro}} + \mu_{\tiny \mbox{micro}}} \, ,
    \qquad
    &\nu_{\tiny \mbox{micro}}
    &=
    \frac{3 \kappa_{\tiny \mbox{micro}} -2 \mu_{\tiny \mbox{micro}} }{2 (3 \kappa_{\tiny \mbox{micro}} +\mu_{\tiny \mbox{micro}} )} \, ,
    \\*
    \lambda_{\tiny \mbox{micro}}
    &=
    \frac{3\kappa_{\tiny \mbox{micro}} - 2\mu_{\tiny \mbox{micro}}}{3} \, ,
    \qquad
    &\kappa_{\tiny \mbox{micro}}
    &=
    \frac{2\mu_{\tiny \mbox{micro}} + 3\lambda_{\tiny \mbox{micro}}}{3} \, ,
    \label{eq:real_parameter_micro_2}
\end{align}
The macroscopic parameters $ \mu_{\tiny \mbox{macro}}$ and  $ \lambda_{\tiny \mbox{macro}}$ do not directly intervene in the formulation of the relaxed micromorphic model  (\ref{eq:energy_RM}), but the connection is necessarily given by the \textit{Reuss-like homogenization formula} \cite{reuss1929berechnung}
\begin{align}
    \mu_{\mbox{\tiny macro}}
    &=
    \frac{\mu_e \, \mu_{\mbox{\tiny micro}}}{\mu_e + \mu_{\mbox{\tiny micro}}}
    \qquad \Longleftrightarrow \qquad
    \mu_{e}
    =
    \frac{\mu_{\mbox{\tiny macro}} \, \mu_{\mbox{\tiny micro}}}{\mu_{\mbox{\tiny micro}} - \mu_{\mbox{\tiny macro}}} \, ,
    \label{eq:para_meso_scale}
    \\*
    \kappa_{\mbox{\tiny macro}}
    &=
    \frac{\kappa_e \, \kappa_{\mbox{\tiny micro}}}{\kappa_e + \kappa_{\mbox{\tiny micro}}}
    \qquad \Longleftrightarrow \qquad
    \kappa_{e}
    =
    \frac{\kappa_{\mbox{\tiny macro}} \, \kappa_{\mbox{\tiny micro}}}{\kappa_{\mbox{\tiny micro}} - \kappa_{\mbox{\tiny macro}}}
    \, .
    \notag
\end{align}
\begin{figure}[H]
	\centering
    \begin{subfigure}{0.4\textwidth}
    	\centering
	    \includegraphics[width=\textwidth]{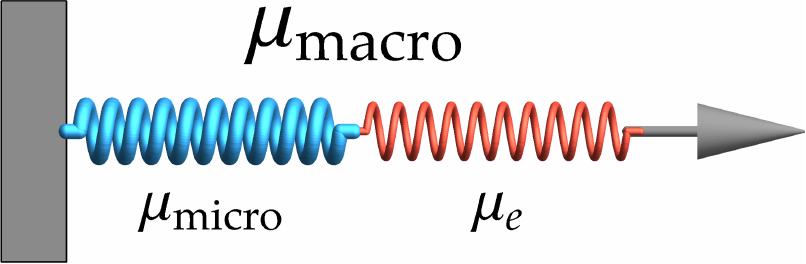}
    \end{subfigure}
    \hspace{1.5cm}
    \begin{subfigure}{0.4\textwidth}
	    \centering
		\includegraphics[width=\textwidth]{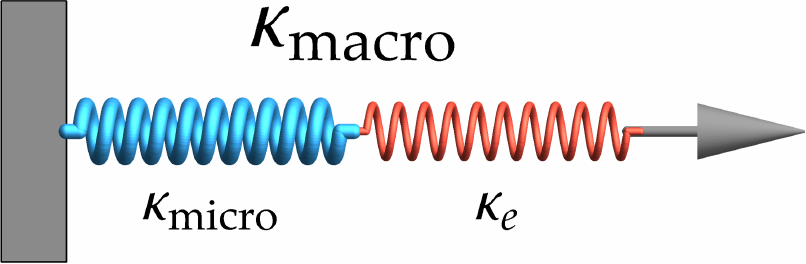}
    \end{subfigure}
\caption{
Macro and micro-scale stiffness governed by two springs in series. If $\mu_{\mbox{\tiny micro}} \to \infty$, this implies that $\mu_{\mbox{\tiny macro}} = \mu_e$.
In all suitable cases for our family of considered generalized continua (depending on the kinematics), we use the same/similar lower order energy expression (the energy without curvature).}
\label{fig:springs_micro_macro}
\end{figure}

Note that the Cosserat couple modulus $\mu_c \geq 0$ is not appearing in the homogenization formulas  (\ref{eq:para_meso_scale}).
As a consequence, both parameter sets (\ref{eq:real_parameter_macro_2})-(\ref{eq:para_meso_scale}) can be identified independently of the scale consideration (being particularly careful with the techniques for the micro-parameters identification) and they uniquely determine the meso-scale parameter set $ \mu_{\tiny \mbox{e}}$, $ \lambda_{\tiny \mbox{e}}$ appearing in  (\ref{eq:para_meso_scale})$_2$.

The  general expression of the strain energy for the isotropic relaxed micromorphic continuum is
\begin{align}
	W \left(\boldsymbol{\mbox{D}u}, \boldsymbol{P},\mbox{Curl}\,\boldsymbol{P}\right)
	=
	&
	\, \mu_{e} \left\lVert \mbox{sym} \left(\boldsymbol{\mbox{D}u} - \boldsymbol{P} \right) \right\rVert^{2}
	+ \frac{\lambda_{e}}{2} \mbox{tr}^2 \left(\boldsymbol{\mbox{D}u} - \boldsymbol{P} \right) 
	+ \mu_{c} \left\lVert \mbox{skew} \left(\boldsymbol{\mbox{D}u} - \boldsymbol{P} \right) \right\rVert^{2}
	\notag
	\\*
	&
	+ \mu_{\tiny \mbox{micro}} \left\lVert \mbox{sym}\,\boldsymbol{P} \right\rVert^{2}
	+ \frac{\lambda_{\tiny \mbox{micro}}}{2} \mbox{tr}^2 \left(\boldsymbol{P} \right)
	\label{eq:energy_RM}
	\\*
	&
	+ \frac{\mu \,L_c^2 }{2} \,
	\left(
	a_1 \, \left\lVert \mbox{dev sym} \, \mbox{Curl} \, \boldsymbol{P}\right\rVert^2 +
	a_2 \, \left\lVert \mbox{skew} \,  \mbox{Curl} \, \boldsymbol{P}\right\rVert^2 +
	\frac{a_3}{3} \, \mbox{tr}^2 \left(\mbox{Curl} \, \boldsymbol{P}\right)
	\right)
	\notag
	\, ,
	\notag
\end{align}
where ($\mu_e$,$\lambda_e$), ($\mu_{\mbox{\tiny micro}}$,$\lambda_{\mbox{\tiny micro}}$), $\mu_c$, $L_c > 0$, and ($a_1$,$a_2$,$a_3$) are the parameters related to the meso-scale, the parameters related to the micro-scale, the Cosserat couple modulus, the characteristic length, and the three general isotropic curvature parameters, respectively.
This energy expression represents the  most general isotropic form possible for the relaxed micromorphic model.
It is important to underline that, given the subsequent ansatz  (\ref{eq:ansatz_RM}), it holds that $\mbox{skew} \,  \mbox{Curl} \, \boldsymbol{P}=0$. This reduces immediately the number of curvature parameters appearing in the torsion solution.
 In the absence of body forces, the equilibrium equations  are then
\begin{align}
	\mbox{Div}\overbrace{\left[2\mu_{e}\,\mbox{sym} \left(\boldsymbol{\mbox{D}u} - \boldsymbol{P} \right) + \lambda_{e} \mbox{tr} \left(\boldsymbol{\mbox{D}u} - \boldsymbol{P} \right) \boldsymbol{\mathbbm{1}}
		+ 2\mu_{c}\,\mbox{skew} \left(\boldsymbol{\mbox{D}u} - \boldsymbol{P} \right)\right]}^{\mathlarger{\widetilde{\sigma}}\coloneqq}
	&= \boldsymbol{0},
	\notag
	\\*
	\widetilde{\sigma}
	- 2 \mu_{\mbox{\tiny micro}}\,\mbox{sym}\,\boldsymbol{P} - \lambda_{\tiny \mbox{micro}} \mbox{tr} \left(\boldsymbol{P}\right) \boldsymbol{\mathbbm{1}}
	- 
	\mu \, L_{c}^{2} \, \mbox{Curl}
	\left(
	a_1 \, \mbox{dev sym} \, \mbox{Curl} \, \boldsymbol{P} +
	\frac{a_3}{3} \, \mbox{tr} \left(\mbox{Curl} \, \boldsymbol{P}\right)\mathbbm{1}
	\right)
	&= \boldsymbol{0}.
	\label{eq:equi_RM}
\end{align}
The boundary conditions at the lateral  free surface are 
\begin{align}
	\boldsymbol{\widetilde{t}}(r = R) &= 
	\boldsymbol{\widetilde{\sigma}}(r) \cdot \boldsymbol{e}_{r} = 
	\boldsymbol{0}_{\mathbb{R}^3} \, ,
	&\mbox{(traction free)} \, ,
	\label{eq:BC_RM_gen}
	\\*
	\boldsymbol{\eta}(r = R) &= 
	\boldsymbol{m} (r) \cdot \boldsymbol{\epsilon} \cdot \boldsymbol{e}_{r} = 
	\boldsymbol{m} (r) \times \boldsymbol{e}_{r} = 
	\boldsymbol{0}_{\mathbb{R}^{3 \times 3}} \, ,
	&\mbox{(moment free)}
	\, ,
	\notag
\end{align}
where
\begin{equation}
\boldsymbol{m} = \mu \, L_c^2 \, \left(
a_1 \, \mbox{dev sym} \, \mbox{Curl} \, \boldsymbol{P} +
\frac{a_3}{3} \, \mbox{tr} \left(\mbox{Curl} \, \boldsymbol{P}\right)\mathbbm{1}
\right)
\end{equation}
is a generalized non-symmetric second order moment tensor,
the (non-symmetric) force-stress tensor $\boldsymbol{\widetilde{\sigma}}$ is given in  (\ref{eq:equi_RM}),
$\boldsymbol{e}_{r}$ is the radial unit vector, and
$\boldsymbol{\epsilon}$ is the Levi-Civita tensor.
 The  vector $\boldsymbol{\widetilde{t}}(r) \in \mathbb{R}^{3}$ is the generalised traction and the tensor $\boldsymbol{\eta}(r) \in \mathbb{R}^{3 \times 3}$ is called the generalized double traction tensor.
According to the cylindrical reference  system shown in Fig.~\ref{fig:intro_1}, the ansatz for the displacement and for the micro-distortion $\boldsymbol{P}$ is
\begin{align}
	\boldsymbol{u}(x_1,x_2,x_3) &=
	\boldsymbol{u}(r,\varphi,z) =
	\boldsymbol{\vartheta}
	\left(
	\begin{array}{c}
		-x_2(r,\varphi) \, x_3(z) \\
		x_1(r,\varphi) \, x_3(z) \\
		0 
	\end{array}
	\right)
	\, ,
	\notag
	\\*
	\boldsymbol{P}(x_1,x_2,x_3) &=
	\boldsymbol{P}(r,\varphi,z) =
	\boldsymbol{\vartheta}
	\left(
	\begin{array}{ccc}
		0   & -x_3(z) & - g_{2}(r) \, x_2(r,\varphi) \\
		x_3(z) &    0 &   g_{2}(r) \, x_1(r,\varphi) \\
		g_{1}(r) \, x_2(r,\varphi) & - g_{1}(r) \, x_1(r,\varphi) & 0 \\
	\end{array}
	\right)
	\, ,
	\label{eq:ansatz_RM}
\end{align}
where $x_1(r,\varphi) = r \, \cos \varphi$, $x_2(r,\varphi) = r \, \sin \varphi$, and $x_3(z) = z$.
The Cartesian $\boldsymbol{\mbox{D}u}$ and the Cartesian $\mbox{Curl} \, \boldsymbol{P}$ expressed in the cylindrical variables $(r,\varphi,z)$ are
\begin{align}
	\boldsymbol{\mbox{D}u}(r,\varphi,z)
	&
	= 
	\boldsymbol{\vartheta}
	\left(
	\begin{array}{ccc}
		0 & -z & -r \, \sin \varphi \\
		z &  0 &  r \, \cos \varphi \\
		0 &  0 & 0 \\
	\end{array}
	\right)
	\, ,
	\label{eq:grad_RM}
	\\*[3mm]
	\mbox{Curl} \, \boldsymbol{P}(r,\varphi,z)
	&
	= 
	\boldsymbol{\vartheta}
	\left(
    \begin{array}{ccc}
        1 - g_{2}(r ) - r \, g_{2}'(r ) \, \sin ^2 \varphi
        &
        r \, g_{2}'(r ) \, \sin \varphi \, \cos \varphi
        &
        0
        \\
        r \, g_{2}'(r ) \, \sin \varphi \, \cos \varphi
        &
        1 - g_{2}(r ) - r \, g_{2}'(r ) \, \cos ^2 \varphi
        &
        0
        \\
        0
        &
        0
        &
        - \left(2 \, g_{1}(r ) + r \, g_{1}'(r ) \right) \\
    \end{array}
    \right) \, .
	\notag
\end{align}
It can be remarked   that $\mbox{Curl} \, \boldsymbol{P}$ is symmetric.

Inserting the ansatz  (\ref{eq:ansatz_RM})-(\ref{eq:grad_RM}) in  (\ref{eq:equi_RM}), the 12 equilibrium equations are reduced to the following 4 ordinary differential equilibrium equations
\begin{align}
	\frac{1}{3} \, \boldsymbol{\vartheta} \, \sin \varphi \left(r  \left(\mu \, L_c^2 \left( (a_{1} - a_{3}) \, g_{1}''(r ) - (2 a_{1} + a_{3}) \, g_{2}''(r )\right) + 3 \mu _c \, (g_{1}(r ) + g_{2}(r )-1)
	\right.\right.
	\qquad
	\qquad
	\qquad
	&
	\notag
	\\*
	\left.\left.
	-3 \left(\mu _e + \mu _{\mbox{\tiny micro}}\right) \, (g_{1}(r ) - g_{2}(r )) - 3 \mu _e\right) + 3 \mu \, L_c^2 \left((a_{1}-a_{3}) \, g_{1}'(r )-(2 a_{1} + a_{3}) \, g_{2}'(r )\right)\right)
	&= 0
	\, ,
	\notag
	\\*
	\pushleft{\frac{1}{3} \, \boldsymbol{\vartheta} \, \cos \varphi \left(r \left(\mu \, L_c^2 \left((a_{3} - a_{1}) \, g_{1}''(r ) + (2 a_{1}+a_{3}) \, g_{2}''(r )\right) - 3 \, \mu _c (g_{1}(r ) + g_{2}(r ) - 1)
	\right.\right.}
	&
	\notag
	\\*
	\left.\left.
	+3 \left(\mu _e + \mu _{\mbox{\tiny micro}}\right) (g_{1}(r ) - g_{2}(r )) + 3 \mu _e\right) + 3 \, \mu \, L_c^2 \left((a_{3} - a_{1}) \, g_{1}'(r ) + (2 a_{1} + a_{3}) \, g_{2}'(r )\right)\right)
	&= 0
	\, ,
	\label{eq:equi_equa_RM}
	\\*
	\pushleft{\frac{1}{3} \, \boldsymbol{\vartheta} \, \sin \varphi \left(r  \left(\mu \, L_c^2 \left((2 a_{1} + a_{3}) \, g_{1}''(r ) + (a_{3} - a_{1}) \, g_{2}''(r )\right) - 3 \left(\mu _c \, (g_{1}(r ) + g_{2}(r ) - 1)
	\right.\right.\right.}
	&
	\notag
	\\*
	\left.\left.\left.
	+\left(\mu _e + \mu _{\mbox{\tiny micro}}\right) (g_{1}(r ) - g_{2}(r ))\right) - 3 \mu _e\right) + 3 \, \mu \, L_c^2 \left((2 a_{1} + a_{3}) \, g_{1}'(r ) + (a_{3} - a_{1}) \, g_{2}'(r )\right)\right)
	&= 0
	\, ,
	\notag
	\\*
	\pushleft{\frac{1}{3} \, \boldsymbol{\vartheta} \, \cos \varphi \left(r  \left(\mu \, L_c^2 \left((a_{1} - a_{3}) \, g_{2}''(r )-(2 a_{1} + a_{3}) \, g_{1}''(r )\right) + 3 \left(\mu _c \, (g_{1}(r ) + g_{2}(r ) - 1)
	\right.\right.\right.}
	&
	\notag
	\\*
	\left.\left.\left.
	+\left(\mu _e + \mu _{\mbox{\tiny micro}}\right) (g_{1}(r ) - g_{2}(r ))\right) + 3 \mu _e\right) + 3 \, \mu \, L_c^2 \left((a_{1} - a_{3}) \, g_{2}'(r )-(2 a_{1} + a_{3}) \, g_{1}'(r )\right)\right)
	&= 0
	\, .
	\notag
\end{align}
It is important to underline that  (\ref{eq:equi_RM})$_1$ is identically satisfied, and that from the entire set  of four equilibrium equations (\ref{eq:equi_equa_RM})  only two are not redundant since  (\ref{eq:equi_equa_RM})$_1=\tan \varphi$  (\ref{eq:equi_equa_RM})$_2$ and  (\ref{eq:equi_equa_RM})$_3=\tan \varphi$  (\ref{eq:equi_equa_RM})$_4$.

It is also pointed out that the two remaining linearly independent equations  (\ref{eq:equi_equa_RM})$_{1,3}$ can be uncoupled and are of the Bessel ODE type (see Appendix \ref{app:bessel}). Indeed, if we \textit{take their sum and difference}, while being careful of substituting $g_p(r)=g_{1}(r) + g_{2}(r)$ and $g_m(r)=g_{1}(r) - g_{2}(r)$ along with their derivatives, we deduce
\begin{align}
	\boldsymbol{\vartheta} \, \sin \varphi \left(a_{1} \, \mu \, L_c^2 \left(3 \, g_{m}'(r ) + r \, g_{m}''(r )\right) - 2 \, r \, \mu _e (g_{m}(r ) + 1) - 2 \, r \, g_{m}(r ) \, \mu _{\mbox{\tiny micro}}\right)
	&= 0
	\, ,
	\label{eq:equi_equa_RM_2}
	\\*
	\frac{1}{3} \, \boldsymbol{\vartheta} \, \sin \varphi \left(6 \, r \, \mu _c \, (g_{p}(r ) - 1) - \mu \, L_c^2 (a_{1} + 2 a_{3}) \left(3 \, g_{p}'(r ) + r \, g_{p}''(r )\right)\right)
	&= 0
	\, .
	\notag
\end{align}
Since $g_1(r) \coloneqq \frac{g_{p}(r) + g_{m}(r)}{2}$ and $g_2(r) \coloneqq \frac{g_{p}(r) - g_{m}(r)}{2}$, the solution in terms of $g_1(r)$ and $g_2(r)$ of  (\ref{eq:equi_equa_RM_2}) is 
\begin{align}
	g_1(r) = \, & 
	\frac{1}{2} \left(
	1
	-\frac{
	 i A_{1} \, I_1\left(      \frac{r \, f_2}{L_c}\right)
	-  A_{2} \, Y_1\left(-i \, \frac{r \, f_2}{L_c}\right)
	+i B_{1} \, I_1\left(      \frac{r \, f_1}{L_c}\right)
	-  B_{2} \, Y_1\left(-i \, \frac{r \, f_1}{L_c}\right)
	}{
	r
	}
	-\frac{\mu _e}{\mu _e+\mu _{\mbox{\tiny micro}}}
	\right) \, ,
	\notag
	\\*
	g_2(r) = \, & 
	\frac{1}{2} \left(
	1
	+\frac{
	 i A_{1} \, I_1\left(      \frac{r \, f_2}{L_c}\right)
	-  A_{2} \, Y_1\left(-i \, \frac{r \, f_2}{L_c}\right)
	-i B_{1} \, I_1\left(      \frac{r \, f_1}{L_c}\right)
	+  B_{2} \, Y_1\left(-i \, \frac{r \, f_1}{L_c}\right)
	}{
	r
	}
	+\frac{\mu _e}{\mu _e+\mu _{\mbox{\tiny micro}}}
	\right) \, ,
	\label{eq:sol_fun_RM}
	\\*
	f_{1} \coloneqq \, &
	\sqrt{\frac{6 \mu _c}{(a_{1} + 2 a_{3}) \, \mu}} \, ,
	\qquad
	f_{2} \coloneqq
	\sqrt{\frac{2(\mu _e + \mu _{\mbox{\tiny micro}})}{a_{1} \, \mu }} \, ,
	\notag
\end{align}
where $I_{n}\left(\cdot\right)$ is the \textit{modified Bessel function of the first kind}, $Y_{n}\left(\cdot\right)$ is the \textit{Bessel function of the second kind} (see appendix \ref{app:bessel} for the formal definitions), and
$A_1$, $B_1$, $A_2$, $B_2$ are integration constants.

The values of $A_1$, $B_1$ are determined from  the boundary conditions  (\ref{eq:BC_RM_gen}), while, due to the divergent behaviour of the Bessel function of the second kind at $r=0$, we have to set $A_2=0$ and $B_2=0$ in order to have a continuous solution.
The fulfilment of the boundary conditions  (\ref{eq:BC_RM_gen}) allows us to find the expressions of the integration constants
\begin{align}
	A_1 = \, &
	\frac{
	i \, L_c
	\left(
	3 f_1 \, R \, z_1 \, I_0\left(\frac{R \, f_1}{L_c}\right)
	-2 L_c \, I_1\left(\frac{R \, f_1}{L_c}\right)
	\right)
	}{
	f_2 \, L_c \, I_0\left(\frac{R \, f_2}{L_c}\right) I_1\left(\frac{R \, f_1}{L_c}\right)
	+f_1 \, z_1 \, I_0\left(\frac{R \, f_1}{L_c}\right) \left(L_c \, I_1\left(\frac{R \, f_2}{L_c}\right)
	-2 f_2 \, R \, I_0\left(\frac{R \, f_2}{L_c}\right)\right)
	}
	\frac{
	\mu _{\mbox{\tiny micro}} 
	}{
	\mu _e+\mu _{\mbox{\tiny micro}}
	}
	\, ,
	\label{eq:BC_RM_3}
	\\
	B_1 = \, & 
	\frac{
	i \, L_c
	\left(
	f_2 \, R \, I_0\left(\frac{R \, f_2}{L_c}\right)
	-2 L_c \, I_1\left(\frac{R \, f_2}{L_c}\right)
	\right)
	}{
	f_2 \, L_c \, I_0\left(\frac{R \, f_2}{L_c}\right) I_1\left(\frac{R \, f_1}{L_c}\right)
	+f_1 \, z_1 \, I_0\left(\frac{R \, f_1}{L_c}\right) \left(L_c \, I_1\left(\frac{R \, f_2}{L_c}\right)
	-2 f_2 \, R \, I_0\left(\frac{R \, f_2}{L_c}\right)\right)
	}
	\frac{
	\mu _{\mbox{\tiny micro}} 
	}{
	\mu _e+\mu _{\mbox{\tiny micro}}
	}
	\, ,
	\notag
	\\
	z_1 \coloneqq \, & \frac{a_1 + 2a_3}{3a_1}
	\, .
	\notag
\end{align}
\allowdisplaybreaks
The classical torque, the higher-order torque, and energy (per unit length d$z$) expressions are
\begin{align}
	M_{\mbox{c}} (\boldsymbol{\vartheta}) 
	\coloneqq&
	\int_{0}^{2\pi}
	\int_{0}^{R}
	\Big[
	\langle
	\widetilde{\boldsymbol{\sigma}} \, \boldsymbol{e}_{z} , \boldsymbol{e}_{\varphi}
	\rangle
	r
	\Big] r
	\, \mbox{d}r \, \mbox{d}\varphi
	\notag
	\\*
	=&
	\left[
	\left(
	\frac{
	4 \mu _c \, f_2 \, I_2\left(\frac{R \, f_1}{L_c}\right) I_2\left(\frac{R \, f_2}{L_c}\right) \left(\frac{L_c}{R}\right)^2
	}{
	\mu _e \, f_1 \left(f_1 \, z_1 \, I_0\left(\frac{R \, f_1}{L_c}\right) \left(2 f_2 \, I_0\left(\frac{R \, f_2}{L_c}\right) 
	- I_1\left(\frac{R \, f_2}{L_c}\right) \frac{L_c}{R} \right)
	- f_2 \, I_0\left(\frac{R \, f_2}{L_c}\right) I_1\left(\frac{R \, f_1}{L_c}\right) \frac{L_c}{R}\right)
	}
	\right.
	\right.
	\notag
	\\*[2mm]
	&+
	\frac{
	f_1 \, z_1 \, I_0\left(\frac{R \, f_1}{L_c}\right) \left(
	24 I_1\left(\frac{R \, f_2}{L_c}\right) \left(\frac{L_c}{R}\right)^3
	-12 f_2 \, I_0\left(\frac{R \, f_2}{L_c}\right) \left(\frac{L_c}{R}\right)^2
	- f_2^2 \, I_1\left(\frac{R \, f_2}{L_c}\right) \frac{L_c}{R}
	+2 f_2^3 \, I_0\left(\frac{R \, f_2}{L_c}\right)\right) 
	}{
	f_2^2 \left(f_1 \, z_1 \, I_0 \left(\frac{R \, f_1}{L_c}\right) \left(2 f_2 \, I_0\left(\frac{R \, f_2}{L_c}\right)
	- I_1\left(\frac{R \, f_2}{L_c}\right) \frac{R}{L_c} \right)
	- f_2 \, I_0\left(\frac{R \, f_2}{L_c}\right) I_1\left(\frac{R \, f_1}{L_c}\right) \frac{R}{L_c}\right)
	}
	\notag
	\\*[2mm]
	&
	\left.
	\left.
	-
	\frac{
	I_1\left(\frac{R \, f_1}{L_c}\right) \left(
	16 I_1\left(\frac{R \, f_2}{L_c}\right) \left(\frac{L_c}{R}\right)^4
	-8 f_2 \, I_0\left(\frac{R \, f_2}{L_c}\right) \left(\frac{L_c}{R}\right)^3
	+ f_2^3 \, I_0\left(\frac{R \, f_2}{L_c}\right) \frac{L_c}{R}
	\right) 
	}{
	f_2^2 \left(f_1 \, z_1 \, I_0 \left(\frac{R \, f_1}{L_c}\right) \left(2 f_2 \, I_0\left(\frac{R \, f_2}{L_c}\right)
	- I_1\left(\frac{R \, f_2}{L_c}\right) \frac{R}{L_c} \right)
	- f_2 \, I_0\left(\frac{R \, f_2}{L_c}\right) I_1\left(\frac{R \, f_1}{L_c}\right) \frac{R}{L_c}\right)
	}
	\right)
	\frac{\mu _e \, \mu _{\mbox{\tiny micro}}}{\mu _e + \mu _{\mbox{\tiny micro}}}
	\right]
	I_{p} \, 
	\boldsymbol{\vartheta}
	\notag
	\\*
	=& 
	\, T_{\mbox{c}} \, \boldsymbol{\vartheta} \, ,
	\notag
	\\[3mm]
	M_{\mbox{m}}(\boldsymbol{\vartheta}) 
	\coloneqq&
	\int_{0}^{2\pi}
	\int_{0}^{R}
	\Big[
	\langle
	\left(\boldsymbol{m} \times \boldsymbol{e}_{z}\right)
	\boldsymbol{e}_{\varphi} ,
	\boldsymbol{e}_{r}
	\rangle
	-
	\langle
	\left(\boldsymbol{m} \times \boldsymbol{e}_{z}\right)
	\boldsymbol{e}_{r} ,
	\boldsymbol{e}_{\varphi}
	\rangle
	\Big]
	\, r
	\, \mbox{d}r \, \mbox{d}\varphi
	\label{eq:torque_stiffness_RM}
	\\*
	=&
	\left[
	\frac{
		I_2\left(\frac{R \, f_2}{L_c}\right)
		\left(
		q_{1} \, I_0 \left(\frac{R \, f_1}{L_c}\right)
		\left(\frac{L_c}{R}\right)^2
		-q_{2} \, I_1 \left(\frac{R \, f_1}{L_c}\right)
		\left(\frac{L_c}{R^3}\right)^3
		\right)
	}{
		f_1 \, z_1 \, I_0\left(\frac{R \, f_1}{L_c}\right) \left(2 f_2 \, I_0\left(\frac{R \, f_2}{L_c}\right) 
		- I_1\left(\frac{R \, f_2}{L_c}\right) \frac{L_c}{R} \right)
		- f_2 \, I_0\left(\frac{R \, f_2}{L_c}\right) I_1\left(\frac{R \, f_1}{L_c}\right) \frac{L_c}{R}
	}
	\frac{4 \mu \, \mu _{\mbox{\tiny micro}}}{3 \left(\mu _e + \mu _{\mbox{\tiny micro}}\right)}
	\right]
	I_{p} \, 
	\boldsymbol{\vartheta}
	\notag
	\\*
	=& 
	\, T_{\mbox{m}} \, \boldsymbol{\vartheta} \, ,
	\notag
	\\[3mm]
	W_{\mbox{tot}} (\boldsymbol{\vartheta}) 
	\coloneqq&
	\int_{0}^{2\pi}
	\int_{0}^{R}
	W \left(\boldsymbol{\mbox{D}u}, \boldsymbol{P}, \mbox{Curl}\boldsymbol{P}\right) \, \, r
	\, \mbox{d}r \, \mbox{d}\varphi
	\notag
	\\*
	=&
	\frac{1}{2}
	\left[
	\left(
	\frac{
		4 \mu _c \, f_2 \, I_2\left(\frac{R \, f_1}{L_c}\right) I_2\left(\frac{R \, f_2}{L_c}\right) \left(\frac{L_c}{R}\right)^2
	}{
		\mu _e \, f_1 \left(f_1 \, z_1 \, I_0\left(\frac{R \, f_1}{L_c}\right) \left(2 f_2 \, I_0\left(\frac{R \, f_2}{L_c}\right) 
		- I_1\left(\frac{R \, f_2}{L_c}\right) \frac{L_c}{R} \right)
		- f_2 \, I_0\left(\frac{R \, f_2}{L_c}\right) I_1\left(\frac{R \, f_1}{L_c}\right) \frac{L_c}{R}\right)
	}
	\right.
	\right.
	\notag
	\\*[2mm]
	&+
	\frac{
		f_1 \, z_1 \, I_0\left(\frac{R \, f_1}{L_c}\right) \left(
		24 I_1\left(\frac{R \, f_2}{L_c}\right) \left(\frac{L_c}{R}\right)^3
		-12 f_2 \, I_0\left(\frac{R \, f_2}{L_c}\right) \left(\frac{L_c}{R}\right)^2
		- f_2^2 \, I_1\left(\frac{R \, f_2}{L_c}\right) \frac{L_c}{R}
		+2 f_2^3 \, I_0\left(\frac{R \, f_2}{L_c}\right)\right) 
	}{
		f_2^2 \left(f_1 \, z_1 \, I_0 \left(\frac{R \, f_1}{L_c}\right) \left(2 f_2 \, I_0\left(\frac{R \, f_2}{L_c}\right)
		- I_1\left(\frac{R \, f_2}{L_c}\right) \frac{R}{L_c} \right)
		- f_2 \, I_0\left(\frac{R \, f_2}{L_c}\right) I_1\left(\frac{R \, f_1}{L_c}\right) \frac{R}{L_c}\right)
	}
	\notag
	\\*[2mm]
	&
	-
	\frac{
		I_1\left(\frac{R \, f_1}{L_c}\right) \left(
		16 I_1\left(\frac{R \, f_2}{L_c}\right) \left(\frac{L_c}{R}\right)^4
		-8 f_2 \, I_0\left(\frac{R \, f_2}{L_c}\right) \left(\frac{L_c}{R}\right)^3
		+ f_2^3 \, I_0\left(\frac{R \, f_2}{L_c}\right) \frac{L_c}{R}
		\right) 
	}{
		f_2^2 \left(f_1 \, z_1 \, I_0 \left(\frac{R \, f_1}{L_c}\right) \left(2 f_2 \, I_0\left(\frac{R \, f_2}{L_c}\right)
		- I_1\left(\frac{R \, f_2}{L_c}\right) \frac{R}{L_c} \right)
		- f_2 \, I_0\left(\frac{R \, f_2}{L_c}\right) I_1\left(\frac{R \, f_1}{L_c}\right) \frac{R}{L_c}\right)
	}
	\notag
	\\*[2mm]
	&
	\left.
	\left.
	+
	\frac{
		4\mu \, I_2\left(\frac{R \, f_2}{L_c}\right)
		\left(
		q_{1} \, I_0 \left(\frac{R \, f_1}{L_c}\right)
		\left(\frac{L_c}{R}\right)^2
		-q_{2} \, I_1 \left(\frac{R \, f_1}{L_c}\right)
		\left(\frac{L_c}{R^3}\right)^3
		\right)
	}{
		3\mu _e \left(f_1 \, z_1 \, I_0\left(\frac{R \, f_1}{L_c}\right) \left(2 f_2 \, I_0\left(\frac{R \, f_2}{L_c}\right) 
		- I_1\left(\frac{R \, f_2}{L_c}\right) \frac{L_c}{R} \right)
		- f_2 \, I_0\left(\frac{R \, f_2}{L_c}\right) I_1\left(\frac{R \, f_1}{L_c}\right) \frac{L_c}{R}\right)
	}
	\right)
	\frac{\mu _e \, \mu _{\mbox{\tiny micro}}}{\mu _e + \mu _{\mbox{\tiny micro}}}
	\right]
	I_{p} \, 
	\boldsymbol{\vartheta}^2
	\notag
	\\*
	=& \frac{1}{2} \, T_{\mbox{w}} \, \boldsymbol{\vartheta}^2
	\, ,
	\notag
	\\[3mm]
	q_{1} \coloneqq& \, 3 a_{1} \, f_1 \, f_2 \, z_1 \, ,
	\qquad
	q_{2} \coloneqq 2 f_2 (a_{1}-a_{3}) \, .
	\notag
\end{align}
Again it holds
\allowdisplaybreaks
\begin{equation}
\frac{\mbox{d}}{\mbox{d}\boldsymbol{\vartheta}}W_{\mbox{tot}}(\boldsymbol{\vartheta}) = M_{\mbox{c}} (\boldsymbol{\vartheta}) + M_{\mbox{m}} (\boldsymbol{\vartheta}) \, ,
\qquad\qquad\qquad
\frac{\mbox{d}^2}{\mbox{d}\boldsymbol{\vartheta}^2}W_{\mbox{tot}}(\boldsymbol{\vartheta})
= T_{\mbox{c}} + T_{\mbox{m}}
= T_{\mbox{w}}
\, .
\end{equation}
Both quantities $M_{\mbox{c}}$ and $W_{\mbox{tot}}$ are immediately accessible in any higher order generalized continuum model.
However, the precise form of $M_{\mbox{m}}$ is difficult to guess.
The latter identity can, therefore, also be seen as an implicit definition of the higher order moment $M_{\mbox{m}}$.
In the Appendix we will provide an independent way of obtaining the notation for $M_{\mbox{m}}$ starting form considerations done on the Cosserat model (see Appendix  \ref{app:class_Coss}).
We provide again the homogenization relations between the macro-parameters, the meso- (with index $(\cdot)_{e}$), and the micro-parameters \cite{d2019effective,neff2014unifying,neff2019identification}
\begin{equation}
\begin{array}{c}
\mu _{\tiny \mbox{macro}} = \frac{\mu_e \, \mu_{\tiny \mbox{micro}}}{\mu_e + \mu_{\tiny \mbox{micro}}} \, ,
\qquad\qquad
\kappa _{\tiny \mbox{macro}} = \frac{\kappa_e \, \kappa_{\tiny \mbox{micro}}}{\kappa_e + \kappa_{\tiny \mbox{micro}}} \, ,
\qquad\qquad
\text{with}
\quad
\left\{
\begin{array}{l}
\kappa_{\tiny \mbox{i}} = \frac{2\mu_{i} + 3\lambda_{i}}{3} \, ,
\\
i = \left\{e, {\tiny \mbox{micro}}, {\tiny \mbox{macro}}\right\} \, ,
\end{array}
\right.
\end{array}
\label{eq:static_homo_relation}
\end{equation}
which can be used to define the following torsional stiffnesses
\begin{equation}
    T_{\tiny \mbox{macro}} =
    \mu_{\tiny \mbox{macro}} \, I_{p} =
    \frac{\mu_{\tiny \mbox{micro}} \, \mu_{e}}{\mu_{\tiny \mbox{micro}} + \mu_{e}} \, I_{p}
    \, ,
    \qquad\qquad
    T_{\tiny \mbox{micro}} =
    \mu_{\tiny \mbox{micro}} \, I_{p}
    \, ,
    \qquad\qquad
    T_{e} =
    \mu_{e} \, I_{p}
    \, .
    \label{eq:diff_tors_stiff}
\end{equation}

\begin{figure}[H]
	\begin{subfigure}{0.32\textwidth}
		\centering
		\includegraphics[width=\textwidth]{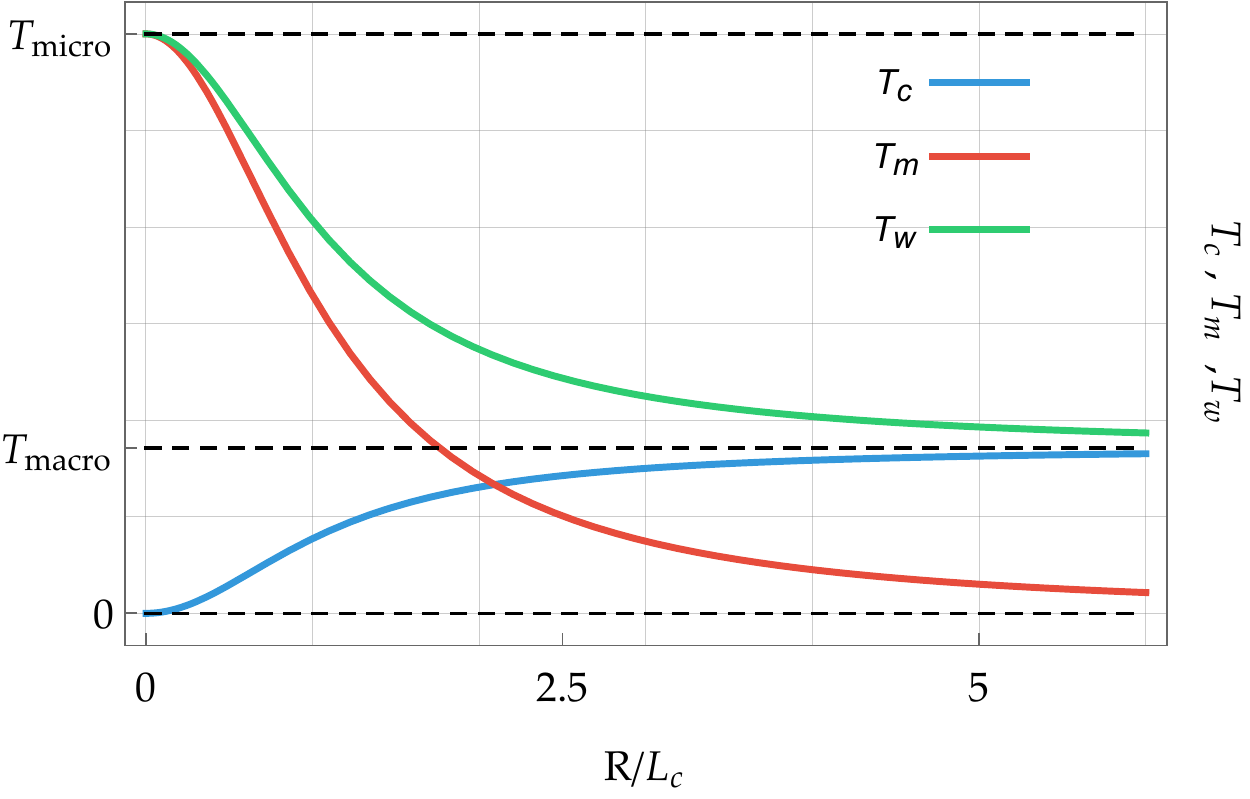}
		\caption{}
	\end{subfigure}
	\hfill
	\begin{subfigure}{0.32\textwidth}
		\centering
		\includegraphics[width=\textwidth]{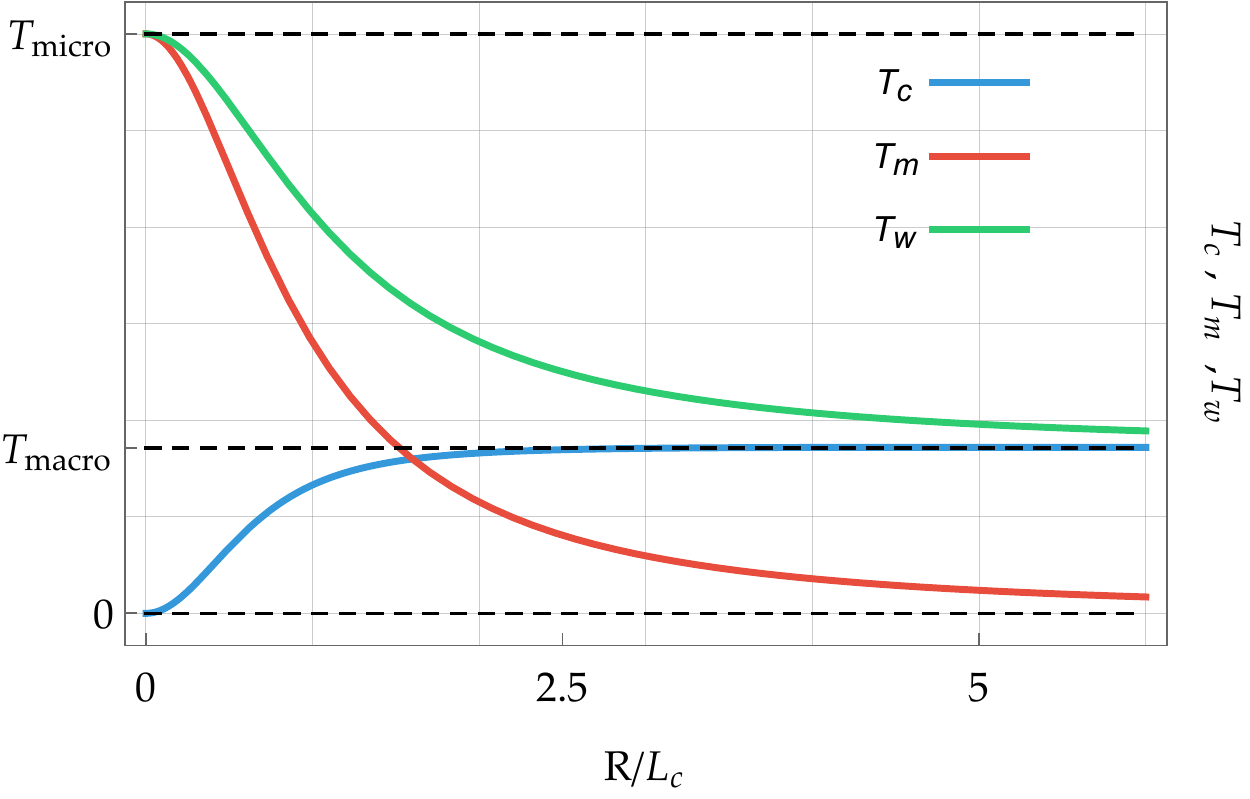}
		\caption{}
	\end{subfigure}
	\hfill
	\begin{subfigure}{0.32\textwidth}
		\centering
		\includegraphics[width=\textwidth]{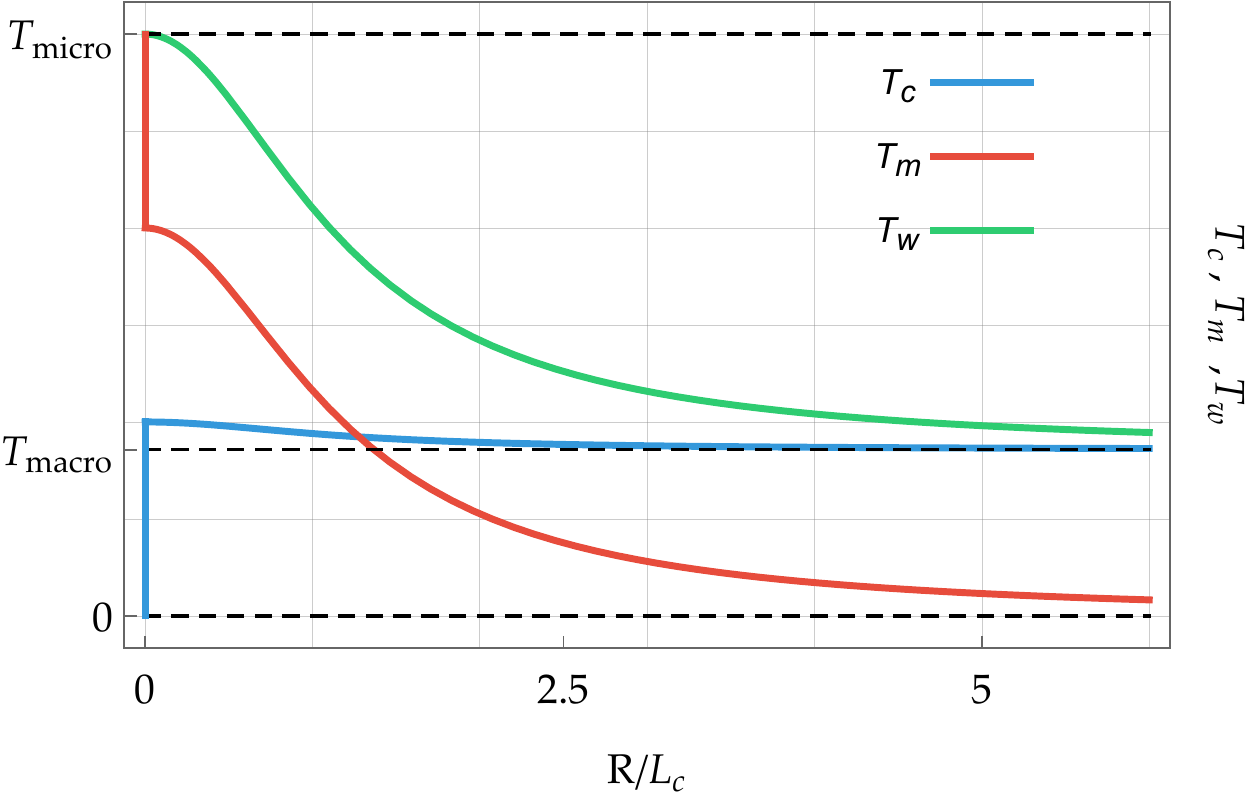}
		\caption{}
	\end{subfigure}
	\caption{(\textbf{Relaxed micromorphic model}) Torsional stiffness for the classical torque $T_{\mbox{c}}$, the higher-order torque $T_{\mbox{m}}$, and the torque energy $T_{\mbox{w}}$ while varying $L_c$ for (a) $\mu_c \to 0$, (b) $\mu_c=1/2$, and (c) $\mu_c \to \infty$. The torsional stiffness remains bounded as $L_c \to \infty$ ($R\to 0$). The values of the parameters used are: $\mu = 1$, $\mu _e= 1/10$, $\mu _{\mbox{\tiny micro}}= 1/4$, $a_1= 1/5$, $a_2= 1/6$, $a_3= 1/7$, $R= 1$.}
	\label{fig:all_plot_RM_3}
\end{figure}

The plots of the torsional stiffness for the classical torque (light blue), the higher-order torque (red), and the torque energy (green) for $\mu_c=\{0,1/2,\infty\}$ while varying $L_c$ is shown in Fig.~\ref{fig:all_plot_RM_3}.

It is here highlighted that the torsional stiffness obtainable from the energy $T_{\mbox{w}}$ is the only stiffness available experimentally.
\subsection{Limits}
\subsubsection{The relaxed micromorphic model with symmetric force stresses ($\mu_c \to 0$)}
\begin{align}
W_{\mbox{tot}} (\boldsymbol{\vartheta}) 
\coloneqq&\notag
\int_{0}^{2\pi}
\int_{0}^{R}
W \left(\boldsymbol{\mbox{D}u}, \boldsymbol{P}, \mbox{Curl}\boldsymbol{P}\right) \, \, r
\, \mbox{d}r \, \mbox{d}\varphi
\label{eq:torque_stiffness_RM_mc}
\\*
=&
\frac{1}{2}
\left[
\frac{
f_2 I_0\left(\frac{R f_2}{L_c}\right) \left(v_{2} \, \mu _{\mbox{\tiny micro}} \frac{L_c^2}{R^2}+\mu _e \left(\mu _e+\mu _{\mbox{\tiny micro}}\right)\right)
}{
f_2 I_0\left(\frac{R f_2}{L_c}\right)
-2 v_{1} I_1\left(\frac{R f_2}{L_c}\right) \frac{L_c}{R}
}
\right.
\\*\notag
&
\qquad\qquad\qquad\qquad
\left.
-
\frac{
	 2 I_1\left(\frac{R f_2}{L_c}\right) \left(v_{2} \, \mu _{\mbox{\tiny micro}}\frac{L_c^2}{R^2}+v_{1} \mu _e \left(\mu _e+\mu _{\mbox{\tiny micro}}\right)\right) \frac{L_c}{R}
}{
	f_2 I_0\left(\frac{R f_2}{L_c}\right)
	-2 v_{1} I_1\left(\frac{R f_2}{L_c}\right) \frac{L_c}{R}
}
\right]
\frac{
	\mu _{\mbox{\tiny micro}}
}{
	\left(\mu _e+\mu _{\mbox{\tiny micro}}\right)^2
}
I_{p} \, 
\boldsymbol{\vartheta}^2
\, ,
\notag
\\[3mm]
v_{1} \coloneqq & \, \frac{a_{1}+2 a_{3}}{a_{1}+8 a_{3}} \, ,
\qquad
v_{2} \coloneqq \frac{24 a_{1} \, a_{3} \, \mu }{a_{1}+8 a_{3}} \, .
\notag
\end{align}
\begin{figure}[H]
	\centering
	\includegraphics[height=5.5cm]{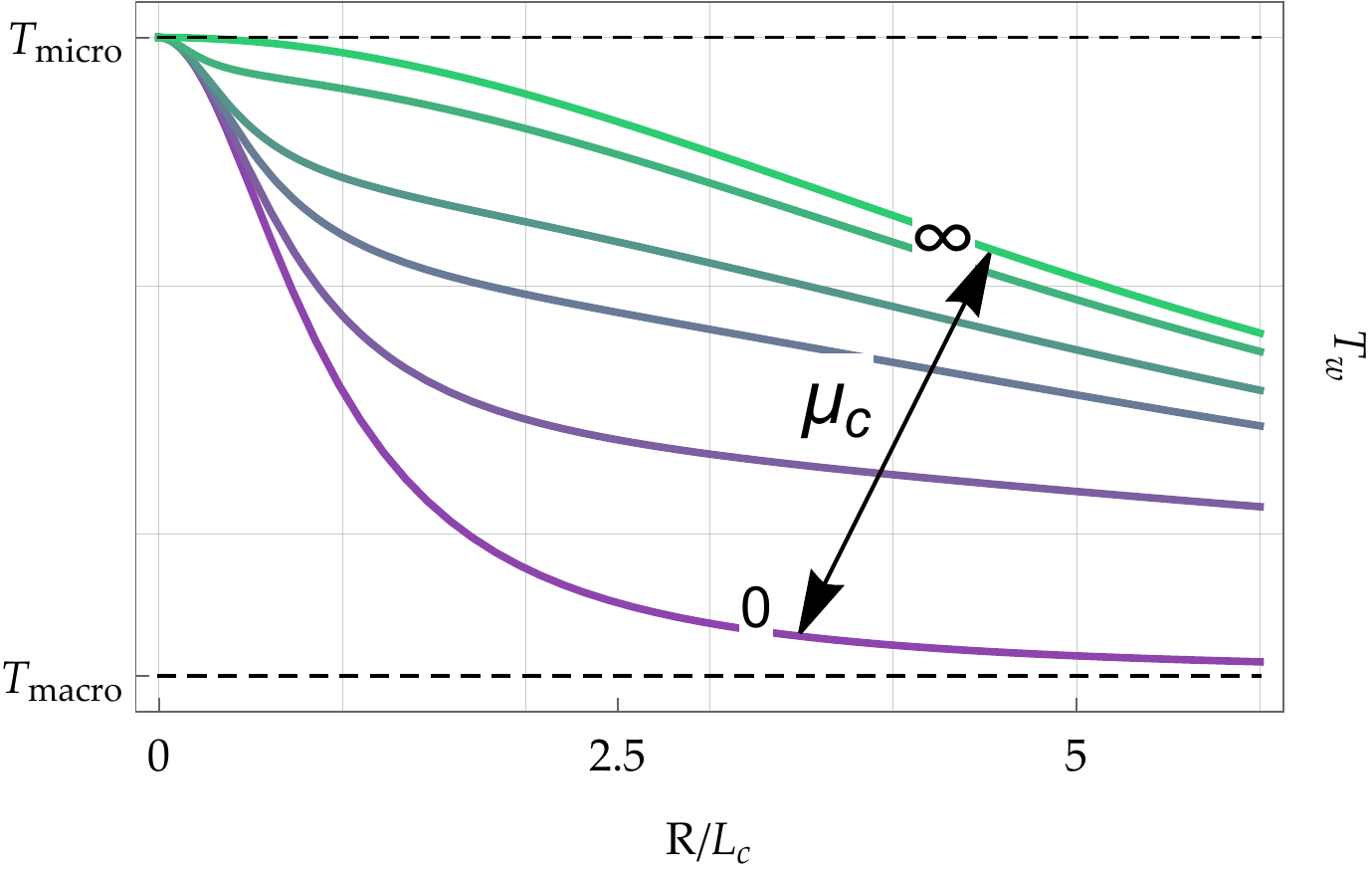}
	\caption{(\textbf{Relaxed micromorphic model})
		Torsional stiffness for the torque energy while varying $L_c$, for different values of $\mu_c=\{0, 1/30, 1/10, 1/5, 1, \infty\}$.
		The torsional stiffness remains bounded as $L_c \to \infty$ ($R\to 0$) and the model does not collapse into a linear elastic one.
		The values of the other parameters used are: $\mu = 1$, $\mu _e= 1/3$, $\mu _{\mbox{\tiny micro}}= 1/4$, $a_1= 10$, $a_3= 1/50$, $R= 1$.
		Here, varying $\mu_c$ does not intervene with $T_{\tiny \text{macro}}$ and $T_{\tiny \text{micro}}$.
		}
	\label{fig:all_plot_RM_mc}
\end{figure}
Note that the torsional stiffness at the micro-scale $T_{\mbox{\tiny micro}}$ is here independent of the Cosserat couple modulus $\mu_c$, see  (\ref{eq:diff_tors_stiff}).
\subsubsection{The relaxed micromorphic model with conformal curvature energy ($a_3 = 0$) while varying the Cosserat couple modulus $\mu_c$}
\label{sec:relax_conformal}
In the particular case for which the parameter $a_3$ is equal to zero the elastic energy turns into
\begin{align}
	W \left(\boldsymbol{\mbox{D}u}, \boldsymbol{P},\mbox{Curl}\,\boldsymbol{P}\right)
	=
	&
	\, \mu_{e} \left\lVert \mbox{sym} \left(\boldsymbol{\mbox{D}u} - \boldsymbol{P} \right) \right\rVert^{2}
	+ \frac{\lambda_{e}}{2} \mbox{tr}^2 \left(\boldsymbol{\mbox{D}u} - \boldsymbol{P} \right) 
	+ \mu_{c} \left\lVert \mbox{skew} \left(\boldsymbol{\mbox{D}u} - \boldsymbol{P} \right) \right\rVert^{2}
	\label{eq:energy_RMM_conformal}
	\\*
	&
	+ \mu_{\tiny \mbox{micro}} \left\lVert \mbox{sym}\,\boldsymbol{P} \right\rVert^{2}
	+ \frac{\lambda_{\tiny \mbox{micro}}}{2} \mbox{tr}^2 \left(\boldsymbol{P} \right)
	+ \frac{\mu \,L_c^2 }{2} \,
	a_1 \, \left\lVert \mbox{dev sym} \, \mbox{Curl} \, \boldsymbol{P}\right\rVert^2
	\, .
	\notag
\end{align}
In this case, the torsional stiffness at the micro scale, namely for $L_c \to \infty$ ($R\to 0$)
\footnote{
Looking at the analytical solution obtained in \eqref{eq:equi_equa_RM} we see that the expression $\frac{R}{L_{\rm c}}$ solely determines the response. Therefore, we can either fix $R>0$ and send $L_{\rm c}\to \infty$, or fix $L_{\rm c}$ and send $R\to 0$, having the same effect.
}
, depends also on $\mu_c$
\begin{equation}
    \widetilde{T}\coloneqq
    \lim_{L_c\to \infty}T_{w} = \frac{\mu _{\tiny \mbox{micro}} \left(9 \mu _c+\mu _e\right)}{\left(9 \mu _c+\mu _e\right) + \mu _{\tiny \mbox{micro}}} I_{p} \, .
    \label{eq:limi_conformal_RM}
\end{equation}
For $\mu_c \to 0$ we obtain a linear elastic model with stiffness $T_{\mbox{\tiny macro}}$, for $\mu_c \to \infty$ it is recovered a model that has $T_{\mbox{\tiny micro}}$ at the micro-scale, while for intermediate values of $0<\mu_c<\infty$ a torsional stiffness between $T_{\mbox{\tiny macro}}$ and $T_{\mbox{\tiny micro}}$ appears.
\begin{figure}[H]
	\centering
	\includegraphics[height=5.5cm]{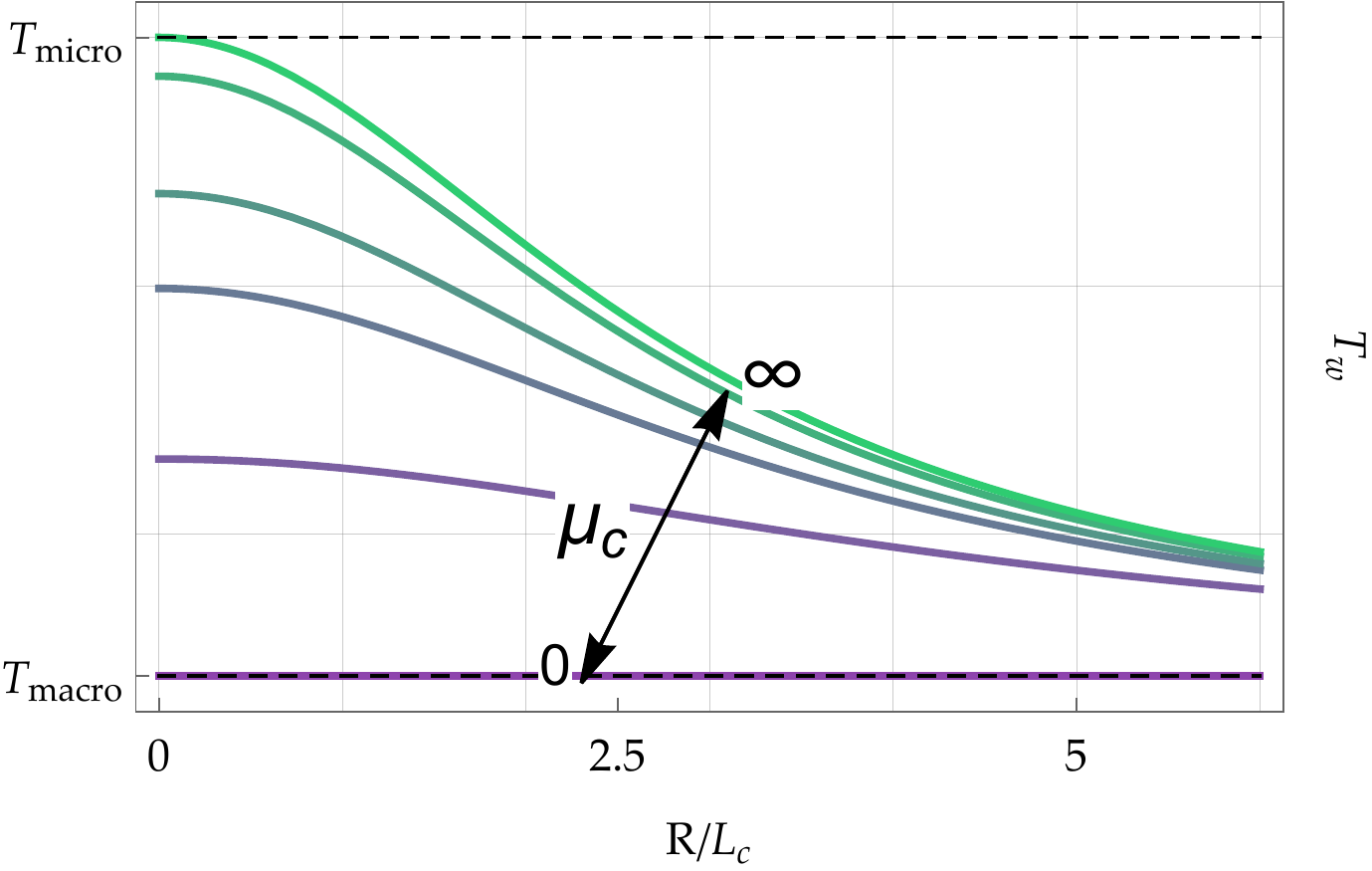}
	\caption{(\textbf{Relaxed micromorphic model with conformal curvature energy})
	Torsional stiffness for the torque energy while varying $L_c$, for different values of $\mu_c=\{0, 1/30, 1/10, 1/5, 1, \infty\}$.
	The torsional stiffness remains bounded as $L_c \to \infty$ ($R\to 0$) and the model does not collapse in a linear elastic one beside the case $\mu_{\mbox{\tiny c}}=0$.
	The values of the other parameters used are: $\mu = 1$, $\mu _e= 1/3$, $\mu _{\mbox{\tiny micro}}= 1/4$, $a_1= 2$, $R= 1$.
	In this case, varying $\mu_c$ influences the torsional stiffness also for small specimen size.
	}
	\label{fig:all_plot_RM_a3_0}
\end{figure}

We may consider  a further limit in  (\ref{eq:limi_conformal_RM}). It holds
\begin{equation}
    \overline{T}
    \coloneqq
    \lim_{\mu_{\mbox{\tiny micro}} \to \infty} \widetilde{T}
    =
    \left( 9 \mu _c+\mu _e \right)I_{p}
    =
    \left( 9 \mu _c+\mu _{\mbox{\tiny macro}} \right)I_{p}
    \, ,
    \label{eq:limi_conformal_RM_2}
\end{equation}
where the last relation for which we have $\mu _e = \mu _{\mbox{\tiny macro}}$ is obtained from (\ref{eq:para_meso_scale})$_2$ taking $\mu _{\mbox{\tiny micro}}\to\infty$.

\subsubsection{The Cosserat model as a limit of the relaxed micromorphic model ($\mu_{\mbox{\tiny micro}} \to \infty$)}
The Cosserat model can be obtained from the relaxed micromorphic model by formally letting $\mu_{\mbox{\tiny micro}} \to \infty$ and $\kappa_{\tiny \mbox{micro}} \to \infty$.
\footnote{
For the torsion problem, $\kappa_{\tiny \mbox{micro}}$ does not intervene.
}
From the homogenization formula  (\ref{eq:static_homo_relation}) it is possible to see that for $\mu_{\mbox{\tiny micro}} \to \mu_{\mbox{\tiny macro}}$ we have $\mu_e \to \infty$, while $\mu_{\mbox{\tiny macro}} = \mu_e$ for $\mu_{\mbox{\tiny micro}} \to \infty$, which is the stiffness at the macro-scale for the Cosserat model.
\begin{figure}[H]
	\begin{subfigure}{0.5\textwidth}
		\centering
		\includegraphics[width=\textwidth]{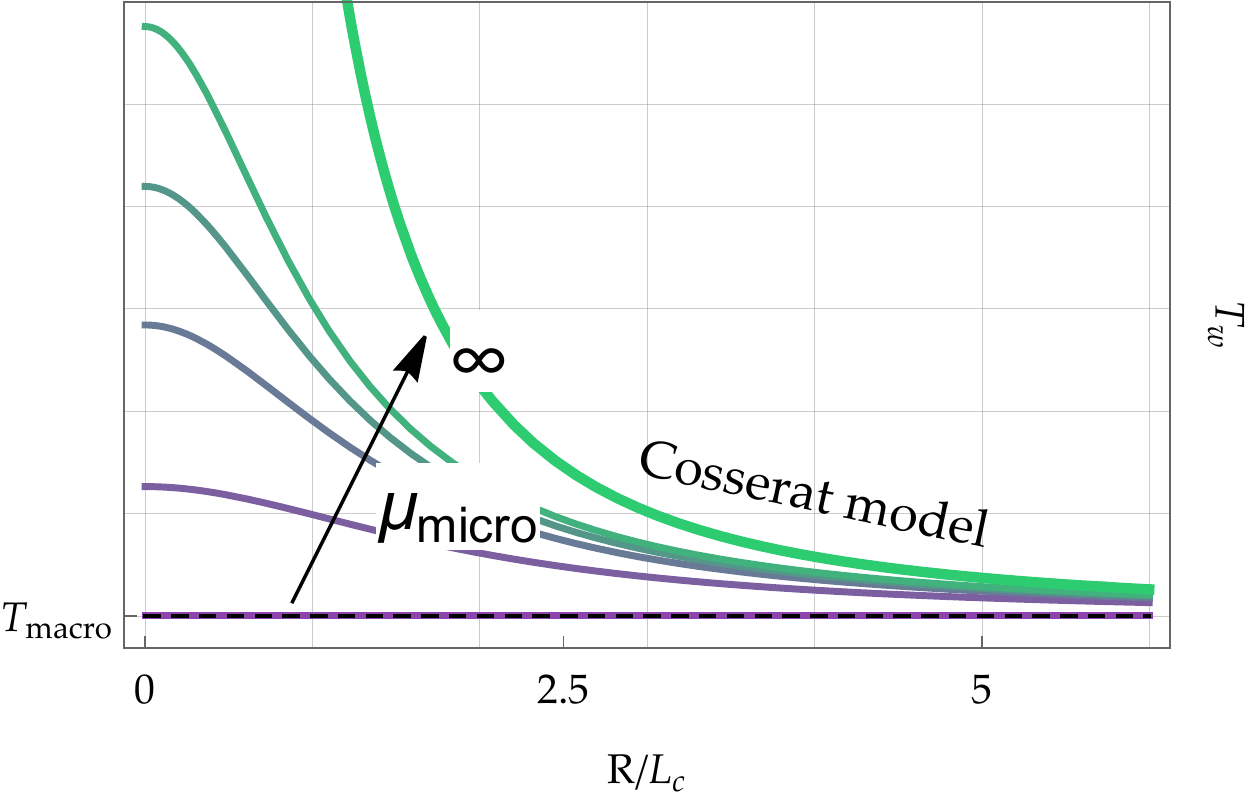}
		\caption{}
	\end{subfigure}
	\hfill
	\begin{subfigure}{0.5\textwidth}
		\centering
		\includegraphics[width=\textwidth]{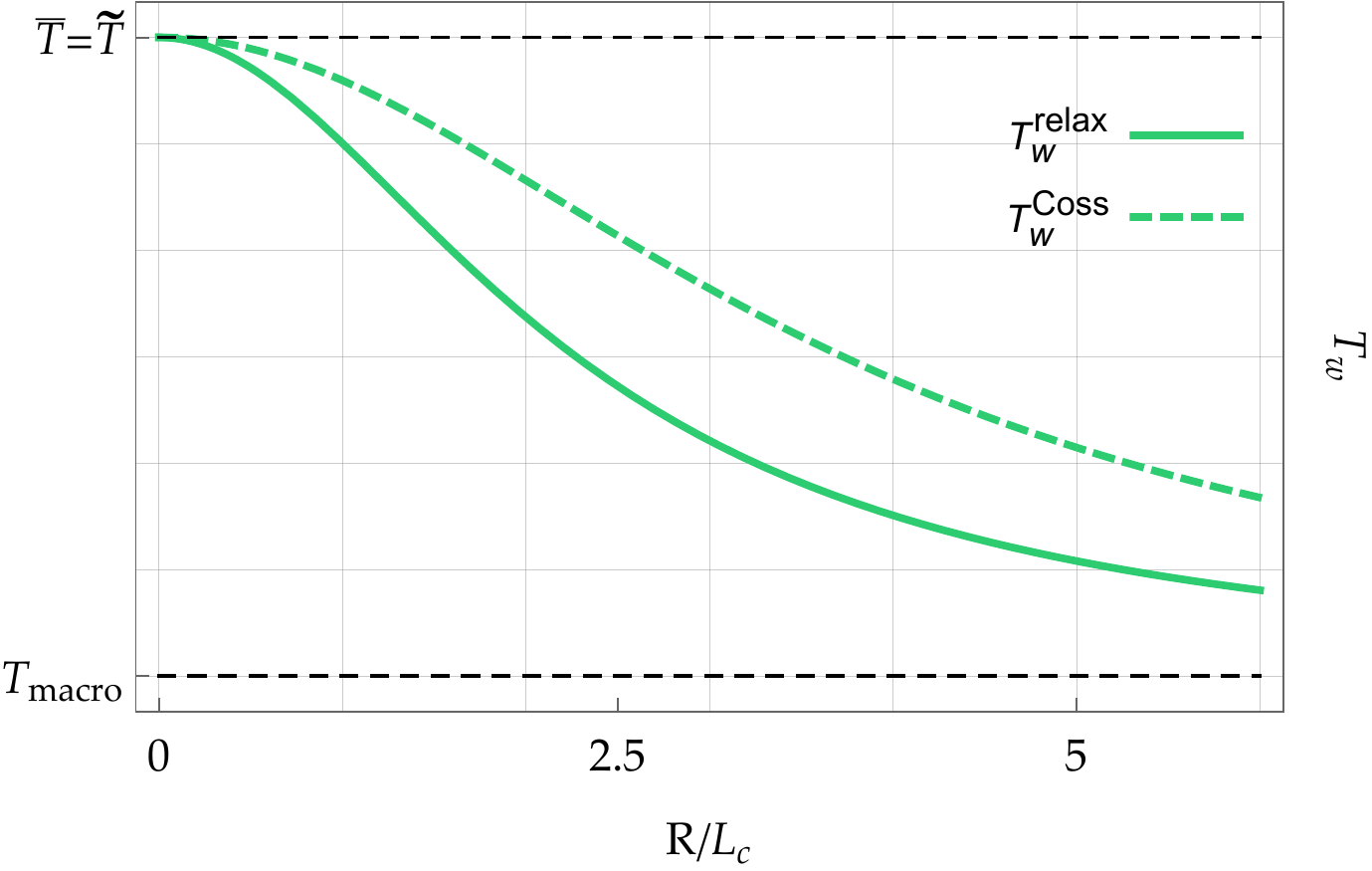}
		\caption{}
	\end{subfigure}
	\caption{
	    (a) (\textbf{Relaxed micromorphic model with full curvature})
		Torsional stiffness for the torque energy while varying $L_c$, for different values of $\mu_{\mbox{\tiny micro}}=\{0, 1/20, 1/7, 1/4, 1/2, \infty\}$.
		The torsional stiffness becomes unbounded as $L_c \to \infty$ ($R\to 0$) when $\mu_{\mbox{\tiny micro}} \to \infty$.
		The values of the other parameters used are: $\mu = 1$, $\mu _{\mbox{\tiny macro}}= 1/10$, $\mu _c= 1/2$, $a_1= 1/5$, $a_3= 1/7$, $R= 1$. The Cosserat solution appears for $\mu_{\mbox{\tiny micro}} \to \infty$.
		(b) (\textbf{Cosserat model and relaxed micromorphic model with conformal curvature}).
		Torsional stiffness for the torque energy while varying $L_c$. The torsional stiffness is bounded as $L_c \to \infty$ ($R\to 0$). For the Cosserat model we chose $\mu _c= 1/9$ while for the relaxed micromorphic model $\mu _c= 1/2$ and $\mu _{\mbox{\tiny micro}}= 3$ in order to have the same upper bound $\overline{T}=\widetilde{T}$. The values of the other parameters used are: $\mu = 1$, $\mu _{\mbox{\tiny macro}}= 1$, $a_1= 5$, $R= 1$.}
	\label{fig:all_plot_RM_mm}
\end{figure}
\subsubsection{Sensitivity of the relaxed micromorphic model with respect to the curvature parameters $a_1$ and $a_3$.}
Sensitivity study for the relaxed micromorphic model while varying $a_1$ and $a_3$ independently.
\begin{figure}[H]
	\begin{subfigure}{0.5\textwidth}
		\centering
		\includegraphics[width=\textwidth]{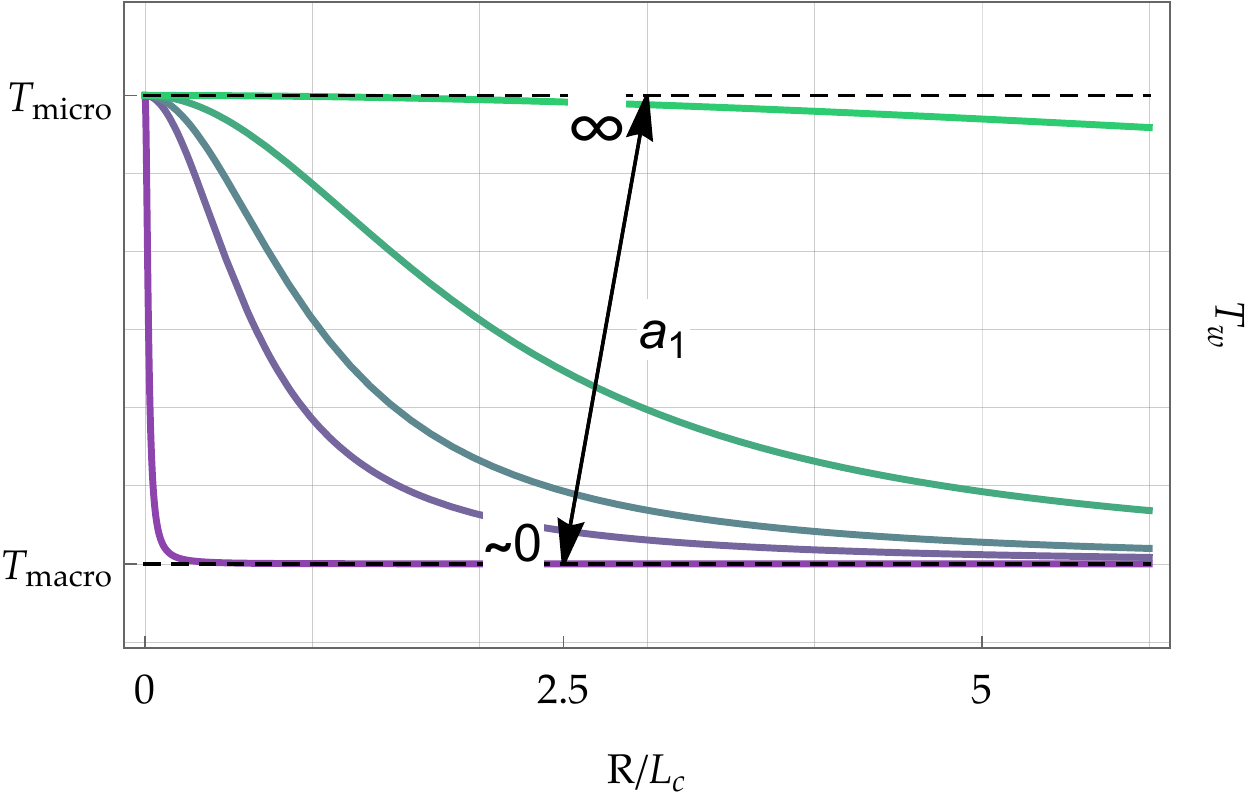}
		\caption{}
	\end{subfigure}
	\hfill
	\begin{subfigure}{0.5\textwidth}
		\centering
		\includegraphics[width=\textwidth]{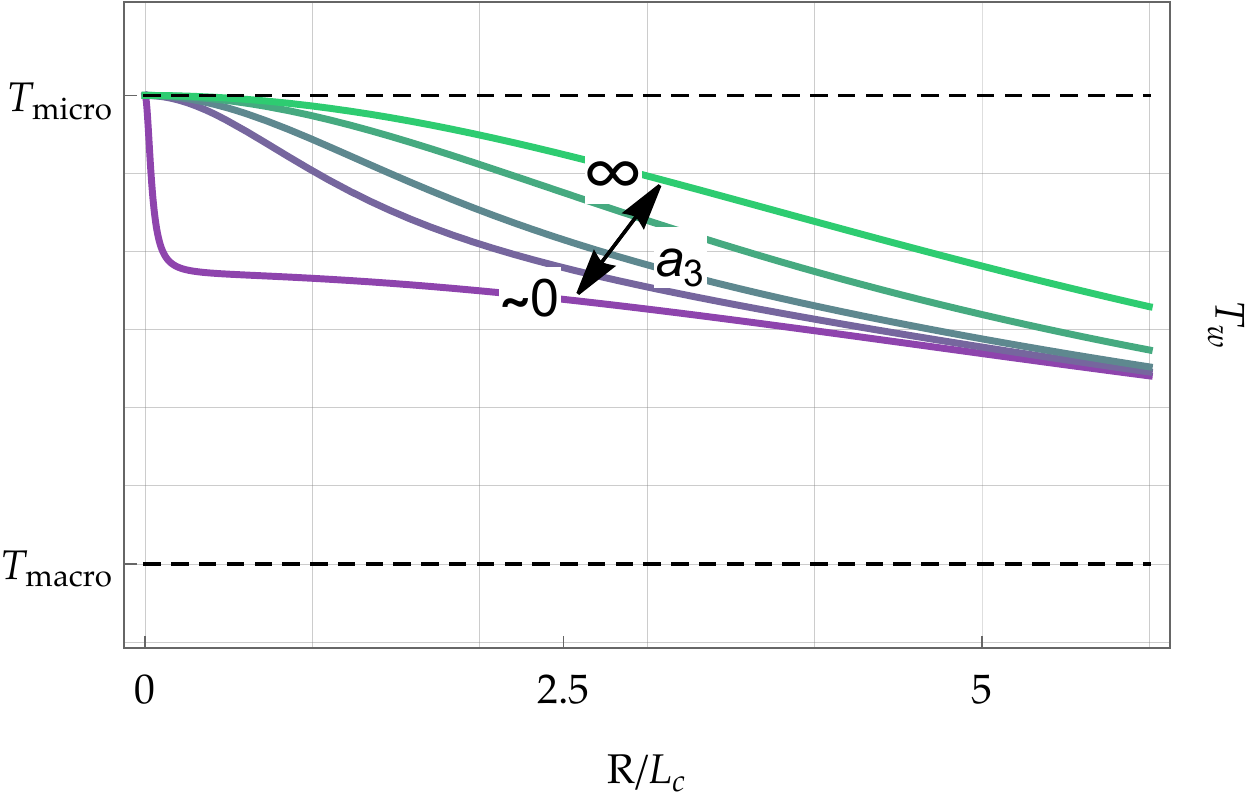}
		\caption{}
	\end{subfigure}
	\caption{(\textbf{Relaxed micromorphic model}) Response of the relaxed micromorphic model while varying (a) the curvature parameter $a_1$ having $a_{3}= 20$ and (b) the curvature parameter $a_3$ having $a_{1}= 20$. The values of the other parameters are $\mu = 1$, $\mu _c= 1/5$, $\mu _e= 1$, $\mu _{\mbox{\tiny micro}}= 1/9$, $R= 1$.}
	\label{fig:para_study_RMM}
\end{figure}
The parametric study represented in Fig. \ref{fig:para_study_RMM} has not been carried out for the limit $a_1 \to 0$ and $a_3 \to 0$ since we would have had an indeterminate form for $L_c \to \infty$, and that is why we used the symbol $\sim 0$.
The solution of the problem while having $a_3 = 0$ a priori is analyzed carefully in  Section  \ref{sec:relax_conformal}, and the solution of the problem while having $a_1 = 0$ a priori make the relaxed model collapse into a classical linear elastic model with torsional stiffness $T_{\tiny \mbox{macro}}$.

\subsection{Finite element simulations}
Using finite element analysis as a tool of comparison, in this section we will \begin{itemize}
	\item [i)]
 test the validity of the solution in terms of the hypothesis of small deformations (i.e., small twist rate);
 \item[ii)] discuss the validity of the St.Venant principle for the relaxed micromorphic model.\end{itemize}

In this analysis we take a finite-size cylindrical rod and we apply opposite and equal finite-rotation at both of its ends ($z=\pm L/2$).
Accordingly, the boundary conditions are
\begin{align}
\boldsymbol{u}(z=\pm L/2)
&= 
\left(
\begin{array}{ccc}
 \cos \,\boldsymbol{\pm\Theta}  & \sin \,\boldsymbol{\pm\Theta}  & 0 \\
 -\sin \,\boldsymbol{\pm\Theta} & \cos \,\boldsymbol{\pm\Theta}  & 0 \\
 0 & 0 & 1 \\
\end{array}
\right)
\left(
\begin{array}{ccc}
 x_1 \\
 x_2 \\
 L \\
\end{array}
\right)
-
\left(
\begin{array}{ccc}
 x_1 \\
 x_2 \\
 L \\
\end{array}
\right)
\, ,
\\*
\boldsymbol{P}(z=\pm L/2) \times \boldsymbol{e}_1
&= 
\left(
\begin{array}{ccc}
 -P_{12} & P_{11} & 0 \\
 -P_{22} & P_{21} & 0 \\
 -P_{32} & P_{31} & 0 \\
\end{array}
\right)
=
\left(
\begin{array}{ccc}
 \sin \,\boldsymbol{\pm\Theta}  & \cos \, \boldsymbol{\pm\Theta} -1 & 0 \\
 1 - \cos \,\boldsymbol{\pm\Theta}  & \sin \,\boldsymbol{\pm\Theta}  & 0 \\
 0 & 0 & 0 \\
\end{array}
\right)
=
\mbox{D}\boldsymbol{u} (z=\pm L/2) \times \boldsymbol{e}_1
\, ,
\notag
\end{align}
where $\boldsymbol{P}(z) \times \boldsymbol{e}_1 = \mbox{D}\boldsymbol{u} (z)  \times \boldsymbol{e}_1$ are the \textbf{consistent  boundary conditions} on the tangential part for the micro-distortion tensor $\boldsymbol{P}$.

\begin{figure}[H]
	\begin{subfigure}{0.5\textwidth}
		\centering
		\includegraphics[width=\textwidth]{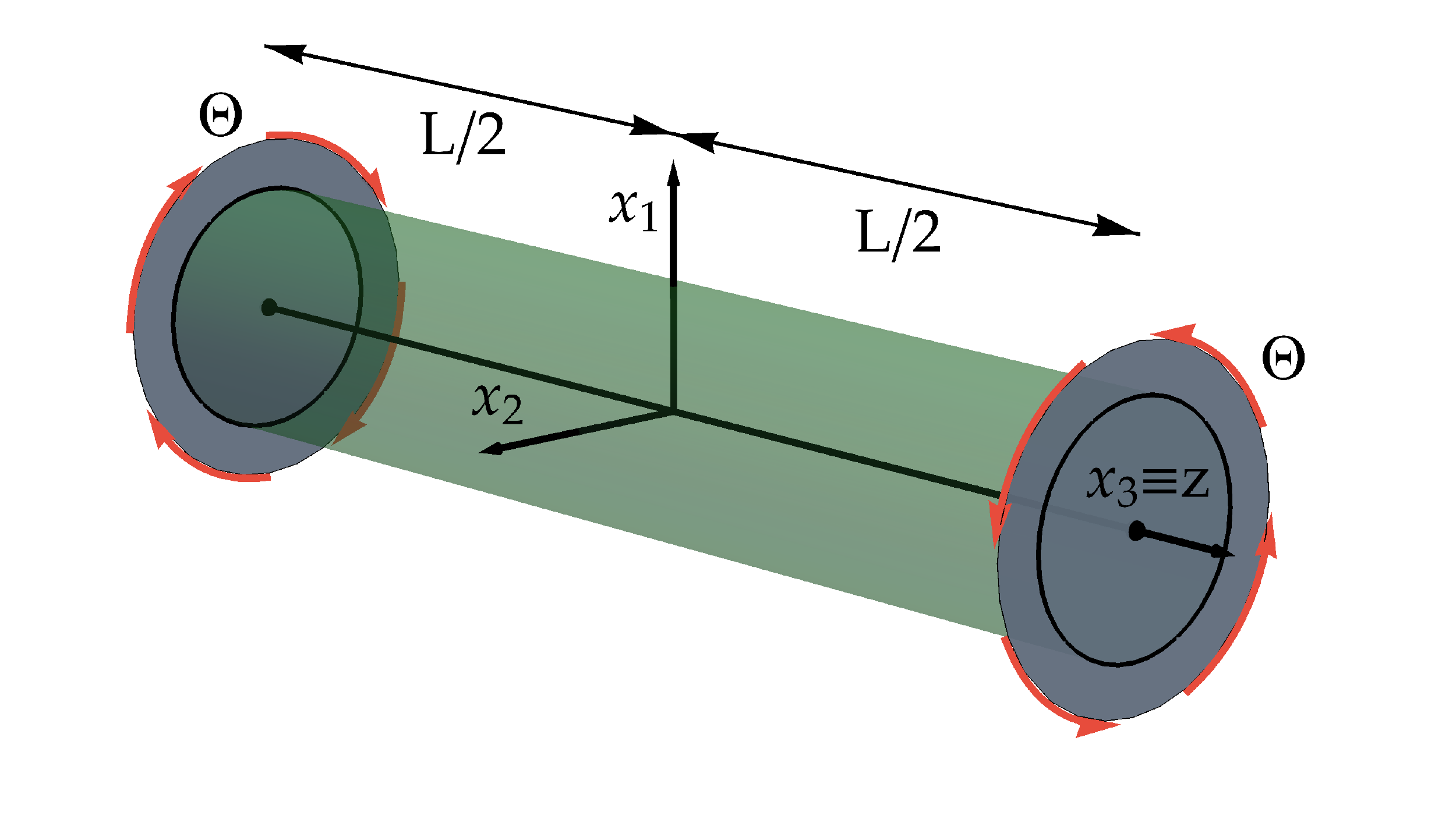}
		\caption{}
	\end{subfigure}
	\hfill
	\begin{subfigure}{0.5\textwidth}
		\centering
		\includegraphics[width=\textwidth]{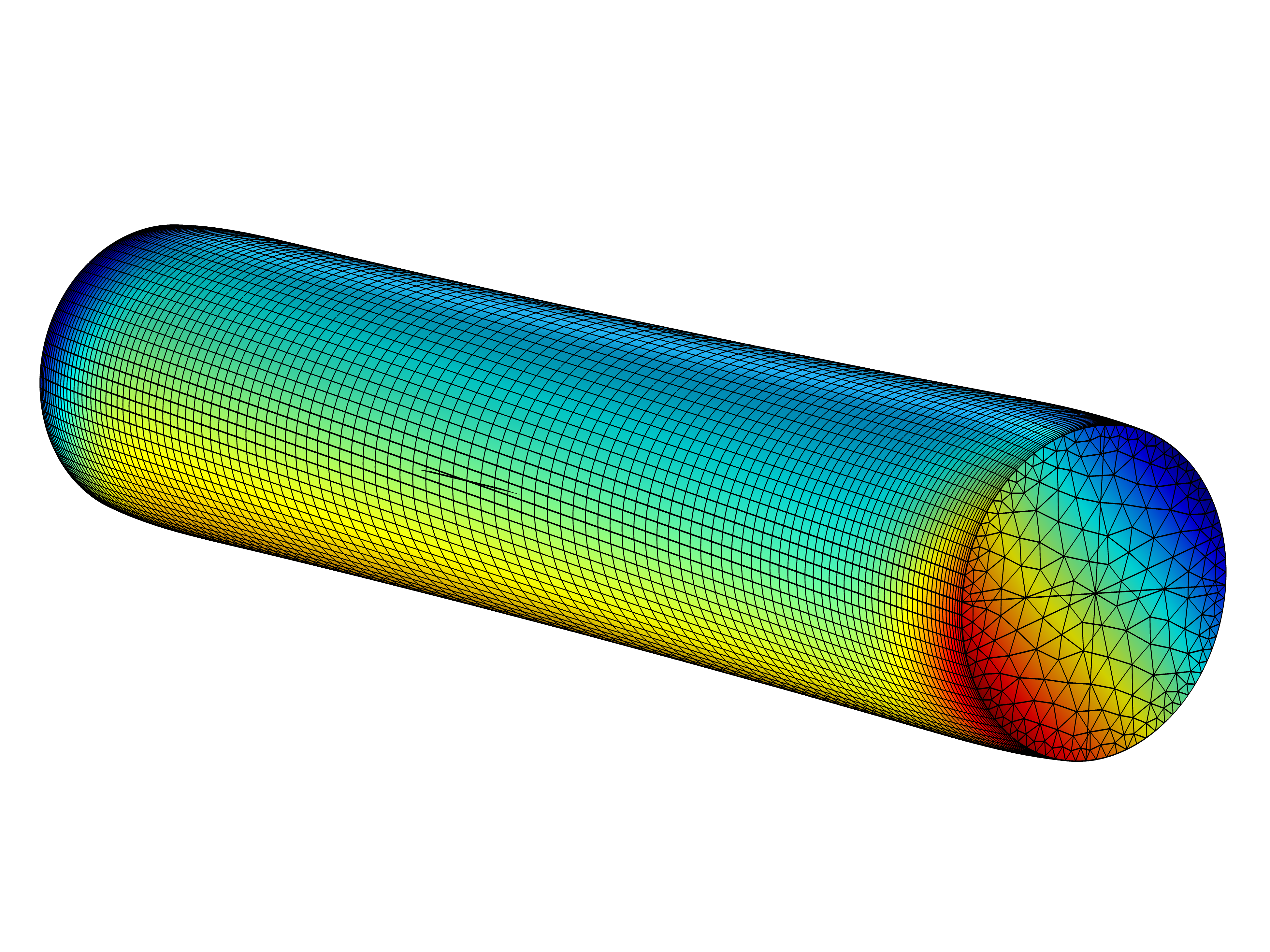}
		\caption{}
	\end{subfigure}
	\caption{(a) boundary conditions scheme for a cylindrical rod of length $L=10$ and radius $R=1$; (b) deformed rod from the finite-element simulation on which it is mapped how the component of the gradient of the displacement $u_{1,3}$ changes.}
	\label{fig:FE_RM_1}
\end{figure}

\begin{figure}[H]
	\begin{subfigure}{0.45\textwidth}
		\centering
		\includegraphics[width=\textwidth]{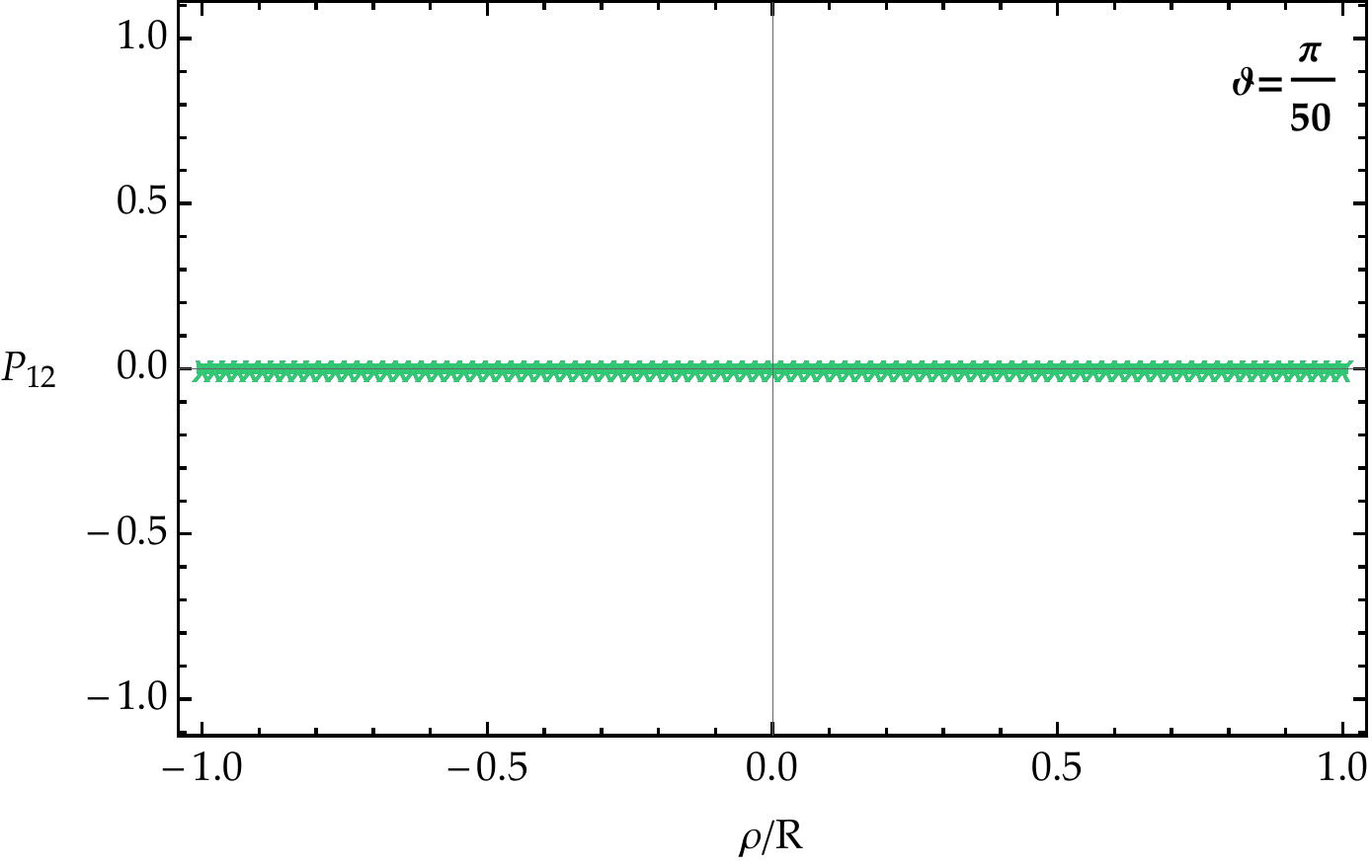}
		\caption{}
	\end{subfigure}
	\hfill
	\begin{subfigure}{0.45\textwidth}
		\centering
		\includegraphics[width=\textwidth]{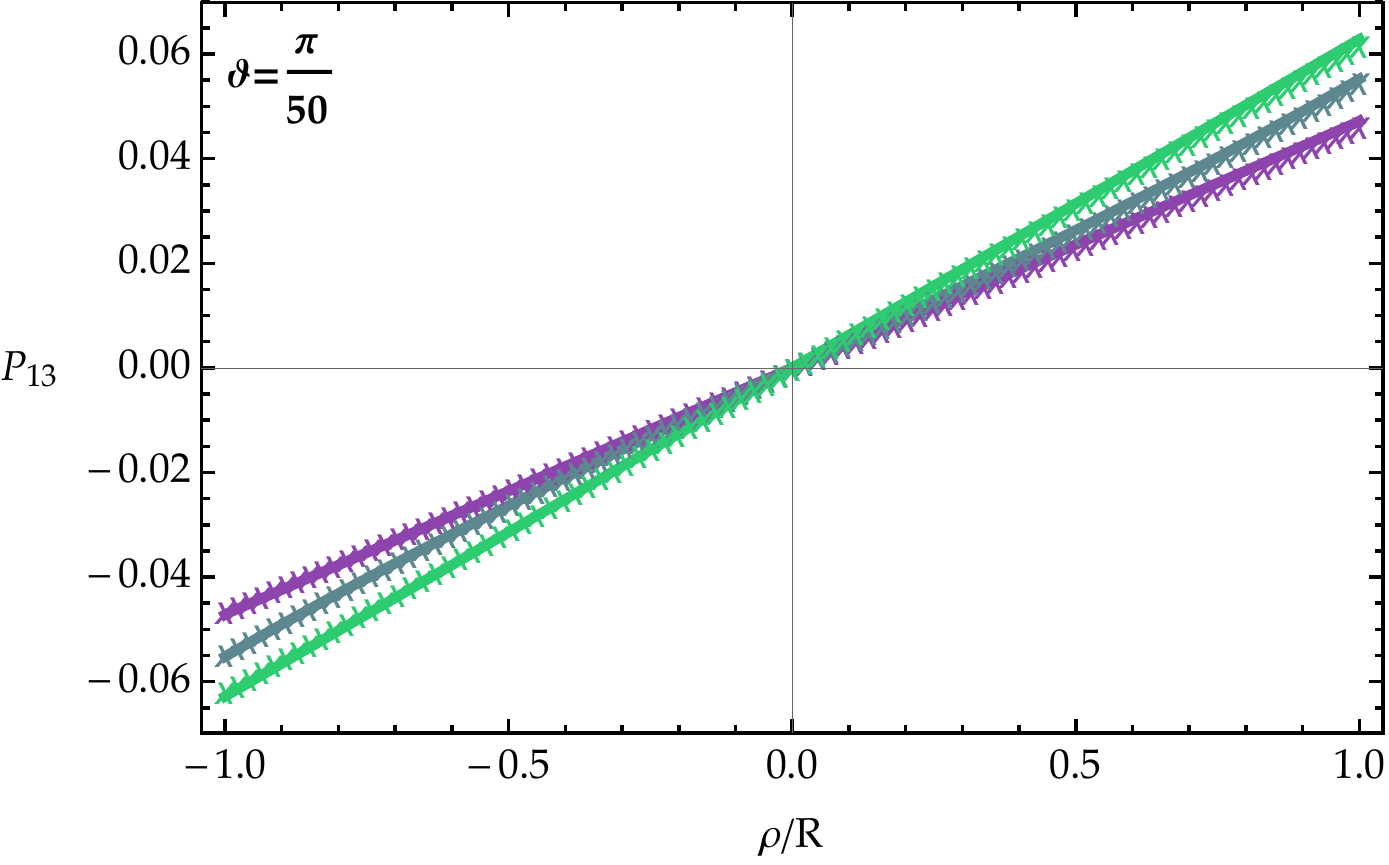}
		\caption{}
	\end{subfigure}
	\hfill
	\begin{subfigure}{0.45\textwidth}
		\centering
		\includegraphics[width=\textwidth]{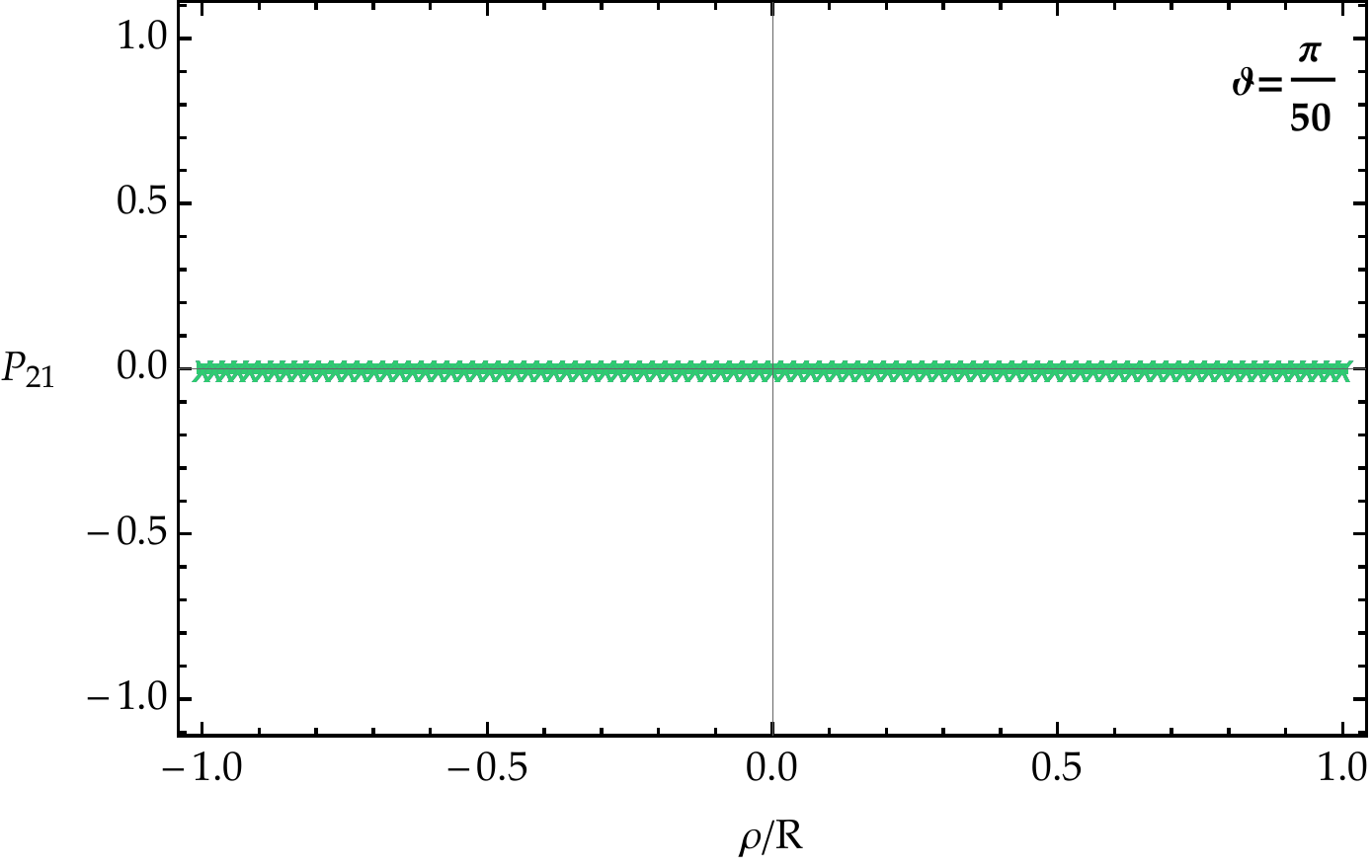}
		\caption{}
	\end{subfigure}
	\hfill
	\begin{subfigure}{0.45\textwidth}
		\centering
		\includegraphics[width=\textwidth]{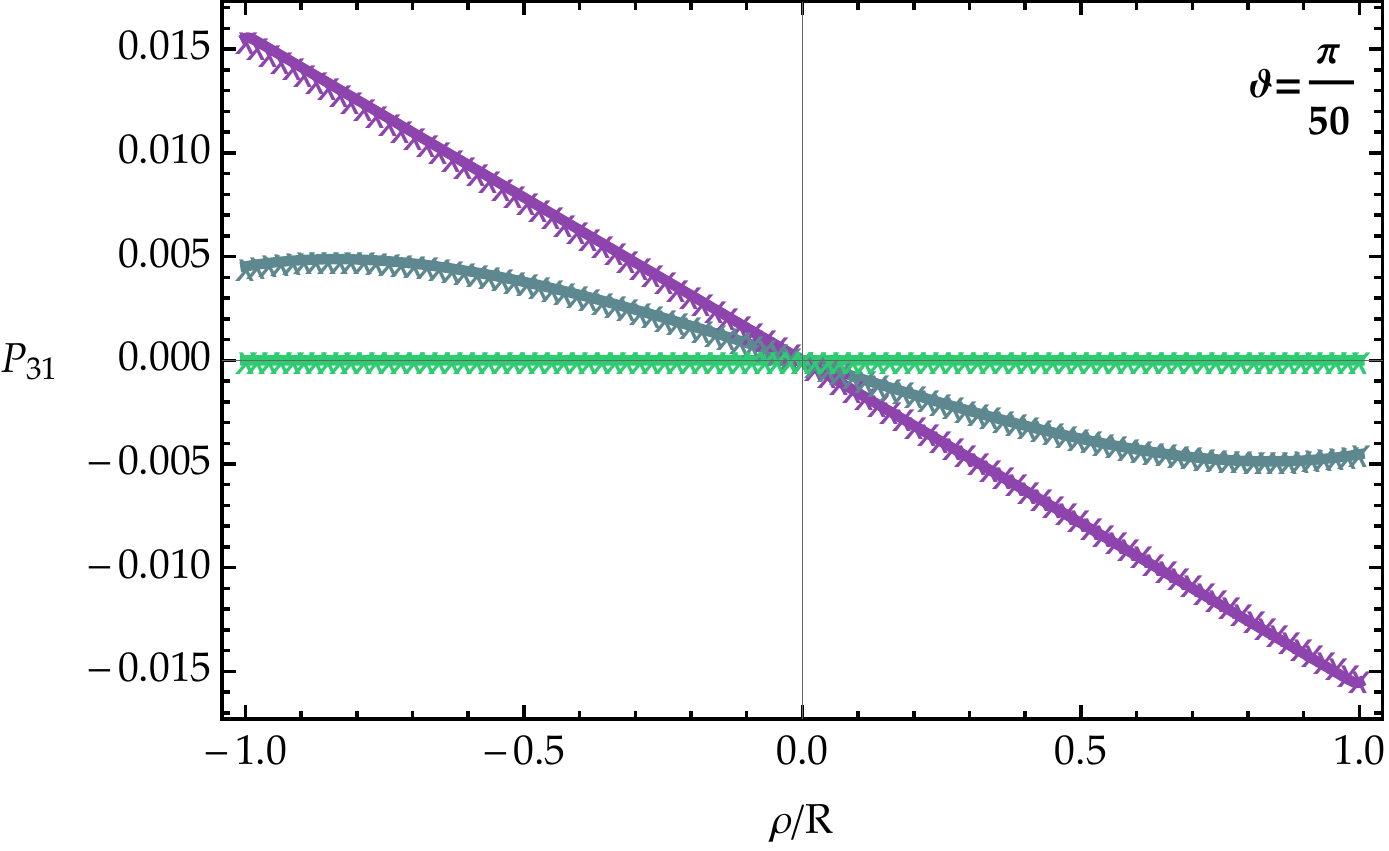}
		\caption{}
	\end{subfigure}
	\caption{Plots of the components (a) $P_{12}$, (b) $P_{13}$, (c) $P_{21}$, and (d) $P_{31}$ of the micro-distortion tensor $\boldsymbol{P}$ at the cross-section $z=0$. The purple line corresponds to $L_c=0$, the gray one to $L_c=1$, and the green one to $L_c=\infty$. The values of the other parameters used are $\mu = 1$, $\mu _c= 1$, $\mu _e= 1$, $\mu _{\mbox{\tiny micro}}= 1$, $a _1= 1$, $a _3= 1$, $R= 1$, $\boldsymbol{\vartheta}= \pi/50$, $L=10$.}
	\label{fig:FE_RM_2}
\end{figure}

In Fig.~\ref{fig:FE_RM_2} it is possible to see how the non identically zero components of the micro-distortion $\boldsymbol{P}$ vary across the diameter aligned to the $x_1$-axis ($\varphi=\pi/2$) of the cross-section placed in the middle of the cylindrical rod ($z=0$).
We chose the middle section in order to study the solution far away enough from the disturbance region on which the boundary conditions have been applied.
The values of the components of $\boldsymbol{P}$ of Fig. \ref{fig:FE_RM_2} (twist rate $\boldsymbol{\vartheta}=\pi/50$) are perfectly in agreement with the analytical solution, confirming the validity of the small-deformation solution obtained in  Section \ref{sec:Relaxed_micro}.

\begin{figure}[H]
	\begin{subfigure}{0.45\textwidth}
		\centering
	\includegraphics[width=\textwidth]{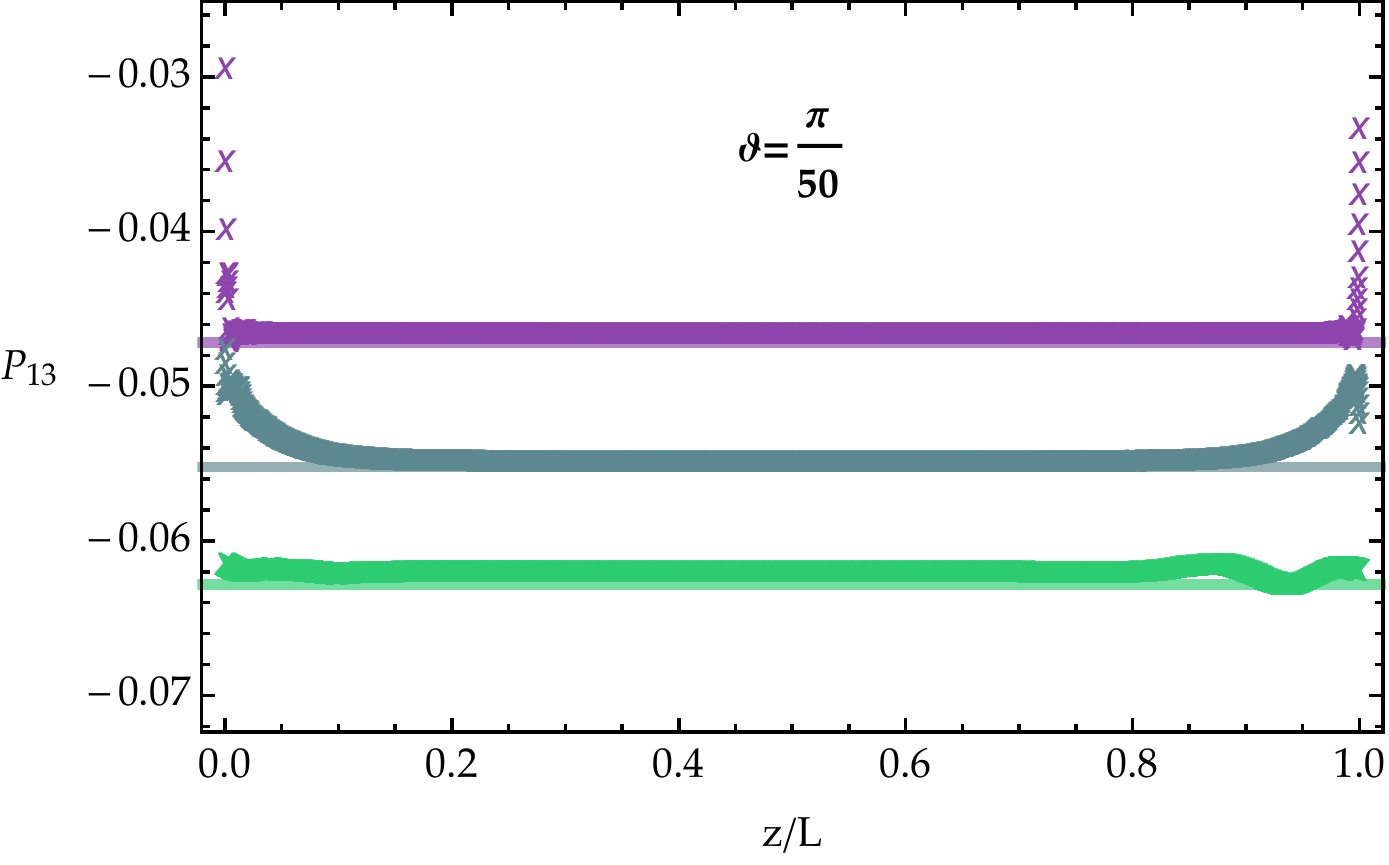}
		\caption{}
	\end{subfigure}
	\hfill
	\begin{subfigure}{0.45\textwidth}
		\centering
	\includegraphics[width=\textwidth]{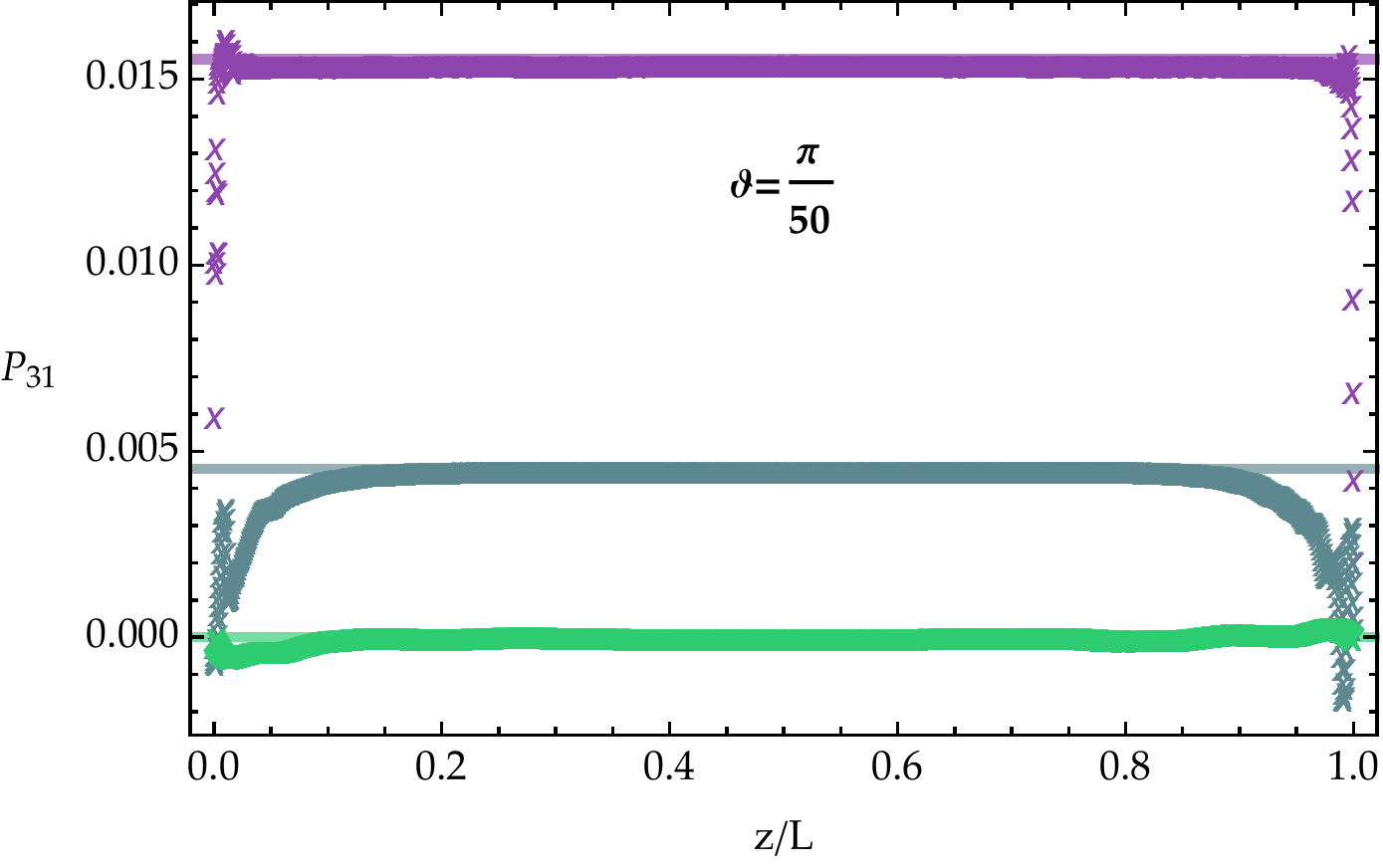}
		\caption{}
	\end{subfigure}
	\caption{
	(a) Plot of how the component $P_{13}$; 
	(b) and $P_{31}$ vary along a line on the external boundary ($\varphi=\pi/2$ and $r=R$): the solid lines are the analytical solution while the marker are the numerical values obtained thanks to a finite-element analysis. The purple line has been obtained for $L_c=0$, the gray one for $L_c=1$, and the green one for $L_c=\infty$. The values of the other parameters used are $\mu = 1$, $\mu _c= 1$, $\mu _e= 1$, $\mu _{\mbox{\tiny micro}}= 1$, $a _1= 1$, $a _3= 1$, $R= 1$, $\varphi=\pi/2$.
	As it can be seen, the solution does not converge stably and not perfectly symmetrically (the mesh is not symmetric) to the analytical one, but nevertheless it converges rapidly.
	}
	\label{fig:FE_RM_3}
\end{figure}

Furthermore, in Fig. \ref{fig:FE_RM_3} we show how the solution obtained while applying consistent boundary conditions converges to the one obtained analytically in a distance from the boundary which is more or less between one radius and one diameter.
This is the pinnacle expression of the Saint Venant principle: we have applied not only a finite-rotation instead of a linearized one, but we have also used consistent boundary conditions for $\boldsymbol{P}$ which we know are different from the correct values that the tangential part of $\boldsymbol{P}$ should have, and we obtained nevertheless the analytical linearized solution after a rather small boundary layer.

\begin{figure}[H]
	\begin{subfigure}{0.45\textwidth}
		\centering
	\includegraphics[width=\textwidth]{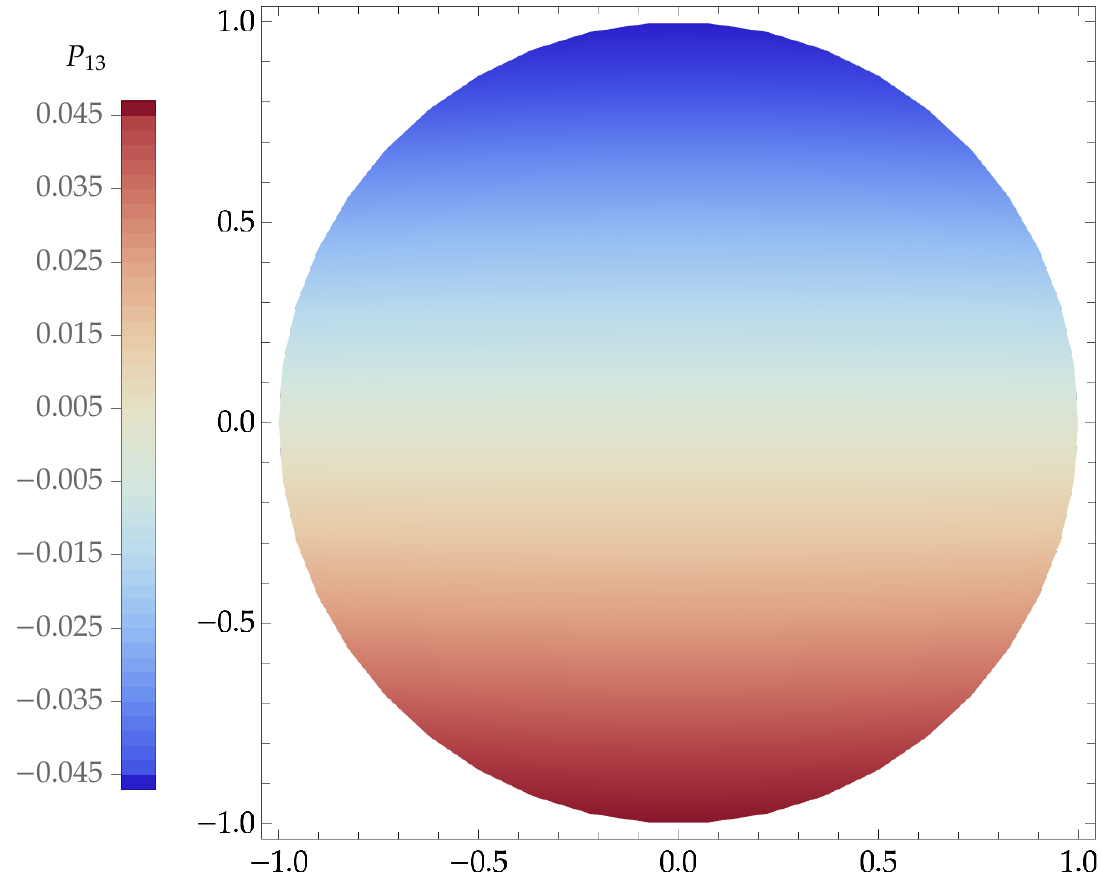}
		\caption{}
	\end{subfigure}
	\hfill
	\begin{subfigure}{0.45\textwidth}
		\centering
	\includegraphics[width=\textwidth]{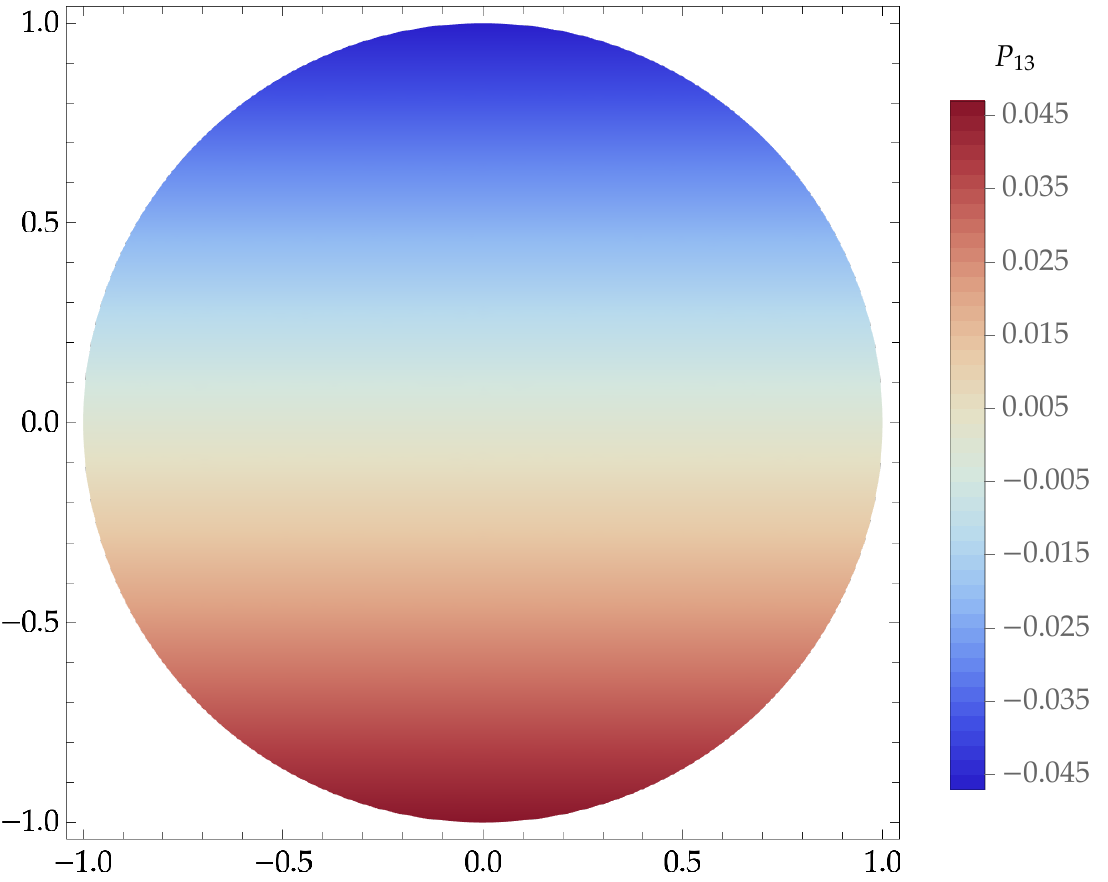}
		\caption{}
	\end{subfigure}
	\caption{
	Plots of the component $P_{13}$ across the section placed at $z=0$ obtained (a) analytically and (b) with the finite-element analysis. The two results are in perfect agreement.
	The values of the other parameters used are $\mu = 1$, $\mu _c= 1$, $\mu _e= 1$, $\mu _{\mbox{\tiny micro}}= 1$, $a _1= 1$, $a _3= 1$, $R= 1$, $\varphi=\pi/2$, $L_c=1$.
	}
	\label{fig:FE_RM_4}
\end{figure}

We describe in particular the component $P_{31}$ (other than the component $P_{13}$) since, due to  the consistent boundary conditions, it is forced to start from zero at the lateral boundary.
In Fig. \ref{fig:FE_RM_3} we plot this component which are evaluated for the length of the rod on the external surface ($\varphi=\pi/2$ and $r=R$).
In Fig. \ref{fig:FE_RM_4} is reported how the component $P_{13}$ vary on the cross section centered in the origin of the reference system ($z=0$).

\textit{The implications of this results are of great value in the context of the identification of the elastic material parameters: it is clear how to apply consistent boundary conditions on a real sample in a laboratory (Dirichlet hard), and thanks to this results, we now know that our analytical solution is taking place far away enough from the boundary layer.}

\section{Torsional problem for the isotropic micro-stretch model in dislocation format}

In the micro-stretch model in dislocation format \cite{neff2014unifying,scalia2000extension,de1997torsion,neff2009mean,kirchner2007mechanics}, contrary to the relaxed micromorphic model, the micro-distortion tensor is devoid  from the deviatoric component $\mbox{dev} \, \mbox{sym} \, \boldsymbol{P} = 0 \Leftrightarrow \boldsymbol{P} = \boldsymbol{A} + \omega \boldsymbol{\mathbbm{1}}$, $\boldsymbol{A} \in \mathfrak{so}(3)$, $\omega \in \mathbb{R}$.
The expression of the strain energy for this model in dislocation format can be written as \cite{neff2014unifying}:
\begin{align}
	W \left(\boldsymbol{\mbox{D}u}, \boldsymbol{A},\omega,\mbox{Curl}\,\left(\boldsymbol{A} - \omega \boldsymbol{\mathbbm{1}}\right)\right) 
	\hspace{-2.5cm}
	&
	\notag
	\\*
	=
	&
	\, \mu_{\tiny \mbox{macro}} \left\lVert \mbox{dev} \, \mbox{sym} \, \boldsymbol{\mbox{D}u} \right\rVert^{2}
	+ \frac{\kappa_{e}}{2} \mbox{tr}^2 \left(\boldsymbol{\mbox{D}u} - \omega \boldsymbol{\mathbbm{1}} \right) 
	+ \mu_{c} \left\lVert \mbox{skew} \left(\boldsymbol{\mbox{D}u} - \boldsymbol{A} \right) \right\rVert^{2}
	+ \frac{9}{2} \, \kappa_{\tiny \mbox{micro}} \, \omega^2
	\label{eq:energy_micro_stretch}
	\\*
	&
	+ \frac{\mu \,L_c^2}{2} \,
	\left(
	a_1 \, \left\lVert \mbox{dev sym} \, \mbox{Curl} \, \boldsymbol{A} \right\rVert^2
	+ a_2 \, \left\lVert \mbox{skew} \,  \mbox{Curl} \, \left(\boldsymbol{A} + \omega \boldsymbol{\mathbbm{1}}\right) \right\rVert^2
	+ \frac{a_3}{3} \, \mbox{tr}^2 \left(\mbox{Curl} \, \boldsymbol{A} \right)
	\right) \, ,
	\notag
\end{align}
since $\mbox{Curl} \left(\omega \boldsymbol{\mathbbm{1}}\right) \in \mathfrak{so}(3)$.
The equilibrium equations, in the absence of body forces,   are then
\begin{align}
	\mbox{Div}\overbrace{\left[
		2\mu_{\tiny \mbox{macro}}\,\mbox{dev}\,\mbox{sym} \, \boldsymbol{\mbox{D}u}
		+ \kappa_{e} \mbox{tr} \left(\boldsymbol{\mbox{D}u} - \omega \boldsymbol{\mathbbm{1}}\right) \boldsymbol{\mathbbm{1}}
		+ 2\mu_{c}\,\mbox{skew} \left(\boldsymbol{\mbox{D}u} - \boldsymbol{A}\right) \right]}^{\mathlarger{\widetilde{\sigma}}\coloneqq}
	&= \boldsymbol{0} \, ,
	\notag
	\\
	2\mu_{c}\,\mbox{skew} \left(\boldsymbol{\mbox{D}u} - \boldsymbol{A}\right)
	\hspace{10.1cm}
	\notag
	\\*
	-\mu \, L_c^2 \, \mbox{skew} \, \mbox{Curl}\,\left(
	a_1 \, \mbox{dev} \, \mbox{sym} \, \mbox{Curl} \, \boldsymbol{A} \, 
	+ a_2 \, \mbox{skew} \, \mbox{Curl} \, \left(\boldsymbol{A} + \, \omega \boldsymbol{\mathbbm{1}}\right) \,
	+ \frac{a_3}{3} \, \mbox{tr} \left(\mbox{Curl} \, \boldsymbol{A} \right)\boldsymbol{\mathbbm{1}} \, 
	\right) &= \boldsymbol{0} \, ,
	\label{eq:equi_micro_stretch}
	\\*
	\mbox{tr}
	\bigg[
	2\mu_{\tiny \mbox{macro}}\,\mbox{dev}\,\mbox{sym} \, \boldsymbol{\mbox{D}u}
	\hspace{10.25cm}
	\notag
	\\*
	+ \kappa_{e} \mbox{tr} \left(\boldsymbol{\mbox{D}u} - \omega \boldsymbol{\mathbbm{1}}\right) \boldsymbol{\mathbbm{1}}
	- \kappa_{\mbox{\tiny micro}} \mbox{tr} \left( \omega \boldsymbol{\mathbbm{1}}\right) \boldsymbol{\mathbbm{1}}
	-\mu \, L_c^2 \,  a_2 \, \mbox{Curl}\,
	\mbox{skew} \, \mbox{Curl} \, \left(\omega \boldsymbol{\mathbbm{1}} + \boldsymbol{A}\right) 
	\bigg]
	&= \boldsymbol{0} \,.
	\notag
\end{align}
The boundary conditions at the free surface are 
\begin{align}
	\boldsymbol{\widetilde{t}}(r = R) &= 
	\boldsymbol{\widetilde{\sigma}}(r) \cdot \boldsymbol{e}_{r} = 
	\boldsymbol{0}_{\mathbb{R}^{3}} \, ,
	\notag
	\\
	\boldsymbol{\eta}(r = R) &= 
	\mbox{skew}\left(\boldsymbol{m} (r) \cdot \boldsymbol{\epsilon} \cdot \boldsymbol{e}_{r} \right) = 
	\mbox{skew}\left(\boldsymbol{m} (r) \times \boldsymbol{e}_{r} \right) = 
	\boldsymbol{0}_{\mathbb{R}^{3\times 3}} \, ,
	\label{eq:BC_micro_stretch}
	\\
	\gamma(r = R) &= 
	\frac{1}{3} \mbox{tr}\left(\boldsymbol{m} (r) \cdot \boldsymbol{\epsilon} \cdot \boldsymbol{e}_{r} \right) = 
	\frac{1}{3} \mbox{tr}\left(\boldsymbol{m} (r) \times \boldsymbol{e}_{r} \right) = 
	0 \, .
	\notag
\end{align}
According with the reference system shown in Fig.~\ref{fig:intro_1}, the ansatz for the displacement and micro-distortion fields is
\begin{align}
	\boldsymbol{u}(x_1,x_2,x_3) &=
	\boldsymbol{u}(r,\varphi,z) =
	\boldsymbol{\vartheta}
	\left(
	\begin{array}{c}
		-x_2(r,\varphi) \, x_3(z) \\
		x_1(r,\varphi) \, x_3(z) \\
		0 
	\end{array}
	\right)
	\, ,
	\qquad
	\omega = 0 \, ,
	\notag
	\\*
	\boldsymbol{A}(x_1,x_2,x_3) &=
	\boldsymbol{A}(r,\varphi,z) = \frac{\boldsymbol{\vartheta}}{2}
	\left(
	\begin{array}{ccc}
		0   & -2x_3(z) & - g_{p}(r) \, x_2(r,\varphi) \\
		2x_3(z) &    0 &   g_{p}(r) \, x_1(r,\varphi) \\
		g_{p}(r) \, x_2(r,\varphi) & - g_{p}(r) \, x_1(r,\varphi) & 0 \\
	\end{array}
	\right)
	\, .
	\label{eq:ansatz_micro_stretch}
\end{align}
Since the ansatz requires $\omega=0$, \textit{the micro-stretch model coincides with the Cosserat model} which will be presented in the next section.

\section{Torsional problem for the isotropic Cosserat continuum}
\label{sec:Cos}

The strain energy for the isotropic Cosserat continuum in dislocation tensor format (curvature expressed in term of ${\rm Curl} A$) can be written as \cite{rizzi2021shear,rizzi2021bending,nano11020453,rueger2018strong,cosserat1909theorie,neff2006cosserat,neff2009new,jeong2009numerical}
\begin{align}
	W \left(\boldsymbol{\mbox{D}u}, \boldsymbol{A},\mbox{Curl}\,\boldsymbol{A}\right) = &
	\, \mu_{\mbox{\tiny macro}} \left\lVert \mbox{sym} \, \boldsymbol{\mbox{D}u} \right\rVert^{2}
	+ \frac{\lambda_{\mbox{\tiny macro}}}{2} \mbox{tr}^2 \left(\boldsymbol{\mbox{D}u} \right) 
	+ \mu_{c} \left\lVert \mbox{skew} \left(\boldsymbol{\mbox{D}u} - \boldsymbol{A} \right) \right\rVert^{2}
	\label{eq:energy_Cos}
	\\
	&
	+ \frac{\mu \, L_c^2}{2}
	\left(
	a_1 \, \left \lVert \mbox{dev} \, \mbox{sym} \, \mbox{Curl} \, \boldsymbol{A}\right \rVert^2 \, 
	+ a_2 \, \left \lVert \mbox{skew} \, \mbox{Curl} \, \boldsymbol{A}\right \rVert^2 \, 
	+ \frac{a_3}{3} \, \mbox{tr}^2 \left(\mbox{Curl} \, \boldsymbol{A} \right)
	\right)  \, ,
	\notag
\end{align}
where $A \in \mathfrak{so}(3)$.
It is underlined that for the ansatz  (\ref{eq:ansatz_Cos}), which will be presented later in this section, it holds that $\mbox{skew} \left(\mbox{Curl} \, \boldsymbol{A} \right)=0$ (see calculation  (\ref{eq:grad_RM})$_{2}$).
The equilibrium equations, in the absence of body forces,   are therefore the following

\begin{align}
	\mbox{Div}\overbrace{\left[2\mu_{\mbox{\tiny macro}}\,\mbox{sym} \, \boldsymbol{\mbox{D}u} + \lambda_{\mbox{\tiny macro}} \mbox{tr} \left(\boldsymbol{\mbox{D}u} \right) \boldsymbol{\mathbbm{1}}
		+ 2\mu_{c}\,\mbox{skew} \left(\boldsymbol{\mbox{D}u} - \boldsymbol{A}\right) \right]}^{\mathlarger{\widetilde{\sigma}}\coloneqq}
	&= \boldsymbol{0} \, ,
	\notag
	\\
	2\mu_{c}\,\mbox{skew} \left(\boldsymbol{\mbox{D}u} - \boldsymbol{A}\right)
	-\mu \, L_c^2 \, \mbox{skew} \, \mbox{Curl}\,
	\left(
	a_1 \, \mbox{dev} \, \mbox{sym} \, \mbox{Curl} \, \boldsymbol{A} \, 
	+ \frac{a_3}{3} \, \mbox{tr} \left(\mbox{Curl} \, \boldsymbol{A} \right)\boldsymbol{\mathbbm{1}} \, 
	\right)
	&= \boldsymbol{0} \,.
	\label{eq:equi_Cos}
\end{align}
The boundary conditions at the free surface are 
\begin{align}
	\boldsymbol{\widetilde{t}}(r = R) &= 
	\boldsymbol{\widetilde{\sigma}}(r) \cdot \boldsymbol{e}_{r} = 
	\boldsymbol{0}_{\mathbb{R}^{3}} \, ,
	\label{eq:BC_Cos_gen}
	\\*
	\boldsymbol{\eta}(r = R) &= 
	\mbox{skew}\left(\boldsymbol{m} (r) \cdot \boldsymbol{\epsilon} \cdot \boldsymbol{e}_{r}\right) = 
	\mbox{skew}\left(\boldsymbol{m} (r) \times \boldsymbol{e}_{r}\right) = 
	\boldsymbol{0}_{\mathbb{R}^{3\times 3}} \, ,
	\notag
\end{align}
where the second order moment stress tensor is now given by
\begin{equation}
\boldsymbol{m} =
\mu \, L_c^2
\left(
a_1 \, \mbox{dev} \, \mbox{sym} \, \mbox{Curl} \, \boldsymbol{A} \, 
+ \frac{a_3}{3} \, \mbox{tr} \left(\mbox{Curl} \, \boldsymbol{A} \right)\boldsymbol{\mathbbm{1}} \, 
\right) \, ,
\end{equation}
the expression of $\boldsymbol{\widetilde{\sigma}}$ is in  (\ref{eq:equi_Cos}), $\boldsymbol{e}_{r}$ is the radial unit vector, and $\boldsymbol{\epsilon}$ is the Levi-Civita tensor.

According to the reference system shown in Fig.~\ref{fig:intro_1}, the ansatz for the displacement field and the micro-rotation is
\begin{align}
	\boldsymbol{u}(x_1,x_2,x_3) &=
	\boldsymbol{u}(r,\varphi,z) = \boldsymbol{\vartheta}
	\left(
	\begin{array}{c}
		-x_2(r,\varphi) \, x_3(z) \\
		x_1(r,\varphi) \, x_3(z) \\
		0 
	\end{array}
	\right)
	\, ,
	\notag
	\\*
	\boldsymbol{A}(x_1,x_2,x_3) &= 
	\boldsymbol{A}(r,\varphi,z) = \frac{\boldsymbol{\vartheta}}{2}
	\left(
	\begin{array}{ccc}
		0   & -2x_3(z) & - g_{p}(r) \, x_2(r,\varphi) \\
		2x_3(z) &    0 &   g_{p}(r) \, x_1(r,\varphi) \\
		g_{p}(r) \, x_2(r,\varphi) & - g_{p}(r) \, x_1(r,\varphi) & 0 \\
	\end{array}
	\right)
	\, ,
	\label{eq:ansatz_Cos}
\end{align}
where, in relation to the ansatz  (\ref{eq:ansatz_RM}), we define $g_{p}(r ) \coloneqq g_{1}(r ) + g_{2}(r )$, so that there is only one unknown function to be determined.
Substituting the ansatz  (\ref{eq:ansatz_Cos}) in  (\ref{eq:equi_Cos}) the 6 equilibrium equations are equivalent to
\begin{align}
	\frac{1}{6} \, \boldsymbol{\vartheta} \, \sin \varphi \left(6 r \, \mu _c (g_{p}(r )-1) - \mu \,  L_c^2 \, (a_{1}+2 a_{3}) \left(3 g'_{p}(r ) + r \, g''_{p}(r )\right)\right) = 0 \, ,
	\label{eq:equi_equa_Cos}
	\\*
	\frac{1}{6} \, \boldsymbol{\vartheta} \, \cos \varphi \left(6 r \, \mu _c (g_{p}(r )-1) - \mu \,  L_c^2 \, (a_{1}+2 a_{3}) \left(3 g'_{p}(r ) + r \, g''_{p}(r )\right)\right) = 0 \, .
	\notag
\end{align}
It is important to underline that  (\ref{eq:equi_Cos})$_1$ is identically satisfied, and that between the two equilibrium equations  (\ref{eq:equi_equa_Cos}) there is only one independent equation since  (\ref{eq:equi_equa_Cos})$_1=\tan \varphi$  (\ref{eq:equi_equa_Cos})$_2$.
The solution of  (\ref{eq:equi_equa_Cos}) is 
\begin{align}
	g_{p}(r ) = \,
	 1
	-\frac{i \, A_{1} I_1\left( \frac{     r \, f_1}{L_c}\right)}{r }
	+\frac{     A_{2} Y_1\left(-\frac{i \, r \, f_1}{L_c}\right)}{r }
	\, ,
	\qquad\qquad\qquad
	f_{1} \coloneqq
	\sqrt{\frac{6 \mu _c}{(a_{1} + 2 a_{3}) \, \mu}} \, ,
	\label{eq:sol_fun_Cos}
\end{align}
where $I_{n}\left(\cdot\right)$ is the \textit{modified Bessel function of the first kind}, $Y_{n}\left(\cdot\right)$ is the \textit{Bessel function of the second kind} (see appendix \ref{app:bessel} for the formal definitions), and
$A_1$, $A_2$ are integration constants.

\begin{figure}[H]
	\begin{subfigure}{0.5\textwidth}
		\centering
		\includegraphics[width=\textwidth]{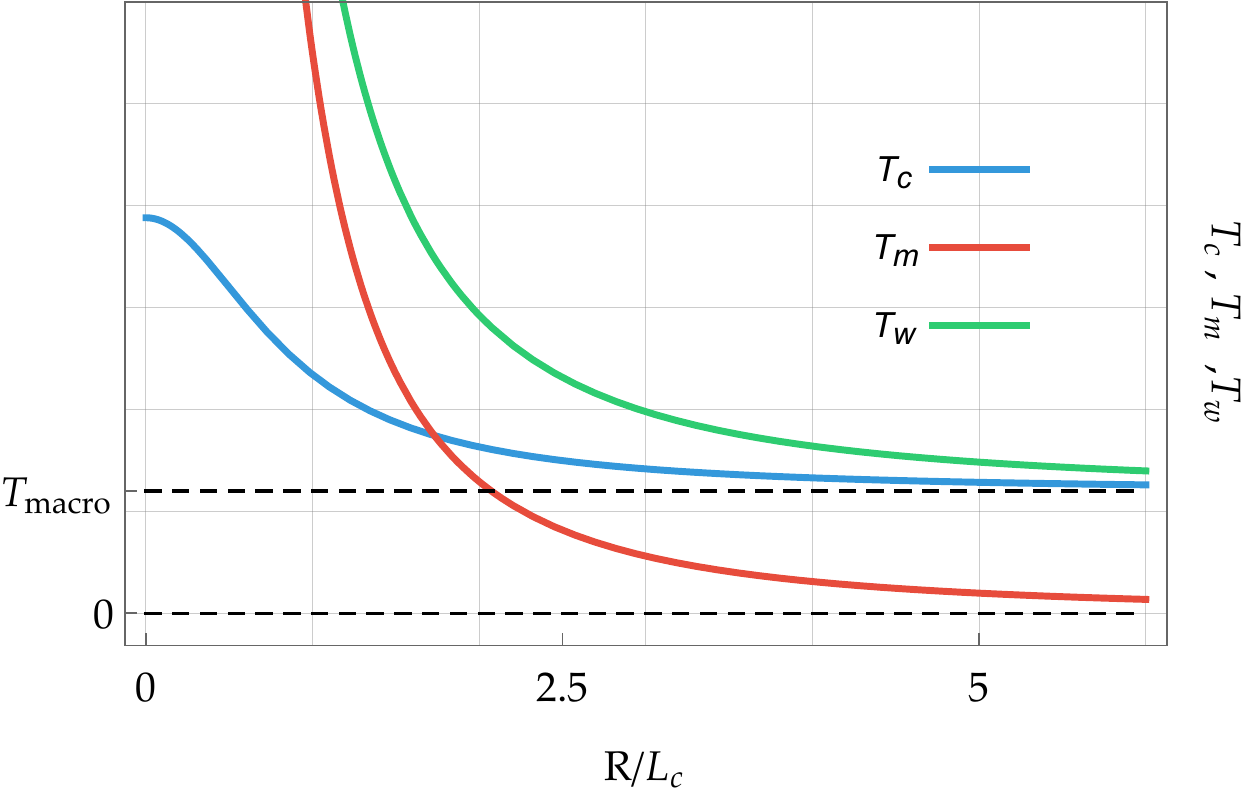}
		\caption{}
	\end{subfigure}
	\hfill
	\begin{subfigure}{0.5\textwidth}
		\centering
		\includegraphics[width=\textwidth]{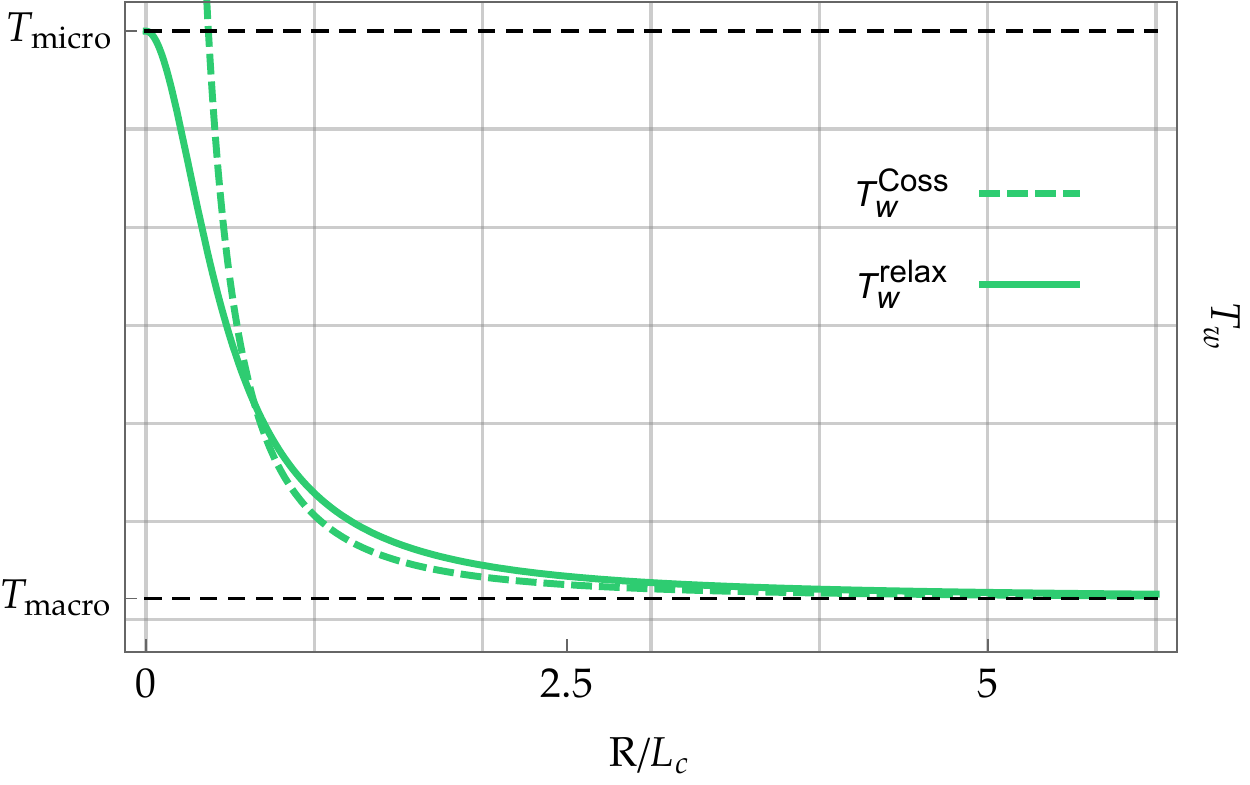}
		\caption{}
	\end{subfigure}
	\caption{
	(a)(\textbf{Cosserat model with full curvature}) Torsional stiffness for the classical torque $T_{\mbox{c}}$, the higher-order torque $T_{\mbox{m}}$, and the torque energy $T_{\mbox{w}}$ while varying $L_c$. The torsional stiffness is unbounded as $L_c \to \infty$ ($R\to 0$).
	(b)(\textbf{Cosserat model with full curvature vs relaxed micromorphic model}) Torsional stiffness for the torque energy ($T_{\mbox{w}}$) while varying $L_c$. Observe that the torsional stiffness remains bounded for the relaxed micromorphic model while it blows up for the Cosserat model as $L_c \to \infty$ ($R\to 0$).
	For best comparison, the characteristic length scale of the Cosserat model has been chosen $L_c^{\text{Coss}}\coloneqq\frac{L_c^{\text{relax}}}{\sqrt{2}}$. The values of the parameters used are: $\mu = 1$, $\mu _c= 1/2$, $\mu _{\mbox{\tiny macro}}= 1/14$, $\mu _{\mbox{\tiny micro}}= 2$ (just for the relaxed micromorphic model), $a_1= 1/5$, $a_3= 1/7$, $R= 1$.}
	\label{fig:all_plot_Cos_3}
\end{figure}

The value of $A_1$ is determined from to the boundary conditions  (\ref{eq:BC_Cos_gen}), where, due to the divergent behaviour of the Bessel function of the second kind at $r=0$, we have to set $A_2=0$ in order to have a continuous solution.
The fulfilment of the boundary conditions  (\ref{eq:BC_Cos_gen}) allows us to find the expressions of the integration constants
\begin{align}
	A_1 = \, 
	-\frac{
	i \, R \, L_c
	}{
	f_1 \, R \, z_1 \left(I_0\left(\frac{R \, f_1}{L_c}\right)
	+ I_2 \left(\frac{R \, f_1}{L_c}\right)\right)
	+ z_2 \, L_c \, I_1\left(\frac{R \, f_1}{L_c}\right)
	}
	\, ,
	\,\,\,\,\,\,\,\,\,\,
	z_1 \coloneqq \,  \frac{a_1 + 2a_3}{3a_1} \, ,
	\,\,\,\,\,\,\,\,\,\,
	z_2 \coloneqq \, \frac{4 a_3 - a_1}{3a_1}
	\, .
	\label{eq:BC_Cos_3}
\end{align}
The classical torque, the higher-order torque, and energy (per unit length d$z$) expressions are
\begin{align}
	M_{\mbox{c}} (\boldsymbol{\vartheta}) 
	\coloneqq&
	\int_{0}^{2\pi}
	\int_{0}^{R}
	\Big[
	\langle
	\widetilde{\boldsymbol{\sigma}} \, \boldsymbol{e}_{z} , \boldsymbol{e}_{\varphi}
	\rangle
	r
	\Big] r
	\, \mbox{d}r \, \mbox{d}\varphi
	\notag
	\\*
	=&
	\left[
	\mu_{\mbox{\tiny macro}}
	+ \frac{
	4 \mu _c \, I_2\left(\frac{R \, f_1}{L_c}\right)
	\frac{L_c^2}{R^2}
	}{
	f_1 \left(2 \, f_1 \, z_1 \, I_0\left(\frac{R \, f_1}{L_c}\right)
	+\left(z_2 - 2 z_1\right) I_1\left(\frac{R \, f_1}{L_c}\right) \frac{L_c}{R}\right)
	}
	\right]
	I_{p} \, 
	\boldsymbol{\vartheta}
	= T_{\mbox{c}} \, \boldsymbol{\vartheta} \, ,
	\notag
	\\[3mm]
	M_{\mbox{m}}(\boldsymbol{\vartheta}) 
	\coloneqq&
	\int_{0}^{2\pi}
	\int_{0}^{R}
	\Big[
	\langle
	\mbox{skew}(\boldsymbol{m} \times \boldsymbol{e}_{z})
	\boldsymbol{e}_{\varphi} ,
	\boldsymbol{e}_{r}
	\rangle
	-
	\langle
	\mbox{skew}(\boldsymbol{m} \times \boldsymbol{e}_{z})
	\boldsymbol{e}_{r} ,
	\boldsymbol{e}_{\varphi}
	\rangle
	\Big]
	\, r
	\, \mbox{d}r \, \mbox{d}\varphi
	\notag
	\\*
	=&
	\left[
	\frac{
	2 \mu \left(
	3 a_{1} \, f_1 \, z_1 \, I_0\left(\frac{R \, f_1}{L_c}\right) \,
	\frac{L_c^2}{R^2}
	-2 (a_{1} - a_{3}) \, I_1\left(\frac{R \, f_1}{L_c}\right)
	\frac{L_c^3}{R^3}
	\right)
	}{
	6 f_1 \, z_1 \, I_0\left(\frac{R \, f_1}{L_c}\right)
	-3 I_1\left(\frac{R \, f_1}{L_c}\right)
	\frac{L_c}{R}
	}
	\right]
	I_{p} \, 
	\boldsymbol{\vartheta}
	= T_{\mbox{m}} \, \boldsymbol{\vartheta} \, ,
	\label{eq:torque_stiffness_Cos}
	\\[3mm]
	W_{\mbox{tot}} (\boldsymbol{\vartheta}) 
	\coloneqq&
	\int_{0}^{2\pi}
	\int_{0}^{R}
	W \left(\boldsymbol{\mbox{D}u}, \boldsymbol{A}, \mbox{Curl}\boldsymbol{A}\right) \, \, r
	\, \mbox{d}r \, \mbox{d}\varphi
	\notag
	\\*
	=&
	\frac{1}{2}
	\left[
	\mu_{\mbox{\tiny macro}}
	+ \frac{
		4 \mu _c \, I_2\left(\frac{R \, f_1}{L_c}\right)
		\frac{L_c^2}{R^2}
	}{
		f_1 \left(2 \, f_1 \, z_1 \, I_0\left(\frac{R \, f_1}{L_c}\right)
		+\left(z_2 - 2 z_1\right) I_1\left(\frac{R \, f_1}{L_c}\right) \frac{L_c}{R}\right)
	}
	\notag
	\right.
	\\*
	& \hspace{3cm}
	\left.
	+
	\frac{
		2 \mu \left(
		3 a_{1} \, f_1 \, z_1 \, I_0\left(\frac{R \, f_1}{L_c}\right) \,
		\frac{L_c^2}{R^2}
		-2 (a_{1} - a_{3}) \, I_1\left(\frac{R \, f_1}{L_c}\right)
		\frac{L_c^3}{R^3}
		\right)
	}{
		6 f_1 \, z_1 \, I_0\left(\frac{R \, f_1}{L_c}\right)
		-3 I_1\left(\frac{R \, f_1}{L_c}\right)
		\frac{L_c}{R}
	}
	\right]
	I_{p} \, 
	\boldsymbol{\vartheta}^2
	\notag
	\\*
	=& \frac{1}{2} \, T_{\mbox{w}} \, \boldsymbol{\vartheta}^2
	\, .
	\notag
\end{align}
The validity of  (\ref{eq:torque_stiffness_Cos})$_2$ for $M_{\mbox{m}}$ will be shown in the Appendix \ref{app:class_Coss}.
The plot of the torsional stiffness for the classical torque (light blue), the higher-order torque (red), and the torque energy (green) while varying $L_c$ is shown in Fig.~\ref{fig:all_plot_Cos_3}.
\subsection{Cosserat conformal curvature case - bounded stiffness in torsion}
\label{sec:Cos_confor}
In the particular case for which the parameter $a_3$ is equal to zero the elastic energy turns into
\begin{align}
	W \left(\boldsymbol{\mbox{D}u}, \boldsymbol{A},\mbox{Curl}\,\boldsymbol{A}\right) = &
	\, \mu_{\mbox{\tiny macro}} \left\lVert \mbox{sym} \, \boldsymbol{\mbox{D}u} \right\rVert^{2}
	+ \frac{\lambda_{\mbox{\tiny macro}}}{2} \mbox{tr}^2 \left(\boldsymbol{\mbox{D}u} \right) 
	+ \mu_{c} \left\lVert \mbox{skew} \left(\boldsymbol{\mbox{D}u} - \boldsymbol{A} \right) \right\rVert^{2}
	\label{eq:energy_Cos_conformal}
	\\
	&
	+ \frac{\mu \, L_c^2}{2}
	\,
	a_1 \, \left \lVert \mbox{dev} \, \mbox{sym} \, \mbox{Curl} \, \boldsymbol{A}\right \rVert^2
	\, .
	\notag
\end{align}
In terms of $\phi=\mbox{axl} (\boldsymbol{A})$, the curvature energy can be written as 
$\frac{\mu \, L_c^2}{2} \, a_1 \, \left \lVert \mbox{dev} \, \mbox{sym} \, \mbox{D} \, \mbox{axl} \left( \boldsymbol{A} \right) \right \rVert^2$ which is the conformal curvature case \cite{neff2010stable}.
In this special case, the torsional stiffness remains bounded as $L_c \to \infty$ ($R\to 0$), namely $\overline{T} \coloneqq \left(9 \mu _c + \mu _{\mbox{\tiny macro}} \right) I_{p}$, which is consistent with the results in  \eqref{eq:limi_conformal_RM_2}.

\begin{figure}[H]
	\centering
	\includegraphics[height=5.5cm]{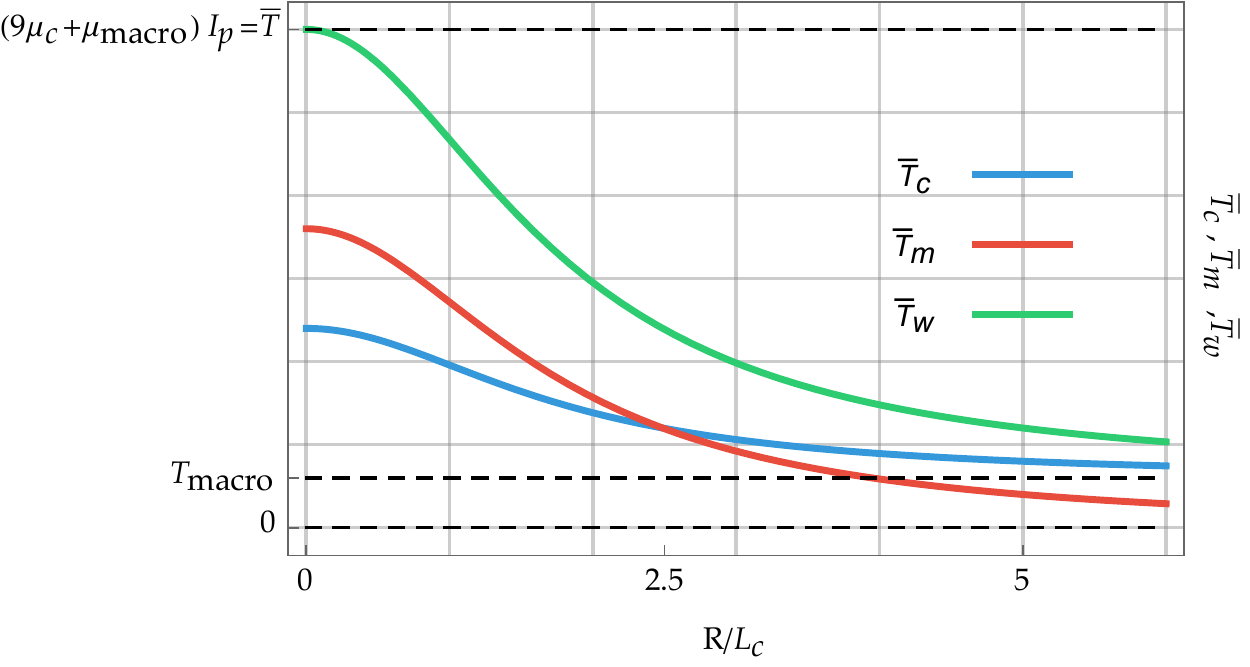}
	\caption{(\textbf{Cosserat model with conformal curvature}) Torsional stiffness for the classical torque $T_{\mbox{c}}$, the higher-order torque $T_{\mbox{m}}$, and the torque energy $T_{\mbox{w}}$ while varying $L_c$. The torsional stiffness is bounded as $L_c \to \infty$ ($R\to 0$). The values of the parameters used are: $\mu = 1$, $\mu _c= 1/2$, $\mu _{\mbox{\tiny macro}}= 1/2$, $a_1= 5$, $R= 1$.
	Here, the Cosserat couple modulus $\mu _c$ is clearly related to the value of the stiffness for small specimen size.}
	\label{fig:all_plot_Cos_4}
\end{figure}
\subsection{Cosserat limit case $\mu_c \to \infty$ (indeterminate couple stress model)}
\begin{align}
\displaystyle\lim_{\mu_c\to \infty} M_{\mbox{c}} (\boldsymbol{\vartheta})
&= \, 
\left[\mu _{\mbox{\tiny macro}} + a_{1} \mu \, \frac{L_c^2}{R^2}\right] I_{p} \, \boldsymbol{\vartheta}
=
T_{\mbox{c}} \, \boldsymbol{\vartheta}
\, ,
\qquad\qquad
\displaystyle\lim_{\mu_c\to \infty} M_{\mbox{m}} (\boldsymbol{\vartheta})
=  \, 
2 a_{1} \, \mu \, \frac{L_c^2}{R^2} \,  I_{p} \, \boldsymbol{\vartheta}
=
T_{\mbox{m}} \, \boldsymbol{\vartheta}
\, ,
\label{eq:torque_stiffness_Cos_Limit}
\\
\displaystyle\lim_{\mu_c\to \infty} W_{\mbox{tot}} (\boldsymbol{\vartheta})
&=  \, 
\frac{1}{2}\left[\mu _{\mbox{\tiny macro}} + 3a_{1} \mu \, \frac{L_c^2}{R^2}\right] I_{p} \, \boldsymbol{\vartheta}^2
=
\frac{1}{2} T_{\mbox{w}} \, \boldsymbol{\vartheta}^2
\, .
\notag
\end{align}
It is highlighted that there is \textit{not a one to one correspondence between the torque} obtained as a limit from the Cosserat model  (\ref{eq:torque_stiffness_Cos_Limit}) and the one obtained using the indeterminate couple stress model from the beginning  (\ref{eq:torque_stiffness_Ind_Coup_Stress}), but of course \textit{the energy (or the sum of the two torques) coincides} and thus the total torque stiffness $T_{\mbox{w}}$ coincides  as well.
\begin{figure}[H]
	\begin{subfigure}{0.5\textwidth}
		\centering
		\includegraphics[width=\textwidth]{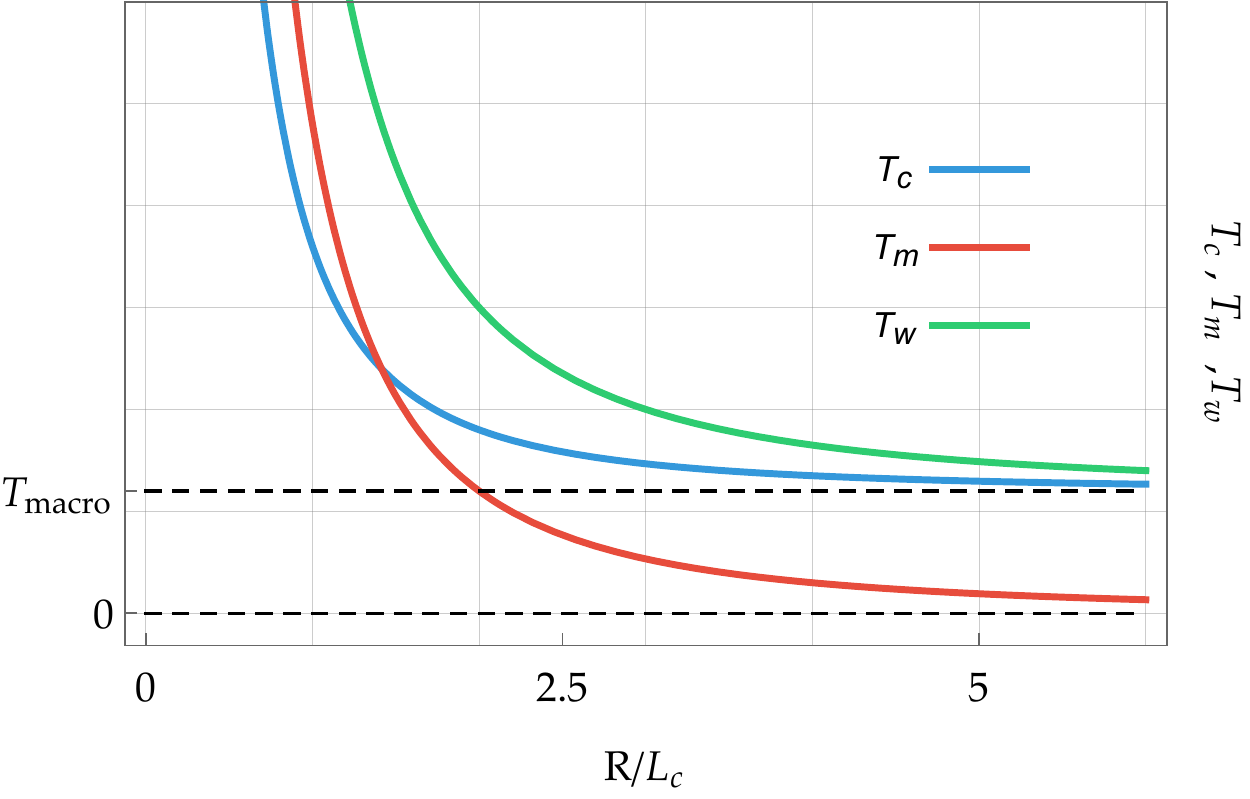}
		\caption{}
	\end{subfigure}
	\begin{subfigure}{0.5\textwidth}
		\centering
		\includegraphics[width=\textwidth]{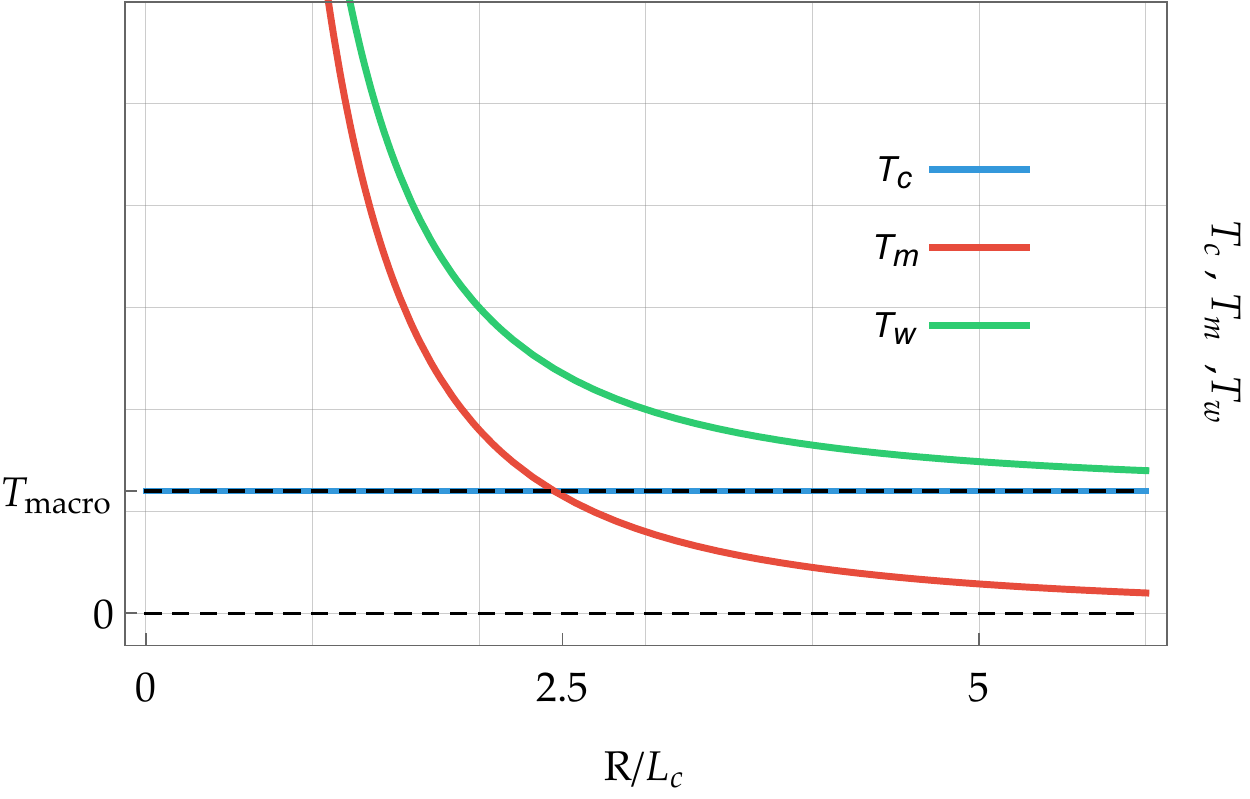}
		\caption{}
	\end{subfigure}
	\caption{(\textbf{Cosserat model vs indeterminate couple stress model})
		Comparison of the torsional stiffness for the classical torque $T_{\mbox{c}}$, the higher-order torque $T_{\mbox{m}}$, and the torque energy $T_{\mbox{w}}$ while varying $L_c$ for (a) the Cosserat model with $\mu_c \to \infty$ and for (b) the indeterminate couple stress model. There is not a one to one correspondence between the respective torque but the energy coincides.}
	\label{fig:torque_stiffness_Coss_vs_Couple_stress}
\end{figure}
\subsection{Cosserat limit case $\mu_c \to 0$.}
\label{sec:muc_zero_Coss}
\begin{align}
\displaystyle\lim_{\mu_c\to 0} M_{\mbox{c}} (\boldsymbol{\vartheta})
&= \, 
\mu _{\mbox{\tiny macro}} I_{p} \, \boldsymbol{\vartheta}
=
T_{\mbox{c}} \, \boldsymbol{\vartheta}
\, ,
\qquad\qquad
\displaystyle\lim_{\mu_c\to 0} M_{\mbox{m}} (\boldsymbol{\vartheta})
=  \, 
\frac{24 \mu \, a_{1} a_{3}}{a_{1}+8 a_{3}}\frac{ L_c^2}{R^2}
\,  I_{p} \, \boldsymbol{\vartheta}
=
T_{\mbox{m}} \, \boldsymbol{\vartheta}
\, ,
\label{eq:torque_stiffness_Cos_Limit_muc_0}
\\
\displaystyle\lim_{\mu_c\to 0} W_{\mbox{tot}} (\boldsymbol{\vartheta})
&=  \, 
\frac{1}{2}
\left[
\mu _{\mbox{\tiny macro}}
+ \frac{24 \mu \, a_{1} a_{3}}{a_{1}+8 a_{3}}\frac{ L_c^2}{R^2}
\right]
I_{p} \, \boldsymbol{\vartheta}^2
=
\frac{1}{2} T_{\mbox{w}} \, \boldsymbol{\vartheta}^2
\, .
\notag
\end{align}
It is highlighted that the Cosserat model does not collapse into a classical linear elastic model for $\mu_c  \to 0$, but it remains proportional to $(L_c/R)^2$ eq.(\ref{eq:torque_stiffness_Cos_Limit_muc_0}).
In this case, the Cosserat model behaves similarly to the indeterminate couple stress model eq.(\ref{eq:torque_stiffness_Cos_Limit}) or eq.(\ref{eq:torque_stiffness_Ind_Coup_Stress}), and it collapses into this model (both the energy and the torques) by formally letting $a_3 \to \infty$ as it can be seen from equations (\ref{eq:torque_stiffness_Cos_Limit}) and (\ref{eq:torque_stiffness_Cos_Limit_muc_0}).
\begin{figure}[H]
    \begin{subfigure}{0.5\textwidth}
		\centering
        \includegraphics[width=\textwidth]{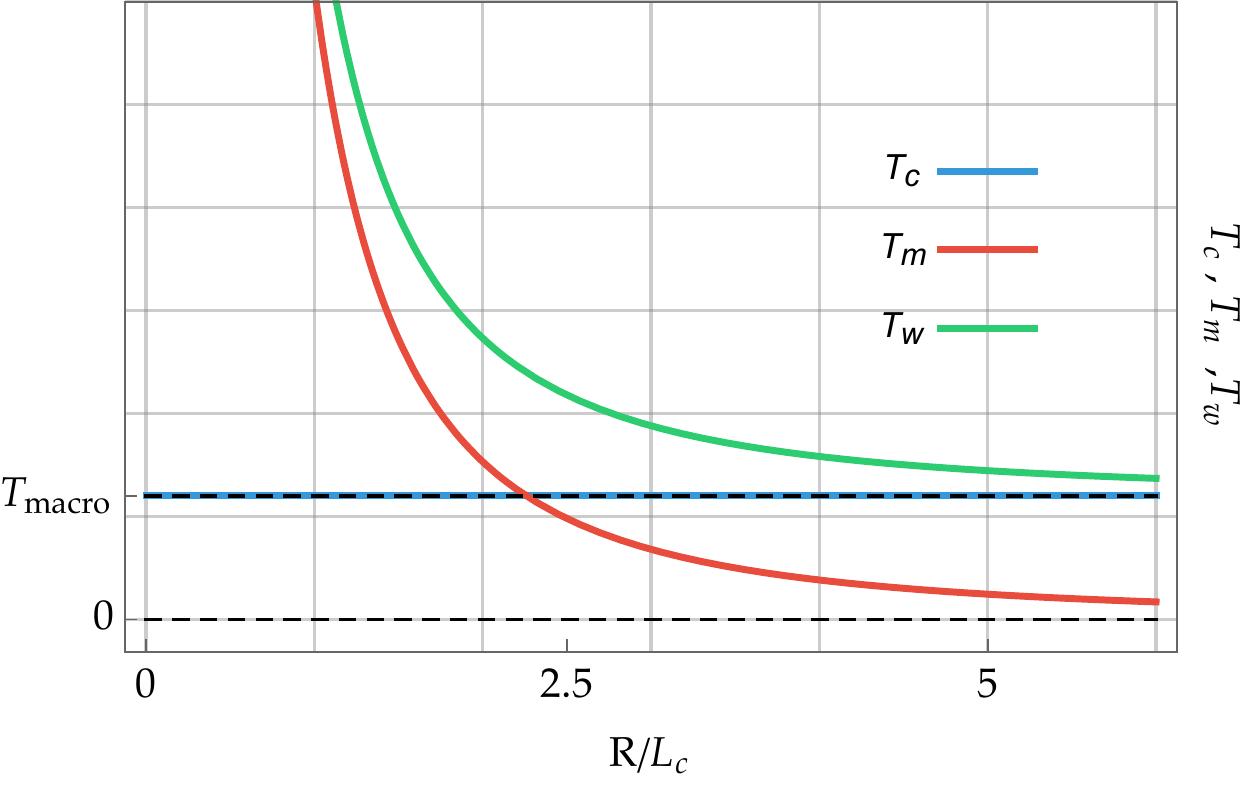}
    \end{subfigure}
	\hfill
    \begin{subfigure}{0.5\textwidth}
		\centering
        \includegraphics[width=\textwidth]{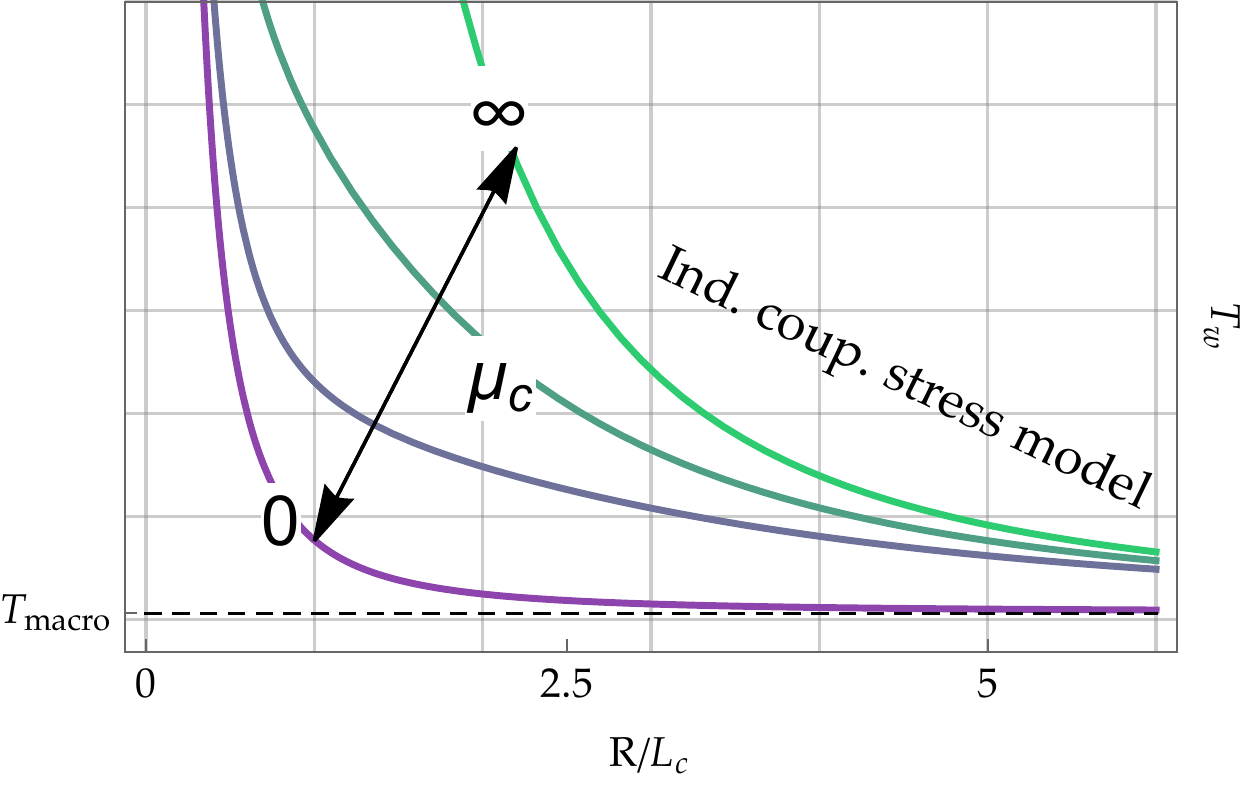}
    \end{subfigure}
	\caption{
	(\textbf{Cosserat model})
	(a) Torsional stiffness for the classical torque $T_{\mbox{c}}$, the higher-order torque $T_{\mbox{m}}$, and the torque energy $T_{\mbox{w}}$ while varying $L_c$ for the limit $\mu_c \to 0$. The model does not collapse into a classical linear elastic one.
	The values of the material parameter used are $\mu=1$, $\mu_e=1/10$, $a_1=1/5$, $a_3=1/7$, $R=1$.
	(b) Sensitivity study on how the Cosserat model behaves while varying $\mu_c=\{0,1/3,1,\infty\}$: for $\mu_c \to \infty$ we recover the indeterminate couple stress model, while for $\mu_c \to 0$ we still have a non linear relation between $T_{\mbox{w}}$ and $R/L_c$ since a classical linear elastic model is not attained (see eq.(\ref{eq:torque_stiffness_Cos_Limit_muc_0})).
	The values of the material parameter used are $\mu=1$, $\mu_e=1/10$, $a_1=12$, $a_3=1/20$, $R=1$.
	}
	\label{fig:torque_stiffness_Coss_muc_0}
\end{figure}
\subsection{Sensitivity of the Cosserat model with respect to the curvature parameters $a_1$ and $a_3$.}
Here, we study the sensitivity for the Cosserat model while varying $a_1$ and $a_3$ independently.
\begin{figure}[H]
	\begin{subfigure}{0.5\textwidth}
		\centering
		\includegraphics[width=\textwidth]{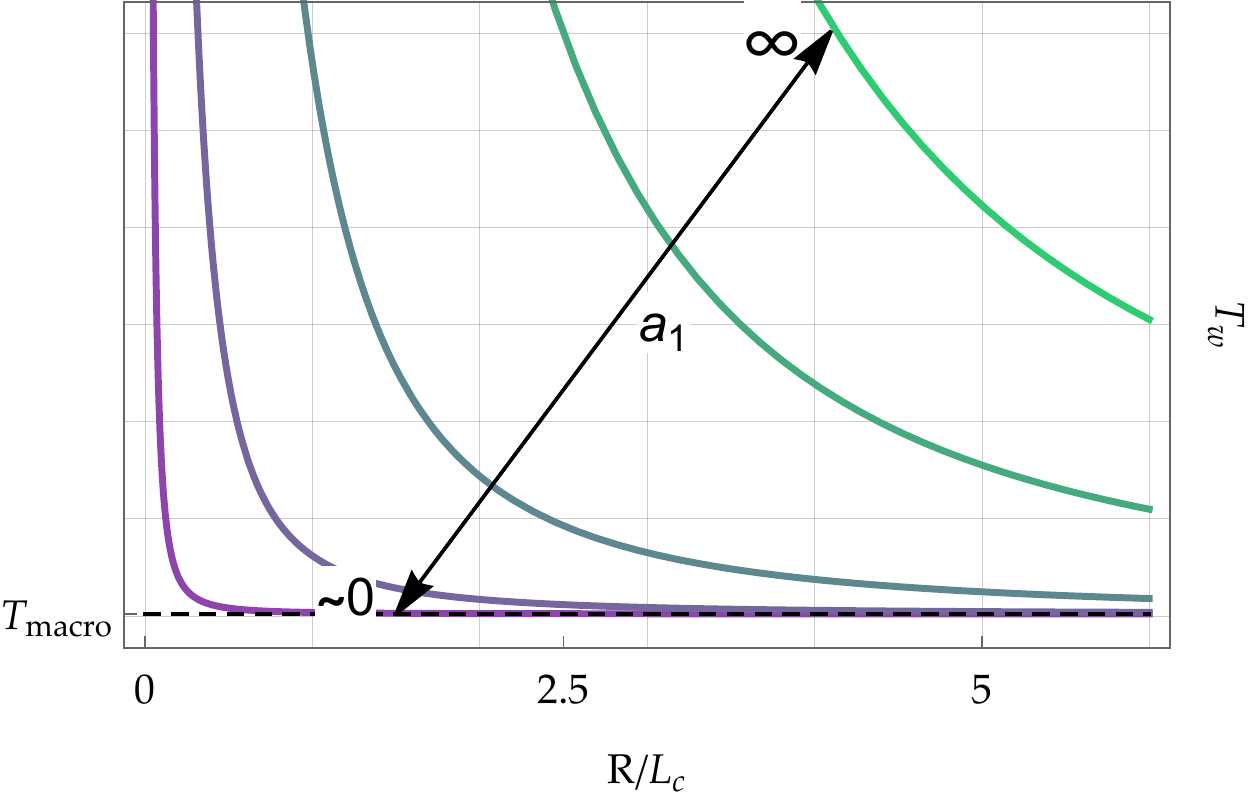}
		\caption{}
	\end{subfigure}
	\hfill
	\begin{subfigure}{0.5\textwidth}
		\centering
		\includegraphics[width=\textwidth]{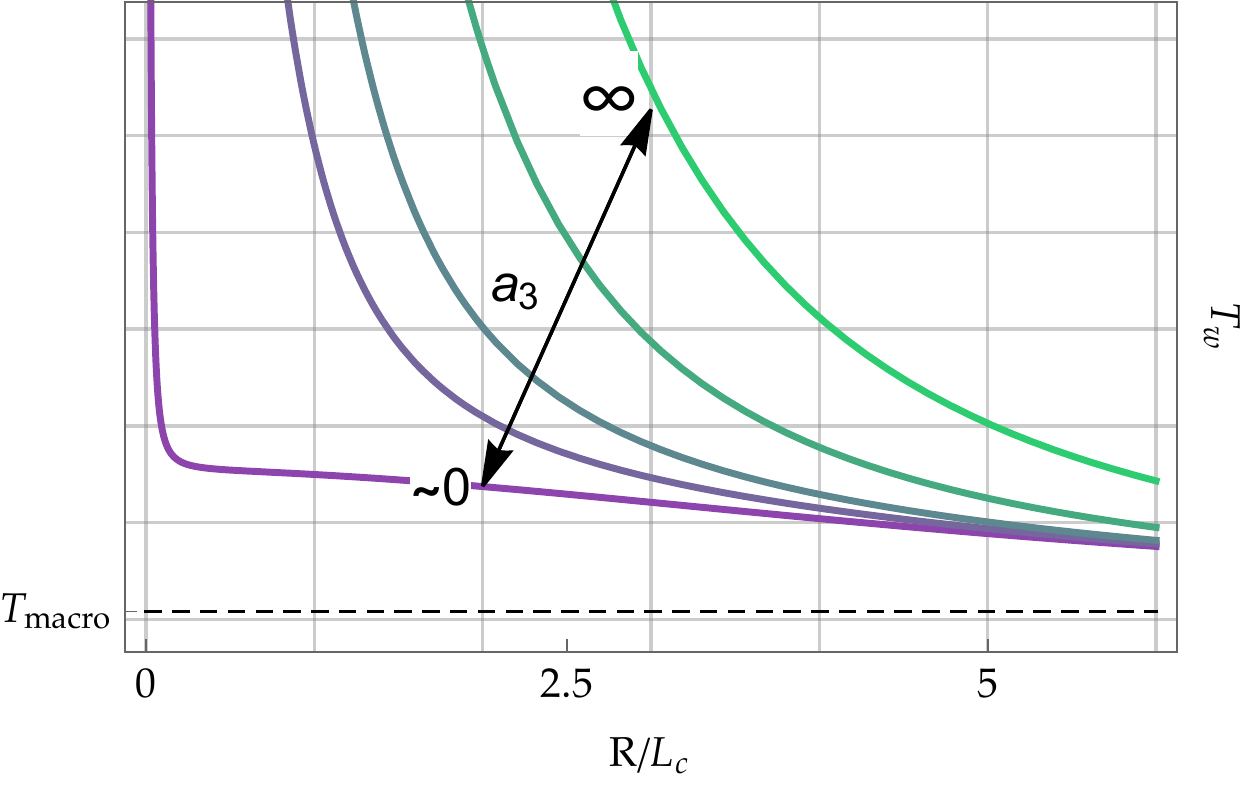}
		\caption{}
	\end{subfigure}
	\caption{(\textbf{Cosserat model}) Response of the Cosserat model while varying (a) the curvature parameter $a_1$ having $a_{3}= 20$ and (b) the curvature parameter $a_3$ having $a_{1}= 20$. The values of the other parameters are $\mu = 1$, $\mu _c= 1/5$, $\mu _{\mbox{\tiny macro}}= 1/10$, $R= 1$.}
	\label{fig:para_study_Cos}
\end{figure}
The parametric study represented in Fig. \ref{fig:para_study_Cos} has not been carried out for the limit $a_1 \to 0$ and $a_3 \to 0$ since we would have had an indeterminate form for $L_c \to \infty$, and that is why we used the symbol $\sim 0$.
The solution of the problem while having $a_3 = 0$ a priori is analyzed carefully in  Section  \ref{sec:Cos_confor}, and the solution of the problem while having $a_1 = 0$ a priori makes the relaxed micromorphic model collapse into a classical linear elastic model with torsional stiffness $T_{\tiny \mbox{macro}}$.
It is also highlighted that the Cosserat model collapses into the indeterminate couple stress model for $a_3\to \infty$ also in this more general case for which $\mu_c$ is arbitrary.
\section{Torsional problem for the isotropic micro-void model in dislocation tensor format}
\label{sec:Micro-void}

The strain energy for the isotropic micro-void continuum in dislocation tensor format can be written as \cite{rizzi2021shear,Cowin1983}
\begin{align}
	W \left(\boldsymbol{\mbox{D}u}, \omega ,\mbox{Curl}\,\left(\omega \boldsymbol{\mathbbm{1}}\right) \right) = &
	\, \mu_{\tiny \mbox{macro}} \left\lVert \mbox{dev} \, \mbox{sym} \, \boldsymbol{\mbox{D}u}\right\rVert^{2}
	+ \frac{\kappa_{e}}{2} \mbox{tr}^2 \left(\boldsymbol{\mbox{D}u} - \omega \boldsymbol{\mathbbm{1}} \right) 
	+ \frac{\kappa_{\tiny \mbox{micro}}}{2} \mbox{tr}^2 \left(\omega \boldsymbol{\mathbbm{1}} \right)
	\label{eq:energy_micro_void}
	\\
	&
	+ \frac{\mu \,L_c^2 }{2} \,
	a_2 \, \left\lVert \mbox{Curl} \, \left(\omega \boldsymbol{\mathbbm{1}}\right)\right\rVert^2 
	\, .
	\notag
\end{align}
Here, $\omega : \mathbb{R}^3 \to \mathbb{R}$ is the additional scalar micro-void degree of freedom \cite{Cowin1983}.
The equilibrium equations, in the absence of body forces,   are\!\!
\footnote{
	Where
	$\kappa_e = \frac{2\mu_e + 3\lambda_e}{3}$ and $\kappa_{\mbox{\tiny micro}} = \frac{2\mu_{\mbox{\tiny micro}} + 3\lambda_{\mbox{\tiny micro}}}{3}$ are the meso- and the micro-scale 3D bulk modulus.
}
\begin{align}
	\mbox{Div}\overbrace{\left[
		2\mu_{\tiny \mbox{macro}} \, \mbox{dev} \, \mbox{sym} \, \boldsymbol{\mbox{D}u}
		+ \kappa_{e} \mbox{tr} \left(\boldsymbol{\mbox{D}u} - \omega \boldsymbol{\mathbbm{1}} \right) \boldsymbol{\mathbbm{1}}
		\right]}^{\mathlarger{\widetilde{\sigma}}\coloneqq}
	&= \boldsymbol{0},
	\label{eq:equi_micro_void}
	\\
	\frac{1}{3}\mbox{tr}\left[\widetilde{\sigma}
	- \kappa_{\tiny \mbox{micro}} \mbox{tr} \left(\omega \boldsymbol{\mathbbm{1}}\right) \boldsymbol{\mathbbm{1}}
	- \mu \, L_{c}^{2} \, a_2 \, \mbox{Curl} \,
	\mbox{Curl} \, \left(\omega \boldsymbol{\mathbbm{1}}\right)
	\right] &= 0.
	\notag
\end{align}
The boundary conditions at the  free surface are
\begin{align}
	\boldsymbol{\widetilde{t}}(r = R) &= 
	\boldsymbol{\widetilde{\sigma}}(r) \cdot \boldsymbol{e}_{r} = 
	\boldsymbol{0}_{\mathbb{R}^{3}} \, ,
	\label{eq:BC_micro_void}
	\\
	\eta(r = R) &= 
	\frac{1}{3}\mbox{tr}\left(\boldsymbol{m} (r) \cdot \boldsymbol{\epsilon} \cdot \boldsymbol{e}_r\right) = 
	\frac{1}{3}\mbox{tr}\left(\boldsymbol{m} (r) \times \boldsymbol{e}_r \right) = 
	0 \, .
	\notag
\end{align}
According with the reference system shown in Fig.~\ref{fig:intro_1}, the ansatz for the displacement field and the function $\omega$ have to be
\begin{equation}
	\boldsymbol{u}(x_1,x_2)=
	\left(
	\begin{array}{c}
		-x_2 \, x_3 \\
		x_1 \, x_3 \\
		0 
	\end{array}
	\right) \, ,
	\qquad
	\omega\left(x_2\right) \boldsymbol{\mathbbm{1}} =
	\left(
	\begin{array}{ccc}
		0 & 0 & 0 \\
		0 & 0 & 0 \\
		0 & 0 & 0 \\
	\end{array}
	\right) \, .
	\label{eq:ansatz_micro_void}
\end{equation}
which clearly reduce the model to a classical linear elastic one.
No further calculation will be carried on and the reader is referred to  Section \ref{sec:Cau}.

\section{Torsional problem for the isotropic couple stress continuum}

The indeterminate couple stress model \cite{neff2015correct,hadjesfandiari2016comparison,tsiatas2011new,neff2016some,koiter1964couple,ghiba2017variant,tong2004size} appears by letting formally the Cosserat couple modulus $\mu_c \to \infty$. This implies the constraint $\boldsymbol{A}=\mbox{skew} \, \boldsymbol{\mbox{D}u} \in \mathfrak{so}(3)$.
It is highlighted that for the torsional problem, we do not have any unknown fields in this model since the displacement $\boldsymbol{u}$ is known a priori.
\!\!\!\footnote{
Since we can show that the classical torsion displacement solution satisfies the external balance equation (\ref{eq:equi_Coup_Stress}) as well as the higher order traction boundary conditions  (\ref{eq:BC_Coup_Stress}).
}

Since $\mbox{tr}(\mbox{Curl} \, \mbox{skew} \, \boldsymbol{\mbox{D}u}) = \left \lVert \mbox{skew} \, \mbox{Curl} \, \mbox{skew} \, \boldsymbol{\mbox{D}u}\right \rVert^2 = 0$ in terms of the ansatz  (\ref{eq:ansatz_Coup_Stress}), the indeterminate couple stress elastic energy for the torsion can be written as
\begin{align}
	W \left(\boldsymbol{\mbox{D}u},\mbox{Curl}\,\mbox{skew} \, \boldsymbol{\mbox{D}u}\right) = &
	\, \mu_{\mbox{\tiny macro}} \left\lVert \mbox{sym} \, \boldsymbol{\mbox{D}u} \right\rVert^{2}
	+ \frac{\lambda_{\mbox{\tiny macro}}}{2} \mbox{tr}^2 \left(\boldsymbol{\mbox{D}u} \right)
	+ \frac{\mu \, L_c^2}{2}
	\, a_1 \, \left \lVert \mbox{sym} \, \mbox{Curl} \, \mbox{skew} \, \boldsymbol{\mbox{D}u}\right \rVert^2 \, 
	\, .
	\label{eq:energy_Coup_Stress}
\end{align}
The equilibrium equations, in the absence of body forces,   are
\!\footnote{
Using Nye's formula \cite{ghiba2017variant}
$\mbox{Curl} \, \boldsymbol{A} = 
\mbox{tr} \left[ \left( \mbox{D} \, \mbox{axl} \boldsymbol{A} \right)^T\right]
- \left( \mbox{D} \, \mbox{axl} \boldsymbol{A} \right)^T
$
for
$
\boldsymbol{A} \in \mathfrak{so}(3)
$
we can rewrite
$
\mbox{Curl} \, \mbox{skew} \, \mbox{D} \boldsymbol{u} =
- \left( \mbox{D} \, \mbox{axl} \left( \mbox{skew} \, \mbox{D} \boldsymbol{u}\right) \right)^T = \frac{1}{2} \mbox{D} \, \mbox{curl} \, \boldsymbol{u}
$, since $\mbox{tr} \left( \mbox{Curl} \, \mbox{skew} \, \mbox{D}\boldsymbol{u} \right) = 0$.
}
\begin{align}
	\mbox{Div}\left[
	2\mu _{\mbox{\tiny macro}} \mbox{sym} \, \boldsymbol{\mbox{D}u} + 
	\lambda _{\mbox{\tiny macro}} \mbox{tr} \left(\boldsymbol{\mbox{D}u}\right) \boldsymbol{\mathbbm{1}}
	+ \mu \, L_c^2 \, \mbox{skew} \, \mbox{Curl}\,
	\left(
	a_1 \, \mbox{dev} \, \mbox{sym} \, \mbox{Curl} \, \mbox{skew} \, \boldsymbol{\mbox{D}u} \,
	\right)
	\right]
	=
	0\, ,
\label{eq:equi_Coup_Stress}
\end{align}
while the (highly non-trivial) boundary traction conditions on the free surface are (for more details see \cite{neff2015correct,hadjesfandiari2016comparison})
\begin{align}
	\pushleft{\boldsymbol{\widetilde{t}}(r = R) =
		\pm \, \bigg\{\left(\boldsymbol{\widetilde{\sigma}} -\frac{1}{2} \mbox{Anti}\left(\mbox{Div} \, \boldsymbol{m} \right)\right)\cdot \boldsymbol{e}_r
		-\frac{1}{2} \boldsymbol{e}_r\times \boldsymbol{\mbox{D}} \left[\left\langle \boldsymbol{e}_r,\mbox{sym} \, \boldsymbol{m}\cdot \boldsymbol{e}_r\right\rangle\right]}
	\notag
	\\*
	-\frac{1}{2}\boldsymbol{\mbox{D}}\left[\mbox{Anti} \, \left(\left(\boldsymbol{\mathbbm{1}} - \boldsymbol{e}_r\otimes \boldsymbol{e}_r\right)\cdot \boldsymbol{m}\cdot \boldsymbol{e}_r\right)\cdot\left(\boldsymbol{\mathbbm{1}} - \boldsymbol{e}_r\otimes \boldsymbol{e}_r\right)\right]:\left(\boldsymbol{\mathbbm{1}} - \boldsymbol{e}_r\otimes \boldsymbol{e}_r\right)\bigg\}
	=&\,
	\boldsymbol{0} \, ,
	\label{eq:BC_Coup_Stress}
	\\*
	\left(\boldsymbol{\mathbbm{1}} - \boldsymbol{e}_r\otimes \boldsymbol{e}_r\right) \cdot \boldsymbol{\eta}(r = R) =
	\pm \,  
	\left(\boldsymbol{\mathbbm{1}} - \boldsymbol{e}_r\otimes \boldsymbol{e}_r\right) \cdot 
	\mbox{Anti} \left[ \left(\boldsymbol{\mathbbm{1}} - \boldsymbol{e}_r\otimes \boldsymbol{e}_r\right)\cdot \boldsymbol{m}\cdot \boldsymbol{e}_r \right]\cdot \boldsymbol{e}_r =& \,
	\boldsymbol{0} \, ,
	\notag
	\\*
	\boldsymbol{\pi}(r = R) = 
	\pm \,  
	\left(
	\mbox{Anti}\left[\left(\boldsymbol{\mathbbm{1}} - \boldsymbol{e}_r\otimes \boldsymbol{e}_r\right)\cdot \boldsymbol{m}\cdot \boldsymbol{e}_r\right]^+
	-
	\mbox{Anti}\left[\left(\boldsymbol{\mathbbm{1}} - \boldsymbol{e}_r\otimes \boldsymbol{e}_r\right)\cdot \boldsymbol{m}\cdot \boldsymbol{e}_r\right]^-
	\right)\cdot \boldsymbol{e}_{\varphi}  =& \,
	\boldsymbol{0} \, ,
	\notag
\end{align}
where $\boldsymbol{\widetilde{\sigma}} = 2\,\mu _e\, \mbox{sym} \, \boldsymbol{\mbox{D}u} + 
\lambda _e\, \mbox{tr} \left(\boldsymbol{\mbox{D}u}\right) \boldsymbol{\mathbbm{1}}$ is the symmetric force stress tensor, $\boldsymbol{e}_r$ is the radial unit vector, and the non-symmetric second order moment stress is
\begin{equation}
    \boldsymbol{m} =
    \mu \, L_c^2 \, \left (
    a_1 \, \mbox{dev} \, \mbox{sym} \, \mbox{Curl} \, \mbox{skew} \, \boldsymbol{\mbox{D}u} \, 
    + a_2 \, \mbox{skew} \, \mbox{Curl} \, \mbox{skew} \, \boldsymbol{\mbox{D}u}
    \right) \, .
\end{equation}
The term
$(
\mbox{Anti}\left[\left(\boldsymbol{\mathbbm{1}} - \boldsymbol{e}_r\otimes \boldsymbol{e}_r\right)\cdot \boldsymbol{m}\cdot \boldsymbol{e}_r\right]^+
-
\mbox{Anti}\left[\left(\boldsymbol{\mathbbm{1}} - \boldsymbol{e}_r\otimes \boldsymbol{e}_r\right)\cdot \boldsymbol{m}\cdot \boldsymbol{e}_r\right]^-
)$ 
is the measure of the discontinuity of 
$ \mbox{Anti}\left[\left(\boldsymbol{\mathbbm{1}} - \boldsymbol{e}_r \otimes \boldsymbol{e}_r\right)\cdot \boldsymbol{m}\cdot \boldsymbol{e}_r\right] $ across the boundary.

According to the reference system shown in Fig.~\ref{fig:intro_1}, the ansatz for the displacement field and consequently the skew-symmetric part of the gradient of the displacement are
\begin{equation}
	\boldsymbol{u}(x_1,x_2)= \boldsymbol{\vartheta}
	\left(
	\begin{array}{c}
		-x_2 \, x_3 \\
		x_1 \, x_3 \\
		0 
	\end{array}
	\right) \, 
	\qquad\Rightarrow\qquad 
	\mbox{skew} \, \boldsymbol{\mbox{D}u} = \frac{\boldsymbol{\vartheta}}{2} 
	\left(
	\begin{array}{ccc}
		0 & -2 x_{3} & -x_{2} \\
		2 x_{3} & 0 & x_{1} \\
		x_{2} & -x_{1} & 0 \\
	\end{array}
	\right) \, .
	\label{eq:ansatz_Coup_Stress}
\end{equation}
Since the ansatz is completely known, it is possible to check that both the equilibrium equations (\ref{eq:equi_Coup_Stress}) and the boundary conditions  (\ref{eq:BC_Coup_Stress}) are identically satisfied\footnote{
	In Hadjesfandiari and Dargush \cite{hadjesfandiari2016comparison} the discussion of higher traction boundary conditions seems to be missing some terms in  (\ref{eq:BC_Coup_Stress}), letting the authors erroneously conclude that the classical displacement pure torsion solution does not satisfy the higher order boundary conditions.
}, 
and it is possible then to evaluate directly the classical torque, the higher-order torque, and the energy.

\begin{figure}[H]
	\centering
	\includegraphics[height=5.5cm]{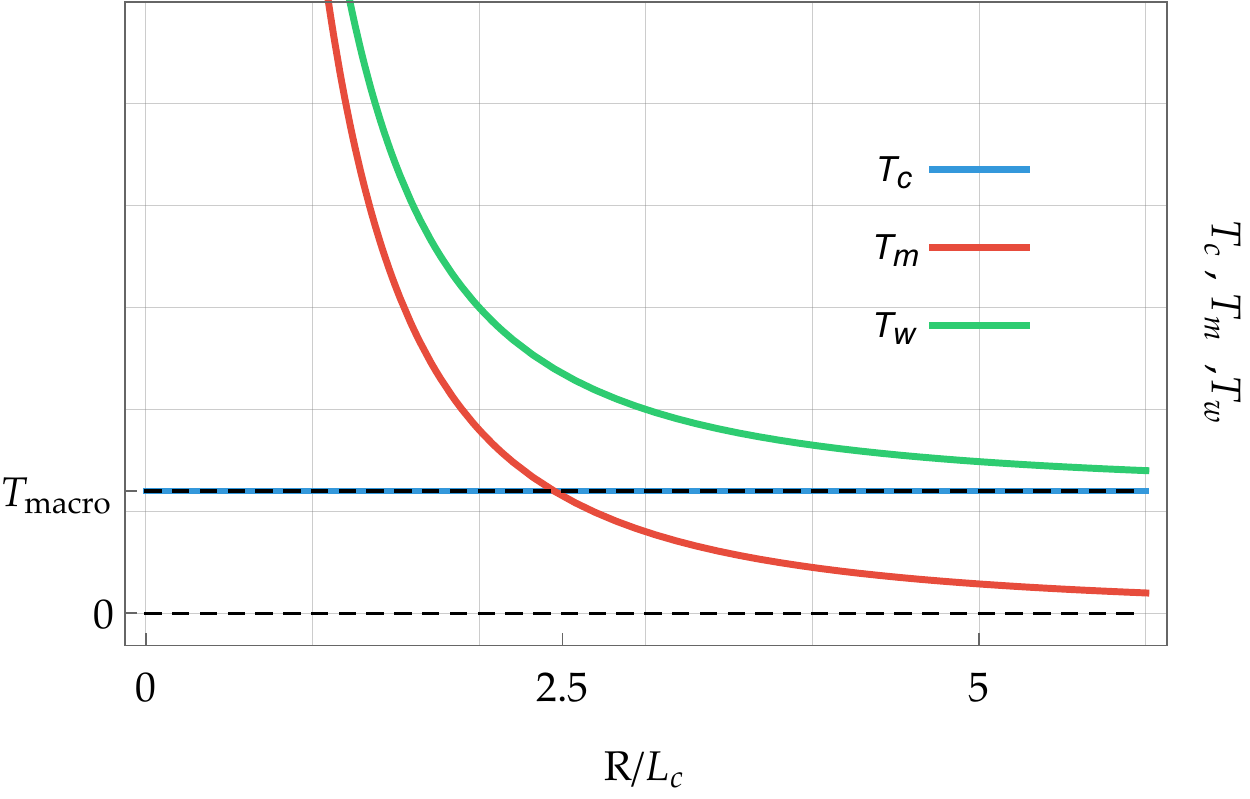}
	\caption{(\textbf{Indeterminate couple stress model}) Torsional stiffness for the classical torque $T_{\mbox{c}}$, the higher-order torque $T_{\mbox{m}}$, and the torque energy $T_{\mbox{w}}$ while varying $L_c$. The torsional stiffness is unbounded as $L_c \to \infty$ ($R\to 0$). The values of the parameters used are: $\mu = 1$, $\mu _e= 1/3$, $a_1= 1/5$, $R= 1$.}
	\label{fig:all_plot_Coup_Stress_3}
\end{figure}

The classical torque, the higher-order torque, and energy (per unit length d$z$) expressions are
\begin{align}
	M_{\mbox{c}} (\boldsymbol{\vartheta}) 
	\coloneqq&
	\int_{0}^{2\pi}
	\int_{0}^{R}
	\Big[
	\langle
	\widetilde{\boldsymbol{\sigma}} \, \boldsymbol{e}_{z} , \boldsymbol{e}_{\varphi}
	\rangle
	r
	\Big] r
	\, \mbox{d}r \, \mbox{d}\varphi
	=
	\mu _e
	I_{p} \, 
	\boldsymbol{\vartheta}
	= T_{\mbox{c}} \, \boldsymbol{\vartheta} \, ,
	\notag
	\\[3mm]
	M_{\mbox{m}}(\boldsymbol{\vartheta}) 
	\coloneqq&
	\int_{0}^{2\pi}
	\int_{0}^{R}
	\Big[
	\langle
	\left(\boldsymbol{m} \times \boldsymbol{e}_{z}\right)
	\boldsymbol{e}_{\varphi} ,
	\boldsymbol{e}_{r}
	\rangle
	-
	\langle
	\left(\boldsymbol{m} \times \boldsymbol{e}_{z}\right)
	\boldsymbol{e}_{r} ,
	\boldsymbol{e}_{\varphi}
	\rangle
	+
	\langle
	\left(\boldsymbol{m} \times \boldsymbol{e}_{r}\right)
	\boldsymbol{e}_{\varphi} ,
	\boldsymbol{e}_{z}
	\rangle
	-
	\langle
	\left(\boldsymbol{m} \times \boldsymbol{e}_{r}\right)
	\boldsymbol{e}_{z} ,
	\boldsymbol{e}_{\varphi}
	\rangle
	\Big]
	\, r
	\, \mbox{d}r \, \mbox{d}\varphi
	\notag
	\\*
	=& \, 
	3 a_{1} \, \mu\, \frac{L_c^2}{R^2} \, 
	I_{p} \, 
	\boldsymbol{\vartheta}
	= T_{\mbox{m}} \, \boldsymbol{\vartheta} \, ,
	\label{eq:torque_stiffness_Ind_Coup_Stress}
	\\[3mm]
	W_{\mbox{tot}} (\boldsymbol{\vartheta}) 
	\coloneqq&
	\int_{0}^{2\pi}
	\int_{0}^{R}
	W \left(\boldsymbol{\mbox{D}u}, \mbox{Curl} \, \mbox{skew} \, \boldsymbol{\mbox{D}u}\right) \, \, r
	\, \mbox{d}r \, \mbox{d}\varphi
	=
	\frac{1}{2}
	\left[
	\mu_{\tiny \mbox{macro}} + 3 a_{1} \, \mu\, \frac{L_c^2}{R^2}
	\right]
	I_{p} \, 
	\boldsymbol{\vartheta}^2
	= \frac{1}{2} \, T_{\mbox{w}} \, \boldsymbol{\vartheta}^2
	\, .
	\notag
\end{align}
The plot of the torsional stiffness for the classical torque, the higher-order torque, and the torque energy while varying $L_c$ is shown in Fig.~\ref{fig:all_plot_Coup_Stress_3}.
\subsection{Torsional problem for the modified and the ``pseudo''-consistent isotropic couple stress model}
The \textbf{modified couple stress model} \cite{munch2017modified,ghiba2017variant,madeo2016new,neff2015correct,neff2009subgrid} consists in choosing $a_1>0$, $a_2=0$ and leads to a symmetric couple stress tensor while the \textbf{(``pseudo'')-consistent couple stress model} \cite{hadjesfandiari2011couple} appears for $a_1=0$, $a_2>0$ and leads to a skew symmetric stress tensor $\boldsymbol{m}$.

Since $\mbox{tr}(\mbox{Curl} \, \mbox{skew} \, \boldsymbol{\mbox{D}u}) = \left \lVert \mbox{skew} \, \mbox{Curl} \, \mbox{skew} \, \boldsymbol{\mbox{D}u}\right \rVert^2 = 0$, the term $\left \lVert \mbox{dev} \,\mbox{sym} \, \mbox{Curl} \, \mbox{skew} \, \boldsymbol{\mbox{D}u}\right \rVert^2$ is the only non zero component in the curvature energy, the form of the energy remains the same.
This implies that, for the torsion problem, the \textbf{modified couple stress model} coincides with the \textbf{indeterminate couple stress model}, and that the \textbf{(``pseudo'')-consistent couple stress model} reduces to a classical linear elastic model without size-effects.
According to the notation \cite{hadjesfandiari2011couple}, the constitutive law can be written as
\!\footnote{
Setting $a_1 \, \mu_{\tiny \mbox{macro}} \, L_c^2 = 8\eta$ we obtain the rigidity as $T_{\mbox{w}} \coloneqq \frac{d^2}{d\boldsymbol{\vartheta}^2}W_{\mbox{tot}}(\boldsymbol{\vartheta}) = \mu_{\tiny \mbox{macro}} \left( 1 + 24 \frac{\eta}{\mu_{\tiny \mbox{macro}}}\frac{1}{R^2} \right) I_p$.
In  (44) of Hadjesfandiari and Dargush \cite{hadjesfandiari2016comparison} we have the relation $\ell^2=\frac{\eta}{\mu_{\tiny \mbox{macro}}}$, while in  (55) we have the formula $T_{\mbox{w}} = \mu_{\tiny \mbox{macro}} \left( 1 + 24\left( \frac{\ell}{R} \right)^2 \right) I_p = \mu_{\tiny \mbox{macro}} \left( 1 + 24\frac{\eta}{\mu_{\tiny \mbox{macro}}}\frac{1}{R^2} \right) I_p$.
}
\begin{equation}
    \widetilde{\boldsymbol{\sigma}}
    =
    2\mu _{\mbox{\tiny macro}} \, \mbox{sym} \, \boldsymbol{\mbox{D}u} + 
	\lambda _{\mbox{\tiny macro}} \, \mbox{tr} \left(\boldsymbol{\mbox{D}u}\right) \boldsymbol{\mathbbm{1}}
	\, ,
	\qquad
	\boldsymbol{m}
	=
	\eta \left( \mbox{D} \, \mbox{curl} \, \boldsymbol{u} \right)^T
	+
	\eta' \, \mbox{D} \, \mbox{curl} \, \boldsymbol{u}
	\, .
\end{equation}
where according to the classical Cosserat notation (see Appendix \ref{app:class_Coss})
\begin{align}
    \eta=\beta=\mu_{\mbox{\tiny macro}} \frac{L_c^2}{2} \frac{a_1-a_2}{2}
    \, ,
    \qquad\qquad\qquad
    \eta'=\gamma=\mu_{\mbox{\tiny macro}} \frac{L_c^2}{2} \frac{a_1+a_2}{2}
    \, .
\end{align}

In this notation, the modified couple stress model appears for $\eta=\eta'$ and the ``pseudo"-consistent couple stress model appears for $\eta=-\eta'$. 
\section{Torsional problem for the classical isotropic micromorphic continuum without mixed terms}
\label{sec:Micro_morphic}

The strain energy for the isotropic micromorphic continuum without mixed terms ($\langle\mbox{sym}\, \boldsymbol{P}, \mbox{sym} \, \left(\boldsymbol{\mbox{D}u} -\boldsymbol{P}\right)\rangle$, etc.) and simplified isotropic curvature can be written as
\begin{align}
	W \left(\boldsymbol{\mbox{D}u}, \boldsymbol{P}, \boldsymbol{\mbox{D}P}\right)
	= &
	\, \mu_{e} \left\lVert \mbox{dev} \, \mbox{sym} \left(\boldsymbol{\mbox{D}u} - \boldsymbol{P} \right) \right\rVert^{2}
	+ \frac{\kappa_{e}}{2} \mbox{tr}^2 \left(\boldsymbol{\mbox{D}u} - \boldsymbol{P} \right)
	+ \mu_{c} \left\lVert \mbox{skew} \left(\boldsymbol{\mbox{D}u} - \boldsymbol{P} \right) \right\rVert^{2}
	\label{eq:energy_MM}
	\\
	&
	+ \mu_{\tiny \mbox{micro}} \left\lVert \mbox{dev} \, \mbox{sym}\,\boldsymbol{P} \right\rVert^{2}
	+ \frac{\kappa_{\tiny \mbox{micro}}}{2} \mbox{tr}^2 \left(\boldsymbol{P} \right)
	\notag
	\\
	&
	+ \frac{\mu \, L_c^2}{2}
	\Bigg(
	a_1 \, \left\lVert \mbox{D}\left(\mbox{dev} \, \mbox{sym} \, \boldsymbol{P}\right) \right\rVert^2
	+ a_2 \, \left\lVert \mbox{D}\left(\mbox{skew} \, \boldsymbol{P}\right) \right\rVert^2
	+ \frac{2}{9} \, a_3 \left\lVert \mbox{D} \left(\mbox{tr} \left(\boldsymbol{P}\right)\boldsymbol{\mathbbm{1}} \right) \right\rVert^2 \Big)
	\Bigg)
	\, .
	\notag
\end{align}
The equilibrium equations, in the absence of body forces,   are the following 
\begin{align}
	\mbox{Div}\overbrace{\left[2\mu_{e} \, \mbox{dev} \, \mbox{sym} \left(\boldsymbol{\mbox{D}u} - \boldsymbol{P} \right) + \kappa_{e} \mbox{tr} \left(\boldsymbol{\mbox{D}u} - \boldsymbol{P} \right) \boldsymbol{\mathbbm{1}}
		+ 2\mu_{c}\,\mbox{skew} \left(\boldsymbol{\mbox{D}u} - \boldsymbol{P} \right)\right]}^{\mathlarger{\widetilde{\sigma}}\coloneqq}
	= \boldsymbol{0} \, ,
	\notag
	\\
	\widetilde{\sigma}
	- 2 \mu_{\mbox{\tiny micro}} \, \mbox{dev} \,\mbox{sym}\,\boldsymbol{P}
	- \kappa_{\tiny \mbox{micro}} \mbox{tr} \left(\boldsymbol{P}\right) \boldsymbol{\mathbbm{1}}
	\hspace{9cm}
	\label{eq:equi_MM}
	\\
	+\mu L_{c}^{2} \,
	\left[
	a_1 \, \mbox{dev} \, \mbox{sym} \, \boldsymbol{\Delta P}
	+ a_2 \, \mbox{skew} \, \boldsymbol{\Delta P}
	+ \frac{2}{9} \, a_3 \, \mbox{tr} \left(\boldsymbol{\Delta P}\right)\boldsymbol{\mathbbm{1}}
	\right]
	= \boldsymbol{0} \, ,
	\notag
\end{align}
where $\boldsymbol{\Delta P} \in \mathbb{R}^{3\times3}$ is taken component-wise.
The boundary conditions at the external surfaces are 
\begin{equation}
	\boldsymbol{\widetilde{t}}(r = R) = 
	\boldsymbol{\widetilde{\sigma}}(r) \cdot \boldsymbol{e}_{r} = 
	\boldsymbol{0}_{\mathbb{R}^{3}}
	\, ,
	\qquad
	\boldsymbol{\eta}(r = R) = 
	\boldsymbol{\mathfrak{m}} (r) \cdot \boldsymbol{e}_{r} =
	\boldsymbol{0}_{\mathbb{R}^{3\times3}}
	\, ,
	\label{eq:BC_MM}
\end{equation}
where 
\begin{equation}
    \boldsymbol{\mathfrak{m}} = \mu \, L_{c}^{2}\,
    \left[
    a_1 \, \mbox{D} \left(\mbox{dev} \, \mbox{sym} \, \boldsymbol{P}\right)
    + a_2 \, \mbox{D} \left(\mbox{skew} \, \boldsymbol{P}\right)
    + \frac{2}{9} \, a_3 \, \mbox{D} \left(\mbox{tr} \left(\boldsymbol{S}\right)\boldsymbol{\mathbbm{1}} \right)
    \right]
\end{equation}
is the third order moment stress tensor,
the expression of $\boldsymbol{\widetilde{\sigma}}$ is given in  (\ref{eq:equi_MM}), and $\boldsymbol{e}_{r}$ is the radial unit vector.
According with the reference system shown in Fig.~\ref{fig:intro_1}, the ansatz for the displacement field and the micro-distortion is
\begin{align}
	\boldsymbol{u}(x_1,x_2,x_3) &=
	\boldsymbol{u}(r,\varphi,z) =
	\boldsymbol{\vartheta}
	\left(
	\begin{array}{c}
		-x_2(r,\varphi) \, x_3(z) \\
		x_1(r,\varphi) \, x_3(z) \\
		0 
	\end{array}
	\right)
	\, ,
	\notag
	\\*
	\boldsymbol{P}(x_1,x_2,x_3) &=
	\boldsymbol{P}(r,\varphi,z) =
	\boldsymbol{\vartheta}
	\left(
	\begin{array}{ccc}
		0   & -x_3(z) & - g_{2}(r) \, x_2(r,\varphi) \\
		x_3(z) &    0 &   g_{2}(r) \, x_1(r,\varphi) \\
		g_{1}(r) \, x_2(r,\varphi) & - g_{1}(r) \, x_1(r,\varphi) & 0 \\
	\end{array}
	\right)
	\, .
	\label{eq:ansatz_MM}
\end{align}
Substituting the ansatz  (\ref{eq:ansatz_MM}) in  (\ref{eq:equi_MM}) the 12 equilibrium equations are equivalent to
\begin{align}
	\frac{1}{2} \boldsymbol{\vartheta}  \sin \varphi  \left(\rho  \left(\mu  L_c^2 \left((a_{1}-a_{2}) g_{1}''(\rho )-(a_{1}+a_{2}) g_{2}''(\rho )\right)+2 \mu _c (g_{1}(\rho )+g_{2}(\rho )-1)
	\right.
	\right.
	\hspace{2cm}
	&
\notag
	\\*
	\left.
	\left.
	-2 \left(\mu _e+\mu _{\mbox{\tiny micro}}\right) (g_{1}(\rho )-g_{2}(\rho ))	-2 \mu _e\right)
	+3 \mu  L_c^2 \left((a_{1}-a_{2}) g_{1}'(\rho )-(a_{1}+a_{2}) g_{2}'(\rho )\right)\right)&= 0
	\, ,
	\notag
	\\*
	\frac{1}{2} \boldsymbol{\vartheta}  \cos \varphi  \left(\rho  \left(\mu  L_c^2 \left((a_{2}-a_{1}) g_{1}''(\rho )+(a_{1}+a_{2}) g_{2}''(\rho )\right)-2 \mu _c (g_{1}(\rho )+g_{2}(\rho )-1)
	\right.
	\right.
	\hspace{2cm}
	&
	\notag
	\\*
	\left.
	\left.
	+2 \left(\mu _e+\mu _{\mbox{\tiny micro}}\right) (g_{1}(\rho )-g_{2}(\rho ))+2 \mu _e\right)
	+3 \mu  L_c^2 \left((a_{2}-a_{1}) g_{1}'(\rho )+(a_{1}+a_{2}) g_{2}'(\rho )\right)\right)&= 0
	\, ,
	\label{eq:equi_equa_MM}
	\\*
	\frac{1}{2} \boldsymbol{\vartheta}  \sin \varphi  \left(\rho  \left(\mu  L_c^2 \left((a_{1}+a_{2}) g_{1}''(\rho )+(a_{2}-a_{1}) g_{2}''(\rho )\right)-2 \left(\mu _c (g_{1}(\rho )+g_{2}(\rho )-1)
	\right.
	\right.
	\right.
	\hspace{2cm}
	&
	\notag
	\\*
	\left.
	\left.
	\left.
	+\left(\mu _e+\mu _{\mbox{\tiny micro}}\right) (g_{1}(\rho )-g_{2}(\rho ))\right)-2 \mu _e\right)
	+3 \mu  L_c^2 \left((a_{1}+a_{2}) g_{1}'(\rho )+(a_{2}-a_{1}) g_{2}'(\rho )\right)\right)&= 0
	\, ,
	\notag
	\\*
	\frac{1}{2} \boldsymbol{\vartheta}  \cos \varphi  \left(\rho  \left(\mu  L_c^2 \left((a_{1}-a_{2}) g_{2}''(\rho )-(a_{1}+a_{2}) g_{1}''(\rho )\right)+2 \left(\mu _c (g_{1}(\rho )+g_{2}(\rho )-1)
	\right.
	\right.
	\right.
	\hspace{2cm}
	&
	\notag
	\\*
	\left.
	\left.
	\left.
	+\left(\mu _e+\mu _{\mbox{\tiny micro}}\right) (g_{1}(\rho )-g_{2}(\rho ))\right)+2 \mu _e\right)
	+3 \mu  L_c^2 \left((a_{1}-a_{2}) g_{2}'(\rho )-(a_{1}+a_{2}) g_{1}'(\rho )\right)\right) &=0
	\, .
	\notag
\end{align}
It is important to underline that  (\ref{eq:equi_MM})$_1$ is identically satisfied, and that between the four equilibrium equations  (\ref{eq:equi_equa_MM}) there are only two that are linearly independent since:  (\ref{eq:equi_equa_MM})$_1=\tan \varphi$  (\ref{eq:equi_equa_MM})$_2$ and  (\ref{eq:equi_equa_MM})$_3=\tan \varphi$  (\ref{eq:equi_equa_MM})$_4$.

It is also pointed out that the two remaining linearly independent equations  (\ref{eq:equi_equa_MM})$_{1,3}$ can be uncoupled\footnote{That this uncoupling takes place at all seems to be connected to the chosen form of the curvature energy. It remains unclear at present whether this feature holds for the most general isotropic curvature expression as well.} and have the form of the Bessel ODE if we take their sum and difference, while being careful of substituting $g_p(r) \coloneqq g_{1}(r) + g_{2}(r)$ and $g_m(r) \coloneqq g_{1}(r) - g_{2}(r)$ along with their derivatives:
\begin{align}
	\boldsymbol{\vartheta}  \sin \varphi  \left(a_{1} \mu  L_c^2 \left(3 g_{m}'(\rho )+\rho \,  g_{m}''(\rho )\right)-2 \rho  \mu _e (g_{m}(\rho )+1)-2 \rho \,  g_{m}(\rho ) \mu _{\mbox{\tiny micro}}\right)
	&= 0
	\, ,
	\label{eq:equi_equa_MM_2}
	\\*
	\boldsymbol{\vartheta}  \sin \varphi  \left(2 \rho \,  \mu _c (g_{p}(\rho )-1)-a_{2} \mu  L_c^2 \left(3 g_{p}'(\rho )+\rho \,  g_{p}''(\rho )\right)\right)
	&= 0
	\, .
	\notag
\end{align}
Since $g_1(r)\coloneqq\frac{g_{p}(r) + g_{m}(r)}{2}$ and $g_2(r)\coloneqq\frac{g_{p}(r) - g_{m}(r)}{2}$, the solution in terms of $g_1(r)$ and $g_2(r)$ of  (\ref{eq:equi_equa_MM_2}) is 
\begin{align}
	g_1(r) = \, & 
	\frac{1}{2} \left(
	1
	-\frac{
		i A_{1} \, I_1\left(      \frac{r \, f_2}{L_c}\right)
		-  A_{2} \, Y_1\left(-i \, \frac{r \, f_2}{L_c}\right)
		+i B_{1} \, I_1\left(      \frac{r \, f_1}{L_c}\right)
		-  B_{2} \, Y_1\left(-i \, \frac{r \, f_1}{L_c}\right)
	}{
		r
	}
	-\frac{\mu _e}{\mu _e+\mu _{\mbox{\tiny micro}}}
	\right) \, ,
	\notag
	\\*
	g_2(r) = \, & 
	\frac{1}{2} \left(
	1
	+\frac{
		i A_{1} \, I_1\left(      \frac{r \, f_2}{L_c}\right)
		-  A_{2} \, Y_1\left(-i \, \frac{r \, f_2}{L_c}\right)
		-i B_{1} \, I_1\left(      \frac{r \, f_1}{L_c}\right)
		+  B_{2} \, Y_1\left(-i \, \frac{r \, f_1}{L_c}\right)
	}{
		r
	}
	+\frac{\mu _e}{\mu _e+\mu _{\mbox{\tiny micro}}}
	\right) \, ,
	\label{eq:sol_fun_MM}
	\\*
	f_{1} \coloneqq \, &
	\sqrt{\frac{2 \mu _c}{a_{2} \, \mu}} \, ,
	\qquad
	f_{2} \coloneqq
	\sqrt{\frac{2(\mu _e + \mu _{\mbox{\tiny micro}})}{a_{1} \, \mu }} \, ,
	\notag
\end{align}
where $I_{n}\left(\cdot\right)$ is the \textit{modified Bessel function of the first kind}, $Y_{n}\left(\cdot\right)$ is the \textit{Bessel function of the second kind} (see appendix \ref{app:bessel} for the formal definitions), and
$A_1$, $B_1$, $A_2$, $B_2$ are integration constants.

\begin{figure}[H]
	\centering
	\includegraphics[height=5.5cm]{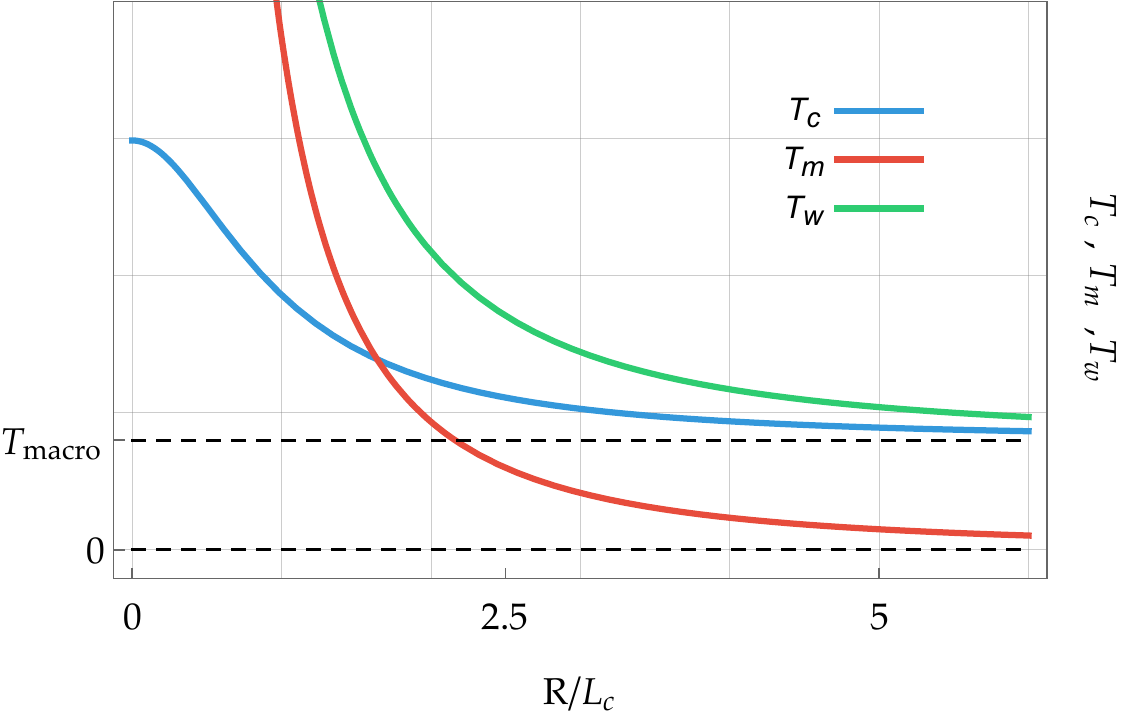}
	\caption{(\textbf{Micromorphic model}, classical case) Torsional stiffness for the classical torque $T_{\mbox{c}}$, the higher-order torque $T_{\mbox{m}}$, and the torque energy $T_{\mbox{w}}$ while varying $L_c$. The torsional stiffness is unbounded as $L_c \to \infty$ ($R\to 0$). The values of the parameters used are: $\mu = 1$, $\mu _e= 1/3$, $\mu_{\mbox{\tiny micro}} = 1/4$, $\mu _c= 1/5$, $a_1= 1/5$, $a_2= 1/6$, $R= 1$.}
	\label{fig:all_plot_MM}
\end{figure}

The values of $A_1$, $B_1$ are determined thanks to the boundary conditions  (\ref{eq:BC_MM}), while, due to the divergent behaviour of the Bessel function of the second kind at $r=0$, we have to set $A_2=0$ and $B_2=0$ in order to have a continuous solution.
The fulfilment of the boundary conditions  (\ref{eq:BC_MM}) allows us to find the expressions of the integration constants
\begin{equation}
	A_1 = \, 
	\frac{2 i L_c \mu _e}{f_2 \left(\mu _e+\mu _{\mbox{\tiny micro}}\right) \left(I_0\left(\frac{R f_2}{L_c}\right)+I_2\left(\frac{R f_2}{L_c}\right)\right)}
	\, ,
	\qquad\qquad 
	B_1 = \,  
	-\frac{2 i L_c}{f_1 \left(I_0\left(\frac{R f_1}{L_c}\right)+I_2\left(\frac{R f_1}{L_c}\right)\right)}
	\, .
	\label{eq:BC_MM_3}
\end{equation}
The classical torque, the higher-order torque, and energy (per unit length d$z$) expressions are
\begin{align}
	M_{\mbox{c}} (\boldsymbol{\vartheta}) 
	\coloneqq&
	\int_{0}^{2\pi}
	\int_{0}^{R}
	\Big[
	\langle
	\widetilde{\boldsymbol{\sigma}} \, \boldsymbol{e}_{z} , \boldsymbol{e}_{\varphi}
	\rangle
	r
	\Big] r
	\, \mbox{d}r \, \mbox{d}\varphi
	\notag
	\\*
	=&
	\left[
	\left(
	\frac{
	8 \mu _c I_2\left(\frac{R f_1}{L_c}\right)
	}{
	f_1^2 I_0\left(\frac{R f_1}{L_c}\right)+f_1^2 I_2\left(\frac{R f_1}{L_c}\right)
	}
	+ \frac{
	8 \mu _e^2 I_2\left(\frac{R f_2}{L_c}\right)
	}{
	\left(\mu _e+\mu _{\mbox{\tiny micro}}\right) \left(f_2^2 I_0\left(\frac{R f_2}{L_c}\right)+f_2^2 I_2\left(\frac{R f_2}{L_c}\right)\right)
	}
	\right) \frac{L_c^2 }{R^2}
	+\frac{\mu _e \, \mu _{\mbox{\tiny micro}}}{\mu _e+\mu _{\mbox{\tiny micro}}}
	\right]
	I_{p} \, 
	\boldsymbol{\vartheta}
	\notag
	\\*
	=& \, T_{\mbox{c}} \, \boldsymbol{\vartheta} \, ,
	\notag
	\\[3mm]
	M_{\mbox{m}}(\boldsymbol{\vartheta}) 
	\coloneqq&
	\int_{0}^{2\pi}
	\int_{0}^{R}
	\Big[
	\langle
	\left(\boldsymbol{\mathfrak{m}} \, \boldsymbol{e}_{z}\right)
	\boldsymbol{e}_{\varphi} ,
	\boldsymbol{e}_{r}
	\rangle
	-
	\langle
	\left(\boldsymbol{\mathfrak{m}} \, \boldsymbol{e}_{z}\right)
	\boldsymbol{e}_{r} ,
	\boldsymbol{e}_{\varphi}
	\rangle
	\Big]
	\, r
	\, \mbox{d}r \, \mbox{d}\varphi
	=
	4 a_{2} \mu \frac{L_c^2}{R^2}
	I_{p} \, 
	\boldsymbol{\vartheta}
	=
	T_{\mbox{m}} \, \boldsymbol{\vartheta} \, ,
	\notag
	\\[3mm]
	W_{\mbox{tot}} (\boldsymbol{\vartheta}) 
	\coloneqq&
	\int_{0}^{2\pi}
	\int_{0}^{R}
	W \left(\boldsymbol{\mbox{D}u}, \boldsymbol{P}, \mbox{D}\boldsymbol{P}\right) \, \, r
	\, \mbox{d}r \, \mbox{d}\varphi
	\notag
	\\*
	=&
	\frac{1}{2}
	\Bigg[
	\Bigg(
	\frac{
		8 \mu _c I_2\left(\frac{R f_1}{L_c}\right)
	}{
		f_1^2 I_0\left(\frac{R f_1}{L_c}\right)+f_1^2 I_2\left(\frac{R f_1}{L_c}\right)
	}
	+ \frac{
		8 \mu _e^2 I_2\left(\frac{R f_2}{L_c}\right)
	}{
		\left(\mu _e+\mu _{\mbox{\tiny micro}}\right) \left(f_2^2 I_0\left(\frac{R f_2}{L_c}\right)+f_2^2 I_2\left(\frac{R f_2}{L_c}\right)\right)
	}
	\Bigg)
	\frac{L_c^2 }{R^2}
	\label{eq:torque_stiffness_MM}
	\\*
	&
	\hspace{0.5cm}
	+
	\underbrace{
	\frac{\mu _e \, \mu _{\mbox{\tiny micro}}}{\mu _e+\mu _{\mbox{\tiny micro}}}
	}_{\mu_{\tiny \mbox{macro}}}
	+4 a_{2} \mu \frac{L_c^2}{R^2}
	\Bigg]
	I_{p} \, 
	\boldsymbol{\vartheta}^2
	=
	\frac{1}{2} \, T_{\mbox{w}} \, \boldsymbol{\vartheta}^2
	\, ,
	\notag
\end{align}
and again it holds,
\begin{equation}
\frac{\mbox{d}}{\mbox{d}\boldsymbol{\vartheta}}W_{\mbox{tot}}(\boldsymbol{\vartheta}) = M_{\mbox{c}} (\boldsymbol{\vartheta}) + M_{\mbox{m}} (\boldsymbol{\vartheta}) \, ,
\qquad\qquad\qquad
\frac{\mbox{d}^2}{\mbox{d}\boldsymbol{\vartheta}^2}W_{\mbox{tot}}(\boldsymbol{\vartheta})
= T_{\mbox{c}} + T_{\mbox{m}}
= T_{\mbox{w}}
\, .
\end{equation}
It is underlined that the boundary conditions for the micromorphic model are consistent with the relaxed micromorphic model's one, being careful of changing $\mathfrak{\boldsymbol{m}}\,e_{z}$ with $\boldsymbol{m} \times e_{z}$.
The plot of the torsional stiffness for the classical torque (light blue), the higher-order torque (red), and the torque energy (green) while varying $L_c$ is shown in Fig.~\ref{fig:all_plot_MM}.
\subsection{Limits}
\subsubsection{The classical micromorphic model with symmetric forces stresses ($\mu_c \to 0$): nothing special}
\begin{figure}[H]
	\centering
	\includegraphics[height=5.5cm]{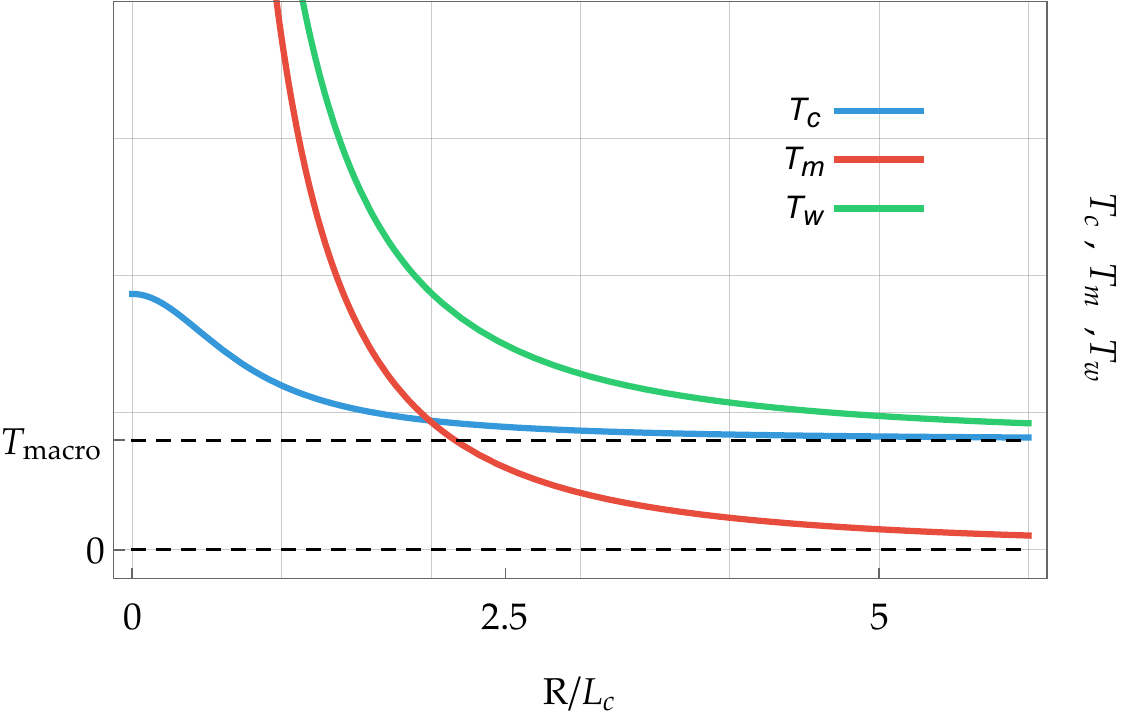}
	\caption{(\textbf{Micromorphic model})
		Torsional stiffness for the torque energy while varying $L_c$, for different values of $\mu_c=\{0, 1/30, 1/10, 1/5, 1, \infty\}$.
		The torsional stiffness remains bounded as $L_c \to \infty$ ($R\to 0$) and the model does not collapse in a linear elastic one.
		The values of the other parameters used are: $\mu = 1$, $\mu _e= 1/3$, $\mu _{\mbox{\tiny micro}}= 1/4$, $a_1= 2$, $a_3= 1/20$, $R= 1$.}
	\label{fig:all_plot_MM_mc_0}
\end{figure}
\subsubsection{The classical micromorphic model with reduced curvature energy ($a_2 = 0$)}
The classical micromorphic model with reduced curvature energy ($a_2 = 0$) collapses into the micro-strain model (Section \ref{sec:Micro_strain} with $\mbox{sym}P$) thus becoming independent with respect to the Cosserat couple modulus $\mu_c$ (see  (\ref{eq:diff_tors_stiff}) for the different stiffnesses expressions).

\begin{figure}[H]
	\centering
	\includegraphics[height=5.5cm]{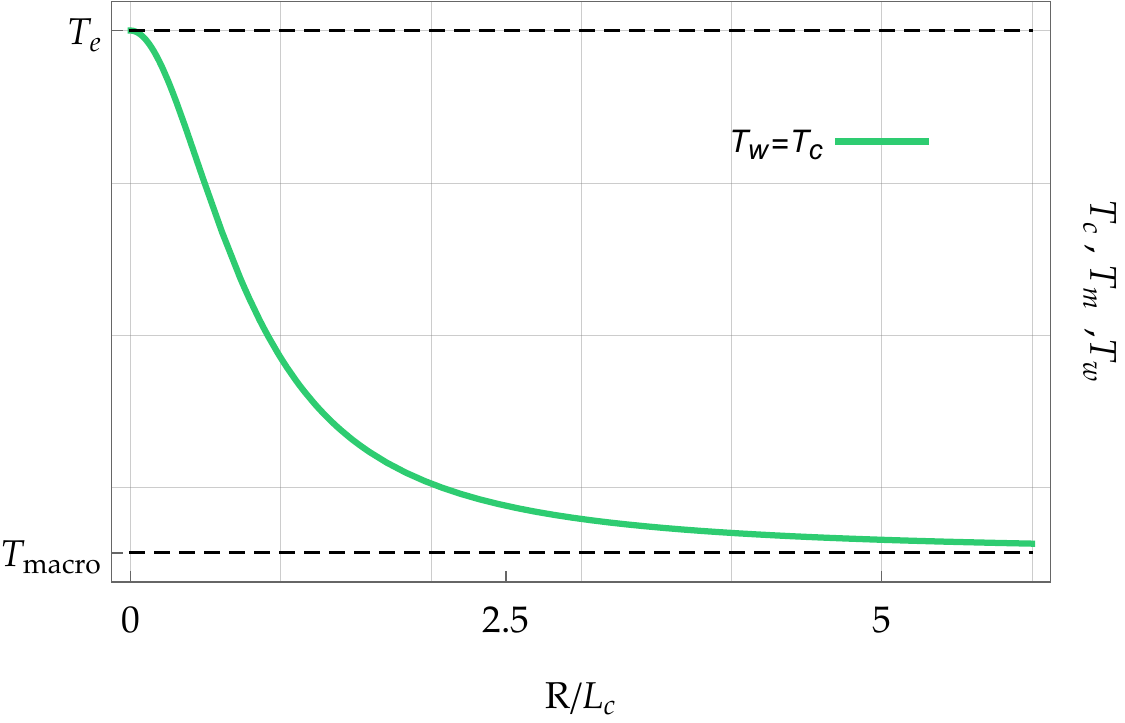}
	\caption{(\textbf{Micromorphic model})
		Torsional stiffness for the torque energy while varying $L_c$, for different values of $\mu_c=\{0, 1/30, 1/10, 1/5, 1, \infty\}$.
		The torsional stiffness remains bounded as $L_c \to \infty$ ($R\to 0$) and the model does not collapse in a linear elastic one.
		The values of the other parameters used are: $\mu = 1$, $\mu _e= 1/3$, $\mu _{\mbox{\tiny micro}}= 1/4$, $a_1= 2$, $a_3= 1/20$, $R= 1$.
		In this case, the stiffness for arbitrary small sample size is governed by $T_e$ and not $T_{\mbox{\tiny micro}}$.
		The reason for this is explained in Appendix \ref{app:const_disto}.
		}
	\label{fig:all_plot_MM_a3_0_mc_0}
\end{figure}
\section{Torsional problem for the micro-strain model without mixed terms}
\label{sec:Micro_strain}

The micro-strain model \cite{forest2006nonlinear,hutter2015micromorphic} is a particular case of the classical Mindlin-Eringen model, in which it is assumed a priori that the micro-distortion remains symmetric, $\boldsymbol{P}=\boldsymbol{S}\in \mbox{Sym}(3)$.
\footnote{
Shaat \cite{shaat2018reduced} uses the micro-strain model with mixed terms and a degenerate curvature expression in $\mbox{D}\boldsymbol{S}$, omitting $S_{11,1}$, $S_{22,2}$, and $S_{33,3}$.
}

A torsion solution for a more general case with mixed terms has been derived in \cite{hutter2016application}, but here we
employ a reduced isotropic curvature expression to make the calculations more manageable.

It is underlined that the micro-strain model cannot be obtained  as a limit case of  the relaxed micromorphic model and vice versa, although there are some similarities.
The strain energy which we consider is
\begin{align}
	W \left(\boldsymbol{\mbox{D}u}, \boldsymbol{S}, \boldsymbol{\mbox{D}S}\right) 
	=&
	\, \mu_{e} \left\lVert \mbox{dev} \left(\mbox{sym} \, \boldsymbol{\mbox{D}u} - \boldsymbol{S}\right) \right\rVert^{2}
	+ \frac{\kappa_{e}}{2} \mbox{tr}^2 \left(\boldsymbol{\mbox{D}u} - \boldsymbol{S} \right)
	+ \mu_{\tiny \mbox{micro}} \left\lVert \mbox{dev} \, \boldsymbol{S} \right\rVert^{2} 
	+ \frac{\kappa_{\tiny \mbox{micro}}}{2} \mbox{tr}^2 \left(\boldsymbol{S} \right)
	\label{eq:energy_MStrain}
	\\*
	&
	+ \frac{\mu \, L_c^2}{2} \,
	\left(
	a_1 \, \left\lVert \boldsymbol{\mbox{D}} \left(\mbox{dev} \, \boldsymbol{S}\right) \right\rVert^2
	+ \frac{2}{9} \, a_3 \, \left\lVert \mbox{D} \left(\mbox{tr} \left(\boldsymbol{S}\right)\boldsymbol{\mathbbm{1}} \right) \right\rVert^2
	\right)
	\, .
	\notag
\end{align}
The chosen 2-parameter curvature expression represents a simplified isotropic curvature (the full isotropic curvature for the micro-strain model would still counts 8 parameters \cite{barbagallo2016transparent}).
In this form, the micro-strain model can be obtained from the classical micromorphic model (Section  \ref{sec:Micro_morphic}), in general, by setting $\mu_c =0$ and $a_2 =0$.
For the torsion problem, the condition $a_2 =0$ alone is sufficient.

It is underlined that for the ansatz  (\ref{eq:ansatz_MStrain}), which will be presented later in this section, it holds $\mbox{tr} \left(\boldsymbol{S}\right)=0$.
The equilibrium equations, in the absence of body forces,   are therefore the following 
\begin{align}
	\mbox{Div}\overbrace{\left[
		2\mu_{e} \, \mbox{dev} \left(\mbox{sym} \, \boldsymbol{\mbox{D}u} - \boldsymbol{S}\right)
		+ \kappa_{e} \, \mbox{tr} \left(\boldsymbol{\mbox{D}u} - \boldsymbol{S} \right) \boldsymbol{\mathbbm{1}}
		\right]}^{\mathlarger{\widetilde{\sigma}}\coloneqq}
	= \boldsymbol{0},
	\notag
	\\*
	2\mu_{e} \, \mbox{dev} \left(\mbox{sym} \, \boldsymbol{\mbox{D}u} - \boldsymbol{S}\right)
	+ \kappa_{e} \, \mbox{tr} \left(\boldsymbol{\mbox{D}u} - \boldsymbol{S} \right) \boldsymbol{\mathbbm{1}}
	- 2 \mu_{\tiny \mbox{micro}} \, \mbox{dev} \,\boldsymbol{S}
	- \kappa_{\tiny \mbox{micro}} \, \mbox{tr} \left(\boldsymbol{S}\right) \boldsymbol{\mathbbm{1}} \, \, 
	\hspace{2.5cm}
	\label{eq:equi_MStrain}
	\\*
	+ \, \mu \, L_{c}^{2}\,
	\mbox{sym} \, \mbox{Div} \, 
	\left[
	a_1 \, \mbox{D} \left(\mbox{dev} \, \boldsymbol{S}\right)
	+ \frac{2}{9} \, a_3 \, \mbox{D} \left(\mbox{tr} \left(\boldsymbol{S}\right)\boldsymbol{\mathbbm{1}} \right)
	\right]
	= \boldsymbol{0} \, .
	\notag
\end{align}
\\
The boundary conditions at the external free surfaces are 
\begin{equation}
	\boldsymbol{\widetilde{t}}(r = R) = 
	\boldsymbol{\widetilde{\sigma}}(r) \cdot \boldsymbol{e}_{r} = 
	\boldsymbol{0}_{\mathbb{R}^3} \, ,
	\qquad\qquad
	\boldsymbol{\eta}(r = R) = 
	\mbox{sym}\left(\boldsymbol{\mathfrak{m}} (r) \cdot \boldsymbol{e}_{r}\right) =
	\boldsymbol{0}_{\mathbb{R}^{3\times 3}}
	\, ,
	\label{eq:BC_MStrain}
\end{equation}
where
\begin{equation}
    \boldsymbol{\mathfrak{m}} = \mu \, L_{c}^{2}\,
    \left[
    a_1 \, \mbox{D} \left(\mbox{dev} \, \boldsymbol{S}\right)
    + \frac{2}{9} \, a_3 \, \mbox{D} \left(\mbox{tr} \left(\boldsymbol{S}\right)\boldsymbol{\mathbbm{1}} \right)
    \right]
    \end{equation}
is the third order moment stress tensor, the expression of $\boldsymbol{\widetilde{\sigma}}$ is in  (\ref{eq:equi_MStrain}), $\boldsymbol{e}_{r}$ is the radial unit vector.
According with the reference system shown in Fig.~\ref{fig:intro_1}, the ansatz for the displacement field and the micro-distortion is
\begin{align}
	\boldsymbol{u}(r,\varphi,z) &= \boldsymbol{\vartheta}
	\left(
	\begin{array}{c}
		-x_2(r,\varphi) \, x_3(z) \\
		x_1(r,\varphi) \, x_3(z) \\
		0 
	\end{array}
	\right)
	\, ,
	\label{eq:ansatz_MStrain}
	\\*
	\boldsymbol{S}(r,\varphi,z) &= \frac{\boldsymbol{\vartheta}}{2}
	\left(
	\begin{array}{ccc}
		0   & 0 & g_{m}(r) \, x_2(r,\varphi) \\
		0 &    0 &  - g_{m}(r) \, x_1(r,\varphi) \\
		g_{m}(r) \, x_2(r,\varphi) & - g_{m}(r) \, x_1(r,\varphi) & 0 \\
	\end{array}
	\right)
	\, ,
	\notag
\end{align}
where, in relation to the ansatz  (\ref{eq:ansatz_RM}), $g_{m}(r ) \coloneqq g_{1}(r ) - g_{2}(r )$.
Substituting the ansatz  (\ref{eq:ansatz_MStrain}) in  (\ref{eq:equi_MStrain}) the 9 equilibrium equation are equivalent to
\begin{align}
	\frac{1}{2} \boldsymbol{\vartheta}  \sin \varphi
	\left(
	a_{1} \, \mu \, L_c^2 \left(3 g_{m}'(r ) + r \, g_{m}''(r )\right)
	-2 r \, \mu _e \, (g_{m}(r ) + 1)
	-2 r \, g_{m}(r ) \, \mu _{\mbox{\tiny micro}}
	\right)
	&= 0 \, ,
	\label{eq:equi_equa_MStrain}
	\\*
	-\frac{1}{2} \boldsymbol{\vartheta}  \cos \varphi
	\left(
	a_{1} \, \mu \, L_c^2 \left(3 g_{m}'(r ) + r \, g_{m}''(r )\right)
	-2 r \, \mu _e (g_{m}(r ) + 1)
	-2 r \, g_{m}(r ) \, \mu _{\mbox{\tiny micro}}
	\right)
	&= 0 \, .
	\notag
\end{align}
Between the two equilibrium equations  (\ref{eq:equi_equa_MStrain}) there is only one independent equation since
\\
 (\ref{eq:equi_equa_MStrain})$_1$ = $-\tan \varphi$  (\ref{eq:equi_equa_MStrain})$_2$.
The solution of  (\ref{eq:equi_equa_MStrain}) is 
\begin{align}
	g_{m}(r) = \,
	\frac{A_{2} Y_1\left(-\frac{i r  f_1}{L_c}\right)-i A_{1} I_1\left(\frac{r  f_1}{L_c}\right)}{r}
	-\frac{\mu _e}{\mu _e+\mu _{\mbox{\tiny micro}}}
	\, ,
	\qquad\qquad\qquad
	f_{1} \coloneqq
	\sqrt{\frac{2(\mu _e+\mu _{\mbox{\tiny micro}})}{a_{1} \mu}} \, ,
	\label{eq:sol_fun_MStrain}
\end{align}
where $I_{n}\left(\cdot\right)$ is the \textit{modified Bessel function of the first kind}, $Y_{n}\left(\cdot\right)$ is the \textit{Bessel function of the second kind} (see appendix \ref{app:bessel} for the formal definitions), and
$A_1$, $A_2$ are integration constants.

The value of $A_1$ is determined thanks to the boundary conditions  (\ref{eq:BC_MStrain}), while, due to the divergent behaviour of the Bessel function of the second kind at $r=0$, we have to set $A_2=0$ in order to have a continuous solution.
The fulfilment of the boundary conditions  (\ref{eq:BC_MStrain}) allows us to find the expressions of the integration constants
\begin{align}
	A_1 = \, 
	\frac{2 i \, L_c }{I_0\left(\frac{R f_1}{L_c}\right)+I_2\left(\frac{R f_1}{L_c}\right)}
	\frac{\mu _e}{f_1 (\mu _e+\mu _{\mbox{\tiny micro}})}
	\, .
	\label{eq:BC_MStrain_3}
\end{align}
The classical torque, the higher-order torque, and the energy (per unit length d$z$) expressions are
\begin{align}
	M_{\mbox{c}} (\boldsymbol{\vartheta}) 
	\coloneqq&
	\int_{0}^{2\pi}
	\int_{0}^{R}
	\Big[
	\langle
	\widetilde{\boldsymbol{\sigma}} \, \boldsymbol{e}_{z} , \boldsymbol{e}_{\varphi}
	\rangle
	r
	\Big] r
	\, \mbox{d}r \, \mbox{d}\varphi
	\notag
	\\*
	=&
	\left[
	\frac{\mu _e \, \mu _{\mbox{\tiny micro}}}{\mu _e+\mu _{\mbox{\tiny micro}}}
	+
	\frac{\mu _e^2 \, \mu \, a_{1}}{\left(\mu _e+\mu _{\mbox{\tiny micro}}\right)^2}
	\frac{
	4 \, I_2\left(\frac{R \, f_1}{L_c}\right)
	}{
	I_0\left(\frac{R \, f_1}{L_c}\right)
	+ I_2\left(\frac{R \, f_1}{L_c}\right) 
	}
	\frac{L_c^2}{R^2}
	\right]
	I_{p} \, 
	\boldsymbol{\vartheta}
	= T_{\mbox{c}} \, \boldsymbol{\vartheta} \, ,
	\notag
	\\[3mm]
	M_{\mbox{m}}(\boldsymbol{\vartheta}) 
	\coloneqq&
	\int_{0}^{2\pi}
	\int_{0}^{R}
	\Big[
	\langle
	\mbox{sym} (\boldsymbol{\mathfrak{m}} \boldsymbol{e}_{z})
	\boldsymbol{e}_{\varphi} ,
	\boldsymbol{e}_{r}
	\rangle
	-
	\langle
	\mbox{sym} (\boldsymbol{\mathfrak{m}} \, \boldsymbol{e}_{z})
	\boldsymbol{e}_{r} ,
	\boldsymbol{e}_{\varphi}
	\rangle
	\Big]
	\, r
	\, \mbox{d}r \, \mbox{d}\varphi	
	=
	0 \, ,
	\label{eq:torque_stiffness_MStrain}
	\\[3mm]
	W_{\mbox{tot}} (\boldsymbol{\vartheta}) 
	\coloneqq&
	\int_{0}^{2\pi}
	\int_{0}^{R}
	W \left(\boldsymbol{\mbox{D}u}, \boldsymbol{S}, \mbox{D}\boldsymbol{S}\right) \, \, r
	\, \mbox{d}r \, \mbox{d}\varphi
	\notag
	\\*
	=&
	\frac{1}{2}
	\left[
	\underbrace{
	\frac{\mu _e \, \mu _{\mbox{\tiny micro}}}{\mu _e+\mu _{\mbox{\tiny micro}}}
	}_{\mu_{\tiny \mbox{macro}}}
	+
	\frac{\mu _e^2 \, \mu \, a_{1}}{\left(\mu _e+\mu _{\mbox{\tiny micro}}\right)^2}
	\frac{
		4 \, I_2\left(\frac{R \, f_1}{L_c}\right)
	}{
		I_0\left(\frac{R \, f_1}{L_c}\right)
		+ I_2\left(\frac{R \, f_1}{L_c}\right) 
	}
	\frac{L_c^2}{R^2}
	\right]
	I_{p} \, 
	\boldsymbol{\vartheta}^2
	= \frac{1}{2} \, T_{\mbox{w}} \, \boldsymbol{\vartheta}^2
	\, .
	\notag
\end{align}

\begin{figure}[H]
	\centering
	\includegraphics[height=5.5cm]{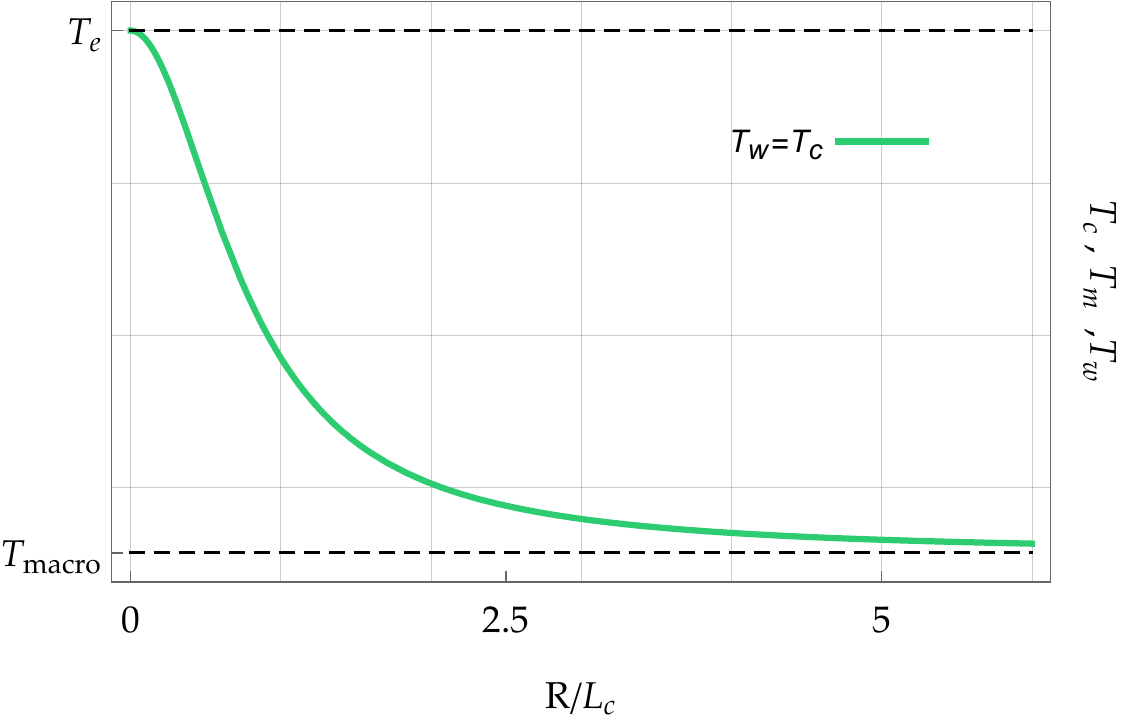}
	\caption{(\textbf{Micro-strain model}) Torque energy $T_{\mbox{w}}$ while varying $L_c$. Observe that the torsional stiffness remains bounded as $L_c \to \infty$ ($R \to 0$). The values of the parameters used are: $\mu _{e} = 1/3$, $\mu _{\tiny \mbox{micro}} = 1/4$, $\mu = 1$, $a _1 = 1/5$.
	In this case, the stiffness for arbitrary small sample size is governed by $T_e$ and not $T_{\mbox{\tiny micro}}$.
	}
	\label{fig:all_plot_MStrain}
\end{figure}

The plot of the torsional stiffness for the classical torque, the higher-order torque, and the torque energy while varying $L_c$ is shown in Fig.~\ref{fig:all_plot_Cos_3}.
Since the higher-order torque is zero and the following relation holds
\begin{equation}
\frac{\mbox{d}}{\mbox{d}\boldsymbol{\vartheta}}W_{\mbox{tot}}(\boldsymbol{\vartheta}) = M_{\mbox{c}} (\boldsymbol{\vartheta}) + M_{\mbox{m}} (\boldsymbol{\vartheta}) = M_{\mbox{c}} (\boldsymbol{\vartheta})
\, ,
\qquad\qquad\qquad
\frac{\mbox{d}^2}{\mbox{d}\boldsymbol{\vartheta}^2}W_{\mbox{tot}}(\boldsymbol{\vartheta})
= T_{\mbox{c}} + T_{\mbox{m}}
= T_{\mbox{c}}
= T_{\mbox{w}}
\, ,
\end{equation}
only the plot of the energy (per unit length d$z$) while changing $L_c$ is shown in Fig.~\ref{fig:all_plot_MStrain}
The energy of the model remains bounded, as for the shear and bending problem \cite{rizzi2021shear,rizzi2021bending}, since for both problems the higher-order moments are zero, and this does not create a conflict with the boundary condition as $L_c \to \infty$ (see  (\ref{eq:diff_tors_stiff}) for the different stiffnesses expressions).
Note that there is no simple way to a-priori guess that the small size torsional response is given by $T_e$ since $\boldsymbol{S} \in \text{Sym}(3)$ is not easily seen to be zero.
In Appendix \ref{app:const_disto} we show that the variational limit for $L_c\to\infty$ is indeed realized by $\boldsymbol{S} (\boldsymbol{x})=\boldsymbol{\overline{S}}=\boldsymbol{0}$ and this shows that the limit stiffness is given by $T_e$.

\section{Torsional problem for the second gradient continuum}
The expression of the most general isotropic strain energy for the second gradient continuum is \cite{mindlin1964micro,neff2009subgrid}
\begin{align}
	W \left(\boldsymbol{\mbox{D}u}, \boldsymbol{\mbox{D}^2 u}\right)
	= &
	\, \mu_{\tiny \mbox{macro}} \left\lVert \mbox{sym} \, \boldsymbol{\mbox{D}u} \right\lVert^2
	+ \frac{\lambda_{\tiny \mbox{macro}}}{2} \, \mbox{tr}^2 \left(\boldsymbol{\mbox{D}u}\right)
	\label{eq:energy_Strain_Grad_gen}
	\\*
	&
	+ \widehat{a}_1 \, \chi_{iik} \, \chi_{kjj}
	+ \widehat{a}_2 \, \chi_{ijj} \, \chi_{ikk}
	+ \widehat{a}_3 \, \chi_{iik} \, \chi_{jjk}
	+ \widehat{a}_4 \, \chi_{ijk} \, \chi_{ijk}
	+ \widehat{a}_5 \, \chi_{ijk} \, \chi_{kji}
	\, ,
	\notag
\end{align}
where $\boldsymbol{\chi} = \boldsymbol{\mbox{D}}^2 \boldsymbol{u}$ ($\chi_{ijk} = \frac{\partial^2 u_k}{\partial x_i \, \partial x_j}$).
The expression we are going to use in the following is a simplified isotropic strain energy with three curvature parameters
\begin{align}
	W \left(\boldsymbol{\mbox{D}u}, \boldsymbol{\mbox{D}^2 u}\right)
	= &
	\, \mu_{\tiny \mbox{macro}} \left\lVert \mbox{sym}\,\boldsymbol{\mbox{D}u} \right\rVert^{2}
	+ \frac{\lambda_{\tiny \mbox{macro}}}{2} \mbox{tr}^2 \left(\boldsymbol{\mbox{D}u} \right)
	\label{eq:energy_Strain_Grad}
	\\*
	&
	+ \frac{\mu \, L_c^2}{2}
	\left(
	a_1 \, \left\lVert \mbox{D} \Big(\mbox{dev} \, \mbox{sym} \, \boldsymbol{\mbox{D} u}\Big) \right\rVert^2
	+ a_2 \, \left\lVert \mbox{D} \Big(\mbox{skew} \, \boldsymbol{\mbox{D} u}\Big) \right\rVert^2
	+  \frac{2}{9} \, a_3 \, \left\lVert \mbox{D}
	\Big( 
	\mbox{tr} \left(\boldsymbol{\mbox{D} u}\right) \, \boldsymbol{\mathbbm{1}}
	\Big) \right\rVert^2
	\right)
	\, .
	\notag
\end{align}
The equilibrium equation, in the absence of body forces, is
\begin{align}
	\mbox{Div}\bigg[
	2 \mu_{\mbox{\tiny macro}} \,\mbox{sym}\,\boldsymbol{\mbox{D}u}
	+ \lambda_{\tiny \mbox{macro}} \mbox{tr} \left(\boldsymbol{\mbox{D}u}\right) \boldsymbol{\mathbbm{1}}
	\hspace{8cm}
	\label{eq:equi_Strain_Grad}
	\\*
	- \mu L_{c}^{2} \,
	\left(
	a_1 \, \mbox{dev} \, \mbox{sym} \, \boldsymbol{\Delta} \left(\boldsymbol{\mbox{D}u}\right)
	+ a_2 \, \mbox{skew} \, \boldsymbol{\Delta} \left(\boldsymbol{\mbox{D}u}\right)
	+ \frac{2}{9} \, a_3 \, \mbox{tr} \left(\boldsymbol{\Delta} \left(\boldsymbol{\mbox{D}u}\right)\right)\boldsymbol{\mathbbm{1}}
	\right)
	\bigg]
	= \boldsymbol{0} \, ,
	\notag
\end{align}
where $\boldsymbol{\Delta} \left(\boldsymbol{\mbox{D}u}\right) \in \mathbb{R}^{3\times3}$ is taken component-wise.
The non-trivial boundary conditions at the free surface are
\begin{align}
	\boldsymbol{\widetilde{t}}(r = R) =
	& \, 
	\widetilde{\boldsymbol{\sigma}} \, \boldsymbol{e}_{r} 
	+ \left[ \left( \boldsymbol{e}_{r} \otimes \boldsymbol{e}_{r} \right) \boldsymbol{:} \boldsymbol{\nabla} \boldsymbol{\mathfrak{m}}\right] \boldsymbol{e}_{r}
	- 2 \left[\left(\boldsymbol{\mathbbm{1}} - \boldsymbol{e}_{r} \otimes \boldsymbol{e}_{r} \right)  \boldsymbol{:} \boldsymbol{\nabla} \boldsymbol{\mathfrak{m}}\right] \, \boldsymbol{e}_{r}
    \label{eq:BC_SG}
	\\*
	& \hspace{2cm}
	+ \left(\left[\left(\boldsymbol{\mathbbm{1}} - \boldsymbol{e}_{r} \otimes \boldsymbol{e}_{r} \right)  \boldsymbol{:} \boldsymbol{\nabla} \boldsymbol{e}_{r}\right]
	\left( \boldsymbol{e}_{r} \otimes \boldsymbol{e}_{r} \right)
	- \left[\left(\boldsymbol{\mathbbm{1}} - \boldsymbol{e}_{r} \otimes \boldsymbol{e}_{r} \right)   \left( \boldsymbol{\nabla} \boldsymbol{e}_{r} \right)^T \right] \right) \boldsymbol{:} \boldsymbol{\mathfrak{m}}
	= \boldsymbol{0}_{\mathbb{R}^3} \, ,
	\notag
	\\
    \boldsymbol{\widetilde{\eta}}(r = R) =
    &
    \,
    \left( \boldsymbol{e}_{r} \otimes \boldsymbol{e}_{r} \right) \boldsymbol{:} \boldsymbol{\mathfrak{m}}
	= \boldsymbol{0}_{\mathbb{R}^{3}} \, ,
	\footnotemark
    \notag
\end{align}
\footnotetext{
In index notation $\left( \boldsymbol{\mathbbm{1}} - \boldsymbol{e}_{r} \otimes \boldsymbol{e}_{r} \right) \boldsymbol{:} \boldsymbol{\nabla} \mathfrak{\boldsymbol{m}} = (\delta_{ip} - n_i n_p) \mathfrak{m}_{ijk,p}$.
}
where, since the boundary surface is smooth, one set of boundary condition is identically satisfied (see \cite{mindlin1964micro,madeo2016new} for all the details).
According to the reference system shown in Fig.~\ref{fig:intro_1}, the ansatz for the displacement field and consequently the gradient of the displacement are
\begin{equation}
	\boldsymbol{u}(x_1,x_2)= \boldsymbol{\vartheta}
	\left(
	\begin{array}{c}
		-x_2 \, x_3 \\
		x_1 \, x_3 \\
		0 
	\end{array}
	\right) \, 
	\qquad \Rightarrow \qquad 
	\boldsymbol{\mbox{D}u} = \frac{\boldsymbol{\vartheta}}{2} 
	\left(
	\begin{array}{ccc}
		0 & -2 x_{3} & -2x_{2} \\
		2 x_{3} & 0 & 2x_{1} \\
		0 & 0 & 0 \\
	\end{array}
	\right) \, .
	\label{eq:ansatz_Strain_Grad}
\end{equation}
Since the ansatz is completely known, it is possible to check that the equilibrium equation  (\ref{eq:equi_Strain_Grad}) and the boundary conditions  (\ref{eq:ansatz_Strain_Grad}) are identically satisfied and it is possible to evaluate directly the classical torque, the higher-order torque, and the energy.

\begin{figure}[H]
	\centering
	\includegraphics[height=5.5cm]{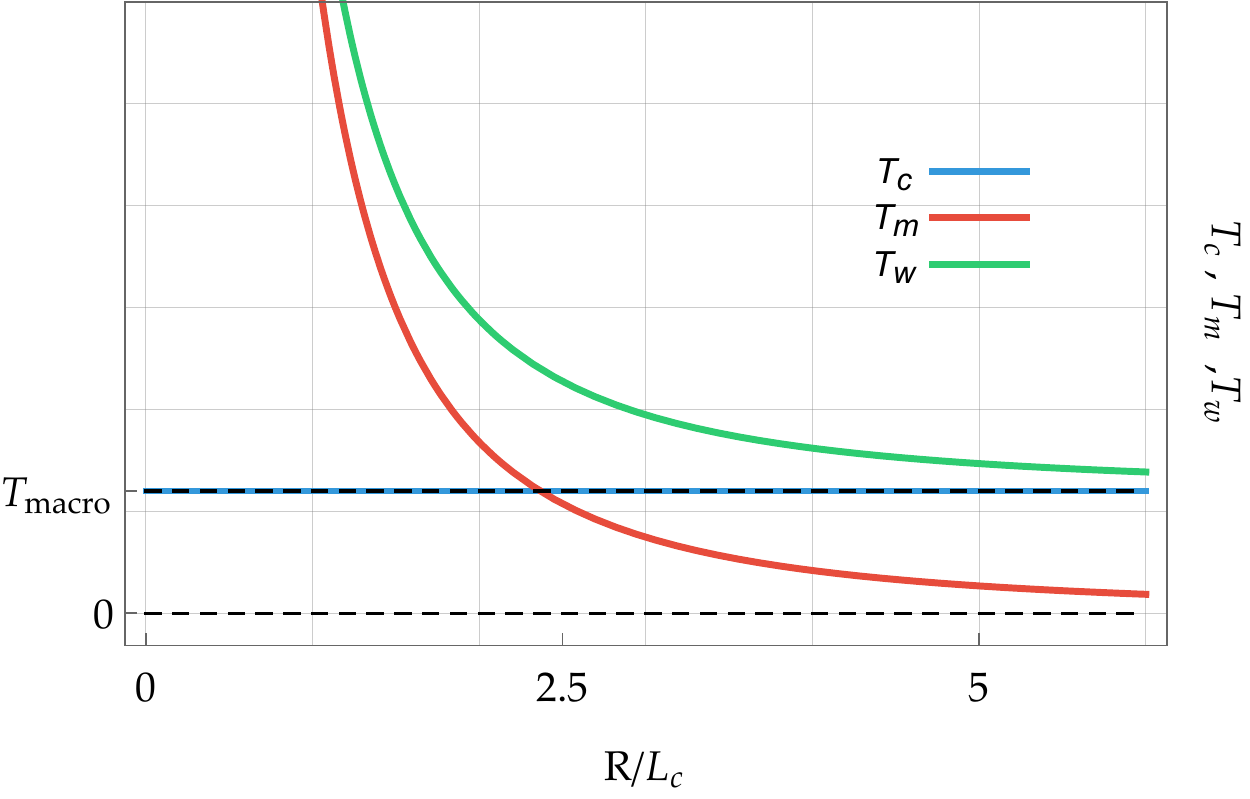}
	\caption{(\textbf{Second gradient model}) Torsional stiffness for the classical torque $T_{\mbox{c}}$, the higher-order torque $T_{\mbox{m}}$, and the torque energy $T_{\mbox{w}}$ while varying $L_c$. The torsional stiffness is unbounded as $L_c \to \infty$ ($R\to 0$). The values of the parameters used are: $\mu = 1$, $\mu _{\mbox{\tiny macro}}= 1/4$, $a_1= 1/5$, $a_3= 1/6$, $R= 1$.}
	\label{fig:all_plot_Strain_Grad_3}
\end{figure}

The classical torque, the higher-order torque, and energy (per unit length d$z$) expressions are
\begin{align}
	M_{\mbox{c}} (\boldsymbol{\vartheta}) 
	\coloneqq&
	\int_{0}^{2\pi}
	\int_{0}^{R}
	\Big[
	\langle
	\widetilde{\boldsymbol{\sigma}} \, \boldsymbol{e}_{z} , \boldsymbol{e}_{\varphi}
	\rangle
	r
	\Big] r
	\, \mbox{d}r \, \mbox{d}\varphi
	=
	\mu_{\mbox{\tiny macro}}
	I_{p} \, 
	\boldsymbol{\vartheta}
	= T_{\mbox{c}} \, \boldsymbol{\vartheta} \, ,
	\notag
	\\*[3mm]
	M_{\mbox{m}}(\boldsymbol{\vartheta}) 
	\coloneqq&
	\int_{0}^{2\pi}
	\int_{0}^{R}
	\Big[
	\langle
	\left(\boldsymbol{\mathfrak{m}} \, \boldsymbol{e}_{z}\right)
	\boldsymbol{e}_{r} ,
	\boldsymbol{e}_{\varphi}
	\rangle
	-
	\langle
	\left(\boldsymbol{\mathfrak{m}} \, \boldsymbol{e}_{z}\right)
	\boldsymbol{e}_{\varphi} ,
	\boldsymbol{e}_{r}
	\rangle
	+
	\langle
	\left(\boldsymbol{\mathfrak{m}} \, \boldsymbol{e}_{r}\right)
	\boldsymbol{e}_{z} ,
	\boldsymbol{e}_{\varphi}
	\rangle
	-
	\langle
	\left(\boldsymbol{\mathfrak{m}} \, \boldsymbol{e}_{\varphi}\right)
	\boldsymbol{e}_{z} ,
	\boldsymbol{e}_{r}
	\rangle
	\Big]
	\, r
	\, \mbox{d}r \, \mbox{d}\varphi 
	\notag
	\\*
	=& \,
	2 \mu (a_{1}+3 a_{2}) \frac{L_c^2}{R^2} \, 
	I_{p} \, 
	\boldsymbol{\vartheta}
	= T_{\mbox{m}} \, \boldsymbol{\vartheta} \, ,
	\notag
	\\*[3mm]
	W_{\mbox{tot}} (\boldsymbol{\vartheta}) 
	\coloneqq&
	\int_{0}^{2\pi}
	\int_{0}^{R}
	W \left(\boldsymbol{\mbox{D}u}, \boldsymbol{\mbox{D}^2 u}\right) \, \, r
	\, \mbox{d}r \, \mbox{d}\varphi
	=
	\frac{1}{2}
	\left[
	\mu_{\mbox{\tiny macro}}
	+ 2 \mu (a_{1}+3 a_{2}) \frac{L_c^2}{R^2}
	\right]
	I_{p} \, 
	\boldsymbol{\vartheta}^2
	= \frac{1}{2} \, T_{\mbox{w}} \, \boldsymbol{\vartheta}^2
	\, .
	\label{eq:torque_stiffness_Strain_Grad}
\end{align}
The plot of the torsional stiffness for the classical torque (light blue), the higher-order torque (red), and the torque energy (green) while varying $L_c$ is shown in Fig.~\ref{fig:all_plot_Strain_Grad_3}.
\subsection{The strain gradient continuum as a limit of the micro-strain model}
If we let $\mu_e, \kappa_e \to \infty$ in the micro-strain model, we obtain in the limit a strain gradient model with elastic energy
\begin{align}
	W \left(\boldsymbol{\mbox{D}u}, \boldsymbol{\mbox{D} \, \mbox{sym} \, \mbox{D} u}\right)
	= &
	\, \mu_{\tiny \mbox{macro}} \left\lVert \mbox{sym}\,\boldsymbol{\mbox{D}u} \right\rVert^{2}
	+ \frac{\lambda_{\tiny \mbox{macro}}}{2} \mbox{tr}^2 \left(\boldsymbol{\mbox{D}u} \right)
	\label{eq:energy_Strain_Grad_z}
	\\*
	&
	+ \frac{\mu \, L_c^2}{2}
	\left(
	a_1 \, \left\lVert \mbox{D} \Big(\mbox{dev} \, \mbox{sym} \, \boldsymbol{\mbox{D} u}\Big) \right\rVert^2
	+a_3 \, \left\lVert \mbox{D} \Big(\mbox{tr} \left( \boldsymbol{\mbox{D} u}
	\right)
	\boldsymbol{\mathbbm{1}}
	\Big) \right\rVert^2
	\right)
	\, .
	\notag
\end{align}
Since $\mbox{tr} \left( \boldsymbol{\mbox{D} u} \right) = 0$ for our ansatz  (\ref{eq:ansatz_disp}), the equilibrium equations, in the absence of body forces,   are
\begin{align}
	\mbox{Div}\bigg[
	2 \mu_{\mbox{\tiny macro}} \,\mbox{sym}\,\boldsymbol{\mbox{D}u}
	+ \lambda_{\tiny \mbox{macro}} \mbox{tr} \left(\boldsymbol{\mbox{D}u}\right) \boldsymbol{\mathbbm{1}}
	- \mu L_{c}^{2} \,
	a_1 \, \mbox{dev} \, \mbox{sym} \, \boldsymbol{\Delta} \left(\boldsymbol{\mbox{D}u}\right)
	\bigg]
	= \boldsymbol{0} \, ,
	\label{eq:equi_Strain_Grad_z}
\end{align}
where $\boldsymbol{\Delta} \left(\boldsymbol{\mbox{D}u}\right) \in \mathbb{R}^{3\times3}$ is taken component-wise.
\begin{figure}[H]
	\centering
	\includegraphics[height=5.5cm]{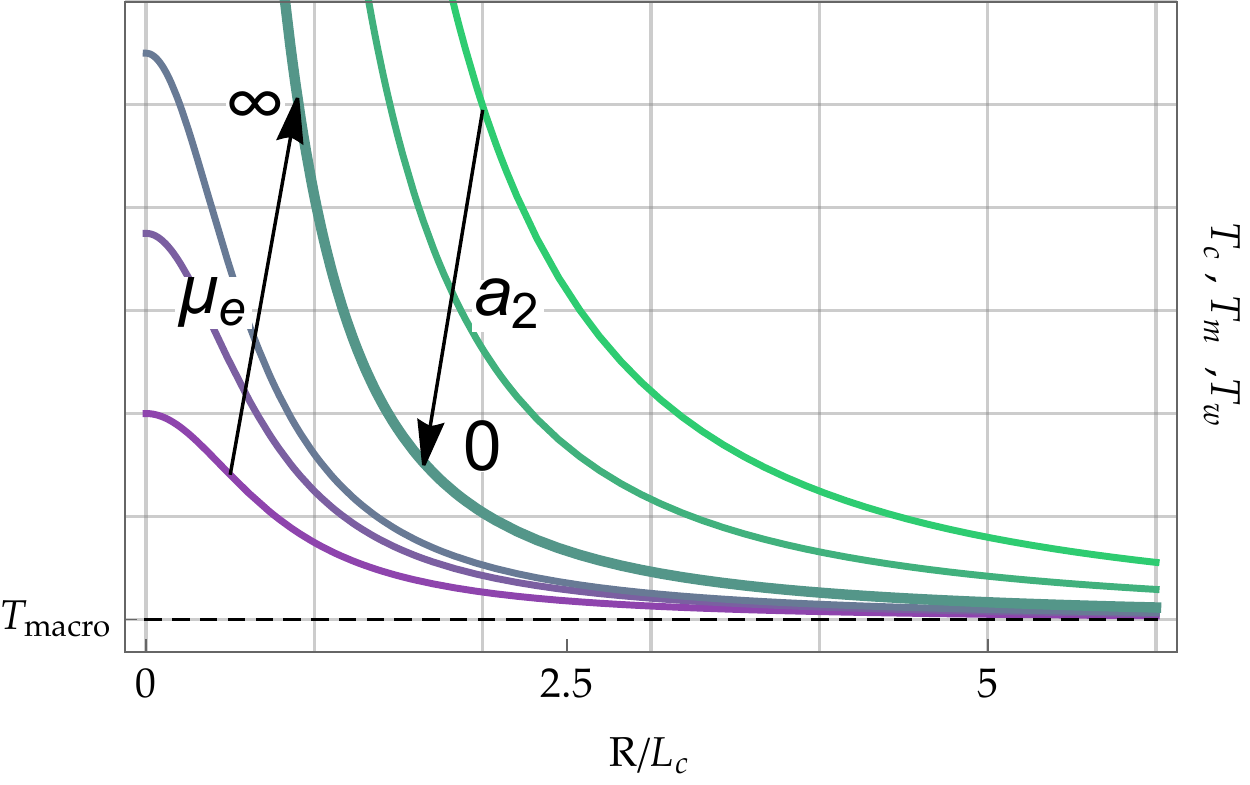}
	\caption{The purplish curves show how the micro-strain model particularises to the strain gradient model for $\mu_e\to\infty$ (the following set has been used $\mu_e=\{ 1/3,1/2,2/3,\infty\}$).
	The greenish curves show how the choice of $a_2=0$ guarantees the formal equivalence (which is however always substantially true) between the second gradient model and the strain gradient model \cite{mindlin1964micro} (the following set has been used $a_2=\{1/4,1/10,0\}$).
	The values of the other parameters used are: $\mu = 1$, $\mu _{\mbox{\tiny macro}}= 1/7$, $a_1= 1/5$, $R= 1$.}
	\label{fig:all_plot_Strain_Grad_Macro_strain}
\end{figure}
\section{Ad-hoc model containing Cosserat and micro-strain effects}
Given $\boldsymbol{S} \in \mbox{Sym}(3)$ and $\boldsymbol{A} \in \mathfrak{so}(3)$, the strain energy which we consider now is
\begin{align}
	W \left(\boldsymbol{\mbox{D}u}, \boldsymbol{A}, \boldsymbol{S}, \text{Curl} \, \boldsymbol{A}, \boldsymbol{\mbox{D}S}\right) 
	=&
	\, \mu_{e} \left\lVert \mbox{sym} \, \boldsymbol{\mbox{D}u} - \boldsymbol{S} \right\rVert^{2}
	+ \frac{\lambda_{e}}{2} \mbox{tr}^2 \left(\boldsymbol{\mbox{D}u} - \boldsymbol{S} \right)
	+ \mu_{c} \left\lVert \mbox{skew} \left(\boldsymbol{\mbox{D}u} - \boldsymbol{A} \right) \right\rVert^{2}
	\notag
	\\*
	&
	+ \mu_{\tiny \mbox{micro}} \left\lVert \mbox{dev} \, \boldsymbol{S} \right\rVert^{2} 
	+ \frac{\kappa_{\tiny \mbox{micro}}}{2} \mbox{tr}^2 \left(\boldsymbol{S} \right)
	\label{eq:energy_AD}
	\\*
	&
	+ \frac{\mu \, L_c^2}{2}
	\left(
	a_1 \, \left \lVert \mbox{dev} \, \mbox{sym} \, \mbox{Curl} \, \boldsymbol{A}\right \rVert^2 \, 
	+ \frac{a_3}{3} \, \mbox{tr}^2 \left(\mbox{Curl} \, \boldsymbol{A} \right)
	+ a_4 \, \left\lVert \boldsymbol{\mbox{D}} \left(\mbox{dev} \, \boldsymbol{S}\right) \right\rVert^2
	\right) \, ,
	\notag
\end{align}
since $\left \lVert \mbox{skew} \, \mbox{Curl} \, \boldsymbol{A}\right \rVert^2 =  \left\lVert \boldsymbol{\mbox{D}} \left(\mbox{skew} \, \boldsymbol{S}\right) \right\rVert^2 = \left\lVert \mbox{D} \left(\mbox{tr} \left(\boldsymbol{S}\right)\boldsymbol{\mathbbm{1}} \right) \right\rVert^2 = 0$ in terms of the ansatz  (\ref{eq:ansatz_AD}).

The equilibrium equations, in the absence of body forces,  are the following 
\begin{align}
	\mbox{Div}\overbrace{\left[
		2\mu_{e} \, \left(\mbox{sym} \, \boldsymbol{\mbox{D}u} - \boldsymbol{S}\right)
		+ \lambda_{e} \, \mbox{tr} \left(\boldsymbol{\mbox{D}u} - \boldsymbol{S} \right) \boldsymbol{\mathbbm{1}}
		+2\mu_{c} \, \left(\mbox{skew} \, \boldsymbol{\mbox{D}u} - \boldsymbol{A}\right)
		\right]}^{\mathlarger{\widetilde{\sigma}}\coloneqq}
	&= \boldsymbol{0},
	\notag
	\\*
	2\mu_{c}\,\mbox{skew} \left(\boldsymbol{\mbox{D}u} - \boldsymbol{A}\right)
	-\mu \, L_c^2 \, \mbox{skew} \, \mbox{Curl}\,
	\left(
	a_1 \, \mbox{dev} \, \mbox{sym} \, \mbox{Curl} \, \boldsymbol{A} \, 
	+ \frac{a_3}{3} \, \mbox{tr} \left(\mbox{Curl} \, \boldsymbol{A} \right)\boldsymbol{\mathbbm{1}} \, 
	\right)
	&= \boldsymbol{0}
	\label{eq:equi_AD}
	\\*
	2\mu_{e} \, \left(\mbox{sym} \, \boldsymbol{\mbox{D}u} - \boldsymbol{S}\right)
	+ \lambda_{e} \, \mbox{tr} \left(\boldsymbol{\mbox{D}u} - \boldsymbol{S} \right) \boldsymbol{\mathbbm{1}}
	- 2 \mu_{\tiny \mbox{micro}} \, \boldsymbol{S}
	- \lambda_{\tiny \mbox{micro}} \, \mbox{tr} \left(\boldsymbol{S}\right) \boldsymbol{\mathbbm{1}}
	+ \, \mu \, L_{c}^{2}\,
	a_4 \, \mbox{sym} \, \boldsymbol{\Delta} \left(\mbox{dev} \, \boldsymbol{S}\right)
	&= \boldsymbol{0} \, .
	\notag
\end{align}
The boundary conditions at the external free surfaces are 
\begin{align}
	\boldsymbol{\widetilde{t}}(r = R) &= 
	\boldsymbol{\widetilde{\sigma}}(r) \cdot \boldsymbol{e}_{r} = 
	\boldsymbol{0}_{\mathbb{R}^3} \, ,
	\label{eq:BC_AD}
	\\*
	\boldsymbol{\eta}_b(r = R) &= 
	\mbox{skew}\left(\boldsymbol{m} (r) \cdot \boldsymbol{\epsilon} \cdot \boldsymbol{e}_{r}\right) = 
	\mbox{skew}\left(\boldsymbol{m} (r) \times \boldsymbol{e}_{r}\right) = 
	\boldsymbol{0}_{\mathbb{R}^{3\times 3}}
	\, ,
	\notag
	\\*
	\boldsymbol{\eta}_a(r = R) &= 
	\mbox{sym}\left(\boldsymbol{\mathfrak{m}} (r) \cdot \boldsymbol{e}_{r}\right) =
	\boldsymbol{0}_{\mathbb{R}^{3\times 3}}
	\, ,
	\notag
\end{align}
where
\begin{align}
\boldsymbol{m} &=
\mu \, L_c^2
\left(
a_1 \, \mbox{dev} \, \mbox{sym} \, \mbox{Curl} \, \boldsymbol{A} \, 
+ \frac{a_3}{3} \, \mbox{tr} \left(\mbox{Curl} \, \boldsymbol{A} \right)\boldsymbol{\mathbbm{1}} \, 
\right)
\, ,
\\*
\boldsymbol{\mathfrak{m}} &=
\mu \, L_{c}^{2}\,
a_4 \, \mbox{D} \left(\mbox{dev} \, \boldsymbol{S}\right)
\, ,
\notag
\end{align}
is the second and third order moment stress tensor respectively, the expression of $\boldsymbol{\widetilde{\sigma}}$ is in  (\ref{eq:equi_AD}), $\boldsymbol{e}_{r}$ is the radial unit vector.
According with the reference system shown in Fig.~\ref{fig:intro_1}, the ansatz for the displacement field and the micro-distortion is
\begin{align}
	\boldsymbol{u}(r,\varphi,z) &= \boldsymbol{\vartheta}
	\left(
	\begin{array}{c}
		-x_2(r,\varphi) \, x_3(z) \\
		x_1(r,\varphi) \, x_3(z) \\
		0 
	\end{array}
	\right)
	\, ,
	\notag
	\\*
	\boldsymbol{A}(r,\varphi,z) &= \frac{\boldsymbol{\vartheta}}{2}
	\left(
	\begin{array}{ccc}
		0   & -2x_3(z) & - g_{p}(r) \, x_2(r,\varphi) \\
		2x_3(z) &    0 &   g_{p}(r) \, x_1(r,\varphi) \\
		g_{p}(r) \, x_2(r,\varphi) & - g_{p}(r) \, x_1(r,\varphi) & 0 \\
	\end{array}
	\right) \, ,
	\label{eq:ansatz_AD}
	\\*
	\boldsymbol{S}(r,\varphi,z) &= \frac{\boldsymbol{\vartheta}}{2}
	\left(
	\begin{array}{ccc}
		0   & 0 & g_{m}(r) \, x_2(r,\varphi) \\
		0 &    0 &  - g_{m}(r) \, x_1(r,\varphi) \\
		g_{m}(r) \, x_2(r,\varphi) & - g_{m}(r) \, x_1(r,\varphi) & 0 \\
	\end{array}
	\right)
	\, ,
	\notag
\end{align}
where, in relation to the ansatz  (\ref{eq:ansatz_RM}), $g_{m}(r ) \coloneqq g_{1}(r ) - g_{2}(r )$ and $g_{p}(r ) \coloneqq g_{1}(r ) + g_{2}(r )$.
Substituting the ansatz  (\ref{eq:ansatz_AD}) in  (\ref{eq:equi_AD}) the 15 equilibrium equation are equivalent to
\begin{align}
	\frac{1}{6} \, \boldsymbol{\vartheta} \, \sin \varphi \left(6 r \, \mu _c (g_{p}(r )-1) - \mu \,  L_c^2 \, (a_{1}+2 a_{3}) \left(3 g'_{p}(r ) + r \, g''_{p}(r )\right)\right) 
	&= 0 \, ,
	\notag
	\\*
	\frac{1}{6} \, \boldsymbol{\vartheta} \, \cos \varphi \left(6 r \, \mu _c (g_{p}(r )-1) - \mu \,  L_c^2 \, (a_{1}+2 a_{3}) \left(3 g'_{p}(r ) + r \, g''_{p}(r )\right)\right) 
	&= 0 \, ,
	\label{eq:equi_equa_AD}
	\\*
	\frac{1}{2} \boldsymbol{\vartheta}  \sin \varphi
	\left(
	a_{4} \, \mu \, L_c^2 \left(3 g_{m}'(r ) + r \, g_{m}''(r )\right)
	-2 r \, \mu _e \, (g_{m}(r ) + 1)
	-2 r \, g_{m}(r ) \, \mu _{\mbox{\tiny micro}}
	\right)
	&= 0 \, ,
	\notag
	\\*
	-\frac{1}{2} \boldsymbol{\vartheta}  \cos \varphi
	\left(
	a_{4} \, \mu \, L_c^2 \left(3 g_{m}'(r ) + r \, g_{m}''(r )\right)
	-2 r \, \mu _e (g_{m}(r ) + 1)
	-2 r \, g_{m}(r ) \, \mu _{\mbox{\tiny micro}}
	\right)
	&= 0 \, .
\notag
\end{align}
Between the two equilibrium equations  (\ref{eq:equi_equa_AD}) there are only two independent equation since
 (\ref{eq:equi_equa_AD})$_1 = -\tan \varphi$~(\ref{eq:equi_equa_AD})$_2$ and
 (\ref{eq:equi_equa_AD})$_3 = \tan \varphi$~(\ref{eq:equi_equa_AD})$_4$.
The solution of  (\ref{eq:equi_equa_AD}) is 
\begin{align}
    g_{p}(r ) &= \,
	 1
	-\frac{i \, A_{1} I_1\left( \frac{     r \, f_1}{L_c}\right)}{r }
	+\frac{     A_{2} Y_1\left(-\frac{i \, r \, f_1}{L_c}\right)}{r }
	\, ,
	&&
	f_{1} \coloneqq
	\sqrt{\frac{6 \mu _c}{(a_{1} + 2 a_{3}) \, \mu}} \, ,
	\\*
	g_{m}(r) &= \,
	\frac{A_{2} Y_1\left(-\frac{i r  f_2}{L_c}\right)-i A_{1} I_1\left(\frac{r  f_2}{L_c}\right)}{r}
	-\frac{\mu _e}{\mu _e+\mu _{\mbox{\tiny micro}}}
	\, ,
	&&
	f_{2} \coloneqq
	\sqrt{\frac{2(\mu _e+\mu _{\mbox{\tiny micro}})}{a_{4} \mu}} \, ,
	\label{eq:sol_fun_MStrainz}
\end{align}
where $I_{n}\left(\cdot\right)$ is the \textit{modified Bessel function of the first kind}, $Y_{n}\left(\cdot\right)$ is the \textit{Bessel function of the second kind} (see appendix \ref{app:bessel} for the formal definitions), and
$A_1$, $A_2$, $A_3$, $A_4$ are integration constants.

\begin{figure}[H]
	\centering
	\includegraphics[height=5.5cm]{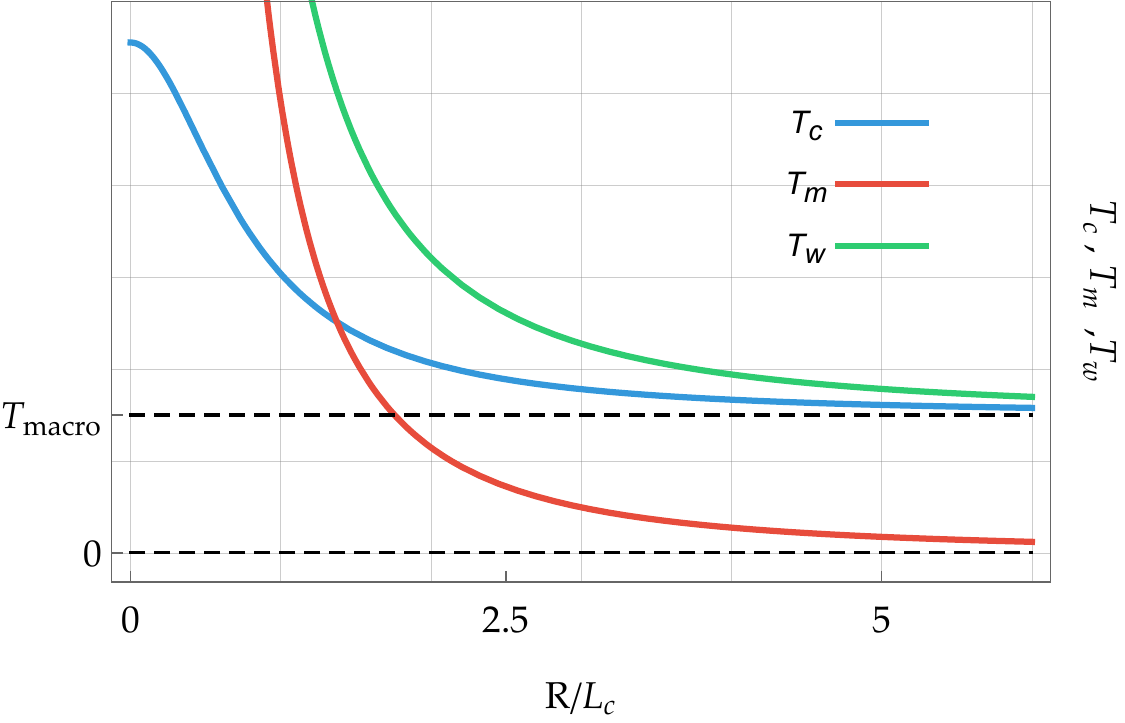}
	\caption{(\textbf{Ad-hoc model}) Torsional stiffness for the classical torque $T_{\mbox{c}}$, the higher-order torque $T_{\mbox{m}}$, and the torque energy $T_{\mbox{w}}$ while varying $L_c$.
	The torsional stiffness is unbounded as $L_c \to \infty$ ($R\to 0$) due to the Cosserat effects.
	The values of the parameters used are: $\mu = 1$, $\mu_c = 1/2$, $\mu _e= 1/3$, $\mu _{\mbox{\tiny micro}}= 1/4$, $a_1= 1/5$, $a_3= 1/6$, $a_4= 1/7$, $R= 1$.}
	\label{fig:all_plot_AD}
\end{figure}

The values of $A_1$ and $A_2$ are determined thanks to the boundary conditions  (\ref{eq:BC_AD}), while, due to the divergent behaviour of the Bessel function of the second kind at $r=0$, we have to set $A_2=A_4=0$ in order to have a continuous solution.
The fulfilment of the boundary conditions  (\ref{eq:BC_AD}) allows us to find the expressions of the integration constants
\begin{align}
    A_1 &= \, 
	-\frac{
	i \, R \, L_c
	}{
	f_1 \, R \, z_1 \left(I_0\left(\frac{R \, f_1}{L_c}\right)
	+ I_2 \left(\frac{R \, f_1}{L_c}\right)\right)
	+ z_2 \, L_c \, I_1\left(\frac{R \, f_1}{L_c}\right)
	}
	\, ,
	&&
	z_1 \coloneqq \,  \frac{a_1 + 2a_3}{3a_1} \, ,
	\notag
	\\
	A_3 &= \, 
	\frac{2 i \, L_c }{I_0\left(\frac{R f_2}{L_c}\right)+I_2\left(\frac{R f_2}{L_c}\right)}
	\frac{\mu _e}{f_2 (\mu _e+\mu _{\mbox{\tiny micro}})}
	\, ,
	&&
	z_2 \coloneqq \, \frac{4 a_3 - a_1}{3a_1}
	\, .
	\label{eq:BC_AD_3}
\end{align}
The classical torque, the higher-order torque, and the energy (per unit length d$z$) expressions are
\begin{align}
	M_{\mbox{c}} (\boldsymbol{\vartheta}) 
	\coloneqq&
	\int_{0}^{2\pi}
	\int_{0}^{R}
	\Big[
	\langle
	\widetilde{\boldsymbol{\sigma}} \, \boldsymbol{e}_{z} , \boldsymbol{e}_{\varphi}
	\rangle
	r
	\Big] r
	\, \mbox{d}r \, \mbox{d}\varphi 
	\notag
	\\*
	=&
	\left[
	\frac{
	4 \mu _c \, I_2\left(\frac{R \, f_1}{L_c}\right)
	\frac{L_c^2}{R^2}
	}{
	f_1 \left(2 \, f_1 \, z_1 \, I_0\left(\frac{R \, f_1}{L_c}\right)
	+\left(z_2 - 2 z_1\right) I_1\left(\frac{R \, f_1}{L_c}\right) \frac{L_c}{R}\right)
	}
	\right.
	\notag
	\\*
	&
	\left.
	\hspace{0.25cm}
	+
	\frac{\mu _e \, \mu _{\mbox{\tiny micro}}}{\mu _e+\mu _{\mbox{\tiny micro}}}
	+
	\frac{\mu _e^2 \, \mu \, a_{4}}{\left(\mu _e+\mu _{\mbox{\tiny micro}}\right)^2}
	\frac{
	4 \, I_2\left(\frac{R \, f_2}{L_c}\right)
	}{
	I_0\left(\frac{R \, f_2}{L_c}\right)
	+ I_2\left(\frac{R \, f_2}{L_c}\right) 
	}
	\frac{L_c^2}{R^2}
	\right]
	I_{p} \, 
	\boldsymbol{\vartheta}
	= T_{\mbox{c}} \, \boldsymbol{\vartheta} \, ,
	\notag
	\\[3mm]
	M_{\mbox{m}}(\boldsymbol{\vartheta}) 
	\coloneqq&
	\int_{0}^{2\pi}
	\int_{0}^{R}
	\Big[
	\overbrace{
	\langle
	\mbox{sym} (\boldsymbol{\mathfrak{m}} \, \boldsymbol{e}_{z})
	\boldsymbol{e}_{\varphi} ,
	\boldsymbol{e}_{r}
	\rangle
	-
	\langle
	\mbox{sym} (\boldsymbol{\mathfrak{m}} \, \boldsymbol{e}_{z})
	\boldsymbol{e}_{r} ,
	\boldsymbol{e}_{\varphi}
	\rangle
	}^{\mbox{micro-strain component}}
	\notag
	\\*
	&
	\hspace{1.2cm}
	+
	\underbrace{
	\langle
	\mbox{skew}(\boldsymbol{m} \times \boldsymbol{e}_{z})
	\boldsymbol{e}_{\varphi} ,
	\boldsymbol{e}_{r}
	\rangle
	-
	\langle
	\mbox{skew}(\boldsymbol{m} \times \boldsymbol{e}_{z})
	\boldsymbol{e}_{r} ,
	\boldsymbol{e}_{\varphi}
	\rangle
	}_{\mbox{Cosserat component}}
	\Big]
	\, r
	\, \mbox{d}r \, \mbox{d}\varphi
	\notag
	\\*
	=&
	\left[
	\frac{
	2 \mu \left(
	3 a_{1} \, f_1 \, z_1 \, I_0\left(\frac{R \, f_1}{L_c}\right) \,
	\frac{L_c^2}{R^2}
	-2 (a_{1} - a_{3}) \, I_1\left(\frac{R \, f_1}{L_c}\right)
	\frac{L_c^3}{R^3}
	\right)
	}{
	6 f_1 \, z_1 \, I_0\left(\frac{R \, f_1}{L_c}\right)
	-3 I_1\left(\frac{R \, f_1}{L_c}\right)
	\frac{L_c}{R}
	}
	\right]
	=
	T_{\mbox{m}} \, \boldsymbol{\vartheta} \, ,
	\label{eq:torque_stiffness_MStrainz}
	\\[3mm]
	W_{\mbox{tot}} (\boldsymbol{\vartheta}) 
	\coloneqq&
	\int_{0}^{2\pi}
	\int_{0}^{R}
	W \left(\boldsymbol{\mbox{D}u}, \boldsymbol{A}, \boldsymbol{S}, \text{Curl} \, \boldsymbol{A}, \boldsymbol{\mbox{D}S}\right)  \, \, r
	\, \mbox{d}r \, \mbox{d}\varphi
	\notag
	\\*
	=&
	\frac{1}{2}
	\left[
	\frac{
	4 \mu _c \, I_2\left(\frac{R \, f_1}{L_c}\right)
	\frac{L_c^2}{R^2}
	}{
	f_1 \left(2 \, f_1 \, z_1 \, I_0\left(\frac{R \, f_1}{L_c}\right)
	+\left(z_2 - 2 z_1\right) I_1\left(\frac{R \, f_1}{L_c}\right) \frac{L_c}{R}\right)
	}
	\right.
	\notag
	\\*
	&
	\hspace{0.5cm}
	+
	\frac{\mu _e \, \mu _{\mbox{\tiny micro}}}{\mu _e+\mu _{\mbox{\tiny micro}}}
	+
	\frac{\mu _e^2 \, \mu \, a_{4}}{\left(\mu _e+\mu _{\mbox{\tiny micro}}\right)^2}
	\frac{
	4 \, I_2\left(\frac{R \, f_2}{L_c}\right)
	}{
	I_0\left(\frac{R \, f_2}{L_c}\right)
	+ I_2\left(\frac{R \, f_2}{L_c}\right) 
	}
	\frac{L_c^2}{R^2}
	\notag
	\\*
	&
	\left.
	\hspace{0.5cm}
	+
	\frac{
	2 \mu \left(
	3 a_{1} \, f_1 \, z_1 \, I_0\left(\frac{R \, f_1}{L_c}\right) \,
	\frac{L_c^2}{R^2}
	-2 (a_{1} - a_{3}) \, I_1\left(\frac{R \, f_1}{L_c}\right)
	\frac{L_c^3}{R^3}
	\right)
	}{
	6 f_1 \, z_1 \, I_0\left(\frac{R \, f_1}{L_c}\right)
	-3 I_1\left(\frac{R \, f_1}{L_c}\right)
	\frac{L_c}{R}
	}
	\right]
	I_{p} \, 
	\boldsymbol{\vartheta}^2
	= \frac{1}{2} \, T_{\mbox{w}} \, \boldsymbol{\vartheta}^2
	\, .
	\notag
\end{align}
It is highlighted that, like for the micro-strain model (Section  \ref{sec:Micro_strain}), the higher order torque contribution
$\langle
\left(\boldsymbol{\mathfrak{m}} \, \boldsymbol{e}_{z}\right)
\boldsymbol{e}_{\varphi} ,
\boldsymbol{e}_{\varphi}
\rangle$ 
is equal to zero.
The plot of the torsional stiffness for the classical torque, the higher-order torque, and the torque energy while varying $L_c$ is shown in Fig.~\ref{fig:all_plot_AD}.
Again, it holds
\begin{equation}
\frac{\mbox{d}}{\mbox{d}\boldsymbol{\vartheta}}W_{\mbox{tot}}(\boldsymbol{\vartheta}) = M_{\mbox{c}} (\boldsymbol{\vartheta}) + M_{\mbox{m}} (\boldsymbol{\vartheta}) \, ,
\qquad\qquad\qquad
\frac{\mbox{d}^2}{\mbox{d}\boldsymbol{\vartheta}^2}W_{\mbox{tot}}(\boldsymbol{\vartheta})
= T_{\mbox{c}} + T_{\mbox{m}}
= T_{\mbox{w}}
\, .
\end{equation}
\section{Summary and conclusions}
We have derived the analytical expressions of the torsional rigidity for a family of generalized continua capable of  modelling size-dependence in the sense that more slender specimens are comparatively stiffer.
We only consider (simplified) isotropic expressions so as to better compare the different models with each other.
For example, a strain gradient continuum, by construction, does not have mixed energy terms.
Therefore, we omitted these terms in all models.
Excluding the mixed terms like $\langle \text{sym} \,  \boldsymbol{\text{D}u}, \text{sym} \,  \boldsymbol{\text{D}u} -\boldsymbol{P} \rangle$ also simplifies considerably the investigation of positive definiteness.
Indeed, all presented models are positive definite if the usual relations
\begin{alignat}{2}
    \mu _{\mbox{\tiny macro}} &> 0 \, ,
    \qquad
    &\kappa _{\mbox{\tiny macro}} &= \frac{2\mu _{\mbox{\tiny macro}} + 3\lambda _{\mbox{\tiny macro}}}{3}  > 0 \, ,
    \notag
    \\*
    \mu _{\mbox{\tiny micro}} &> 0 \, ,
    \qquad
    &\kappa _{\mbox{\tiny micro}} &= \frac{2\mu _{\mbox{\tiny micro}} + 3\lambda _{\mbox{\tiny micro}}}{3}  > 0 \, ,
    \\* 
    \mu _{\mbox{\tiny micro}} > \mu _{\mbox{\tiny macro}} \Longrightarrow \mu _{e} &> 0 \, ,
    \qquad
    &\kappa _{e} &= \frac{2\mu _{e} + 3\lambda _{e}}{3}   > 0 \, .
    \notag
\end{alignat}
are satisfied together with individual positivity of all curvature parameters.
In all cases, the displacement follows the classical pure torsion solution.
Despite the conceptual simplicity of the models, we observe already a delicate interplay between the used kinematics and the assumed curvature energy expression.
For example, let us compare the relaxed micromorphic model with the micro-strain model (Section  \ref{sec:Micro_strain}).
Both models have a similar looking lower order energy term (if the Cosserat couple modulus $\mu_c \equiv 0$), but different kinematics and different curvature energy.
For arbitrary slender specimens, the torsional stiffness of the micro-strain model is governed by $\mu _e$, whereas the torsional stiffness of the relaxed micromorphic model is determined by $\mu _{\mbox{\tiny micro}}$.
Thus, the physical interpretation of the material parameters in both models is completely different.
This is surprising at first glance but the reason for this response is finally explained in Appendix \ref{app:const_disto}.

In the end, the more restricted the used kinematics, the less viable a model may become.
In this respect, only the full micromorphic kinematics degree of freedom (12 DOFS) can be advised.
In addition, the curvature energy should not intervene too strongly.
For example, penalizing a full gradient $D\boldsymbol{P}$ in the curvature energy of the classical micromorphic model leads to a stiffness singularity for arbitrary slender specimens, while penalizing only $\text{Curl} \, \boldsymbol{P}$ in the relaxed micromorphic model does not show the same singular response.
Moreover, in the relaxed micromorphic model the interpretation of the lower order material parameters ($\mu _e$, $\mu _{\mbox{\tiny micro}}$, $\mu _{\mbox{\tiny micro}}$, etc.) does not in principle change when different curvature energies are considered.
In the end, it is therefore the relaxed micromorphic model that produces sensible and consistent response in all considered cases.
It remains to be investigated if, together with the previously developed solution for bending and shear \cite{rizzi2021shear,rizzi2021bending}, the present analytical solution allows to identify the complete set of micromorphic parameters of a material from bending, shear and torsion experiments at specimens with different diameters.

{\scriptsize
	\paragraph{{\scriptsize Acknowledgements.}}
	Angela Madeo acknowledges support from the European Commission through the funding of the ERC Consolidator Grant META-LEGO, N° 101001759.
	Angela Madeo and Gianluca Rizzi acknowledge funding from the French Research Agency ANR, “METASMART” (ANR-17CE08-0006).
	Angela Madeo and Gianluca Rizzi acknowledge support from IDEXLYON in the framework of the “Programme Investissement d’Avenir” ANR-16-IDEX-0005.
	Hassam Khan acknowledges the  support of the German Academic Exchange Service (DAAD) and the Higher Education Commission of Pakistan (HEC).
	Ionel Dumitrel Ghiba acknowledges support from a grant of the Romanian Ministry of Research and Innovation, CNCS--UEFISCDI, project number PN-III-P1-1.1-TE-2019-0348, Contract No. TE 8/2020, within PNCDI III.
	Patrizio Neff acknowledges support in the framework of the DFG-Priority Programme 2256 ``Variational Methods for Predicting Complex Phenomena in Engineering Structures and Materials", Neff 902/10-1, Project-No. 440935806.
}

\let\oldbibliography\thebibliography
\renewcommand{\thebibliography}[1]{%
	\oldbibliography{#1}%
	\setlength{\itemsep}{2.5pt}}


\begin{scriptsize}
	\bibliographystyle{plain}
	\bibliography{Torsion_Long_Biblio}
\end{scriptsize}

%
%

\begin{footnotesize}
\appendix
\section{Cylindrical coordinates}
\label{app:cyli_coordi}
The relation between the derivatives for a vector field $\boldsymbol{u}(\boldsymbol{x})$, $\widehat{\boldsymbol{u}}(\boldsymbol{r})$ and a second order tensor $\boldsymbol{P}(\boldsymbol{x})$, $\widehat{\boldsymbol{P}}(\boldsymbol{r})$ with respect an orthogonal set of coordinates ${\boldsymbol{x} = \{x_1,x_2,x_3\}}$ and a cylindrical set  of coordinates ${\boldsymbol{r} = \{r,\varphi,z\}}$ are the following

\begin{align}
	\boldsymbol{\nabla}_{\boldsymbol{x}} \boldsymbol{u}(\boldsymbol{x})
	&=
	\left(\boldsymbol{\nabla}_{\boldsymbol{r}} \widehat{\boldsymbol{u}}(\boldsymbol{r})\right) \boldsymbol{Q}^{(1)} \, ,
	\qquad
	&\boldsymbol{\nabla}_{\boldsymbol{x}} \boldsymbol{P}(\boldsymbol{x})
	&=
	\left(\boldsymbol{\nabla}_{\boldsymbol{r}} \widehat{\boldsymbol{P}}(\boldsymbol{r})\right) \boldsymbol{Q}^{(1)} \, ,
    \\*
	u_{i,j}
	&=
	\widehat{u}_{i,\alpha} Q^{(1)}_{\alpha j} \, ,
	\qquad
	&P_{ij,k}
	&=
	\widehat{P}_{ij,\alpha} Q^{(1)}_{\alpha k} \, ,
\end{align}
where
\begin{equation}
\boldsymbol{Q}^{(1)}=
\left(\boldsymbol{\nabla}_{\boldsymbol{r}} \boldsymbol{x}(\boldsymbol{r})\right)^{-1} =
\left(
\begin{array}{ccc}
	\cos \varphi & \sin \varphi & 0 \\
	-\frac{\sin \varphi}{r } & \frac{\cos \varphi}{r } & 0 \\
	0 & 0 & 1 \\
\end{array}
\right)\, ,
\end{equation}
and 
\begin{align}
	\boldsymbol{\nabla}^2_{\boldsymbol{x}} \boldsymbol{u}(\boldsymbol{x})
	&=
	\left(\boldsymbol{\nabla}^2_{\boldsymbol{r}} \widehat{\boldsymbol{u}}(\boldsymbol{r})\right) \boldsymbol{:} \boldsymbol{Q}^{(2)} 
	+\left(\boldsymbol{\nabla}_{\boldsymbol{r}} \widehat{\boldsymbol{u}}(\boldsymbol{r})\right) \boldsymbol{Q}^{(3)} \, ,
	\qquad
	&\boldsymbol{\nabla}_{\boldsymbol{x}} \boldsymbol{P}(\boldsymbol{x})
	&=
	\left(\boldsymbol{\nabla}^2_{\boldsymbol{r}} \widehat{\boldsymbol{P}}(\boldsymbol{r})\right) \boldsymbol{:} \boldsymbol{Q}^{(2)} 
	+\left(\boldsymbol{\nabla}_{\boldsymbol{r}} \widehat{\boldsymbol{P}}(\boldsymbol{r})\right) \boldsymbol{Q}^{(3)} \, .
    \\*
	u_{i,jk}
	&=
	\widehat{u}_{i,\alpha\beta} Q_{\alpha \beta j k} 
	+ \widehat{u}_{i,\alpha} Q_{\alpha j k}\, ,
	\qquad
	&P_{ij,kl}
	&=
	\widehat{P}_{ij,\alpha\beta} Q_{\alpha \beta k l} 
	+ \widehat{P}_{ij,\alpha} Q_{\alpha k l}\, .
\end{align}
\begin{align}
	\boldsymbol{Q}^{(2)} &=
	\left(
	\begin{array}{ccc}
		\left(
		\begin{array}{ccc}
			\cos ^2 \varphi  & \sin \varphi \cos \varphi & 0 \\
			\sin \varphi \cos \varphi & \sin ^2 \varphi  & 0 \\
			0 & 0 & 0 \\
		\end{array}
		\right) & \left(
		\begin{array}{ccc}
			-\frac{\sin \varphi \cos \varphi}{r } & \frac{\cos (2 \varphi )}{2 r } & 0 \\
			\frac{\cos (2 \varphi )}{2 r } & \frac{\sin \varphi \cos \varphi}{r } & 0 \\
			0 & 0 & 0 \\
		\end{array}
		\right) & \left(
		\begin{array}{ccc}
			0 & 0 & \frac{\cos \varphi}{2} \\
			0 & 0 & \frac{\sin \varphi}{2} \\
			\frac{\cos \varphi}{2} & \frac{\sin \varphi}{2} & 0 \\
		\end{array}
		\right) \\
		\left(
		\begin{array}{ccc}
			-\frac{\sin \varphi \cos \varphi}{r } & \frac{\cos (2 \varphi )}{2 r } & 0 \\
			\frac{\cos (2 \varphi )}{2 r } & \frac{\sin \varphi \cos \varphi}{r } & 0 \\
			0 & 0 & 0 \\
		\end{array}
		\right) & \left(
		\begin{array}{ccc}
			\frac{\sin ^2 \varphi }{r ^2} & -\frac{\sin \varphi \cos \varphi}{r ^2} & 0 \\
			-\frac{\sin \varphi \cos \varphi}{r ^2} & \frac{\cos ^2 \varphi }{r ^2} & 0 \\
			0 & 0 & 0 \\
		\end{array}
		\right) & \left(
		\begin{array}{ccc}
			0 & 0 & -\frac{\sin \varphi}{2 r } \\
			0 & 0 & \frac{\cos \varphi}{2 r } \\
			-\frac{\sin \varphi}{2 r } & \frac{\cos \varphi}{2 r } & 0 \\
		\end{array}
		\right) \\
		\left(
		\begin{array}{ccc}
			0 & 0 & \frac{\cos \varphi}{2} \\
			0 & 0 & \frac{\sin \varphi}{2} \\
			\frac{\cos \varphi}{2} & \frac{\sin \varphi}{2} & 0 \\
		\end{array}
		\right) & \left(
		\begin{array}{ccc}
			0 & 0 & -\frac{\sin \varphi}{2 r } \\
			0 & 0 & \frac{\cos \varphi}{2 r } \\
			-\frac{\sin \varphi}{2 r } & \frac{\cos \varphi}{2 r } & 0 \\
		\end{array}
		\right) & \left(
		\begin{array}{ccc}
			0 & 0 & 0 \\
			0 & 0 & 0 \\
			0 & 0 & 1 \\
		\end{array}
		\right) \\
	\end{array}
	\right) \, ,
	\notag
	\\*
	\boldsymbol{Q}^{(3)} &=
	\left(
	\begin{array}{ccc}
		\left(
		\begin{array}{c}
			\frac{\sin ^2 \varphi }{r } \\
			-\frac{\sin \varphi \cos \varphi}{r } \\
			0 \\
		\end{array}
		\right) & \left(
		\begin{array}{c}
			-\frac{\sin \varphi \cos \varphi}{r } \\
			\frac{\cos ^2 \varphi }{r } \\
			0 \\
		\end{array}
		\right) & \left(
		\begin{array}{c}
			0 \\
			0 \\
			0 \\
		\end{array}
		\right) \\
		\left(
		\begin{array}{c}
			\frac{\sin (2 \varphi )}{r ^2} \\
			-\frac{\cos (2 \varphi )}{r ^2} \\
			0 \\
		\end{array}
		\right) & \left(
		\begin{array}{c}
			-\frac{\cos (2 \varphi )}{r ^2} \\
			-\frac{\sin (2 \varphi )}{r ^2} \\
			0 \\
		\end{array}
		\right) & \left(
		\begin{array}{c}
			0 \\
			0 \\
			0 \\
		\end{array}
		\right) \\
		\left(
		\begin{array}{c}
			0 \\
			0 \\
			0 \\
		\end{array}
		\right) & \left(
		\begin{array}{c}
			0 \\
			0 \\
			0 \\
		\end{array}
		\right) & \left(
		\begin{array}{c}
			0 \\
			0 \\
			0 \\
		\end{array}
		\right) \\
	\end{array}
	\right)\, .
\end{align}
\section{Bessel functions}
\label{app:bessel}
The \textit{Bessel functions} are the solutions $y(x)$ of the Bessel differential equation \cite{watson1995treatise}
\begin{equation}
	x^2 \, \frac{d^2 y}{dx^2} + x \, \frac{d y}{dx} + (x^2 - \alpha^2)y = 0 \, .
	\label{eq:bessel_diff_equa}
\end{equation}
For the particular case in which $\alpha=n$ is an integer, the solution of  (\ref{eq:bessel_diff_equa}) can be expressed as a linear combination of the Bessel function of the \textit{first} $J_{n}(x)$ and \textit{second} $Y_{n}(x)$ kind
\begin{equation}
	y(x) = A_1 \, J_{n}(x) + A_2 \, Y_{n}(x) \, ,
	\label{eq:sol_bessel_diff_equa}
\end{equation}
whose definitions are
\begin{equation}
	J_{n}(x)= \int_{0}^{\pi} \cos(n \, \tau - x \, \sin(\tau)) \, \mbox{d}\tau \, ,
	\qquad\qquad 
	Y_{n}(x)= \frac{J_{n}(x) \, \cos(n \pi) - J_{-n}(x)}{\sin(n \pi)} \, .
	\label{eq:bessel_funcs}
\end{equation}
Moreover, the \textit{modified Bessel functions of the first kind} is defined as
$
I_{n}(x)= i^{-n} J_{n}(i x)
$ .
\section{Classical Cosserat formulation in micro-rotation vector format}
\label{app:class_Coss}
 An overview of the different classical notations for the Cosserat model has been given  in \cite{hassanpour2017micropolar}. In \cite{neff2009new} we have presented the correspondence between the Cosserat model expressed in dislocation format (Curl of the skew symmetric micro-distortion tensor) and in its classical formulation (gradient of the micro-rotation vector $\boldsymbol{\phi}$).
The relation between the coefficients in the two notations is
\begin{align}
\alpha &=\frac{\mu \, L_c^2}{2} \frac{1}{3}\left(4a_3-a_1\right) \, ,
&\beta &=\frac{\mu \, L_c^2}{2} \frac{a_1-a_2}{2} \, ,
&\gamma &=\frac{\mu \, L_c^2}{2} \frac{a_1+a_2}{2} \, ,
\label{eq:coeff_coss_class}
\\*
a_1 &=\frac{\gamma + \beta }{\mu \, L_c^2} \, ,
&a_2 &=\frac{\gamma - \beta }{\mu \, L_c^2} \, ,
&a_3 &=\frac{3\alpha + \beta +\gamma}{4 \mu \, L_c^2} \, .
\notag
\end{align}
Setting $\boldsymbol{\phi} \coloneqq \mbox{axl} (\boldsymbol{A})$ and taking into account  (\ref{eq:coeff_coss_class}), the expression of the strain energy for the isotropic Cosserat continuum can be equivalently expressed as
\begin{align}
W \left(\boldsymbol{\mbox{D}u}, \boldsymbol{A},\mbox{Curl}\,\boldsymbol{A}\right) = &
\, \mu_{\mbox{\tiny macro}} \left\lVert \mbox{sym} \, \boldsymbol{\mbox{D}u} \right\rVert^{2}
+ \frac{\lambda_{\mbox{\tiny macro}}}{2} \mbox{tr}^2 \left(\boldsymbol{\mbox{D}u} \right) 
+ \mu_{c} \left\lVert \mbox{skew} \, \boldsymbol{\mbox{D}u} - \boldsymbol{A} \right\rVert^{2}
\notag
\\*
&+ \frac{\mu \, L_c^2}{2}
\underbrace{
	\left(
	a_1 \, \left \lVert \mbox{dev} \, \mbox{sym} \, \mbox{Curl} \, \boldsymbol{A}\right \rVert^2 \, 
	+ a_2 \, \left \lVert \mbox{skew} \, \mbox{Curl} \, \boldsymbol{A}\right \rVert^2 \, 
	+ \frac{a_3}{3} \, \mbox{tr}^2 \left(\mbox{Curl} \, \boldsymbol{A} \right)
	\right)
}_{\mbox{dislocation tensor format}}
\label{eq:energy_Cos_classic}
\\*
= W \left(\boldsymbol{\mbox{D}u}, \boldsymbol{\phi},\boldsymbol{\mbox{D} \phi}\right) = &
\, \mu_{\mbox{\tiny macro}} \left\lVert \mbox{sym} \, \boldsymbol{\mbox{D}u} \right\rVert^{2}
+ \frac{\lambda_{\mbox{\tiny macro}}}{2} \mbox{tr}^2 \left(\boldsymbol{\mbox{D}u} \right) 
+ \frac{\mu_{c}}{2} \left\lVert \mbox{curl} \boldsymbol{u} - 2\boldsymbol{\phi} \right\rVert^{2}
\notag
\\*
&\underbrace{
	\frac{1}{2}
	\left(
	\alpha \, \mbox{tr}^2 \left(\boldsymbol{\mbox{D} \phi} \right)
	+ \beta \, \langle \boldsymbol{\mbox{D} \phi}^{T},\boldsymbol{\mbox{D} \phi} \rangle 
	+ \gamma \, \left \lVert \boldsymbol{\mbox{D} \phi} \right \rVert^2
	\right)
}_{\mbox{classical micro-rotation vector format}}  \, ,
\notag
\end{align}
since 
\begin{equation}
\begin{split}
\left\lVert \mbox{skew} \, \boldsymbol{\mbox{D}u} - \boldsymbol{A} \right\rVert^{2} =
2\left\lVert \mbox{axl}(\mbox{skew} \, \boldsymbol{\mbox{D}u} - \mbox{Anti}(\boldsymbol{\phi})) \right\rVert^{2} =
2\lVert \frac{1}{2}\mbox{curl} \boldsymbol{u} - \boldsymbol{\phi} \rVert^{2} =
\frac{1}{2}\left\lVert \mbox{curl} \boldsymbol{u} - 2\boldsymbol{\phi} \right\rVert^{2} \, .
\end{split}
\end{equation}
The equilibrium equations, in the absence of body forces,  in the classical notation are
\begin{align}
\mbox{Div}\left[
2\mu_{\mbox{\tiny macro}}\,\mbox{sym} \, \boldsymbol{\mbox{D}u} + \lambda_{\mbox{\tiny macro}} \mbox{tr} \left(\boldsymbol{\mbox{D}u} \right) \boldsymbol{\mathbbm{1}}
\right]
-\mu_{c} \, \mbox{curl} \left[\mbox{curl} \, \boldsymbol{u} - 2 \boldsymbol{\phi}\right]  \,  
&= \boldsymbol{0} \, ,
\label{eq:equi_Cos_classic}
\\*
\mbox{Div}\left[
\alpha \, \mbox{tr} \left(\boldsymbol{\mbox{D} \phi} \right) \, \boldsymbol{\mathbbm{1}}
+ \beta \, \left(\boldsymbol{\mbox{D} \phi}\right)^{T} 
+ \gamma \, \boldsymbol{\mbox{D} \phi}  \, 
\right] 
+2\mu_c \, \left(\mbox{curl} \, \boldsymbol{u} - 2 \boldsymbol{\phi}\right)
&= \boldsymbol{0} \, .
\notag
\end{align}
The boundary conditions at the free surface are 
\begin{align}
\boldsymbol{\widetilde{t}}(r = R) = 
\boldsymbol{\widetilde{\sigma}}(r) \cdot \boldsymbol{e}_{r} = 
\boldsymbol{0}_{\mathbb{R}^{3}} \, ,
\qquad\qquad
\overline{\boldsymbol{\eta}}(r = R) = 
\overline{\boldsymbol{m}} (r) \cdot \boldsymbol{e}_{r} = 
\boldsymbol{0}_{\mathbb{R}^{3\times 3}} \, ,
\label{eq:BC_Cos_classic_gen}
\end{align}
where 
\begin{align}\boldsymbol{\widetilde{\sigma}} = 2\mu_{\mbox{\tiny macro}}\,\mbox{sym} \, \boldsymbol{\mbox{D}u} + \lambda_{\mbox{\tiny macro}} \mbox{tr} \left(\boldsymbol{\mbox{D}u} \right) \boldsymbol{\mathbbm{1}} + 2\mu_c \, \left(\mbox{skew} \, \boldsymbol{\mbox{D}u} - \mbox{Anti}(\boldsymbol{\phi})\right),\end{align}
 $\boldsymbol{e}_{r}$ is the radial unit vector, and the second-order moment stress tensor
\begin{equation}
    \overline{\boldsymbol{m}} = 
    \alpha \, \mbox{tr} \left(\boldsymbol{\mbox{D} \phi} \right) \, \boldsymbol{\mathbbm{1}}+ \beta \, \left(\boldsymbol{\mbox{D} \phi}\right)^{T} + \gamma \, \boldsymbol{\mbox{D} \phi}
    \, .
\end{equation}

According to the reference system shown in Fig.~\ref{fig:intro_1}, the ansatz for the displacement field and the micro-rotation vector turns into
\begin{equation}
	\boldsymbol{u}(x_1,x_2,x_3) =
	\boldsymbol{u}(r,\varphi,z) = \boldsymbol{\vartheta}
	\left(
	\begin{array}{c}
		-x_2(r,\varphi) \, x_3(z) \\
		x_1(r,\varphi) \, x_3(z) \\
		0 
	\end{array}
	\right)
	\, ,
	\qquad
	\boldsymbol{\phi}(x_1,x_2,x_3) = 
	\boldsymbol{\phi}(r,\varphi,z) = \frac{\boldsymbol{\vartheta}}{2}
	\left(
	\begin{array}{ccc}
		- g_{p}(r) \, x_1(r,\varphi)   \\
		- g_{p}(r) \, x_2(r,\varphi)\\
		2x_3(z)
	\end{array}
	\right)
	\, ,
    \label{eq:ansatz_Cos_classic}
\end{equation}
Substituting the ansatz  (\ref{eq:ansatz_Cos_classic}) in  (\ref{eq:equi_Cos_classic}) the equilibrium equations are equivalent to
\begin{align}
	-\frac{1}{2} \boldsymbol{\vartheta} \, \cos \varphi  \left(4 \rho \, \mu _c \, (g(\rho )-1) - (\alpha +\beta +\gamma ) \left(3 g'(\rho )+\rho  g''(\rho )\right)\right) = 0 \, ,
	\label{eq:equi_equa_Cos_Classic}
	\\*
	-\frac{1}{2} \boldsymbol{\vartheta} \, \sin \varphi  \left(4 \rho \, \mu _c \, (g(\rho )-1) - (\alpha +\beta +\gamma ) \left(3 g'(\rho )+\rho  g''(\rho )\right)\right) = 0 \, ,
	\notag
\end{align}
which are completely equivalent to  (\ref{eq:equi_equa_Cos}) in Section~\ref{sec:Cos} once used the relations  (\ref{eq:coeff_coss_class}). Since also the boundary conditions  (\ref{eq:BC_Cos_classic_gen}) are equivalent to the boundary condition  (\ref{eq:BC_Cos_gen}) in Section~\ref{sec:Cos}, further calculations are avoided.
Here, we recall the relations between the two moment stress tensors expressed in the classical format ($\overline{\boldsymbol{m}}$) and in the dislocation format ($\boldsymbol{m}$)
\begin{gather}
    \mbox{dev} \, \mbox{sym} \, \boldsymbol{m}
    =
    - \mbox{dev} \, \mbox{sym} \, \overline{\boldsymbol{m}}
    \, ,
    \qquad\qquad\qquad
    \mbox{skew} \, \boldsymbol{m}
    =
    \mbox{skew} \, \overline{\boldsymbol{m}}
    \, ,
    \qquad\qquad\qquad
    \mbox{tr} \, \boldsymbol{m}
    =
    \frac{1}{2} \mbox{tr} \, \overline{\boldsymbol{m}}
    \, ,
    \label{eq:relation_m_over_m}
\end{gather}
where $\mbox{skew} \, \boldsymbol{m}= \mbox{skew} \, \overline{\boldsymbol{m}} = \boldsymbol{0}$ for the torsional problem, and
\begin{gather}
    \boldsymbol{m} =
    \mu \, L_c^2
    \left(
    a_1 \, \mbox{dev} \, \mbox{sym} \, \mbox{Curl} \, \boldsymbol{A} \, 
    + \frac{a_3}{3} \, \mbox{tr} \left(\mbox{Curl} \, \boldsymbol{A} \right)\boldsymbol{\mathbbm{1}} \, 
    \right) \, ,
    \qquad
    \overline{\boldsymbol{m}} = 
    \alpha \, \mbox{tr} \left(\boldsymbol{\mbox{D} (\mbox{axl} (\boldsymbol{A}))} \right) \, \boldsymbol{\mathbbm{1}}+ \beta \, \left(\boldsymbol{\mbox{D} (\mbox{axl} (\boldsymbol{A}))}\right)^{T} + \gamma \, \boldsymbol{\mbox{D} (\mbox{axl} (\boldsymbol{A}))}
    \, ,
\end{gather}
It is also interesting to show the relation between the two higher-order torques expressed in terms of $\boldsymbol{m}$ and $\overline{\boldsymbol{m}}$, respectively.
First, we observe
\begin{align}
	\langle
	\mbox{skew}\left(\boldsymbol{m} \times \boldsymbol{e}_{z}\right)
	\boldsymbol{e}_{\varphi} ,
	\boldsymbol{e}_{r}
	\rangle
	-
	\langle
	\mbox{skew} (\boldsymbol{m} \times \boldsymbol{e}_{z} )
	\boldsymbol{e}_{r} ,
	\boldsymbol{e}_{\varphi}
	\rangle
	=
	\langle
	\left(\boldsymbol{m} \times \boldsymbol{e}_{z}\right)
	\boldsymbol{e}_{\varphi} ,
	\boldsymbol{e}_{r}
	\rangle
	-
	\langle
	\left(\boldsymbol{m} \times \boldsymbol{e}_{z}\right)
	\boldsymbol{e}_{r} ,
	\boldsymbol{e}_{\varphi}
	\rangle
	\, ,
\end{align}
which does not holds component-wise.
\!\!\!\footnote{
The values of the four terms for a generic second order tensor $\boldsymbol{m}$ are:
$
\langle
\mbox{skew} (\boldsymbol{m} \times \boldsymbol{e}_{z} )
\boldsymbol{e}_{\varphi} ,
\boldsymbol{e}_{r}
\rangle
=
\frac{1}{2}\left( m_{11} + m_{22} \right)
$
, 
$
\langle
\mbox{skew} (\boldsymbol{m} \times \boldsymbol{e}_{z} )
\boldsymbol{e}_{r} ,
\boldsymbol{e}_{\varphi}
\rangle
=
-\frac{1}{2}\left( m_{11} + m_{22} \right)
$
, 
$
\langle
\left(\boldsymbol{m} \times \boldsymbol{e}_{z}\right)
\boldsymbol{e}_{\varphi} ,
\boldsymbol{e}_{r}
\rangle
=
m_{22} \sin ^2 \varphi + m_{11} \cos ^2 \varphi + (m_{12}+m_{21}) \sin \varphi \, \cos \varphi 
$
, 
$
\langle
\left(\boldsymbol{m} \times \boldsymbol{e}_{z}\right)
\boldsymbol{e}_{r} ,
\boldsymbol{e}_{\varphi}
\rangle
=
-m_{22} \cos ^2 \varphi - m_{11} \sin ^2 \varphi +(m_{12}+m_{21}) \sin \varphi  \, \cos \,\varphi
$.
Note that
$
\langle
\mbox{skew}\left(\boldsymbol{m} \times \boldsymbol{e}_{z}\right)
\boldsymbol{e}_{\varphi} ,
\boldsymbol{e}_{r}
\rangle
\neq
\langle
\left(\boldsymbol{m} \times \boldsymbol{e}_{z}\right)
\boldsymbol{e}_{\varphi} ,
\boldsymbol{e}_{r}
\rangle
$
and 
$
\langle
\mbox{skew}\left(\boldsymbol{m} \times \boldsymbol{e}_{z}\right)
\boldsymbol{e}_{r} ,
\boldsymbol{e}_{\varphi}
\rangle
\neq
\langle
\left(\boldsymbol{m} \times \boldsymbol{e}_{z}\right)
\boldsymbol{e}_{r} ,
\boldsymbol{e}_{\varphi}
\rangle
$.
}

Using that the cross product between two unit vectors gives the third one, and
\begin{align}
    \left( \boldsymbol{m} \times \boldsymbol{v} \right)\boldsymbol{w}
    =
    \boldsymbol{m} \left( \boldsymbol{v} \times \boldsymbol{w} \right)
    \qquad
    \forall\boldsymbol{v},\boldsymbol{w}\in\mathbb{R}^3
    \quad
    \mbox{and}
    \quad
    \forall\boldsymbol{m}\in\mathbb{R}^{3\times 3}
    \, ,
\end{align}
it is possible to write
\begin{align}
    \langle
	\left(\boldsymbol{m} \times \boldsymbol{e}_{z}\right)
	\boldsymbol{e}_{\varphi} ,
	\boldsymbol{e}_{r}
	\rangle
	-
	\langle
	\left(\boldsymbol{m} \times \boldsymbol{e}_{z}\right)
	\boldsymbol{e}_{r} ,
	\boldsymbol{e}_{\varphi}
	\rangle
	=
	\langle
	\boldsymbol{m}
	\left(
	\boldsymbol{e}_{\varphi}
	\times
	\boldsymbol{e}_{z}
	\right)
	,
	\boldsymbol{e}_{r}
	\rangle
	-
	\langle
	\boldsymbol{m}
	\left(
	\boldsymbol{e}_{r}
	\times
	\boldsymbol{e}_{z}
	\right)
	,
	\boldsymbol{e}_{\varphi}
	\rangle
	=
	-
	\left[
	\langle
	\boldsymbol{m}
	\,
	\boldsymbol{e}_{r}
	,
	\boldsymbol{e}_{r}
	\rangle
	+
	\langle
	\boldsymbol{m}
	\,
	\boldsymbol{e}_{\varphi}
	,
	\boldsymbol{e}_{\varphi}
	\rangle
	\right]
	\, .
\end{align}
Since $\left( \boldsymbol{e}_{r}\otimes\boldsymbol{e}_{r} + \boldsymbol{e}_{\varphi}\otimes\boldsymbol{e}_{\varphi} + \boldsymbol{e}_{z}\otimes\boldsymbol{e}_{z} \right)= \boldsymbol{\mathbbm{1}}$
we may convert the double dot-product into a dyadic product as follows
\begin{align}
	&
	-
	\left[
	\langle
	\boldsymbol{m}
	\,
	\boldsymbol{e}_{r}
	,
	\boldsymbol{e}_{r}
	\rangle
	+
	\langle
	\boldsymbol{m}
	\,
	\boldsymbol{e}_{\varphi}
	,
	\boldsymbol{e}_{\varphi}
	\rangle
	\right]
	=
	-
	\langle
	\boldsymbol{m}
	,
	\left(
	\boldsymbol{e}_{r}
	\otimes
	\boldsymbol{e}_{r}
	+
	\boldsymbol{e}_{\varphi}
	\otimes
	\boldsymbol{e}_{\varphi}
	\right)
	\rangle
	=
	-
	\langle
	\boldsymbol{m}
	,
	\left(
	\boldsymbol{\mathbbm{1}}
	-
	\boldsymbol{e}_{z}
	\otimes
	\boldsymbol{e}_{z}
	\right)
	\rangle
	\, .
\end{align}
Substituting the relation  (\ref{eq:relation_m_over_m}) between $\boldsymbol{m}$ and $\overline{\boldsymbol{m}}$ we have
\begin{align}
	-
	\langle
	\boldsymbol{m}
	,
	\left(
	\boldsymbol{\mathbbm{1}}
	-
	\boldsymbol{e}_{z}
	\otimes
	\boldsymbol{e}_{z}
	\right)
	\rangle
	=
	-
	\langle \left(
	- \mbox{dev} \, \overline{\boldsymbol{m}}
	+ \mbox{skew} \, \overline{\boldsymbol{m}}
	+ \frac{1}{3} \frac{\mbox{tr}(\overline{\boldsymbol{m}})}{2} \boldsymbol{\mathbbm{1}}
	\right)
	,
	\left(
	\boldsymbol{\mathbbm{1}}
	-
	\boldsymbol{e}_{z}
	\otimes
	\boldsymbol{e}_{z}
	\right)
	\rangle
	\, .
\end{align}
Since $\langle\boldsymbol{\mathbbm{1}},\boldsymbol{\mathbbm{1}}\rangle = 3$, $\langle \boldsymbol{\mathbbm{1}},\left(\boldsymbol{e}_{z}\otimes\boldsymbol{e}_{z}\right)\rangle = 1$, and $\overline{\boldsymbol{m}}$ is decomposed into its three orthogonal components (except for multiplying factors) we can write
\begin{align}
    &
    -
	\langle \left(
	- \mbox{dev} \, \overline{\boldsymbol{m}}
	+ \mbox{skew} \, \overline{\boldsymbol{m}}
	+ \frac{1}{3} \frac{\mbox{tr}(\overline{\boldsymbol{m}})}{2} \boldsymbol{\mathbbm{1}}
	\right)
	,
	\underbrace{
	\left(
	\boldsymbol{\mathbbm{1}}
	-
	\boldsymbol{e}_{z}
	\otimes
	\boldsymbol{e}_{z}
	\right)
	}_{\in \, \mbox{Sym (3)}}
	\rangle
	=
	-
	\langle
	\left(
	- \mbox{dev} \, \overline{\boldsymbol{m}}
	+ \frac{1}{3} \frac{\mbox{tr}(\overline{\boldsymbol{m}})}{2} \boldsymbol{\mathbbm{1}}
	\right)
	,
	\underbrace{
	\left(
	\boldsymbol{\mathbbm{1}}
	-
	\boldsymbol{e}_{z}
	\otimes
	\boldsymbol{e}_{z}
	\right)
	}_{\in \, \mbox{Sym (3)}}
	\rangle
	=
	\\*
	&
	\langle
	\mbox{dev} \, \overline{\boldsymbol{m}}
	,
	\left(
	\boldsymbol{\mathbbm{1}}
	-
	\boldsymbol{e}_{z}
	\otimes
	\boldsymbol{e}_{z}
	\right)
	\rangle
	-
	\frac{1}{3} \frac{\mbox{tr}(\overline{\boldsymbol{m}})}{2}
	\left[
	3-1
	\right]
	=
	-
	\langle
	\mbox{dev} \, \overline{\boldsymbol{m}}
	,
	\boldsymbol{e}_{z}
	\otimes
	\boldsymbol{e}_{z}
	\rangle
	-
	\frac{1}{3} \mbox{tr}(\overline{\boldsymbol{m}})
	=
	\\*
	&
	-
	\langle
	\mbox{dev} \, \overline{\boldsymbol{m}}
	,
	\boldsymbol{e}_{z}
	\otimes
	\boldsymbol{e}_{z}
	\rangle
	-
	\frac{1}{3} \mbox{tr}(\overline{\boldsymbol{m}})
	\langle
	\boldsymbol{\mathbbm{1}}
	,
	\boldsymbol{e}_{z}
	\otimes
	\boldsymbol{e}_{z}
	\rangle
	=
	-
	\langle
	\mbox{dev} \, \overline{\boldsymbol{m}}
	+
	\frac{1}{3} \mbox{tr}(\overline{\boldsymbol{m}}) \boldsymbol{\mathbbm{1}}
	,
	\boldsymbol{e}_{z}
	\otimes
	\boldsymbol{e}_{z}
	\rangle
	=
	\langle
	\overline{\boldsymbol{m}}
	,
	\boldsymbol{e}_{z}
	\otimes
	\boldsymbol{e}_{z}
	\rangle
	=
	-
	\langle
	\overline{\boldsymbol{m}} \, 
	\boldsymbol{e}_{z} ,
	\boldsymbol{e}_{z}
	\rangle
	= - \overline{\boldsymbol{m}}_{zz}
	\, .
	\label{eq:equivalence_Cos_clas}
\end{align}

The last relation
is a pure algebraic relation valid for all $\boldsymbol{m},\overline{\boldsymbol{m}}$ related by  (\ref{eq:ansatz_Cos_classic}) and
$\boldsymbol{e}_{r},\boldsymbol{e}_{\varphi},\boldsymbol{e}_{z}$ are given in  (\ref{eq:polar_unit_vector}).
Thus, we have shown that 
\begin{align}
	\langle
	\mbox{skew}\left(\boldsymbol{m} \times \boldsymbol{e}_{z}\right)
	\boldsymbol{e}_{\varphi} ,
	\boldsymbol{e}_{r}
	\rangle
	-
	\langle
	\mbox{skew}\left(\boldsymbol{m} \times \boldsymbol{e}_{z}\right)
	\boldsymbol{e}_{r} ,
	\boldsymbol{e}_{\varphi}
	\rangle
	=
	\overline{\boldsymbol{m}}_{zz}
    \, .
\end{align}
This solution is valid for a generic second order tensor and for a generic vector triplet.
Using  (\ref{eq:equivalence_Cos_clas}), we finally see that
\begin{align}
    \int_{\Gamma}
    \overline{\boldsymbol{m}}_{zz} \, r \, 
    \mbox{d}r \, \mbox{d}\varphi
    =
    \int_{\Gamma}
    -
    \Big[
	\langle
	\mbox{skew}\left(\boldsymbol{m} \times \boldsymbol{e}_{z}\right)
	\boldsymbol{e}_{\varphi} ,
	\boldsymbol{e}_{r}
	\rangle
	-
	\langle
	\mbox{skew}\left(\boldsymbol{m} \times \boldsymbol{e}_{z}\right)
	\boldsymbol{e}_{r} ,
	\boldsymbol{e}_{\varphi}
	\rangle
	\Big]\, r \,
    \mbox{d}r \, \mbox{d}\varphi
    \, .
\end{align}

The ratio $\Omega$ between the Cosserat torsional stiffness and the classical value that can be found e.g. in \cite{lakes1995experimental,gauthier1975quest,anderson1994size} is
\begin{align}
    \Omega &= 1 + 6 \left( \frac{\ell_t}{R} \right)^2
    \left[
    \frac{1-4/3 \Psi \, \chi}{1 - \Psi \, \chi}
    \right]
    \, ,
    \qquad\qquad\qquad
    \ell_t^2 = \frac{\beta + \gamma}{2\mu _{\tiny \mbox{macro}}}
    \, ,
    \qquad\qquad\qquad
    \Psi = \frac{\beta + \gamma}{\alpha + \beta + \gamma}
    \, ,
    \label{eq:torque_stiffness_Cos_Lakes}
    \\*
    \chi &= \frac{I_1 (p \, R)}{p \, R \, I_0 (p \, R)}
    \, ,
    \qquad\qquad\qquad
    p^2 = \frac{4 \, \mu_c}{\alpha + \beta + \gamma}
    \, ,
    \qquad\qquad\qquad
    N^2 = \frac{\mu_c}{\mu _{\tiny \text{macro}} + \mu _c }
    \, .
    \notag
\end{align}
where $\ell_t$ is the \textit{characteristic length for torsion}, $\Psi$ is the polar ratio, $N$ is the \textit{Cosserat coupling number}, $\alpha$, $\beta$, and $\gamma$ are the curvature coefficients in the classical Cosserat formulation, $\mu_{\tiny \text{macro}}$ is the classical Cauchy shear modulus, $\mu_c$ is the Cosserat couple modulus, and $I_n$ is the modified Bessel function of the first kind of order $n$.
\begin{figure}[H]
		\centering
		\includegraphics[width=0.5\textwidth]{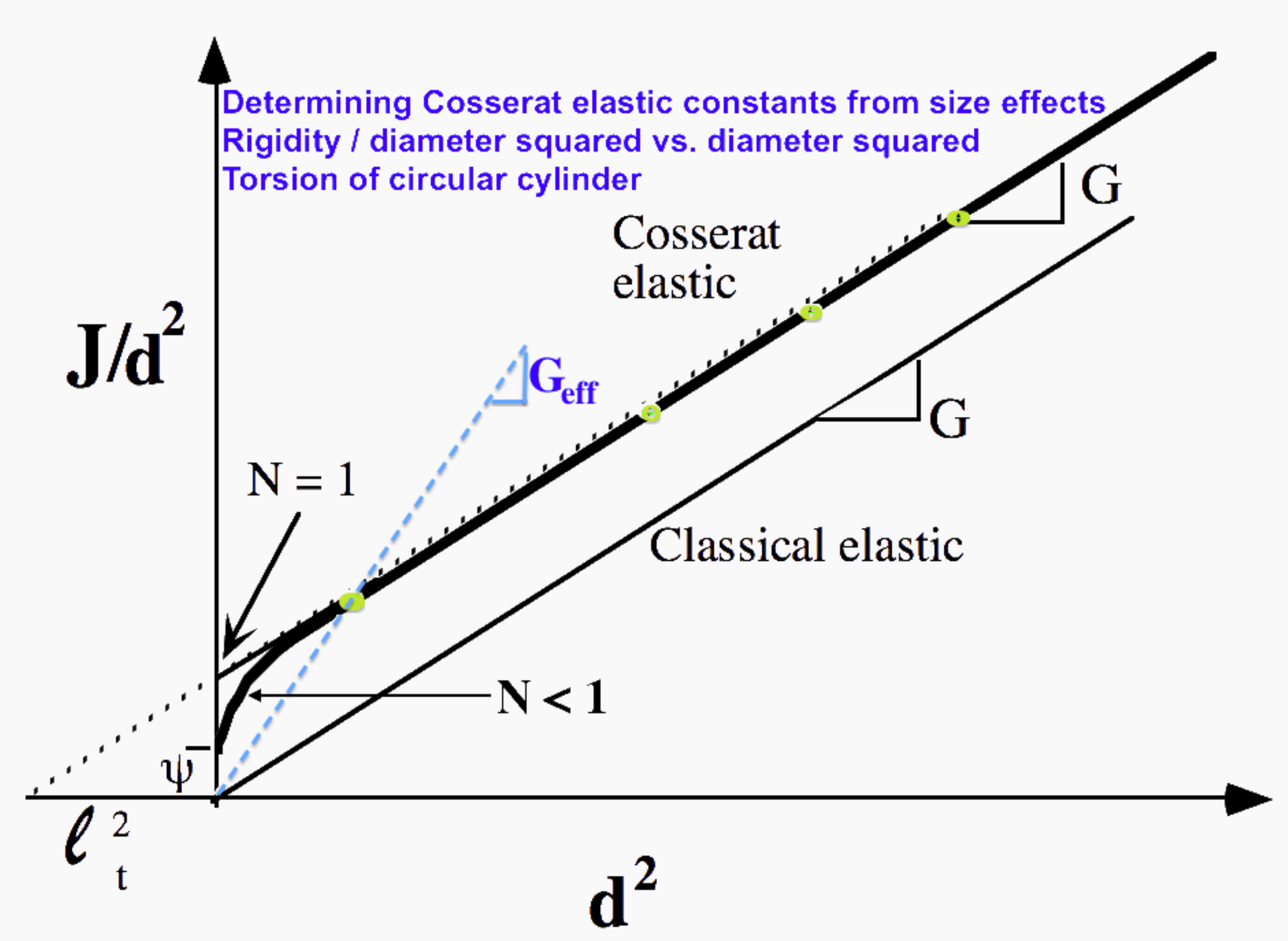}
		\caption{For comparison, the figure has been taken from Lakes \cite{lakes1995experimental}. Here, $G=\mu_{\tiny \text{macro}}$ denotes the classical shear modulus and $N\to1$ corresponds to $\mu_c\to\infty$.}
	\label{fig:Lakes_figure}
\end{figure}

To go from  (\ref{eq:torque_stiffness_Cos}) to  (\ref{eq:torque_stiffness_Cos_Lakes}) we have to use the relations  (\ref{eq:coeff_coss_class}), while remembering to incorporate the term $\mu \, L_c^2$ (the terms not reported do not change between the two notations)
\begin{align}
    \Omega &= 1 + 6 \left( \frac{\ell_t}{R} \right)^2
    \left[
    \frac{1-4/3 \Psi \, \chi}{1 - \Psi \, \chi}
    \right]
    \, ,
    \qquad\qquad
    \ell_t^2 = \frac{a_1}{2 \mu_{\tiny \text{macro}}}
    \, ,
    \qquad\qquad
    \Psi = \frac{3a_1}{2a_1 + 4a_3}
    \, ,
    \qquad\qquad
    p^2 = \frac{6\mu_c}{a_1 + 2a_3}
    \, .
    \label{eq:torque_stiffness_Cos_Lakes_2}
\end{align}
\begin{figure}[H]
	\begin{subfigure}{0.45\textwidth}
		\centering
		\includegraphics[width=\textwidth]{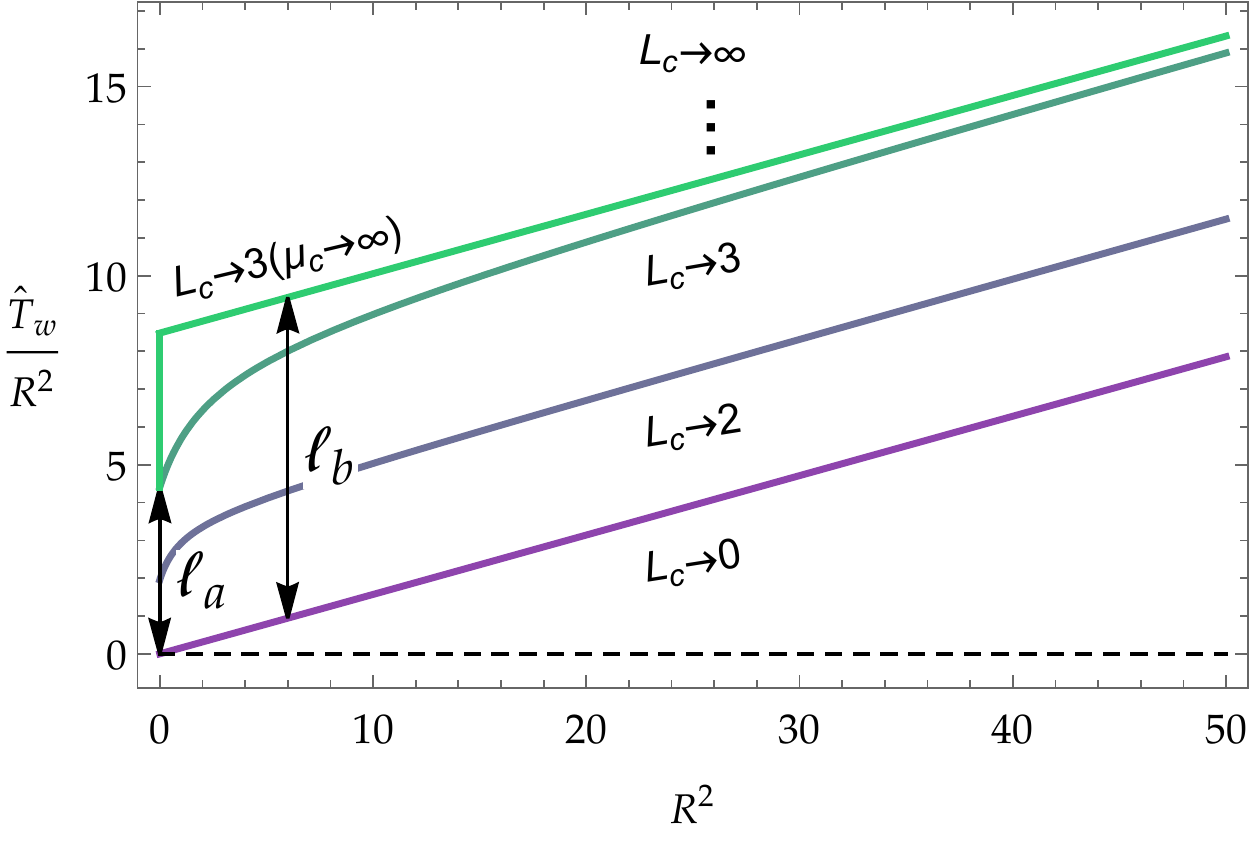}
		\caption{}
	\end{subfigure}
	\hfill
	\begin{subfigure}{0.45\textwidth}
		\centering
		\includegraphics[width=\textwidth]{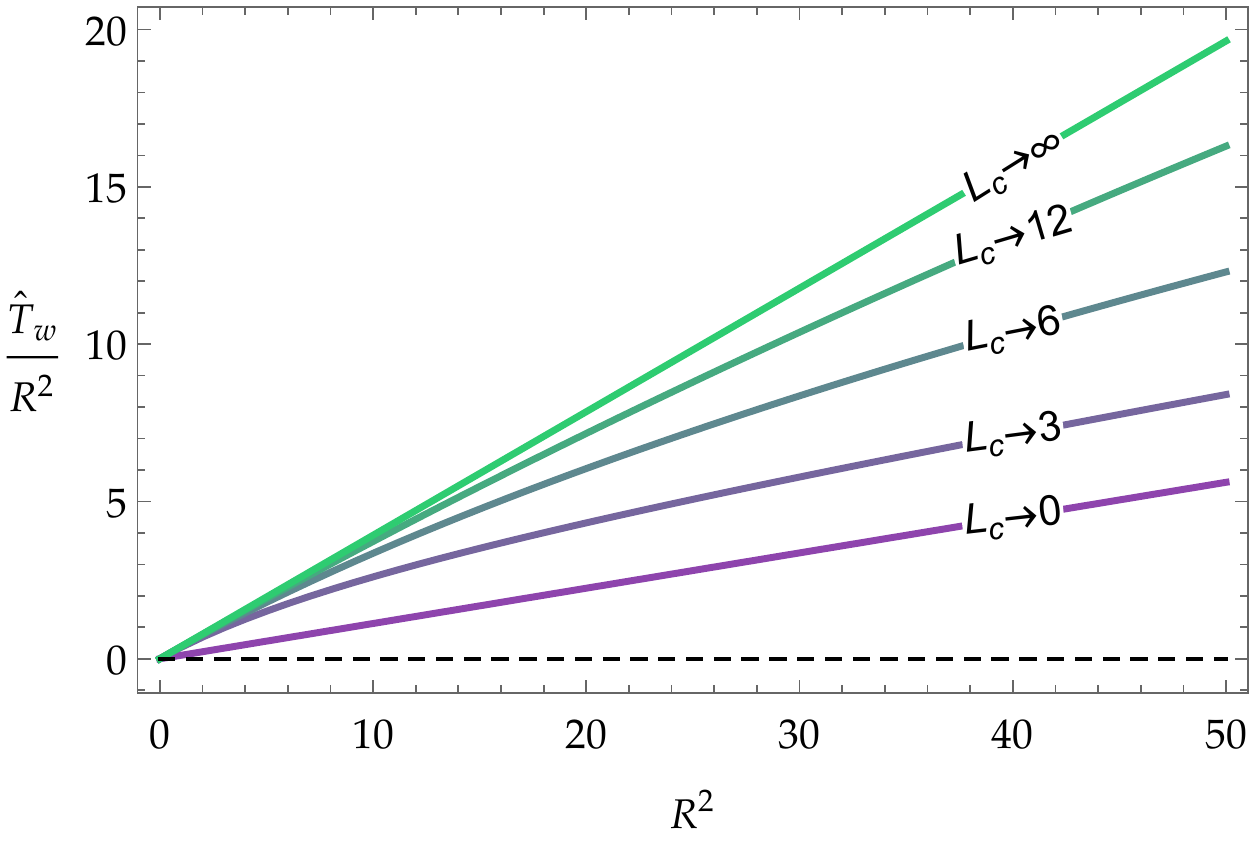}
		\caption{}
	\end{subfigure}
	\caption{(a) Cosserat model; (b) relaxed micromorphic model.
	The values of the coefficient used are:
	$\mu = 1$, $\mu _c = 1/2$, $\mu _{\tiny \text{macro}}=1/14$, $\mu _{\tiny \text{micro}} = 1/4$, $a_1= 1/5$, $a_3=1/37$.}
	\label{fig:torque_stiffness_Lakes_Cos_RMM}
\end{figure}
In Fig. \ref{fig:torque_stiffness_Lakes_Cos_RMM} we report how the torsional stiffness divided by the radius of the cylindrical rod squared ($\widehat{T}_w/R^2$) vary with respect to the radius squared $R^2$ for the Cosserat model and the relaxed micromorphic model
where 
\begin{equation}
    \ell_a =  \mu  \, L_c^{2} \, 12 \pi \frac{a_{1} \, a_{3}}{a_{1} + 8 a_{3}} =
    \frac{\pi  (\beta +\gamma ) (3 \alpha +\beta +\gamma )}{2 \alpha +\beta +\gamma }
     \, ,
    \qquad\qquad\qquad
    \ell_b =  \mu  \, L_c^{2} \, \frac{3}{2} \pi \, a_{1} = \frac{3}{2} \pi  (\beta +\gamma )
    \, .
\end{equation}
It is highlighted that the Cosserat model do not tent to a classical linear elastic model for $\mu_c \to 0$ as it can bee seen from eq.(\ref{eq:torque_stiffness_Cos_Lakes_2}) or eq.(\ref{eq:torque_stiffness_Cos_Limit_muc_0}).
\begin{figure}[H]
	\begin{subfigure}{0.45\textwidth}
		\centering
		\includegraphics[width=\textwidth]{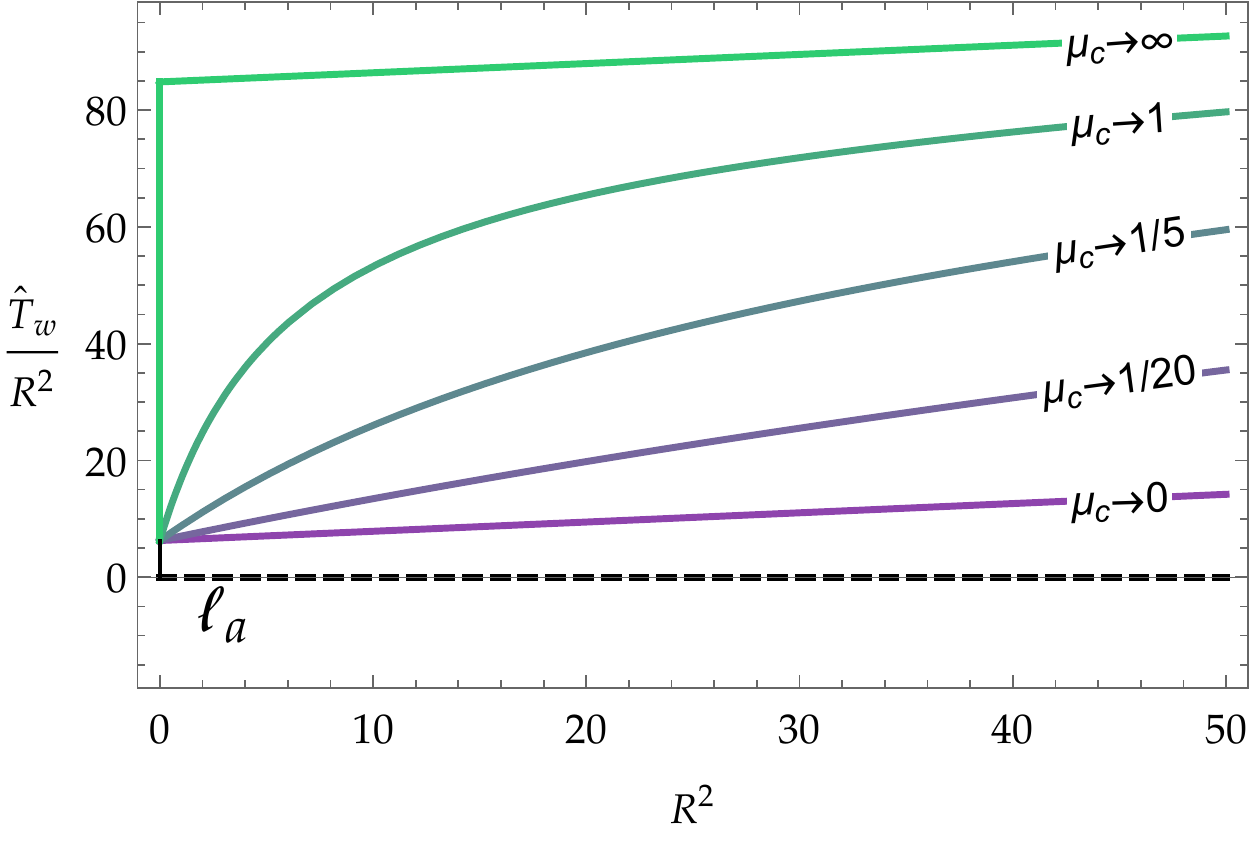}
		\caption{}
	\end{subfigure}
	\hfill
	\begin{subfigure}{0.45\textwidth}
		\centering
		\includegraphics[width=\textwidth]{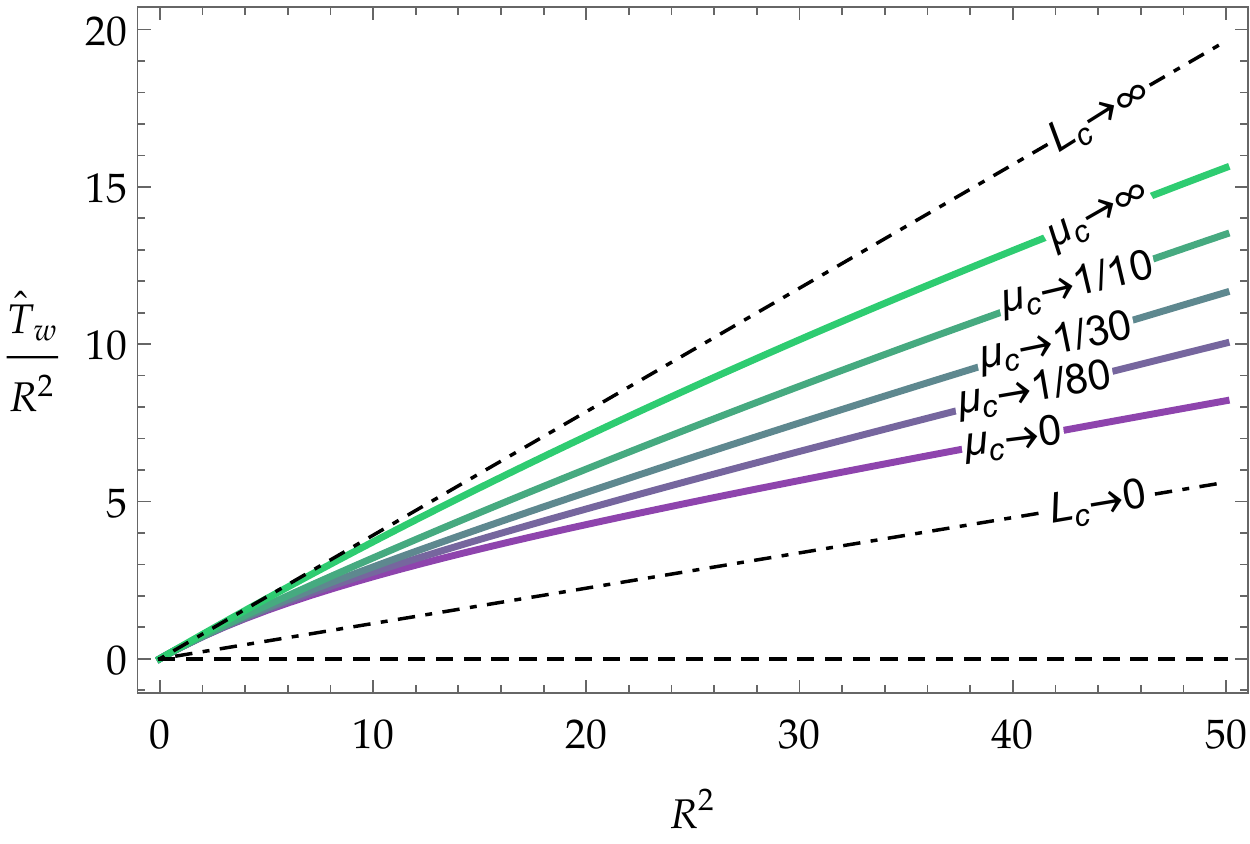}
		\caption{}
	\end{subfigure}
	\caption{
	(a) Cosserat model; (b) relaxed micromorphic model.
	The values of the coefficient used are:
	$\mu =1$, $\mu _{\tiny \text{macro}}=1/14$, $\mu _{\tiny \text{micro}}=1/4$, $a_1=2$, $a_3=1/50$, $L_c=3$.
	In the Cosserat model, the solution for $\mu_c \to \infty$ (the indeterminate couple stress model) shows a jump.
	}
	\label{fig:torque_stiffness_Lakes_Cos_RMM_mc}
\end{figure}
It is underlined that for the relaxed micromorphic model the stiffness is bounded by the one obtained for $L_c \to 0$ (macro) and $L_c \to \infty$ (micro): the macro-stiffness ($L_c \to 0$) is the limit to which all curves with finite $L_c$ tend asymptotically to for $R^2 \to \infty$ (this limit has been cut in order make possible to distinguish all the curves), while the micro-stiffness ($L_c \to \infty$) is the limit to which all the curves tend asymptotically to for $R^2 \to 0$.

\section{Ad-hoc minimization for $L_c \to \infty$ in the full micromorphic model and in the micro-strain model}
\label{app:const_disto}
Looking at the curvature energy of the full micromorphic model (or the micro-strain model) it is clear that for $L_c \to \infty$ the micro-distortion tensor field $\boldsymbol{P}$ must be constant $\boldsymbol{P} = \overline{\boldsymbol{P}}$, provided all curvature coefficients are strictly positive.
We calculate this constant in the following.
Thus we consider
\begin{align}
	\min\limits_{\boldsymbol{u},\overline{\boldsymbol{P}}}
	\left[
	\int_{\Omega}
	\, \mu_{e} \left\lVert \mbox{dev} \, \mbox{sym} \left(\boldsymbol{\mbox{D}u} - \overline{\boldsymbol{P}}\right) \right\rVert^{2}
	+ \frac{\kappa_{e}}{2} \mbox{tr}^2 \left(\boldsymbol{\mbox{D}u} - \overline{\boldsymbol{P}} \right)
	+ \mu_{c} \left\lVert \mbox{skew} \left( \boldsymbol{\mbox{D}u} - \overline{\boldsymbol{P}}\right) \right\rVert^{2}
	\right.
	\\*
	&
	\hspace{-4cm}
	\left.
	+ \mu_{\tiny \mbox{micro}} \left\lVert \mbox{dev} \, \mbox{sym} \, \overline{\boldsymbol{P}} \right\rVert^{2} 
	+ \frac{\kappa_{\tiny \mbox{micro}}}{2} \mbox{tr}^2 \left(\overline{\boldsymbol{P}} \right)
	\mbox{d}V
	\right]
	 \, .
	\notag
\end{align}
The weak form is given by
\begin{align}
	\int_{\Omega}
	2\mu_{e} \left\langle \mbox{dev} \, \mbox{sym} \left(\boldsymbol{\mbox{D}u} - \overline{\boldsymbol{P}}\right), -\delta \overline{\boldsymbol{P}}\right\rangle
	+ \kappa_{e} \, \mbox{tr} \left(\boldsymbol{\mbox{D}u} - \overline{\boldsymbol{P}} \right)
	\left\langle \boldsymbol{\mathbbm{1}} , -\delta \overline{\boldsymbol{P}}\right\rangle
	+ 2\mu_{c} \, \left\langle \mbox{skew} \left( \overline{\boldsymbol{P}}\right), -\delta \overline{\boldsymbol{P}}\right\rangle
	\hspace{+3.8cm}
	\\*
	+ 2\mu_{\tiny \mbox{micro}} \, \left\langle \mbox{dev} \, \mbox{sym} \left( \overline{\boldsymbol{P}}\right), \delta \overline{\boldsymbol{P}}\right\rangle
	+ \kappa_{\tiny \mbox{micro}} \, \mbox{tr} \left(\overline{\boldsymbol{P}} \right)
	\left\langle \boldsymbol{\mathbbm{1}} , \delta \overline{\boldsymbol{P}}\right\rangle
	\mbox{d}V
	=
	0 \qquad \forall \, \delta \overline{\boldsymbol{P}} \, .
    \notag
\end{align}

\begin{align}
	\int_{\Omega}
	\left\langle
	2\mu_{e} \, \mbox{dev} \, \mbox{sym} \left(\boldsymbol{\mbox{D}u} - \overline{\boldsymbol{P}}\right)
	+ \kappa_{e} \, \mbox{tr} \left(\boldsymbol{\mbox{D}u} - \overline{\boldsymbol{P}} \right)
	\boldsymbol{\mathbbm{1}}
	+ 2\mu_{c} \, \mbox{skew} \left(\boldsymbol{\mbox{D}u} - \overline{\boldsymbol{P}}\right)
	\right.
	\hspace{+6cm}
	\\*
	\left.
	- 2\mu_{\tiny \mbox{micro}} \, \mbox{dev} \, \mbox{sym} \left( \overline{\boldsymbol{P}}\right)
	- \kappa_{\tiny \mbox{micro}} \, \mbox{tr} \left(\overline{\boldsymbol{P}} \right) \boldsymbol{\mathbbm{1}}
	 , \delta \overline{\boldsymbol{P}}
	 \right\rangle
	\mbox{d}V
	=
	0 \qquad \forall \, \delta \overline{\boldsymbol{P}}
	\, .
	\notag
\end{align}
For constant $\delta \overline{\boldsymbol{P}}$ this can be rewritten as
\begin{align}
	\left\langle 
	\int_{\Omega}
	2\mu_{e} \, \mbox{dev} \, \mbox{sym} \, \left(\boldsymbol{\mbox{D}u} - \overline{\boldsymbol{P}}\right)
	+ \kappa_{e} \, \mbox{tr} \left(\boldsymbol{\mbox{D}u} - \overline{\boldsymbol{P}} \right)
	\boldsymbol{\mathbbm{1}}
	+ 2\mu_{c} \, \mbox{skew} \left(\boldsymbol{\mbox{D}u} - \overline{\boldsymbol{P}}\right)
	\right.
	\hspace{+6cm}
	\\*
	\left.
	- 2\mu_{\tiny \mbox{micro}} \, \mbox{dev} \, \mbox{sym} \left( \overline{\boldsymbol{P}}\right)
	- \kappa_{\tiny \mbox{micro}} \, \mbox{tr} \left(\overline{\boldsymbol{P}} \right) \boldsymbol{\mathbbm{1}} 
	\, \mbox{d}V
    , \delta \overline{\boldsymbol{P}}\right\rangle
	=
	0 \qquad \forall \, \delta \overline{\boldsymbol{P}}
	\, .
	\notag
\end{align}
Since $\delta \overline{\boldsymbol{P}}$ is arbitrary, this implies that 
\begin{align}
	\int_{\Omega}
	2\mu_{e} \, \mbox{dev} \, \mbox{sym} \left(\boldsymbol{\mbox{D}u} - \overline{\boldsymbol{P}}\right)
	+ \kappa_{e} \, \mbox{tr} \left(\boldsymbol{\mbox{D}u} - \overline{\boldsymbol{P}} \right)
	\boldsymbol{\mathbbm{1}}
	+ 2\mu_{c} \, \mbox{skew} \left(\boldsymbol{\mbox{D}u} - \overline{\boldsymbol{P}}\right)
	- 2\mu_{\tiny \mbox{micro}} \, \mbox{dev}  \, \mbox{sym} \left( \overline{\boldsymbol{P}}\right)
	- \kappa_{\tiny \mbox{micro}} \, \mbox{tr} \left(\overline{\boldsymbol{P}} \right) \boldsymbol{\mathbbm{1}} 
	\, \mbox{d}V
	=
	0
	\, ,
\end{align}
or
\begin{align}
	\int_{\Omega}
	2\mu_{e} \, \mbox{dev} \, \mbox{sym} \, \boldsymbol{\mbox{D}u}
	+ \kappa_{e} \, \mbox{tr} \left(\boldsymbol{\mbox{D}u} \right)
	\boldsymbol{\mathbbm{1}}
	+ 2\mu_{c} \, \mbox{skew} \, \boldsymbol{\mbox{D}u}
	\, \mbox{d}V
	=
	\hspace{+8.3cm}
	\\*
	=
	\int_{\Omega}
	2\mu_{e} \, \mbox{dev} \, \mbox{sym} \, \overline{\boldsymbol{P}}
	+ \kappa_{e} \, \mbox{tr} \left( \overline{\boldsymbol{P}} \right)
	\boldsymbol{\mathbbm{1}}
	+ 2\mu_{c} \, \mbox{skew} \, \overline{\boldsymbol{P}}
	+ 2\mu_{\tiny \mbox{micro}} \, \mbox{dev} \, \mbox{sym} \, \overline{\boldsymbol{P}}
	+ \kappa_{\tiny \mbox{micro}} \, \mbox{tr} \left(\overline{\boldsymbol{P}} \right) \boldsymbol{\mathbbm{1}} 
	\, \mbox{d}V
	\, .
	\notag
\end{align}
Using the orthogonality of ${\rm dev\, sym \cdot},\, {\rm skew} \cdot$ and $\tr(\cdot)\,\boldsymbol{\mathbbm{1}}$ we obtain
\begin{gather}
	\int_{\Omega}
	2\mu_{e} \, \mbox{dev} \, \mbox{sym} \, \boldsymbol{\mbox{D}u}
	\, \mbox{d}V
	=
	\int_{\Omega}
	2\mu_{e} \, \mbox{dev} \, \mbox{sym} \, \overline{\boldsymbol{P}}
	+ 2\mu_{\tiny \mbox{micro}} \, \mbox{dev} \, \mbox{sym} \, \overline{\boldsymbol{P}}
	\, \mbox{d}V
	\, ,
	\\
	\int_{\Omega}
	\kappa_{e} \, \mbox{tr} \left(\boldsymbol{\mbox{D}u} \right)
	\, \mbox{d}V
	=
	\int_{\Omega}
	\kappa_{e} \, \mbox{tr} \left( \overline{\boldsymbol{P}} \right)
	+ \kappa_{\tiny \mbox{micro}} \, \mbox{tr} \left(\overline{\boldsymbol{P}} \right) 
	\, \mbox{d}V
	\, ,
	\qquad
	\int_{\Omega}
	2\mu_{c} \, \mbox{skew} \, \boldsymbol{\mbox{D}u}
	\, \mbox{d}V
	=
	\int_{\Omega}
	2\mu_{c} \, \mbox{skew} \, \overline{\boldsymbol{P}}
	\, \mbox{d}V
	\, ,
	\notag
\end{gather}
and since $\overline{\boldsymbol{P}}$ is constant we can write 
\begin{gather}
	\mbox{dev} \, \mbox{sym} \, \overline{\boldsymbol{P}}
	=
	\frac{1}{\left | \Omega \right |}
	\int_{\Omega}
	\frac{\mu_{e}}{\mu_{e}+\mu_{\tiny \mbox{micro}}} \, \mbox{dev} \, \mbox{sym} \, \boldsymbol{\mbox{D}u}
	\, \mbox{d}V
	\, ,
	\qquad
	\mbox{tr} \left( \overline{\boldsymbol{P}}\right)
	=
	\frac{1}{\left | \Omega \right |}
	\int_{\Omega}
	\frac{\kappa_{e}}{\kappa_{e}+\kappa_{\tiny \mbox{micro}}} \, \mbox{tr} \left(\boldsymbol{\mbox{D}u}\right)
	\, \mbox{d}V
	\, ,
	\label{eq:microdist_constant_MM_0}
	\\*
	\mbox{skew} \, \overline{\boldsymbol{P}}
	=
	\frac{1}{\left | \Omega \right |}
	\int_{\Omega}
	\mbox{skew} \, \boldsymbol{\mbox{D}u}
	\, \mbox{d}V
	\, .
	\notag
\end{gather}
Since  dev sym, skew, and tr are linear operators, we obtain equivalently 
\begin{gather}
	\mbox{dev} \, \mbox{sym} \, \overline{\boldsymbol{P}}
	=
	\frac{\mu_{e}}{\mu_{e}+\mu_{\tiny \mbox{micro}}} \, 
	\mbox{dev} \, \mbox{sym} \, 
	\left(
	\frac{1}{\left | \Omega \right |} \,
	\int_{\Omega}
	\boldsymbol{\mbox{D}u}
	\, \mbox{d}V
	\right)
	\, ,
	\qquad
	\mbox{tr} \left( \overline{\boldsymbol{P}}\right)
	=
	\frac{\kappa_{e}}{\kappa_{e}+\kappa_{\tiny \mbox{micro}}} \, 
	\mbox{tr} \left(
	\frac{1}{\left | \Omega \right |}
	\,
	\int_{\Omega}
	\boldsymbol{\mbox{D}u}
	\, \mbox{d}V
	\right)
	\, ,
	\label{eq:microdist_constant_MM}
	\\*
	\mbox{skew} \, \overline{\boldsymbol{P}}
	=
	\mbox{skew} \, 
	\left(
	\frac{1}{\left | \Omega \right |} \,
	\int_{\Omega}
	\boldsymbol{\mbox{D}u}
	\, \mbox{d}V
	\right)
	\, .
	\notag
\end{gather}

Substituting the ansatz  (\ref{eq:ansatz_MM}) into  (\ref{eq:microdist_constant_MM}) we obtain $\overline{\boldsymbol{P}}=\boldsymbol{0}$.
Analogous calculations can be carried out for the micro-strain model for which $\mbox{skew} \, \overline{\boldsymbol{P}}=0$ and $\mu_c=0$
\begin{gather}
	\mbox{dev} \, \overline{\boldsymbol{S}}
	=
	\frac{\mu_{e}}{\mu_{e}+\mu_{\tiny \mbox{micro}}}
	\mbox{dev} \, \mbox{sym}
	\left(
	\frac{1}{\left | \Omega \right |}
	\, 
	\int_{\Omega}
	\boldsymbol{\mbox{D}u}
	\, \mbox{d}V
	\right)
	\, ,
	\qquad
	\mbox{tr} \left( \overline{\boldsymbol{S}}\right)
	=
	\frac{\kappa_{e}}{\kappa_{e}+\kappa_{\tiny \mbox{micro}}}
	\mbox{tr} \left(
	\frac{1}{\left | \Omega \right |}
	\, 
	\int_{\Omega}
	\boldsymbol{\mbox{D}u}
	\, \mbox{d}V
	\right)
	\, .
	\label{eq:microdist_constant_MStrain}
\end{gather}
Substituting the ansatz  (\ref{eq:ansatz_MStrain}) into  (\ref{eq:microdist_constant_MStrain}) we obtain $\overline{\boldsymbol{S}}=\boldsymbol{0}$.

The integral on the circular cross-section $\Gamma$ of the gradient of the displacement is
\begin{equation}
    \int_{\Gamma} \mbox{D}\boldsymbol{u} (x) \mbox{d}V
    =
    \int_{0}^{2\pi}\int_{0}^{R}
    \mbox{D}\boldsymbol{u} (r,\varphi,z) \, r \, \mbox{d}r \mbox{d}\varphi
    =
    \left(
    \begin{array}{ccc}
    0 & - \pi R^2 \boldsymbol{\vartheta} z & 0 \\
    \pi R^2 \boldsymbol{\vartheta} z & 0 & 0 \\
    0 & 0 & 0 \\
    \end{array}
    \right)
    =
    \pi R^2 \boldsymbol{\Theta}(z)
    \left(
    \begin{array}{ccc}
    0 & -1 & 0 \\
    1 & 0 & 0 \\
    0 & 0 & 0 \\
    \end{array}
    \right)
    \, .
    \label{eq:int_grad}
\end{equation}

From  (\ref{eq:int_grad}) it is possible to see that the symmetric part of the integral of $\mbox{D}\boldsymbol{u}$ on the circular cross-section is zero, while the skew-symmetric part is zero only if the domain is symmetric with respect $z$.

For the Cosserat model, letting $L_c \to \infty$ still implies that $\boldsymbol{A}(\boldsymbol{x})=\overline{\boldsymbol{A}}=\mbox{const.}$ must be constant.
The same calculations as before yield
\begin{align}
	\overline{\boldsymbol{A}}
	=
	\mbox{skew} \, 
	\left(
	\frac{1}{\left | \Omega \right |} \,
	\int_{\Omega}
	\boldsymbol{\mbox{D}u}
	\, \mbox{d}V
	\right)
	=
	\frac{1}{\left | \Omega \right |} \,
	\int_{\Omega}
	\boldsymbol{\mbox{D}u}
	\, \mbox{d}V
	\, ,
	\qquad\qquad
	\mbox{since}
	\quad
	\mbox{D}\boldsymbol{u}\in \mathfrak{so}(3) \, .
	\label{eq:microdist_constant_Cos}
\end{align}
For $\Omega = \left[ 0,L \right] \times \Gamma$ we have
\begin{align}
    \frac{1}{\left | \Omega \right |} \,
	\int_{\Omega}
	\boldsymbol{\mbox{D}u}
	\, \mbox{d}V
	=
	\frac{1}{L \left( \pi R^2 \right)}
	\int_{0}^{L}
	\int_{0}^{2\pi}
	\int_{0}^{R}
	\mbox{D} \boldsymbol{u} \left( r,\varphi,z \right) \, r \, \mbox{d}r\mbox{d}\varphi\mbox{d}z
	=
	\frac{1}{L \left( \pi R^2 \right)}
	\left(
	\left.
	\boldsymbol{\vartheta} \, \pi \, R^2 \frac{z^2}{2}
	\right |_{0}^{L}
	\right)
	=
	\frac{\boldsymbol{\vartheta}}{2} L
	=
	\frac{1}{2} \boldsymbol{\Theta}(L)
	\, .
\end{align}
We remark that the same limit $L_c \to \infty$ in the relaxed micromorphic model yields a linear elastic response with stiffness $\mathbbm{C}_{\tiny \mbox{micro}}$ since $\mbox{Curl} \boldsymbol{P=0}$ does not imply that $P=$ const. but $P=\boldsymbol{\nabla \zeta}$ for some $\zeta:\Omega \in \mathbb{R}^3 \to \mathbb{R}^3$, see \cite{neff2019identification}.

\end{footnotesize}
\end{document}